\RequirePackage{fix-cm}
\documentclass[smallcondensed]{svjour3}                     % onecolumn (standard format)
\smartqed  % flush right qed marks, e.g. at end of proof

\usepackage{listings}
\usepackage{color,xcolor}
\usepackage[square,numbers]{natbib}
\usepackage{tikz}
\usetikzlibrary{positioning,chains,fit,shapes,calc}
\definecolor{grey1}{rgb}{0.0,0.0,0.0}
\definecolor{grey2}{rgb}{0.2,0.2,0.2}
\definecolor{grey3}{rgb}{0.4,0.4,0.4}
\definecolor{grey4}{rgb}{0.5,0.5,0.5}

\definecolor{myblue}{RGB}{100,100,160}
\definecolor{mygreen}{RGB}{80,160,80}

\lstdefinestyle{highlight-fonts}{
	basicstyle=\ttfamily\mdseries\footnotesize,
	commentstyle=\ttfamily\mdseries\itshape\footnotesize,
	moredelim=[l][\color{\hiddenA}\ttfamily\mdseries\itshape\footnotesize]{\%\#},
	moredelim=[l][\color{\hiddenB}\ttfamily\mdseries\itshape\footnotesize]{\%!},
	keywordstyle=\ttfamily\mdseries\footnotesize,
	keywordstyle={[2]\ttfamily\bfseries\footnotesize},
}

\lstdefinestyle{highlight-colors}{
	basicstyle=\color{grey1}\ttfamily\mdseries\footnotesize,
	commentstyle=\color{grey2}\ttfamily\mdseries\itshape\footnotesize,
	moredelim=[l][\color{\hiddenA}\ttfamily\mdseries\itshape\footnotesize]{\%\#},
	moredelim=[l][\color{\hiddenB}\ttfamily\mdseries\itshape\footnotesize]{\%!},
	keywordstyle=\color{grey2}\ttfamily\mdseries\footnotesize,
	keywordstyle={[2]\color{grey2}\ttfamily\bfseries\footnotesize},
}

\lstdefinestyle{gobble1}{
	xleftmargin=0.3em,
}

\lstdefinestyle{gobble2}{
	xleftmargin=-0.9em,
}

\lstdefinestyle{gobble3}{
	xleftmargin=-2.1em,
}

\lstdefinelanguage{coco}[]{Matlab}{
	keywords={},
	morekeywords={},
	morekeywords=[2]{classdef,properties,private,protected,public,%
Access,Static,methods,if,function,end,for,while,else,elseif,switch,%
case,otherwise,do,repeat,until},
	otherkeywords={{,end},:end,(end,end),\{end,end\},[end,end],/private,private/},
	morecomment=[l]{},
	basicstyle=\ttfamily\mdseries\footnotesize,
	commentstyle=\ttfamily\mdseries\footnotesize,
	moredelim=[l][\color{\hiddenA}\ttfamily\mdseries\footnotesize]{\%\#},
	moredelim=[l][\color{\hiddenB}\ttfamily\mdseries\footnotesize]{\%!},
	keywordstyle=\ttfamily\mdseries\footnotesize,
	keywordstyle={[2]\ttfamily\mdseries\footnotesize},
	numbers=none,
	numberstyle=\scriptsize,
	numbersep=1em,
	breaklines=false,
	breakatwhitespace=true,
	breakindent=2.5em,
	showlines=false,
	lineskip=-0.2ex,
	frame=none,
	fontadjust=true,
	columns=[c]fixed,
	basewidth={0.575em,0.45em},
	fontadjust=true,
	tabsize=3,
	showstringspaces=false,
	aboveskip=1.5\medskipamount,
	belowskip=1.5\medskipamount,
	xleftmargin=0.925em,
	xrightmargin=0em,
	rangeprefix={\%!},
	includerangemarker=false,
	belowcaptionskip=\bigskipamount
}

\lstdefinelanguage{coco-highlight-fonts}[]{coco}{
	style=highlight-fonts,
}

\lstdefinelanguage{coco-highlight-colors}[]{coco}{
	style=highlight-colors,
}

\lstdefinelanguage{coco-highlight}[]{coco-highlight-colors}{}

\newcommand{\mcode}[1]{{\lstinline[language=coco-highlight,basewidth={0.6em,0.45em}]|#1|}}

\renewcommand{\d}{\mathrm{d}}

\definecolor{dblue}{rgb}{0,0,0.5}

\renewcommand{\d}{\mathrm{d}}
\newcommand{\R}{\mathbb{R}}
\usepackage{amsmath,dsfont}

\newcommand{\pcn}{P_\mathrm{cn}}
\newcommand{\pbp}{P_\mathrm{bp}}
\newcommand{\pcont}{P_\mathrm{cont}}

\DeclareMathOperator{\prc}{PRC}

\tikzset{
  -|-/.style={
    to path={
      (\tikztostart) -| ($(\tikztostart)!#1!(\tikztotarget)$) |- (\tikztotarget)
      \tikztonodes
    }
  },
  -|-/.default=0.5,
  |-|/.style={
    to path={
      (\tikztostart) |- ($(\tikztostart)!#1!(\tikztotarget)$) -| (\tikztotarget)
      \tikztonodes
    }
  },
  |-|/.default=0.5,
}

\tikzstyle{startstop} = [rectangle, rounded corners, minimum width=3cm, minimum height=1cm,text centered, draw=black, fill=red!30]
\tikzstyle{io} = [trapezium, trapezium left angle=70, trapezium right angle=110, minimum width=3cm, minimum height=1cm, text centered, draw=black, fill=blue!30]
\tikzstyle{process} = [rectangle, minimum width=3cm, minimum height=1cm, text centered, draw=black, fill=orange!30]
\tikzstyle{process-phi} = [rectangle, minimum width=3cm, minimum height=1cm, text centered, draw=black, fill=blue!30]
\tikzstyle{process-psi} = [rectangle, minimum width=3cm, minimum height=1cm, text centered, draw=black, fill=green!30]
\tikzstyle{process-lambda} = [rectangle, minimum width=3cm, minimum height=1cm, text centered, draw=black, fill=orange!30]
\tikzstyle{decision} = [diamond, minimum width=3cm, minimum height=0.5cm, text centered, draw=black, fill=red!30]
\tikzstyle{arrow} = [thick,->,>=stealth]

\usepackage{graphicx}
\usepackage{amsmath,amsfonts,amssymb}
 \usepackage{mathptmx}      % use Times fonts if available on your TeX system
\usepackage{subfig}
\usepackage{stmaryrd}
\usepackage{paralist,multirow}
\usepackage{isomath}
\usepackage{mathtools}
\usepackage{float}
\usepackage{algorithmic}
\usepackage{marvosym} %for \Letter

\usepackage
[
        a4paper,% other options: a3paper, a5paper, etc
        left=1.5cm,
        right=2cm,
        bottom = 1.44in,
]
{geometry}

\usepackage[bottom]{footmisc}

\journalname{Nonlinear Dynamics}
\begin{document}
%\tolerance 1414
%\emergencystretch
\title{Methods of Continuation and their Implementation in the \textsc{COCO} Software Platform with Application to Delay Differential Equations}

\titlerunning{Methods of Continuation and their Implementation in \textsc{coco}}        % if too long for running head

\author{Zaid Ahsan$^{*}$ \and Harry Dankowicz \and Mingwu Li \and Jan Sieber
}


\institute{ Z. Ahsan \and H. Dankowicz \at Department of Mechanical Science and Engineering\\University of Illinois at Urbana-Champaign\\Urbana, Illinois, USA\\\email{zaid2@illinois.edu, danko@illinois.edu} \and M. Li \at Institute for Mechanical Systems\\ETH Z\"{u}rich\\ Z\"{u}rich,  Switzerland\\
\email{mingwli@ethz.ch}\\ \and	J. Sieber \at College of Engineering, Mathematics and Physical Sciences\\University of Exeter\\ Exeter, United Kingdom\\\email{j.sieber@exeter.ac.uk}
}

\date{Received: date / Accepted: date}

\maketitle

\begin{abstract}
This paper treats comprehensively the construction of problems from nonlinear dynamics and constrained optimization amenable to parameter continuation techniques and with particular emphasis on multi-segment boundary-value problems with delay. The discussion is grounded in the context of the \textsc{coco} software package and its explicit support for community-driven development. To this end, the paper first formalizes the \textsc{coco} construction paradigm for augmented continuation problems compatible with simultaneous analysis of implicitly defined manifolds of solutions to nonlinear equations and the corresponding adjoint variables associated with optimization of scalar objective functions along such manifolds. The paper uses applications to data assimilation from finite time histories and phase response analysis of periodic orbits to identify a universal paradigm of construction that permits abstraction and generalization. It then details the theoretical framework for a \textsc{coco}-compatible toolbox able to support the analysis of a large family of delay-coupled multi-segment boundary-value problems, including periodic orbits, quasiperiodic orbits, connecting orbits, initial-value problems, and optimal control problems, as illustrated in a suite of numerical examples. The paper aims to present a pedagogical treatment that is accessible to the novice and inspiring to the expert by appealing to the many senses of the applied nonlinear dynamicist. Sprinkled among a systematic discussion of problem construction, graph representations of delay-coupled problems, and vectorized formulas for problem discretization, the paper includes an original derivation using Lagrangian sensitivity analysis of phase-response functionals for periodic-orbit problems in abstract Banach spaces, as well as a demonstration of the regularizing benefits of multi-dimensional manifold continuation for near-singular problems analyzed using real-time experimental data.

\keywords{delay differential equations \and multi-segment boundary-value problems \and  implicitly defined manifolds \and problem regularization \and constrained optimization \and Lagrange multipliers \and adjoint equations \and graph representations \and phase response curves}
\let\thefootnote\relax\footnote{ \large{\textbf{Note}: This preprint has not undergone peer review (when
applicable) or any post-submission improvements or corrections. The Version of Record of this
article is published in \textbf{Nonlinear Dynamics}, and is available online at https://doi.org/10.1007/s11071-021-06841-1}}
%\begin{table}[h!] \centering
%\begin{tabular}{c} 
%Note: Manuscript has been submitted to Nonlinear Dynamics aaaaaaaaaaaaaaaaaaaaaaaaaaaaaaaaaaaaaaaaaaaaaaaaaaaaaaaaaa
%\end{tabular}
%\end{table}
%\skip\footins=\bigskipamount

\end{abstract}

\newgeometry{bottom=1.44in,right=2cm,left=1.5cm}

\section{Introduction}
\label{sec: Introduction}
We use this section to describe the objectives of this paper and to place its content in the context of the ongoing development of software tools for continuation-based analysis of nonlinear dynamical systems. An overview of the original contributions along with section-wise descriptions of the content is also provided.

\subsection{Motivation}
\label{sec: motivation}
While the possibility of closed-form analysis is fortuitous and perhaps career-changing, computation is the bread and butter of applied research in nonlinear dynamics. Computational techniques derive their power from rigorous mathematical analysis, but extend far beyond the reach of theoretical tools (see, e.g., the analysis of global bifurcations of the Lorenz manifold in~\cite{doedel2006global,guckenheimer2015invariant}).  With their aid, systematic exploration becomes possible, e.g., of the dependence of system responses on model parameters~\cite{abbas2011parametric}, the sensitivity of these responses to parameter uncertainty~\cite{kewlani2012polynomial}, and the determination of optimal selections of parameter values~\cite{amandio2014stochastic,koh2016optimizing}. Such exploration inspires further theoretical advances, including of methods for projecting the dynamics of large-scale systems onto reduced-order models~\cite{haller2016nonlinear,szalai2019model,szalai2020invariant,touze2006nonlinear}, amenable to efficient computation and powerful visualization.

Continuation methods are a class of deterministic computational techniques for exploring smooth manifolds of solutions to nonlinear equations~\cite{dankowicz2013recipes,krauskopf2007numerical}. By now classical algorithms convert common questions of interest to applied dynamicists into such nonlinear equations, enabling their analysis using continuation. Prominent among such uses are bifurcation analyses of equilibria~\cite{govaerts2000numerical}, periodic orbits~\cite{munoz2003continuation}, connecting orbits~\cite{beyn1990numerical}, quasiperiodic invariant tori~\cite{schilder2005continuation}, and stable and unstable manifolds~\cite{dellnitz1996computation} for smooth and piecewise-smooth vector fields, including in problems with delay~\cite{barton2009stability,chavez2020numerical,luzyanina1997computation,roose2007continuation}. By their versatility, continuation methods are an invaluable tool in the researcher's arsenal.

It is the aim of this review to invite new generations of dynamicists to the world of continuation methods while also giving the seasoned researcher plenty of original fodder for thought. The paper purposely avoids stepping over well-trodden ground dealing with the specific algorithms used to cover solution manifolds or with examples of bifurcation analysis, as these topics have been discussed in great detail in a number of key sources~\cite{allgower2003introduction,kuznetsov2013elements}. Instead, completely new content develops a formalism for problem construction, inspired by functionality available in the \textsc{coco} software platform~\cite{COCO}, illustrates its use on problems from data assimilation and phase response analysis, and applies its principles to the detailed construction of a toolbox for analyzing a large class of (possibly non-autonomous) multi-segment boundary-value problems with delay.

Given the ambitious scope and the anticipated range of reader expertise, the paper is intentionally self-contained and designed to address both deductive and inductive learning styles. Emphasis is placed on a formalism that translates directly to computational encoding, say, in the \textsc{coco} framework. Examples are drawn from the existing literature but retrofitted to the universal abstractions proposed in this manuscript. Where appropriate, the text highlights opportunities for practice or more substantial further development. For that is ultimately the measure of this paper's value; the extent to which it spurs original creativity and innovation.

\subsection{Software design}
\label{sec: software design}
There is a natural tension in both theoretical and computational research between the particular and the general. One of us (HD) spent half of his professional career implementing continuation methods on a just-in-time basis. With each new problem, old code scripts were dusted off, debugged, and redeployed. The investment to do anything beyond solving particular problems seemed overwhelming and perhaps not professionally rewarding. But the general beckoned. And others had long before chosen that road, aiming to build translational tools that would bring nonlinear dynamics to the scientific masses~\cite{back1992dstool,kuznetsov1997content}. In so doing, a balance again had to be struck between the particular and the general, between utility and universality. To this day, several of the outcomes of this effort continue to provide invaluable access to the insights of the qualitative theory of dynamical system also to non-experts. These include \textsc{auto}~\cite{doedel2007auto} and the wrapper \textsc{xppaut}~\cite{ermentrout2002simulating}, a popular choice for continuation-based analysis of ordinary differential equations (ODEs), \textsc{matcont}~\cite{dhooge2003matcont} for ODEs/maps, \textsc{dde-biftool}~\cite{engelborghs2002numerical,ddebiftoolmanual} and \textsc{knut}~\cite{Kunt} for delay differential equations (DDEs), \textsc{pde2path}~\cite{kuehn2015efficient,uecker2014pde2path} for partial differential equations (PDEs), and \textsc{hompack}~\cite{watson1987algorithm} for globally-convergent homotopy analysis of arbitrary nonlinear equations. 

With an emphasis on utility, these packages were designed to address specific problem classes/types, while leaving open the possibility of additional creative uses (e.g., the development of special purpose wrappers for \textsc{auto} for computing invariant manifolds~\textsc{manbvp}~\cite{england2005computing}, bifurcation analysis of Filippov systems~\textsc{slidecont}~\cite{dercole2005slidecont}, bifurcation analysis of periodic orbits in hybrid dynamical systems~$\widehat{\text{{\sc tc}}}$~\cite{thota2008tc}, and the computation of global isochrons in~\cite{osinga2010continuation}). In contrast, an emphasis on universality was the guiding principle behind the creation of the \textsc{coco} software package. Instead of building solutions, build a tool that others could use to build solutions. Define the platform, the language of discourse, and the paradigm of problem construction and analysis. Make it easy to pursue the particular and yet worthwhile to support the general. Reward a tight coupling between rigorous mathematics and computational encoding. Invite the community to innovate and substitute, benefiting from a prescribed interface.

An important step in this direction was the decoupling of problem construction from problem analysis. After all, the nonlinear equations analyzed using continuation methods could also be the target of iterative and even stochastic techniques that made no use of the manifold nature of their solutions~\cite{boender1982stochastic,byrd2000trust}. And while other packages emphasized analysis of one-dimensional solution manifolds (with some exceptions, e.g., \textsc{manpak}~\cite{rheinboldt1996manpak} and \textsc{multifario}~\cite{henderson2002multiple}), there was really no good reason to impose this restriction at the stage of problem construction. The real breakthrough, however, was in conceiving of a modular and staged construction paradigm, respecting an oft-occurring arrangement of the problem unknowns into densely coupled communities with sparse coupling to other communities (per the terminology of network science~\cite{porter2009communities}). A powerful application was found in multi-segment boundary-value problems (BVPs), e.g., periodic orbits in hybrid dynamical systems, in which boundary conditions between segments represent such sparse coupling between groups of unknowns separately parameterizing individual segments~\cite{dankowicz2011extended}.

As implemented in \textsc{coco}, this construction paradigm leveraged an object-oriented perspective, conceiving of a system of equations as decomposed into multiple object instances, describing subsets of equations and unknowns with inherent meaning, coupled together through appropriate gluing conditions (cf.\ the terminology used in multibody systems~\cite{otter1996modeling,schiehlen2013advanced}). With the recognition of common examples of mathematical objects (e.g., equilibria, trajectory segments, and periodic orbits) as constituting abstract classes of equations and unknowns, there resulted a hierarchy of problem construction whereby new abstract classes could be constructed from the composition of existing ones, and different versions of existing abstract classes could be substituted at will. As an example, problems involving the simultaneous analysis of an equilibrium (E), a periodic orbit (P), and an E-to-P connecting orbit were constructed with ease by leveraging existing abstract classes for each of these objects, glued together with a sparse set of boundary conditions~\cite{dankowicz2013recipes,dankowicz2011continuation,krauskopf2008lin}.

With the recognition of problem construction as distinct from problem analysis, more emphasis could also be placed on developing alternative approaches to continuation along solution manifolds, including the possibility of analysis along multi-dimensional manifolds. Where such algorithms in other packages were more tightly connected to the particular defining problems, the implementation in \textsc{coco} sought to remain agnostic as to the origin of the system of governing equations. This level of generality, of course, came at a cost as the particular solutions specific to a problem class could not be anticipated a priori. A satisfactory solution to this tension, also generalizable to the multi-dimensional case, was arrived at only in the past few years~\cite{dankowicz2020multidimensional,yuqing2018thesis}.

With these observations in mind, it is clear that software design has become a matter worthy of independent study, also to the community of nonlinear dynamics researchers. Moreover, with the appropriate attention to its theoretical underpinnings, such study also comes with scholarly reward. One example is the recent expansion to the original staged construction paradigm of \textsc{coco} in support of the parallel staged construction of (a critical subset of) the adjoint necessary conditions for extrema along constraint manifolds~\cite{li2017staged,li2020optimization}. This expansion reflects the decomposition of a problem Lagrangian into a sum of individual constraints linearly paired with corresponding adjoint variables (also called dual variables or Lagrange multipliers) that measure the sensitivities of an objective function to constraint violations at stationary points of the Lagrangian. Since the Lagrangian is linear in the adjoint variables, the contributions to the adjoint conditions from each term of the Lagrangian are also linear in the adjoint variables~\cite{ahsan2020optimization}. The complete set of adjoint conditions may therefore again be built in stages in a one-to-one mapping to the stages used to construct the full set of constraints.

It is one aim of this review to describe in detail this staged construction paradigm in a manner compatible with the implementation in \textsc{coco} but sufficiently abstract to allow for independent implementation. The effort involved in such independent development may be appreciated by reference to the history of \textsc{coco}.

\subsection{A brief history of \textsc{coco}} 
\label{sec: history of coco}
The software package \textsc{coco} is the result of joint development since 2007 by Harry Dankowicz and Frank Schilder, and, since 2016, Mingwu Li~\cite{mingwu2020thesis}, with additional contributions from Michael E.~Henderson, Erika Fotsch~\cite{fotsch2016thesis}, and Yuqing Wang~\cite{yuqing2018thesis}. Helpful feedback and contributions are also acknowledged from Jan Sieber and David Barton, and a growing user community~\cite{barton2017control,cao2019nonlinear,heitmann2021arrhythmogenic,liu2017controlling,ponsioen2018automated,zhong2020global}. Extensive discourse on the original design philosophy and mathematical underpinnings of the \textsc{coco} platform is available in the textbook~\cite{dankowicz2013recipes}, which includes a large collection of template toolboxes and example problems.

The first official release of \textsc{coco} coincided with the publication of~\cite{dankowicz2013recipes} in 2013. This included code documentation and detailed demos from~\cite{dankowicz2013recipes}. The November 2015 release introduced fully documented, production-ready toolboxes for common forms of bifurcation analysis of equilibria and periodic orbits in dynamical systems. These provided support for continuation of
\begin{compactitem}
    \item equilibria in smooth dynamical systems using the \mcode{ep} toolbox;
    \item constrained trajectory segments with independent and adaptive discretizations in autonomous and non-autonomous dynamical systems using the \mcode{coll} toolbox; and
    \item single-segment periodic orbits in smooth, autonomous or non-autonomous dynamical systems, and multi-segment periodic orbits in hybrid, autonomous dynamical systems using the \mcode{po} toolbox.
\end{compactitem} 
The November 2017 release made significant updates to the \textsc{coco} core and library of toolboxes and demos to provide support for constrained design optimization on integro-differential boundary-value problems~\cite{li2017staged}. These updates enabled the staged construction of the adjoint equations associated with equality-constrained optimization problems, and provided support for adaptive remeshing of these equations in parallel with updates to the problem discretization of the corresponding boundary-value problems. The March 2020 release of \textsc{coco} extended this functionality to also allow for finite-dimensional inequality constraints, bounding the locus of extrema to an implicitly-defined feasible region~\cite{li2020optimization}.

The original release of \textsc{coco} included the default atlas algorithm \mcode{atlas_1d} for one-dimensional solution manifolds. This was accompanied by a discussion in Parts III and IV of~\cite{dankowicz2013recipes} that described a general methodology for building atlas algorithms and also included an example of a two-dimensional atlas algorithm with fixed step size for non-adaptive continuation problems, inspired by Henderson's \textsc{multifario} package~\cite{henderson2002multiple}. A fully step-size-adaptive implementation of \textsc{multifario} as a \textsc{coco}-compatible atlas algorithm for multi-dimensional manifolds of solutions to non-adaptive continuation problems was included as an alpha version in the November 2017 release. The March 2020 release of \textsc{coco} included the updated atlas algorithm \mcode{atlas\_kd} for multi-dimensional solution manifolds for \emph{adaptive} continuation problems with varying embedding dimension and interpretation of problem unknowns~\cite{dankowicz2020multidimensional}. Usage of the \mcode{atlas_1d} and \mcode{atlas\_kd} atlas algorithms, as well as the basic \textsc{coco} constructors and utilities and those particular to the \mcode{ep}, \mcode{coll}, and \mcode{po} toolboxes is illustrated in numerous examples in tutorial documents included with the \textsc{coco} release. Each example corresponds to fully documented code in the release. 

One of the purported strengths of the \textsc{coco} package vis-\`{a}-vis its peers is its extensibility~\cite{blyth2020}. For example, it has not been the intent of the \textsc{coco} development to build graphical user interfaces to the methods and data invoked and processed during analysis of a continuation problem, although some low-level data processing and visualization routines are included with the \textsc{coco} core. Instead, support for run-time access to data is available in \textsc{coco}, for example, through a signal-and-slot mechanism as described in~\cite{dankowicz2013recipes}. Such a communication protocol allows independent development of user interfaces without modifications to the core. An example of such independent development is the analysis of hybrid dynamical systems described in~\cite{chong2016numerical} using a graphical user interface to the \textsc{coco} core and the \mcode{po} toolbox. New classes of problems may also be analyzed using \textsc{coco} without a preexisting toolbox for this purpose. Examples include the coupling in~\cite{formica2013coupling} of \textsc{coco} and the \mcode{po} toolbox with the \textsc{comsol multiphysics} finite element software, the analysis of quasiperiodic invariant tori in~\cite{li2020tor} using the \mcode{coll} toolbox, and the integration of the construction of spectral submanifolds in~\cite{ssmtool2} with frequency response analysis using the \mcode{ep} and \mcode{po} toolboxes. With the help of suitable wrappers, the data structures generated by the \textsc{coco} construction methodology may also be used by non-\textsc{coco} computational algorithms. An example is the application in~\cite{coco-fmincon} of the \textsc{matlab} function \texttt{foptim} to the continuation problem constructed using \textsc{coco}.

Sophisticated users may also wish to build new toolboxes for others to use. Advanced techniques for bifurcation detection, normal-form analysis, and so on, can be implemented using well-defined interface functions. An example is the implementation in~\cite{coco-shoot} of shooting techniques for continuation of periodic orbits using \textsc{coco}. Another example is the toolbox \mcode{ddebiftool\_coco} that provides \textsc{coco}-compatible access to the defining systems and monitoring functions created by \textsc{dde-biftool} for bifurcation analysis of DDEs \cite{ddebiftoolmanual}, thereby benefiting from the atlas algorithms and nonlinear solvers of \textsc{coco}. A further example is the work by Schilder \emph{et al.} \cite{schilder2015experimental} to develop a \textsc{coco}-compatible toolbox for noise-contaminated zero problems as occur when performing continuation in experiments~\cite{barton2012control,renson2016robust,renson2019application,renson2019numerical}. Shipped with the \textsc{coco} release, their \mcode{continex} toolbox includes an original atlas algorithm and nonlinear solver designed to track one-dimensional solution manifolds given low-precision numerics and high costs for evaluating residuals and their sensitivities.

%In this case, one observes high sensitivity of the frequency response curves with respect to nonlinearities or damping~\citet{barton2012control,detroux2018experimental,renson2019application}. In all of these cases, however, systematic multiple continuations with varying forcing amplitudes established that the underlying response surface in the forcing amplitude-forcing frequency plane is relatively regular (displaying no sharp ``corners'').

\subsection{Contributions of this paper}
\label{sec: contributions}
Rather than a mere review of the state-of-the-art in continuation methods and their applications, this paper makes several original contributions that are not covered elsewhere. These contributions are conceptual and structural and point to innovations in software design of the sort discussed above. They expand access to known solutions, rather than offer new solutions to known problems.

Foremost among these contributions is a detailed guide for the construction of a \textsc{coco}-compatible toolbox for analyzing families of solutions to multi-segment boundary-value problems with discrete delays, as well as for finding stationary points of scalar-value objective functions along such families (see~\cite{bartoszewski2011solving,calver2017numerical,chai2013unified,engelborghs2001collocation,gollmann2009optimal,khasawneh2011multi,shinohara2007boundary,TraversoMagri2019}, especially~\cite{gollmann2009optimal,TraversoMagri2019} for a similar usage of auxiliary variables to represent time-delayed terms). Not only does such a toolbox not exist previously for \textsc{coco}, but is also not available through other packages. The \mcode{ddebiftool\_coco} toolbox mentioned above, for example, is not designed to couple multiple trajectory segments and lacks tools for automatic construction of the adjoint contributions. By adhering to the object-oriented construction paradigm, the treatment in this paper demonstrates how very general classes of boundary-value problems may be addressed within a single framework, avoiding the need to develop individual solutions for periodic orbits, quasiperiodic orbits, connecting orbits, initial-value problems, and optimal control problems.

Prominent among the additional contributions of this paper is a detailed discussion of a data assimilation problem with delay inspired by the analysis in~\cite{TraversoMagri2019}. Here, a constraint Lagrangian is used to generate explicit adjoint conditions in a form amenable to a continuation-based analysis per the successive continuation framework in~\cite{ahsan2020optimization,li2017staged,li2020optimization} as an alternative to the gradient-based optimization approach of~\cite{TraversoMagri2019}. In contrast to~\cite{TraversoMagri2019}, the discussion highlights the natural decomposition of the governing constraints into a multi-segment boundary-value problem, the linear dependence on suitably defined adjoint variables, and the way in which time delay in the governing differential constraints translates into time-advanced terms in the adjoint differential equations.

A contribution of unexpected importance is the original derivation using a Lagrangian formalism of the governing equations for computing phase response curves associated with limit cycles in a general Banach-space setting (in contrast, e.g., to~\cite{izhikevich1997weakly} where derivation is based on the adjoint equation of the reduced phase models). This treatment demonstrates how a phase response functional may be constructed from the adjoint variables associated with the sensitivities of the orbital period to violations of the differential constraints and periodic boundary conditions. As the adjoint conditions may again be constructed automatically from variations of a constraint Lagrangian, the discussion points to the immediate use of \textsc{coco}-compatible toolboxes that provide such support without the need for further development. This is in contrast to support for phase response analysis in other software packages, for example \textsc{matcont}~\cite{dhooge2003matcont}, which implement a reduced set of adjoint differential equations, boundary conditions, and normalization conditions that must be derived separately for each class of problem. Importantly, the Lagrangian foundation developed in this paper also suggests that out-of-the-box use of existing \textsc{coco} toolboxes for limit cycles in hybrid dynamical systems (e.g., piecewise-smooth vector fields) would permit such phase response analysis.

Finally, of notable mention is an original discussion of the benefits of multi-dimensional continuation for managing uncertainty in singularly-perturbed or noise-contaminated problems, for example applications involving experimental data. In such cases, uncertain input data or randomly disturbed residuals (caused by measurement errors) may result in a dramatic degree of output uncertainty even if single-parameter continuation were computationally feasible. As shown in this paper, however, this singular behavior may be regulated or entirely eliminated using multi-dimensional continuation independently of the value of the damping.

\subsection{Organization of this paper}
\label{sec: organization}
The body of this paper is divided into four sections book-ended by the present introduction and a concluding discussion in Section~\ref{sec: Conclusions}. Section~\ref{sec: Problem formulation} develops the principles of staged problem construction for the so-called augmented continuation problem, amenable to analysis of constraint manifolds and optimization along such manifolds. Several examples are used to first motivate this framework and, subsequently, illustrate its application to advanced analysis of problems with delay. A pattern of universality uncovered by the treatment in Section~\ref{sec: Problem formulation} is converted into a rigorous mathematical formalism for a  \textsc{coco}-compatible toolbox for multi-segment boundary-value problems with discrete delays in Section~\ref{sec: Toolbox construction}. Straightforward generalizations of the toolbox and applications to the computation of phase response curves for limit cycles, homoclinic connections, quasiperiodic invariant tori, and optimal control inputs are considered in Section~\ref{sec: Numerical examples}. The text then turns briefly to opportunities for further development of the basic toolbox functionality in Section~\ref{sec: Further development}.

Several parts of the discussion in Section~\ref{sec: Problem formulation} may be read independently from the remainder of the text, although clearly at some loss to the continuity of the flow. This certainly applies to the description of the \textsc{coco} formalism in Section~\ref{sec: The coco formalism} and to the applications to data assimilation in Section~\ref{sec: data assimilation} (except for Section~\ref{sec: data assimilation problem construction} and parts of Section~\ref{sec: data assimilation problem analysis}) and phase response analysis in Section~\ref{sec: phase response curves} (except for Section~\ref{sec: phase response curves problem construction}). The discussion of problem discretization in Section~\ref{sec: discretization} may be skipped on a first reading. For the reader interested in the detailed implementation or considering an independent development, this section stresses the importance of systematic notation and rigor also in the encoding of a problem in order to ensure code verifiability. Finally, while the ordering of the text implies a natural flow, at times a nonlinear approach to reading this review may be appropriate. The reader may wish to skip ahead to anticipate the implications of the design decisions or return to an earlier section to better appreciate its purpose. That is encouraged.

\section{Problem formulation}
\label{sec: Problem formulation}

% \textcolor{red}{It is customary in treatments of continuation methods to begin with a discussion of the implicit function theorem~\cite{doedel2007lecture} (or, more generally, the implicit manifold theorem), as the theoretical foundation for analyzing solutions of abstract nonlinear problems. Such a discussion naturally concerns itself with a decomposition of the unknowns into independent and dependent variables, and establishes conditions under which such a decomposition makes (local) sense. These conditions are then leveraged to give meaning to the notion of continuation: the local and continuous expansion of the known universe of solutions.}~\textcolor{blue}{The lecture notes by Doedel~\cite{doedel2007lecture} follow this custom and give comprehensive discussions on the implicit-function theorem and then the numerical continuation of solutions to nonlinear problems}. \js{Dankowicz \emph{et al}.\ \cite{dankowicz2020multidimensional} describe how the \textsc{coco} multi-dimensional atlas algorithm \mcode{atlas\_kd} permits the user to choose the embedding space for the implicitly defined manifold to ensure applicability of the implicit manifold theorem for adaptive nonlinear problems (for example, problems with adaptive discretization or phase condition.)}

It is customary in treatments of continuation methods (e.g., \cite{doedel2007lecture}) to begin with a discussion of the implicit function theorem, as the theoretical foundation for analyzing solutions of abstract nonlinear problems. Such a discussion naturally concerns itself with a decomposition of the unknowns into independent and dependent variables, and establishes conditions under which such a decomposition makes (local) sense. These conditions are then leveraged to give meaning to the notion of \emph{continuation}: the local and continuous expansion of the known universe of solutions along implicitly defined manifolds.

Here, we largely depart from such a focus on solutions and their geometry by instead emphasizing the process of problem construction. Our concern is not principally with the techniques used to \emph{perform} continuation, but with a systematic approach to formulating problems \emph{amenable} to continuation, without imposing any preferred decompositions among the problem unknowns. As we show in this section, such a problem-oriented focus may yield benefits also to the process of continuation, as different formulations are more or less well-conditioned. Nevertheless, our primary aim is to identify patterns and structure in the way common problems arise in the study of dynamical systems, and to build useful abstractions around such patterns.

It is instructive to begin this journey into methods of continuation and their implementations in software within the realm of problems amenable to closed-form analysis. Such analysis removes consideration of various numerical approximations, inevitable in a computational implementation, and offers an opportunity for code verification. For the particular examples considered in this section, it points to generalizations to nonlinear problems without closed-form solutions. More importantly, it illustrates principles of intuitive and flexible problem construction, partially agnostic to the final objectives of the analysis. We argue that such flexibility should take precedence in the engineering of general-purpose software for continuation problems.

%\textcolor{red}{We would like to emphasize that the focus of this section and the overall study in general is on different problem formulations. The systematic exploration of solutions to these problems generally involves the application of implicit function theorem (IFT). This theorem provides suitable conditions on system Jacobians for existence of local solution manifolds and also provides information to overcome singularities near singular solution points. While the IFT provides the mathematics for finding solutions of nonlinear problems~\cite{doedel2007lecture}, this review explores systematic ways to build up these problems.}
%%%%%%%%%%%%%%%%%%%
% Inflection points
%%%%%%%%%%%%%%%%%%%
\subsection{Looking for inflection points}
\label{sec: inflection points}
Many problems of interest in the analysis and control of nonlinear dynamical systems may be formulated as problems of constrained design optimization (see, e.g., the study of periodically forced bioreactors in \cite{d2010choice} or bubble motion driven by acoustic forcing in \cite{toilliez2008optimized,wyczalkowski2003optimization}). In this section, we consider the search for optimal points along manifolds of solutions to algebraic and/or differential constraints in terms of objective functions characterizing the local manifold geometry (for an applied context, see~\cite{acharya2020non} for a recent study of non-monotonic dependence of the response dynamics of premixed flames on forcing amplitude).

\begin{figure}[ht]
\centering
\includegraphics[width=0.45\textwidth]{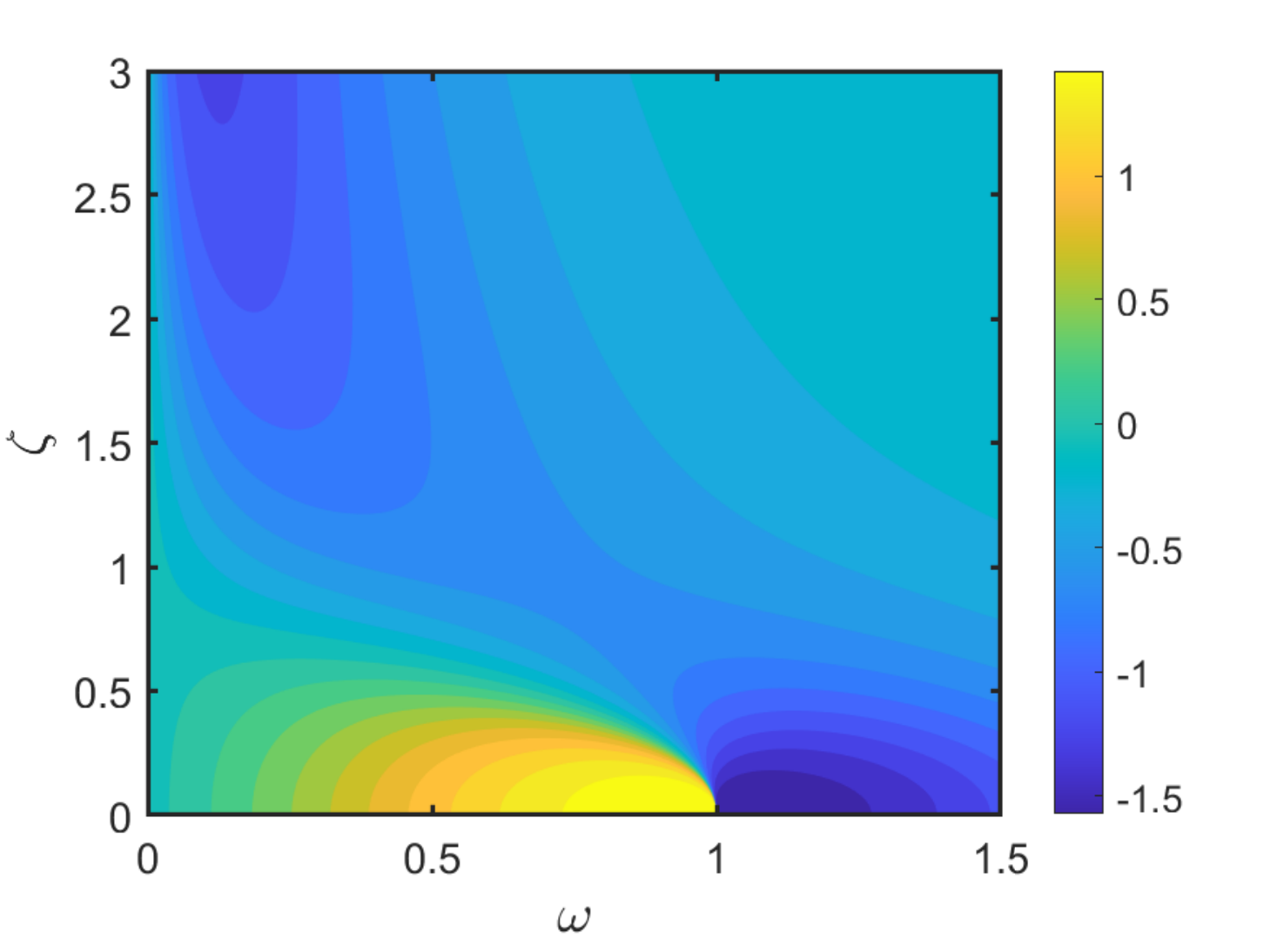}
\includegraphics[width=0.45\textwidth]{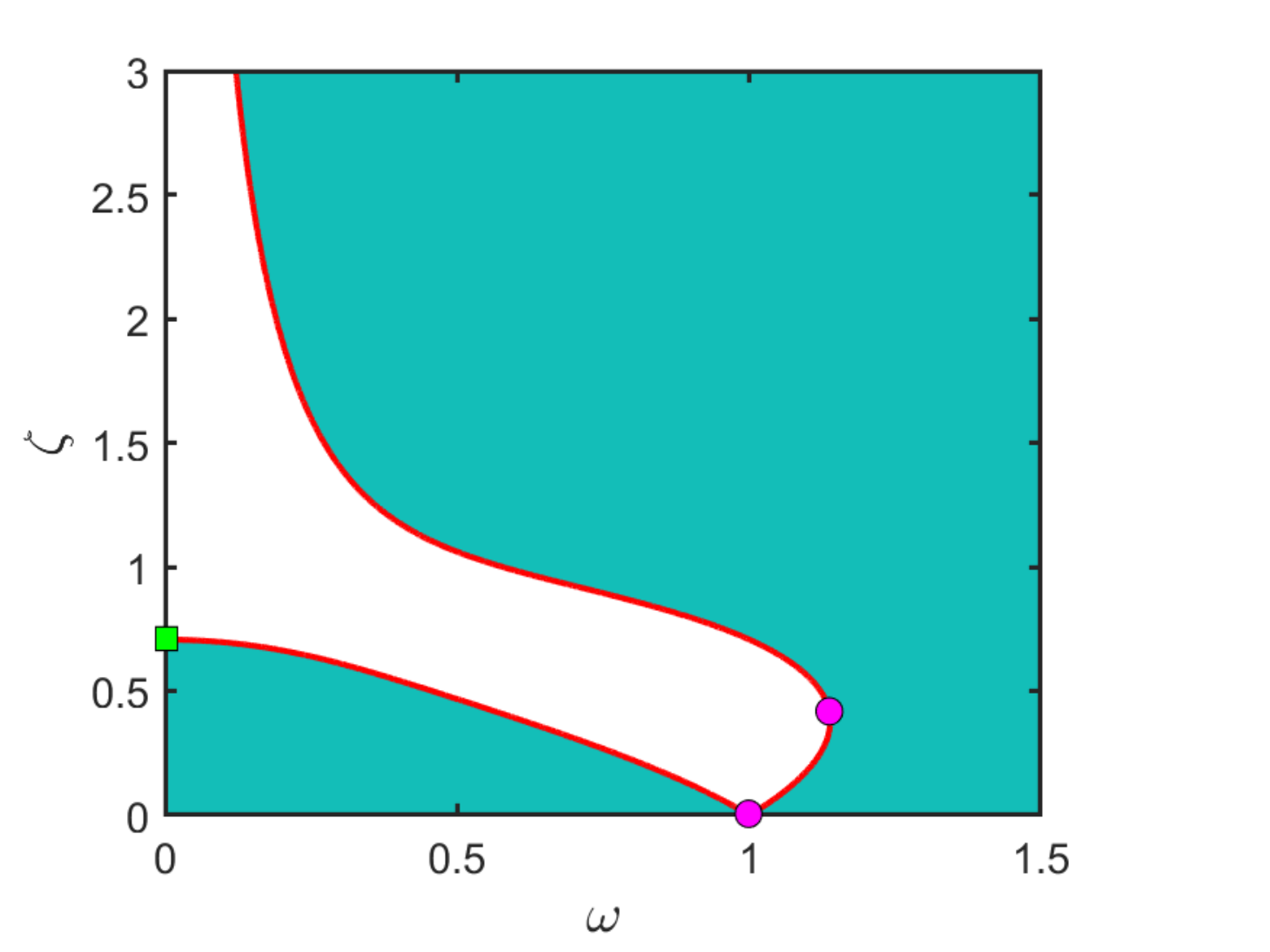}
\caption{(left) Contour plot of the inverse tangent ($\arctan$) of the partial derivative of the response amplitude $C$ in \eqref{eq:explamp} with respect to $\omega$. (The inverse tangent operator is used to handle the singularity of the partial derivative when $(\zeta,\omega)\to(0,1).$) (right) The zero level sets of the polynomial in \eqref{eq:zetaomega} (red lines) coincide with the zero contour of the second order partial derivative of the response amplitude $C$ with respect to $\omega$ (bounding the dark green region). The filled circles (magenta) are singular points when the curve is parameterized by $\zeta(\omega)$, and the filled box (light green) is a singular point if the curve is parameterized by $\omega(\zeta)$.}
\label{fig:amp_om}
\end{figure}

\begin{figure}[ht]
\centering
\includegraphics[width=0.45\textwidth]{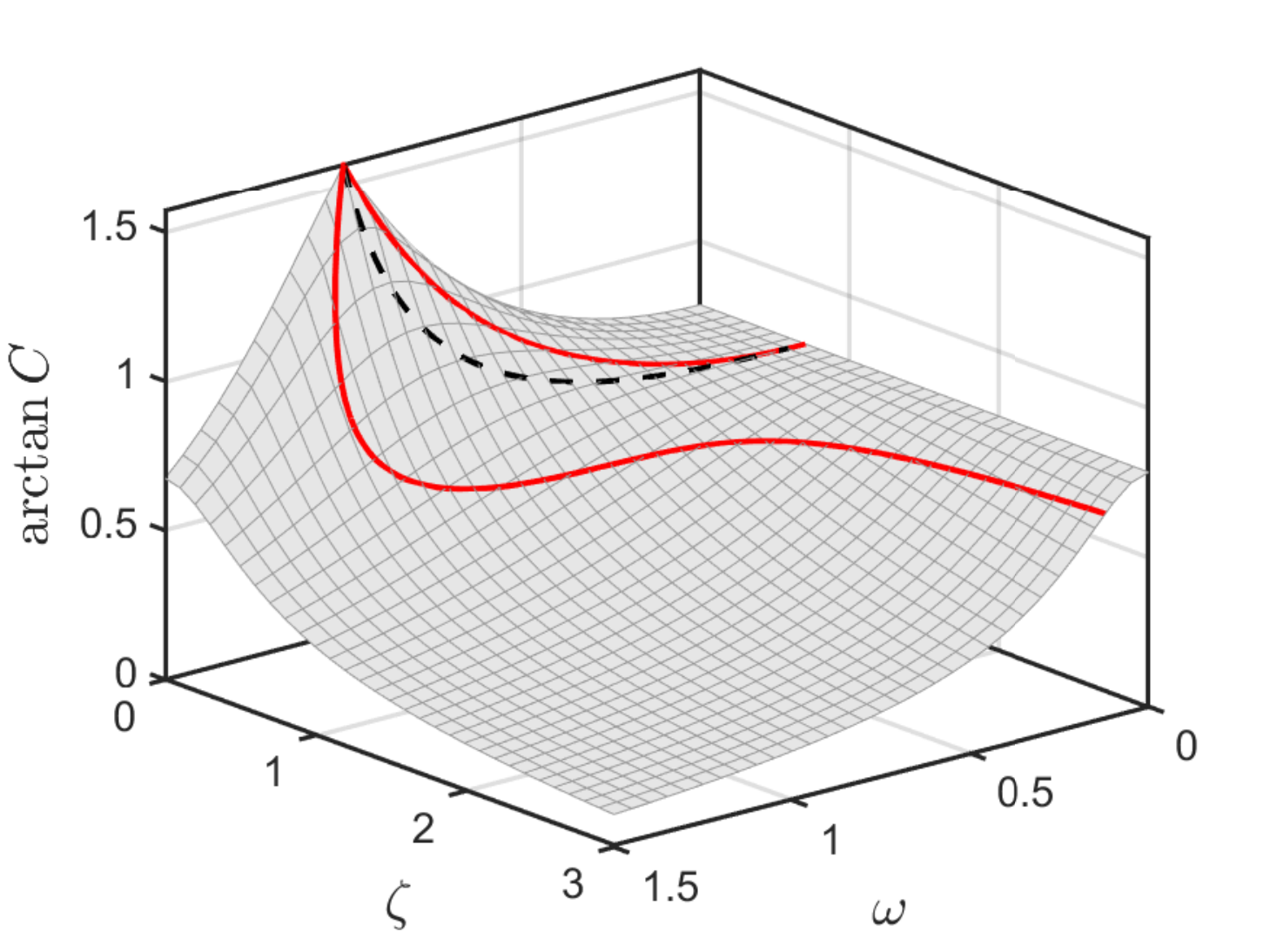}
\caption{Frequency response surface of the harmonically forced linear oscillator. Here, the surface plot is based on the explicit expression for $C$ in \eqref{eq:explamp}, the solid (red) lines are the sought extrema in the rate of change of the response amplitude $C$ with respect to $\omega$ based on~\eqref{eq:zetaomega} and~\eqref{eq:extrema-amp-omega}, and the dashed line (black) locates the global maximum of the response amplitude for $\zeta\leq1/\sqrt{2}$.}
\label{fig:amp}
\end{figure}

Specifically, along the family of steady-state periodic responses of a harmonically-excited, linear oscillator parameterized by the excitation frequency $\omega$, at most two values of $\omega$ correspond to local extrema in the rate of change of the response amplitude with respect to $\omega$, as shown in the left panel of Fig.~\ref{fig:amp_om}. To locate these values, we write the governing equation in the normalized form
\begin{equation}
\label{eq:linosc}
    \ddot{x}+2\zeta\dot{x}+x=\cos\omega t,\,\zeta,\omega>0,
\end{equation}
make the ansatz $x(t)=C\cos(\omega t-\theta)$ for $C>0$, and obtain
\begin{equation}
\label{eq:explamp}
    C=\frac{1}{\sqrt{(1-\omega^2)^2+4\zeta^2\omega^2}}.
\end{equation} 
Differentiation twice with respect to $\omega$ then yields inflection points at the roots of the polynomial
\begin{equation}
\label{eq:zetaomega}
    3\omega^6+5(2\zeta^2-1)\omega^4+(16\zeta^4-16\zeta^2+1)\omega^2+1-2\zeta^2
\end{equation}
or, equivalently, at points $(\zeta,\omega)$ with
\begin{equation}
\label{eq:extrema-amp-omega}
    \zeta=\frac{1}{4\omega}\sqrt{1+8\omega^2-5\omega^4\pm\sqrt{1+38\omega^4-23\omega^8}}
\end{equation}
as illustrated in the right panel of Fig.~\ref{fig:amp_om}. It follows that only one such root exists for $\zeta>1/\sqrt{2}$, whereas two roots bracket the global maximum of the response amplitude at $\omega=\sqrt{1-2\zeta^2}$ for $\zeta\le 1/\sqrt{2}$. At these points, the response amplitude is given by
\begin{equation}
    \frac{2}{\sqrt{5-\omega^4\pm\sqrt{1+38\omega^4-23\omega^8}}}
\end{equation}
and $1/\sqrt{1-\omega^4}$, respectively, as shown in Fig.~\ref{fig:amp}.

In lieu of the analysis afforded by the explicit expression \eqref{eq:explamp} for the response amplitude, consider the equations
\begin{equation}
    (1-\omega^2)A+2\zeta\omega B-1=0,\,(1-\omega^2)B-2\zeta\omega A=0
\end{equation}
obtained by substitution of the ansatz $A\cos\omega t+B\sin\omega t$ in \eqref{eq:linosc}. To locate inflection points in the response amplitude $C=\sqrt{A^2+B^2}$, we may directly constrain a finite-difference approximation of its second derivative per the collection of polynomial constraints
\begin{gather}
     C_1^2-A_1^2-B_1^2=0,\,C_2^2-A_2^2-B_2^2=0,\,C_3^2-A_3^2-B_3^2=0,\\
   (1-\omega_1^2)A_1+2\zeta\omega_1B_1-1=0,\,(1-\omega_1^2)B_1-2\zeta\omega_1A_1=0,\\
    (1-\omega_2^2)A_2+2\zeta\omega_2B_2-1=0,\,(1-\omega_2^2)B_2-2\zeta\omega_2A_2=0,\\
    (1-\omega_3^2)A_3+2\zeta\omega_3B_3-1=0,\,(1-\omega_3^2)B_3-2\zeta\omega_3A_3=0,\\
    \omega_1-\omega_2-\epsilon=0,\,\omega_2-\omega_3-\epsilon=0,\,C_1-2C_2+C_3=0
\end{gather}
in the limit as $\epsilon\rightarrow 0$. As an alternative, consider instead \emph{the constrained optimization} of the \emph{objective function} $C_1-C_2$ with respect to $\omega_1$ in the limit as $\epsilon\rightarrow 0$, given the polynomial \emph{constraints}
\begin{gather}
    C_1^2-A_1^2-B_1^2=0,\,C_2^2-A_2^2-B_2^2=0,\,\omega_1-\omega_2-\epsilon=0,\label{eq:const1}\\
    (1-\omega_1^2)A_1+2\zeta\omega_1B_1-1=0,\,(1-\omega_1^2)B_1-2\zeta\omega_1A_1=0,\\
    (1-\omega_2^2)A_2+2\zeta\omega_2B_2-1=0,\,(1-\omega_2^2)B_2-2\zeta\omega_2A_2=0.\label{eq:const3}
\end{gather}
By the calculus of variations~\cite{gelfand2000calculus,liberzon2011calculus}, we obtain necessary conditions for such loci of optimality by considering vanishing variations of a suitably constructed constraint Lagrangian. Here, such an analysis results in the constraints \eqref{eq:const1}-\eqref{eq:const3} coupled with the adjoint conditions
\begin{gather}
    1+2C_1\lambda_1=-1+2C_2\lambda_2=0,\label{eq:adjobj}\\
    -2A_1\lambda_1+(1-\omega_1^2)\lambda_4-2\zeta\omega_1\lambda_5=0,\,
    -2A_2\lambda_2+(1-\omega_2^2)\lambda_6-2\zeta\omega_2\lambda_7=0,\\
    -2B_1\lambda_1+2\zeta\omega_1\lambda_4+(1-\omega_1^2)\lambda_5=0,\,
    -2B_2\lambda_2+2\zeta\omega_2\lambda_6+(1-\omega_2^2)\lambda_7=0,\\
    \lambda_3+2(\zeta B_1-\omega_1A_1)\lambda_4-2(\zeta A_1+\omega_1B_1)\lambda_5=0,\\
    -\lambda_3+2(\zeta B_2-\omega_2A_2)\lambda_6-2(\zeta A_2+\omega_2B_2)\lambda_7=0\label{eq:adjobjlast}
\end{gather}
in terms of the Lagrange multipliers $\lambda_1$ through $\lambda_7$ that describe the sensitivity of the objective function at stationary points to violations of each of the constraints~\eqref{eq:const1}-\eqref{eq:const3}. Solutions are obtained only for $\omega_1$ and $\omega_2=\omega_1-\epsilon$ that satisfy the equation
\begin{equation}
\label{eq:epsdiffcond}
    \frac{\omega_1(1-2\zeta^2-\omega_1^2)}{((1-\omega_1^2)^2+4\zeta^2\omega_1^2)^{3/2}}-\frac{\omega_2(1-2\zeta^2-\omega_2^2)}{((1-\omega_2^2)^2+4\zeta^2\omega_2^2)^{3/2}}=0,
\end{equation}
or, equivalently,
\begin{equation}
    0=\frac{3\omega_1^6+5(2\zeta^2-1)\omega_1^4+(16\zeta^4-16\zeta^2+1)\omega_1^2+1-2\zeta^2}{((1-\omega_1^2)^2+4\zeta^2\omega_1^2)^{5/2}}\epsilon+\mathcal{O}(\epsilon^2),
\end{equation}
again yielding the condition \eqref{eq:zetaomega} on $\omega_1$ from the previous paragraph. In this case, $C_1,C_2>0$ imply that
\begin{gather}
    \lambda_1=-\frac{\sqrt{(1-\omega_1^2)^2+4\zeta^2\omega_1^2}}{2},\,\lambda_2=\frac{\sqrt{(1-\omega_2^2)^2+4\zeta^2\omega_2^2}}{2},\label{eq:lagrange1}\\
    \lambda_3=-\frac{2\omega_1(1-2\zeta^2-\omega_1^2)}{((1-\omega_1^2)^2+4\zeta^2\omega_1^2)^{3/2}},\,\lambda_4=-\frac{1}{\sqrt{(1-\omega_1^2)^2+4\zeta^2\omega_1^2}},\label{eq:lagrange2}\\
    \lambda_5=0,\,\lambda_6=\frac{1}{\sqrt{(1-\omega_2^2)^2+4\zeta^2\omega_2^2}},\,\lambda_7=0.\label{eq:lagrange3}
\end{gather}
The reader is encouraged to verify the claim regarding the relationship between these values of the Lagrange multipliers and the sensitivity of the objective function at stationary points to constraint violations. In contrast to the discussion that led directly to \eqref{eq:zetaomega} in the first part of this section, we do not presuppose an explicit expression for the response amplitude, one that can be differentiated arbitrarily with respect to $\omega$. Instead, we use a finite-difference approximation in terms of a fixed change $\epsilon$ in the independent variable and show that the predicted extremum converges to the desired solution when $\epsilon\rightarrow 0$.

We may take a further step back from an explicit analysis by considering the constrained optimization with respect to $\omega_1$ of the objective function $x_{1}(0)-x_{2}(0)$ for $x_1(0),x_2(0)>0$, given the differential constraints
\begin{equation}
\label{eq:diffconst}
    \ddot{x}_1+2\zeta\dot{x}_1+x_1=\cos(\omega_1 t-\theta_1),\,
    \ddot{x}_2+2\zeta\dot{x}_2+x_2=\cos(\omega_2 t-\theta_2),
\end{equation}
the boundary conditions
\begin{gather}
\label{eq:boundconst}
    x_1(0)=x_1(T_1),\,\dot{x}_1(0)=\dot{x}_1(T_1)=0,\,x_2(0)=x_2(T_2),\,\dot{x}_2(0)=\dot{x}_2(T_2)=0
\end{gather}
with $T_1=2\pi/\omega_1$ and $T_2=2\pi/\omega_2$, and the algebraic constraint $\omega_1-\omega_2=\epsilon$ in the limit as $\epsilon\rightarrow 0$. The boundary conditions ensure that solutions are periodic with local extrema at $t=0$. In this case, the necessary conditions for optimality append to these constraints the adjoint conditions
\begin{gather}
\label{eq:diffadjt}
    \ddot{\lambda}_1-2\zeta\dot{\lambda}_1+\lambda_1=0,\,\ddot{\lambda}_2-2\zeta\dot{\lambda}_2+\lambda_2=0,\\
   1-2\zeta\lambda_1(0)+\dot{\lambda}_1(0)+\lambda_3=0,\,-1-2\zeta\lambda_2(0)+\dot{\lambda}_2(0)+\lambda_6=0,\label{eq:diffeqsub}\\
   2\zeta\lambda_1(T_1)-\dot{\lambda}_1(T_1)-\lambda_3=0,\,2\zeta\lambda_2(T_2)-\dot{\lambda}_2(T_2)-\lambda_6=0,\\
   -\lambda_1(0)+\lambda_4=0,\,\lambda_1(T_1)+\lambda_5=0,\,-\lambda_2(0)+\lambda_7=0,\,\lambda_2(T_1)+\lambda_8=0,\\
   -\int_0^{T_1}\lambda_1\sin(\omega_1t-\theta_1)\,\d t=0,\,-\int_0^{T_2}\lambda_2\sin(\omega_2t-\theta_2)\,\d t=0,\\
   \int_0^{T_1}\lambda_1 t\sin(\omega_1 t-\theta_1)\,\d t+2\pi\lambda_3\dot{x}_1(T_1)/\omega_1^2-2\pi\lambda_5\ddot{x}_1(T_1)/\omega_1^2+\lambda_9=0,\\
   \int_0^{T_2}\lambda_2 t\sin(\omega_2 t-\theta_2)\,\d t+2\pi\lambda_6\dot{x}_2(T_2)/\omega_2^2-2\pi\lambda_8\ddot{x}_2(T_2)/\omega_2^2-\lambda_9=0\label{eq:diffadjtlast}
\end{gather}
in terms of the Lagrange multipliers $\lambda_1$ through $\lambda_9$ that describe the sensitivity of the objective function at stationary points to violations of the differential constraints~\eqref{eq:diffconst}, boundary conditions~\eqref{eq:boundconst}, or algebraic constraint $\omega_1-\omega_2-\epsilon=0$, respectively. We again find that solutions exist only for $\omega_1$ and $\omega_2=\omega_1-\epsilon$ that satisfy \eqref{eq:epsdiffcond}, in which case, for example,
\begin{gather}
    \lambda_1(t)=-\frac{1}{2\sqrt{\zeta^2-1}}\left(\frac{e^{t(\zeta-\sqrt{\zeta^2-1})}}{e^{2\pi(\zeta-\sqrt{\zeta^2-1})/\omega_1}-1}-\frac{e^{t(\zeta+\sqrt{\zeta^2-1})}}{e^{2\pi(\zeta+\sqrt{\zeta^2-1})/\omega_1}-1}\right),\\
    \lambda_2(t)=-\frac{1}{2\sqrt{\zeta^2-1}}\left(\frac{e^{t(\zeta+\sqrt{\zeta^2-1})}}{e^{2\pi(\zeta+\sqrt{\zeta^2-1})/\omega_2}-1}-\frac{e^{t(\zeta-\sqrt{\zeta^2-1})}}{e^{2\pi(\zeta-\sqrt{\zeta^2-1})/\omega_2}-1}\right),
\end{gather}
and
\begin{equation}
    \lambda_9=-\frac{2\omega_1(1-2\zeta^2-\omega_1^2)}{((1-\omega_1^2)^2+4\zeta^2\omega_1^2)^{3/2}}.
\end{equation}
The reader is again encouraged to verify the claim regarding the relationship between these values of the Lagrange multipliers and the sensitivity of the objective function at stationary points to constraint violations. In contrast to the previous two approaches, we neither presuppose an explicit expression for the response amplitude nor for the form of the periodic response. Instead, the corresponding adjoint conditions \eqref{eq:diffadjt}--\eqref{eq:diffadjtlast} are here derived directly from the governing differential constraints and boundary conditions in a step that immediately generalizes to nonlinear problems for which closed-form solutions would not be available. As before, the finite-difference approximation in terms of a fixed change $\epsilon$ in the independent variable again approximates the loci of the inflection points to lowest order in $\epsilon$.

For practice, it may be worthwhile to repeat this discussion in the simpler search for a local extremum in the response amplitude under variations in $\omega$, known to exist at $\omega=\sqrt{1-2\zeta^2}$ for $\zeta<1/\sqrt{2}$. In this case, we might consider optimization of $C$ with respect to $\omega$ given the polynomial constraints
\begin{equation}
    C^2-A^2-B^2=0,\,1-A-2B\zeta\omega+A\omega^2=0,\,B-2A\zeta\omega-B\omega^2=0,
\end{equation}
or optimization of $x(0)>0$ with respect to $\omega$ given the boundary-value problem
\begin{equation}
    \ddot{x}+2\zeta\dot{x}+x=\cos(\omega t-\theta),\,
    x(0)=x(2\pi/\omega),\,\dot{x}(0)=\dot{x}(2\pi/\omega)=0.
\end{equation}
Alternatively, we could consider imposition of the additional constraint $C_1=C_2$ to the polynomial constraints \eqref{eq:const1}-\eqref{eq:const3} in the limit as $\epsilon\rightarrow 0$, or imposition of the additional constraint $x_1(0)=x_2(0)$ to the differential constraints \eqref{eq:diffconst}, boundary conditions \eqref{eq:boundconst}, and algebraic constraint $\omega_1-\omega_2=\epsilon$ in the limit as $\epsilon\rightarrow 0$. In doing so, one should reasonably ask which of these approaches generalize to nonlinear boundary-value problems and to other objective functions.

%%%%%%%%%%%%%%%%%%%%%%%%%%
% Lessons and inspirations
%%%%%%%%%%%%%%%%%%%%%%%%%%
\subsection{Lessons and inspirations}
\label{sec: lessons and inspirations}

The examples in the previous section are notably concerned not with a singular excitation response in isolation, but with a property of such a response in relation to nearby responses along a continuous (and locally differentiable) family of responses. Although we held $\zeta$ fixed in our analysis, the implicit relationship in \eqref{eq:zetaomega} further defines continuous families of inflection points and corresponding values of $\zeta$. We are inevitably drawn to a methodology for charting such continuous families and for monitoring the values of one or several objective functions along such families. 

As we approach this task, a count of degrees of freedom proves useful. We generically reduce the number of degrees of freedom by one for every algebraic constraint imposed on an \emph{a priori} unknown algebraic variable. Similarly, for every \emph{a priori} unknown solution to a differential constraint, we generically append as many degrees of freedom as the number of required initial conditions. As an example, Eq.~\eqref{eq:zetaomega} imposes one algebraic constraint on two \emph{a priori} unknown algebraic variables, yielding a problem with (generically) a single degree of freedom. Similarly, the seven constraints \eqref{eq:const1}-\eqref{eq:const3} constrain the ten \emph{a priori} unknown algebraic variables $A_1$, $B_1$, $C_1$, $\omega_1$, $A_2$, $B_2$, $C_2$, $\omega_2$, $\zeta$, and $\epsilon$ to yield a problem with (generically) three degrees of freedom. The eight adjoint conditions \eqref{eq:adjobj}-\eqref{eq:adjobjlast} add seven more \emph{a priori} unknown algebraic variables for a net of (generically) two degrees of freedom. Generically, the differential constraints \eqref{eq:diffconst}, boundary conditions \eqref{eq:boundconst}, and algebraic constraint $\omega_1-\omega_2=\epsilon$ on the \emph{a priori} unknown variables $x_1(\cdot)$, $\omega_1$, $\theta_1$, $x_2(\cdot)$, $\omega_2$, $\theta_2$, $\zeta$, and $\epsilon$ result in a problem with three degrees of freedom. The adjoint conditions \eqref{eq:diffadjt}-\eqref{eq:diffadjtlast} add nine more \emph{a priori} unknown variables for a net of (generically) two degrees of freedom.

The number of degrees of freedom of a differentiable continuation problem characterizes the dimension of a local manifold of solutions through any regular (in the sense of the implicit-function theorem~\cite{krantz2012implicit}) solution point. This dimension represents a deficit of constraints relative to the number of \emph{a priori} unknown variables, and so we often speak of the dimensional deficit of a continuation problem. For all the continuation problems of interest here, the dimensional deficit is a finite number, even as the problem domain may be infinite dimensional.

Problems with zero dimensional deficit generically have at most isolated solutions~\cite{doedel2007lecture}. For example, by inspection of the partial derivative with respect to $\zeta$ and $\omega$, respectively, the roots of the multivariable polynomial in \eqref{eq:zetaomega} are found to be locally unique with respect to $\zeta$ for all positive $\zeta\ne 1/\sqrt{2}$ (cf.~the green square at the right panel of Fig.~\ref{fig:amp_om}) and locally unique with respect to $\omega$ for all positive $\omega\ne 1$ or $((19+8\sqrt{6})/23)^{1/4}$ (cf.~the two magenta circles at the right panel of Fig.~\ref{fig:amp_om}). By inspection of the Jacobian with respect to $A_1$, $B_1$, $C_1$, $\omega_1$, $A_2$, $B_2$, $C_2$, $\omega_2$, and $\lambda_1$ through $\lambda_7$, solutions of the polynomial constraints \eqref{eq:const1}-\eqref{eq:const3} and the corresponding adjoint conditions \eqref{eq:adjobj}-\eqref{eq:adjobjlast} are locally unique with respect to $\zeta$ and $\epsilon$ for all positive $\zeta\ne 1/\sqrt{2}$ and sufficiently small $\epsilon$. Similarly, by inspection of the Jacobian with respect to $A_1$, $B_1$, $C_1$, $A_2$, $B_2$, $C_2$, $\omega_2$, $\zeta$, and $\lambda_1$ through $\lambda_7$, solutions are found to be locally unique with respect to $\omega_1$ and $\epsilon$ for all positive $\omega_1\ne 1$ or $((19+8\sqrt{6})/23)^{1/4}$ and sufficiently small $\epsilon$. For solutions to the differential constraints \eqref{eq:diffconst}, boundary conditions \eqref{eq:boundconst}, algebraic constraint $\omega_1-\omega_2=\epsilon$ and the corresponding adjoint conditions \eqref{eq:diffadjt}-\eqref{eq:diffadjtlast}, the same conclusions would be theoretically available by showing the invertibility of the linearization with respect to $x_1(\cdot)$, $\omega_1$, $\theta_1$, $x_2(\cdot)$, $\omega_2$, $\theta_2$, $\lambda_1(\cdot)$, $\lambda_2(\cdot)$, and $\lambda_3$ through $\lambda_9$ or $x_1(\cdot)$, $\theta_1$, $x_2(\cdot)$, $\omega_2$, $\theta_2$, $\zeta$, $\lambda_1(\cdot)$, $\lambda_2(\cdot)$, and $\lambda_3$ through $\lambda_9$, respectively. This undertaking is left to the reader.

Local uniqueness affords us confidence that an approximate algorithm to locate a solution to a problem with zero dimensional deficit will not be distracted by other nearby solutions. Provided that we initialize a search with an initial solution guess in the vicinity of the sought solution, we trust that a well-designed solver, e.g., based on Newton's or Broyden's methods~\cite{kelley1995iterative}, will rapidly converge to this solution. For the first two formulations of the inflection point problem in Section~\ref{sec: inflection points}, we apply such a solver directly to the system of nonlinear equations. For the formulation in terms of differential boundary-value problems, some form of discretization must first be employed.

Inspired by these observations, a general continuation methodology for a problem $\mathbf{P}$ with nonzero dimensional deficit may be obtained by iteratively
\begin{itemize}
    \item constructing auxiliary constraints~\cite{crisfield1983arc,dankowicz2013recipes,henderson2002multiple} that when appended to $\mathbf{P}$ result in a problem $\mathbf{P}_0$ with zero dimensional deficit;
    \item constructing an initial solution guess for $\mathbf{P}_0$ using a previously found solution to $\mathbf{P}$~\cite{dankowicz2013recipes,govaerts2000numerical,seydel2009practical}; and
    \item solving $\mathbf{P}_0$ using an iterative algorithm based at the initial solution guess.
\end{itemize}
By definition, a solution to $\mathbf{P}_0$ also solves $\mathbf{P}$. The success of such a methodology thus depends on its ability to ensure that solutions to $\mathbf{P}_0$ are locally unique; that the iterative solver is able to converge to such a solution; and that the succession of such solutions suitably captures the geometry of the manifold of solutions to $\mathbf{P}$~\cite{dankowicz2013recipes,guddat1990parametric}. 

\begin{figure}[ht]
\centering
\includegraphics[width=0.45\textwidth]{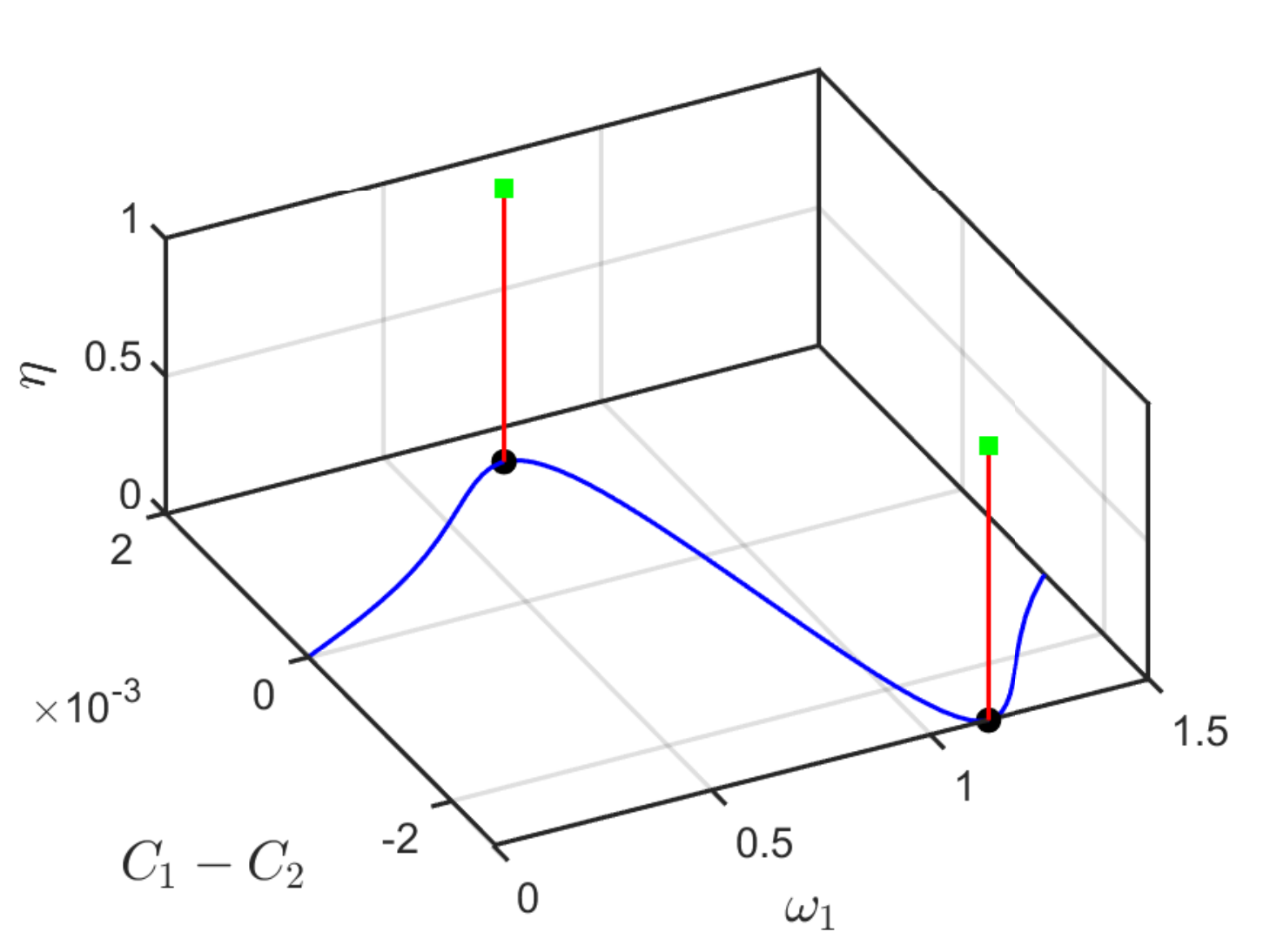}
\includegraphics[width=0.45\textwidth]{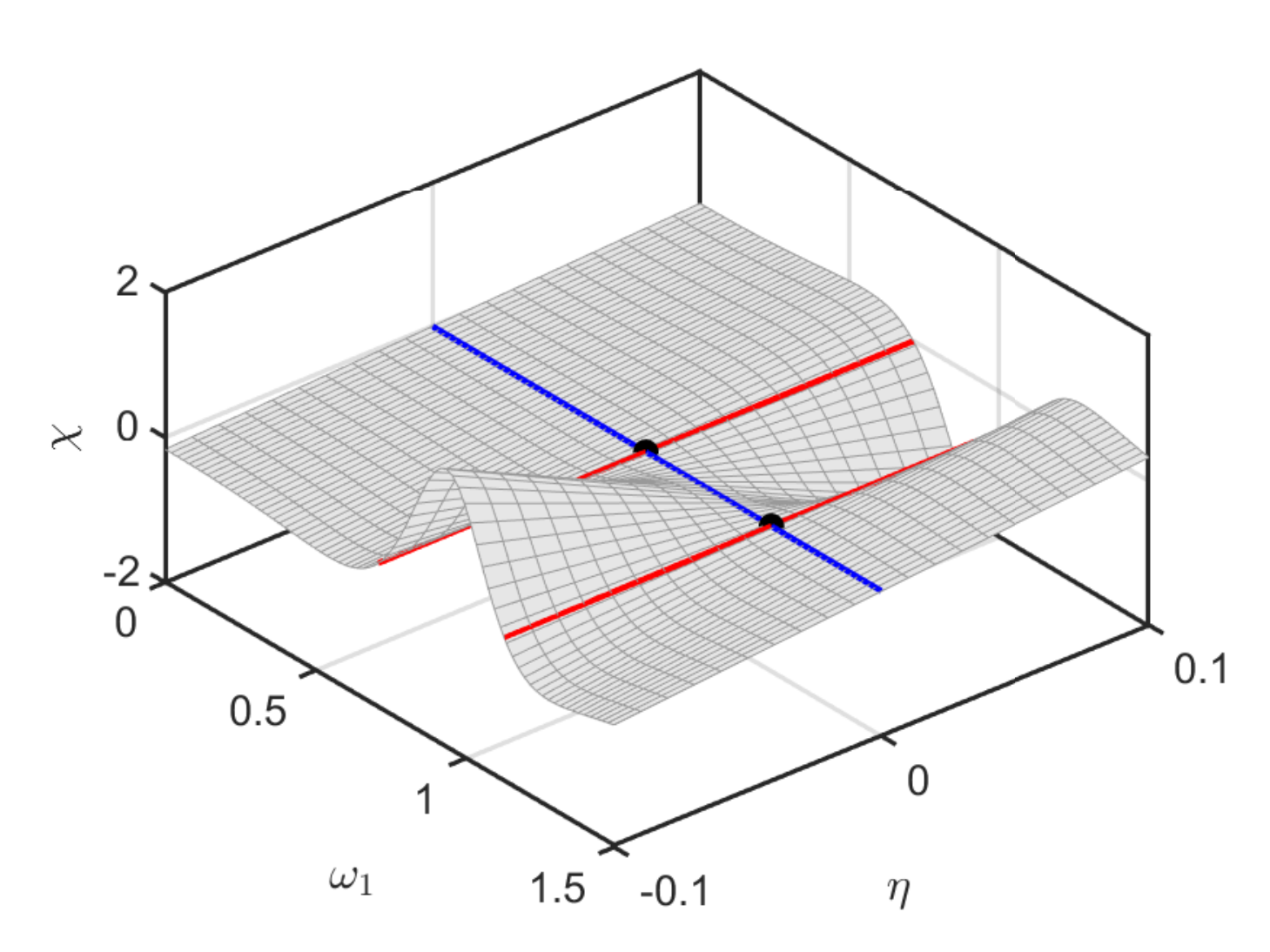}
\caption{(left) Projections of solution branches with vanishing (solid blue) and non-vanishing (solid red) Lagrange multipliers for the problem $\mathbf{P}$ for locating extrema of $C_1-C_2$ with fixed $(\epsilon,\zeta)=(0.001,0.3)$. The intersections (filled circles) correspond to local extrema of $C_1-C_2$ and are singular points of $\mathbf{P}$. The green squares represent solutions with $\eta=1$. (right) Projection of the two-dimensional solution manifold of the problem $\mathbf{P}^\ast$ with fixed $(\epsilon,\zeta)=(0.001,0.3)$. The blue and red straight lines correspond to the identically-colored curves in the left panel and lie in the zero level set of $\chi$ on the solution manifold. The intersections (black filled circles) are regular points of $\mathbf{P}^\ast$.}
\label{fig:inflection_opt_finite}
\end{figure}

Consider, for example, the problem obtained by replacing \eqref{eq:adjobj} in the necessary conditions for an extremum of $C_1-C_2$ under the polynomial constraints \eqref{eq:const1}-\eqref{eq:const3} with
\begin{equation}
    \eta-2C_1\lambda_1=-\eta-2C_2\lambda_2=0.
\end{equation}
For fixed $\zeta$ and $\epsilon$, we obtain a problem $\mathbf{P}$ with nominal dimensional deficit equal to one, generically resulting in the existence of a unique one-dimensional solution curve through any regular solution point. In fact, by linearity and homogeneity of the adjoint conditions \eqref{eq:adjobj}-\eqref{eq:adjobjlast} with respect to $\eta$ and the Lagrange multipliers, one such curve is obtained from solutions $(A_1,B_1,C_1,\omega_1,A_2,B_2,C_2,\omega_2)$ to \eqref{eq:const1}-\eqref{eq:const3} together with $\eta=\lambda_1=\cdots=\lambda_7=0$. For the same reason, all solutions with nonzero $\eta$ lie on a straight line with $\omega_1$ and $\omega_2=\omega_1-\epsilon$ that satisfy \eqref{eq:epsdiffcond} and Lagrange multipliers given by the right-hand sides of \eqref{eq:lagrange1}-\eqref{eq:lagrange3} multiplied by $\eta$. Curiously, but not accidentally~\cite{kernevez1987optimization,li2020optimization}, the two curves intersect precisely at a local extremum of $C_1-C_2$ along the first curve, at a singular point of $\mathbf{P}$, as illustrated in the left panel of Fig.~\ref{fig:inflection_opt_finite}.

For this problem, at each iterate of the continuation methodology we construct $\mathbf{P}_0$ by appending a single auxiliary constraint to $\mathbf{P}$. It comes as no surprise that trouble brews on a vicinity of the singular point as local uniqueness fails there for the sought solution to $\mathbf{P}_0$. With some luck, we may be able to step past the singularity along one of the curves, detect such a crossing, and then switch to the other curve. Such a branch-switching strategy~\cite{kuznetsov2013elements,seydel2009practical} may allow us to locate the sought inflection points starting from an arbitrary solution to \eqref{eq:const1}-\eqref{eq:const3} together with $\eta=\lambda_1=\cdots=\lambda_7=0$.

As an alternative, we seek to construct an augmented continuation problem $\mathbf{P}^\ast$ by introducing one additional \emph{a priori} unknown, say $\chi$, such that the two solution curves to $\mathbf{P}$ satisfy $\mathbf{P}^\ast$ for $\chi=0$. With a bit of care, all solutions of $\mathbf{P}^\ast$ near the singular point of $\mathbf{P}$ are regular points of $\mathbf{P}^\ast$. Here, we simply subtract $\epsilon\chi$ from the left-hand side of \eqref{eq:adjobjlast} such that solutions to $\mathbf{P}^\ast$ are obtained only for $\omega_1$, $\omega_2=\omega_1-\epsilon$, $\eta$, and $\chi$ that satisfy the equation
\begin{equation}
    0=\left(\chi-2\eta\frac{3\omega_1^6+5(2\zeta^2-1)\omega_1^4+(16\zeta^4-16\zeta^2+1)\omega_1^2+1-2\zeta^2}{((1-\omega_1^2)^2+4\zeta^2\omega_1^2)^{5/2}}\right)\epsilon+\mathcal{O}(\epsilon^2).
\end{equation}
For sufficiently small $\epsilon$, it follows that the dimensional deficit of $\mathbf{P}^\ast$ (two) equals the dimension of the solution manifold and all solutions near (and including) the singular point of $\mathbf{P}$ are regular points of $\mathbf{P}^\ast$, as demonstrated in the right panel of Fig.~\ref{fig:inflection_opt_finite}. In this case, at each iterate of the continuation methodology we construct a problem $\mathbf{P}^\ast_0$ with zero dimensional deficit by appending two auxiliary scalar constraints to $\mathbf{P}^\ast$.

As the reader may verify, an equivalent set of observations follows from
\begin{itemize}
    \item the substitution of 
\begin{equation}
    \eta-\lambda_3+2\zeta\lambda_1(0)-\dot{\lambda}_1(0)=0,\,-\eta-\lambda_6+2\zeta\lambda_2(0)-\dot{\lambda}_2(0)=0
\end{equation}
in lieu of \eqref{eq:diffeqsub} to generate a problem $\mathbf{P}$ with nominal dimensional deficit equal to one, but with a singular point at the intersection of two one-dimensional curves of solutions $(x_1(\cdot),\omega_1,\theta_1,x_2(\cdot),\omega_2,\theta_2,\eta,\lambda_1(\cdot),\lambda_2(\cdot),\lambda_3,\ldots,\lambda_9)$; followed by \item subtraction of $\epsilon\chi$ from the left-hand side of \eqref{eq:diffadjtlast} to obtain a problem $\mathbf{P}^\ast$ with dimensional deficit equal to two and with all regular points on the corresponding solution manifold near (and including) the singular point of $\mathbf{P}$ (obtained when $\chi=0$).
\end{itemize}

%%%%%%%%%%%%%%%%%%%%%%%%
% Low-precision numerics
%%%%%%%%%%%%%%%%%%%%%%%%
\subsection{Regularizing nearly singular problems for low-precision numerics}
\label{sec: condition numbers}
Section~\ref{sec: lessons and inspirations} refers to general \emph{families} of solutions, resulting from an arbitrary dimensional deficit, instead of just curves as is common in the literature. In the two examples, the construction of $\mathbf{P}^*$ regularizes the continuation problem $\mathbf{P}$ on a neighborhood of the singular point at the intersection of the solution curves to $\mathbf{P}$. The minimal dimension deficit required to remove the singularity is called the \emph{degree of degeneracy} (or \emph{codimension}) of the singularity~\cite{kuznetsov2013elements}. With a dimension deficit equal to two, we may continue through and past the singular point along a two-dimensional manifold of solutions to $\mathbf{P}^*$. As discussed in this section, continuation along multi-dimensional manifolds may also help address the increased demands on computational robustness in nearly singular problems in the presence of significant numerical uncertainty \cite{arnold1972lectures}. This discussion is motivated by applications in which experimental data is used to evaluate the corresponding constraint residuals and their sensitivities with respect to the problem unknowns. When continuation is performed directly on physical experiments~\cite{barton2012control,barton2017control,gonzalez2017assessing}, these quantities are available only with low precision.

\begin{figure}[ht]
  \centering
  \includegraphics[width=0.95\textwidth]{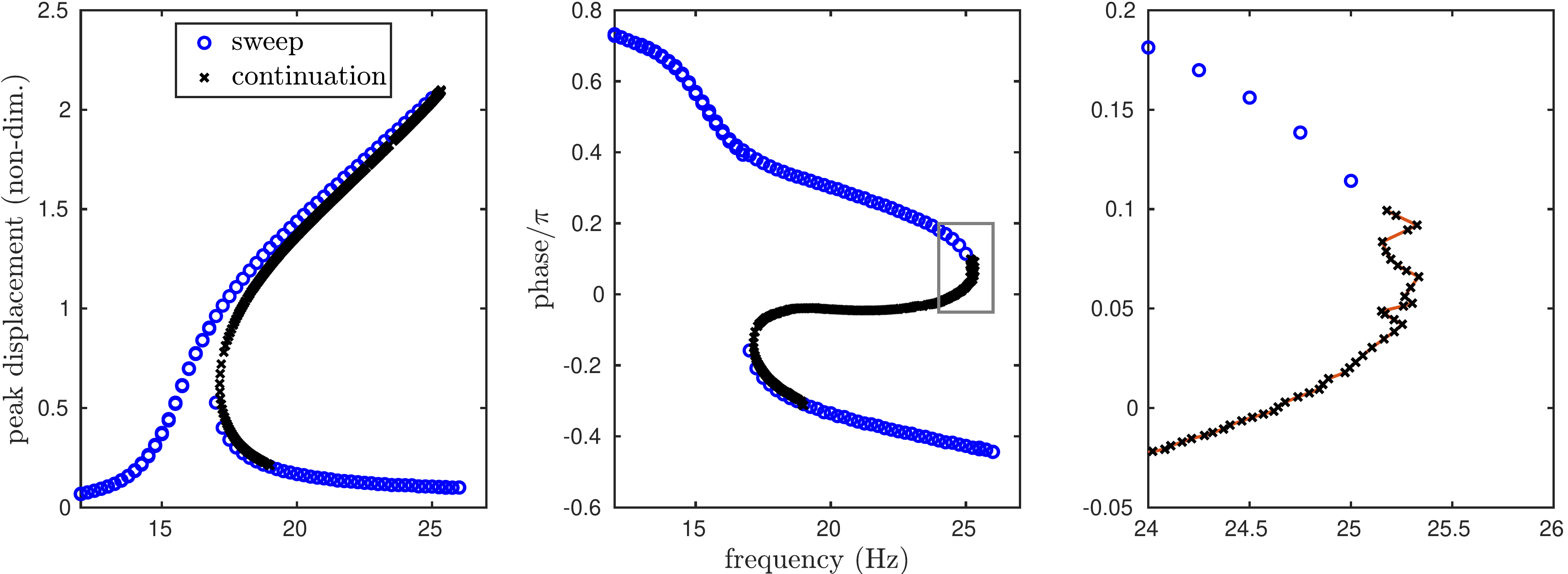}
  \caption[Experimental continuation by \cite{barton2012control}]{The steady-state response of a piezo-electric energy harvesting device to harmonic base excitation obtained using simple parameter sweeps (without feedback control) and control-based continuation: maximal peak displacement (left) and phase (middle and right) of resonator, depending on forcing frequency. The phase displayed is relative to the phase of the forcing and shows increased uncertainty near the peak frequency (right). Data obtained from private communication with David Barton (Fig.~5 in \cite{barton2012control} corresponds to left panel).}
  \label{fig:barton}
\end{figure}

An example application is the analysis of the frequency response of weakly forced mechanical structures in the presence of small damping, e.g., the experimental data in Fig.~\ref{fig:barton} for a piezo-electric energy harvester obtained using control-based continuation in~\cite{barton2012control}. Here, small damping and $\mathcal{O}(1)$ dynamics conspire with low-precision numerics~\cite{barton2017control} to result in an apparent loss of smoothness of the recorded solution branch near the top of the resonance peak. For these operating conditions, we anticipate a resonance curve that is a small perturbation of the (zero-damping) singular limit represented by the backbone curve (obtained using continuation in experiments on purely mechanical structures in~\cite{renson2016robust}). For finite, but small, damping $\zeta$ and fixed forcing phase and amplitude we expect problem condition numbers proportional to $\zeta^{-1}$ near the top of the resonance peak and $\zeta^{-\kappa}$ for $\kappa\in(0,1)$ near its base. 

To illustrate the regularizing benefits of multi-dimensional continuation in a case where they can be demonstrated explicitly, we consider again the linear harmonically forced oscillator. The equations
\begin{align}
    (1-\omega^2)A+2\omega\zeta B-a=0,\,(1-\omega^2)B-2\omega\zeta A-b=0
    \label{eq:condnumbeqs}
\end{align}
for fixed, small $\zeta$ and under variations in $\omega$ are obtained from the ansatz $A\cos\omega t+B\sin\omega t$ substituted into the differential constraint
\begin{equation}\label{condnum:ode}
    \ddot{x}+2\zeta\dot{x}+x=a\cos\omega t +b\sin\omega t.
\end{equation}
In contrast to Section~\ref{sec: inflection points}, we here anticipate small values of the forcing amplitude $\sqrt{a^2+b^2}$ preventing a reduction to the normal form in \eqref{eq:linosc} using an $\mathcal{O}(1)$ scaling. We also keep the phase of the forcing free in \eqref{condnum:ode}.

\begin{figure}[ht]
    \centering
    \includegraphics[width=0.9\textwidth]{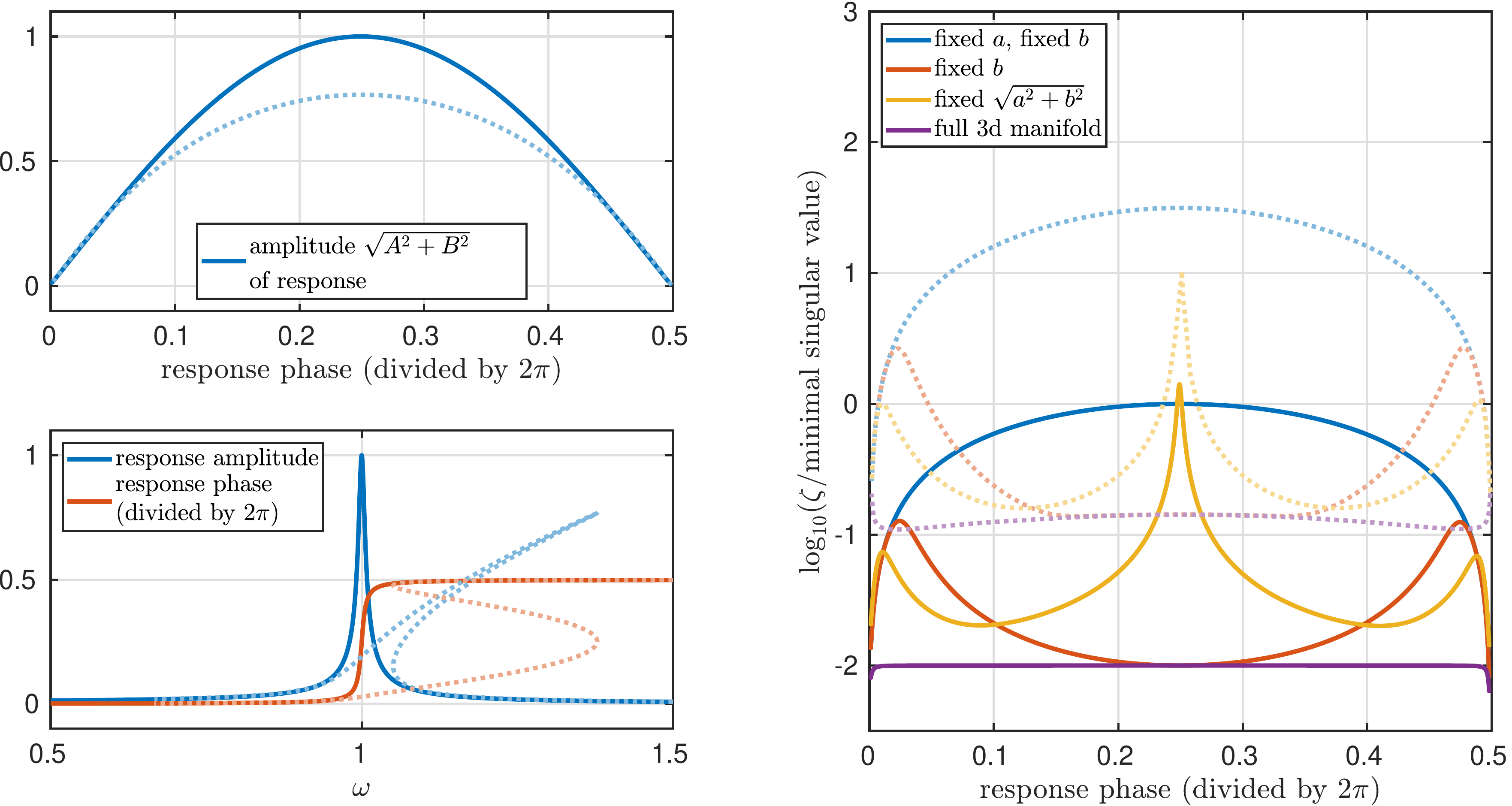}
    \caption{(Top left, solid curves) Illustration of phase-amplitude relation for response $A\cos\omega t+B\sin\omega t$; phase of response is always relative to phase of forcing $a\cos\omega t+b\sin\omega t$. (Bottom left, solid curves) Illustration of amplitude and phase of response depending on forcing frequency relative to the resonance, $\omega-1$. (Right, solid curves) Minimal singular values for all four continuation scenarios in Section~\ref{sec: condition numbers}, depending on phase of response (which parameterizes the resonance peak according to bottom left panel). Dotted curves in the background show corresponding analysis for nonlinear oscillator $\ddot x+2\zeta\dot x+x+x^3=a\cos\omega t+b\sin\omega t$ (with $\sqrt{\smash[b]{a^2+b^2}}=1.5\times10^{-2}$) obtained using \textsc{coco} from the Jacobian of the discretization of the periodic boundary-value problem implemented in the \mcode{coll} toolbox. Here, $\zeta=5\times10^{-3}$.}
    \label{fig:sens:linear}
\end{figure}

As before, we may parameterize $A$ and $B$ explicitly in terms of $a$, $b$, and $\omega$:
\begin{equation}\label{condnum:sol}
    A=\frac{(1-\omega^2)a-2\zeta\omega b}{(1-\omega^2)^2+4\zeta^2\omega^2},\,B=\frac{(1-\omega^2)b+2\zeta\omega a}{(1-\omega^2)^2+4\zeta^2\omega^2}.
\end{equation}
As a consequence, we obtain the rotational symmetry
\begin{equation}
    \begin{pmatrix}a\\b\end{pmatrix}\mapsto R(\phi)\begin{pmatrix}a\\b\end{pmatrix}\Rightarrow\begin{pmatrix}A\\B\end{pmatrix}\mapsto R(\phi)\begin{pmatrix}A\\B\end{pmatrix},\text{ for }R(\phi)=\begin{pmatrix}\cos\phi & -\sin\phi\\\sin\phi & \cos\phi\end{pmatrix}.
\end{equation}
In particular, for fixed $a^2+b^2$, solutions lie on concentric circles in the $(A,B)$ plane of radii
\begin{equation}
    \sqrt{A^2+B^2}=\frac{\sqrt{a^2+b^2}}{\sqrt{(1-\omega^2)^2+4\zeta^2\omega^2}}\label{eq:ABradii}
\end{equation}
and centered on the origin. These are $\mathcal{O}(1)$ when $\sqrt{a^2+b^2},\omega-1=\mathcal{O}(\zeta)$ corresponding to a sharp resonance peak in the frequency-response diagram (see Fig.~\ref{fig:sens:linear}, bottom left panel).

We may explore the full 3-dimensional manifold of solutions to \eqref{eq:condnumbeqs} by considering embedded submanifolds obtained, for example, by holding subsets of variables fixed. Consider, for example, the one-dimensional solution manifold obtained by holding both $a$ and $b$ fixed. We obtain the tangent space at each solution point from the nullspace of the Jacobian matrix of \eqref{eq:condnumbeqs} with respect to $A$, $B$, and $\omega$:
\begin{equation}
    \begin{pmatrix}
    1-\omega^2 & 2\zeta\omega &2\zeta B-2A\omega \\-2\zeta\omega &1-\omega^2 &-2\zeta A-2B\omega
    \end{pmatrix},
\end{equation}
whose singular values equal
\begin{equation}
    \sqrt{(1-\omega^2)^2+4\zeta^2\omega^2},\,\sqrt{(1-\omega^2)^2+4\zeta^2\omega^2+4(A^2+B^2)(\zeta^2+\omega^2)}.
\end{equation}
For $\sqrt{a^2+b^2},\omega-1=\mathcal{O}(\zeta)$, these are $\mathcal{O}(\zeta)$ and $\mathcal{O}(1)$, respectively. The smaller singular value of order $\mathcal{O}(\zeta)$ will cause a great degree of uncertainty of the tangent direction if low-precision computations are performed in this asymptotic limit, even as the radius in \eqref{eq:ABradii} is $\mathcal{O}(1)$. The right panel in Fig.~\ref{fig:sens:linear} quantifies this sensitivity along the resonance peak in terms of the minimal singular value, parameterized using the phase of the response $A\cos\omega t+B\sin\omega t$ relative to the forcing $a\cos\omega t+b\sin\omega t$ (curve for fixed $a$, $b$).

Consider, instead, the two-dimensional solution manifold obtained by holding the forcing amplitude fixed while allowing the phase of the forcing to vary. We may obtain the tangent space at each solution point from the nullspace of the Jacobian matrix of \eqref{eq:condnumbeqs}, augmented by a scaled equation keeping $a^2+b^2$ fixed, with respect to $A$, $B$, $\omega$, $a$, and $b$:
\begin{equation}
    \begin{pmatrix}
    1-\omega^2\quad & 2\zeta\omega\quad  &\phantom{--}2\zeta B-2A\omega\quad & -1\  & \phantom{-}0 \\
    -2\zeta\omega\quad &1-\omega^2 &-2\zeta A-2B\omega  & \phantom{-}0 & -1\\0 & 0 & 0 & \phantom{-}a & \phantom{-}b
    \end{pmatrix}.
    \label{eq:mat1}
\end{equation}
By the rotational symmetry, without loss of generality, we may consider a point on the solution manifold with $b=0$. At this point, the condition number of the matrix in \eqref{eq:mat1} behaves as the inverse of the 2-norm of the $(A,B,\omega)$ columns of the top row of \eqref{eq:mat1} which evaluates to
\begin{equation}
    \sqrt{(1-\omega^2)^2+4\zeta^2\omega^2+\frac{4a^2\omega^2(2\zeta^2-1+\omega^2)^2}{((1-\omega^2)^2+4\zeta^2\omega^2)^2}}.
    \label{eq:matnorm1}
\end{equation}
For $a,\omega-1=\mathcal{O}(\zeta)$ this is $\mathcal{O}(1)$. On the other hand, with the ansatz $a=\mathcal{O}(\zeta)$ and $\omega-1=\mathcal{O}(\zeta^n)$, we find in the asymptotic limit that \eqref{eq:matnorm1} has a local minimum of $\mathcal{O}(\zeta)$ at $\omega-1=-\zeta^2$, local maxima of $\mathcal{O}(1)$ at $\omega-1=\pm\zeta$, and another pair of local minima of $\mathcal{O}(\sqrt{\smash[b]{\zeta}})$ at $\omega-1=\pm\sqrt{\smash[b]{|a|}}/2$. In Fig.~\ref{fig:sens:linear}, the curve for fixed $\sqrt{\smash[b]{a^2+b^2}}$ shows the scaled inverse of \eqref{eq:matnorm1}.
We thus see an improvement in the condition number across the resonance peak where the response amplitude is $\mathcal{O}(1)$, except for the asymptotic limit $a,\sqrt{\omega-1}=\mathcal{O}(\zeta)$, where low-precision computations again would result in a great degree of uncertainty of the corresponding tangent directions. To a lesser degree, increased uncertainty is exhibited also for the asymptotic limit $a,(\omega-1)^2=\mathcal{O}(\zeta)$.

As a third alternative, consider the two-dimensional manifold obtained by holding $b$ fixed at $0$ corresponding to fixing the phase of the forcing. Here, we obtain the tangent space at each solution point from the nullspace of the Jacobian matrix of \eqref{eq:condnumbeqs}, augmented by the equation $b=0$, with respect to $A$, $B$, $\omega$, $a$, and $b$:
\begin{equation}
    \begin{pmatrix}
    1-\omega^2\quad & 2\zeta\omega\quad & \phantom{-}2\zeta B-2A\omega\quad & -1 & \phantom{-}0\\
    -2\zeta\omega\quad & 1-\omega^2\quad & -2\zeta A-2B\omega\quad & \phantom{-}0 & -1\\
    0 & 0 & 0 & \phantom{-}0 & \phantom{-}1
    \end{pmatrix}.\label{cond:mat3}
\end{equation}
The condition number of this matrix behaves as the inverse of the 2-norm of the $(A,B,\omega)$ columns of the second row of \eqref{cond:mat3}, which evaluates to
\begin{equation}
    \sqrt{(1-\omega^2)^2+4\zeta^2\omega^2+\frac{4a^2\zeta^2(1+\omega^2)^2}{((1-\omega^2)^2+4\zeta^2\omega^2)^2}}.
    \label{eq:matnorm2}
\end{equation}
This is again $\mathcal{O}(1)$ for $a,\omega-1=\mathcal{O}(\zeta)$. On the other hand, the ansatz $a=\mathcal{O}(1)$ and $\omega-1=\mathcal{O}(\zeta^n)$ yields local maxima and minima, respectively, in the asymptotic limit at
\begin{equation}
    \omega-1=-\frac{\zeta^2}{2},\,\pm\frac{(2a\zeta)^{1/3}}{\sqrt{2}},
\end{equation}
where the norm in \eqref{eq:matnorm2} behaves as $\mathcal{O}(1)$ and $\mathcal{O}(\zeta^{2/3})$, respectively. With this choice of fixed variables, we have entirely eliminated the poor performance with low-precision computations near the resonant peak, but retain increased uncertainty for the asymptotic limit $a,(\omega-1)^{3/2}=\mathcal{O}(\zeta)$. The right panel in Fig.~\ref{fig:sens:linear} shows that the sensitivity curve for fixed $b$ is lower than the previous two approaches over most of the resonance peak, but still goes to infinity for $\zeta\to0$.

We obtain the full three-dimensional solution manifold by allowing simultaneous variations in $A$, $B$, $\omega$, $a$, and $b$. In this case, the tangent space at each solution point is obtained from the nullspace of the Jacobian matrix
\begin{equation}
    \begin{pmatrix}
    1-\omega^2\quad & 2\zeta\omega\quad & 2\zeta B-2A\omega\quad & -1 & \phantom{-}0\\
    -2\zeta\omega\quad & 1-\omega^2\quad & -2\zeta A-2B\omega\quad & \phantom{-}0 & -1
    \end{pmatrix},
\end{equation}
whose condition number is $\mathcal{O}(1)$ everywhere as $a$ and $b$ may be obtained from the other three variables without the possibility of singularities.  The right panel in Fig.~\ref{fig:sens:linear} shows that the curve for sensitivity when continuing the full manifold is uniformly bounded for small $\zeta$.

While the above analysis was performed for the linear oscillator to permit explicit expressions for the computational sensitivity to low-precision numerics, its predictions remain qualitatively true also for nonlinear oscillators. This is illustrated in Fig.~\ref{fig:sens:linear} for the forced Duffing oscillator $\ddot x+2\zeta \dot x+x+x^3=a \cos\omega t+b\sin\omega t$ (dotted curves obtained using the singular values of the Jacobian of a discretization of the nonlinear periodic-orbit problem) which exhibit similar asymptotic behavior in the limit $\zeta\to0$ as for the linear oscillator. The near constant ratio between the norm of the inverse of the Jacobian in the numerical implementation and the idealized analysis is $\sim N$ (the number of collocation intervals used by the \mcode{coll} toolbox). We leave it to the reader to derive analytical predictions like those found in this section by applying an appropriate perturbation method to this nonlinear problem. We conclude that formulations with dimensional deficits greater than one should be considered whenever one-dimensional bifurcation diagrams appear sensitive to small changes in the problem.

%%%%%%%%%%%%%%%%%%%%
% The COCO formalism
%%%%%%%%%%%%%%%%%%%%
\subsection{The \textsc{coco} formalism}
\label{sec: The coco formalism}

The discussion in Sections~\ref{sec: lessons and inspirations} and \ref{sec: condition numbers} highlights the merits of considering problem construction separately from problem analysis. First decide what you want to do. Then figure out how to do it. The description of the \textsc{coco} construction framework in this section continues in this spirit.

In the general case, we consider continuation problems $\mathbf{P}$ of the form
\begin{equation}
\label{eq: general continuation}
    \Phi(u)=0
\end{equation}
for some Frech\'{e}t differentiable operator $\Phi:\mathcal{U}_\Phi\rightarrow\mathcal{R}_\Phi$ with Banach space domain $\mathcal{U}_\Phi$ and range $\mathcal{R}_\Phi$. At this level of abstraction, there are no distinguishing features of either domain or range. We do not unnecessarily presuppose a dimensional deficit nor assume a particular decomposition of $\mathcal{U}_\Phi$. Instead, we design a general continuation methodology that is accommodating of different dimensional deficits and independent of any substructure of $\mathcal{U}_\Phi$. 

A specialized form of the continuation problem $\mathbf{P}$ in \eqref{eq: general continuation} that self-referentially contains the form in \eqref{eq: general continuation} is given by the extended continuation problem $\mathbf{E}$~\cite{dankowicz2011extended,dankowicz2013recipes} of the form
\begin{equation}
    \label{eq: extended continuation}
    \begin{pmatrix}\Phi(u)\\\Psi(u)-\mu\end{pmatrix}=0
\end{equation}
in terms of the \textit{zero functions} $\Phi$, \textit{continuation variables} $u\in\mathcal{U}_\Phi$, \textit{monitor functions} $\Psi:\mathcal{U}_\Phi\rightarrow\mathbb{R}^{n_\Psi}$, and \textit{continuation parameters} $\mu\in\mathbb{R}^{n_{\Psi}}$. In the special case that $\Psi$ projects onto a finite subspace of $\mathcal{U}_\Phi$, the corresponding $\mu$ amount only to a labeling of these components. More generally, $\mu$ tracks a finite number of solution metrics and, when fixed, restricts attention to a subset of solutions to the \textit{zero problem} $\Phi(u)=0$. The restriction $\mathbf{E}\big|_\mathbb{I}$ obtained by fixing a subset $\mathbb{I}$ of continuation parameters is equivalent to a reduced continuation problem $\mathbf{R}$ in terms of the continuation variables and the remaining continuation parameters. Assuming a dimensional deficit of the zero problem equal to $d$, the number of possible reduced continuation problems equals $2^{\min(d,n_{\Psi})}$.

The examples in Section~\ref{sec: inflection points} illustrate these principles. Each fits the form of \eqref{eq: general continuation} given some association of unknowns with a space $\mathcal{U}_\Phi$ and constraints with $\Phi$. Of course, every such choice, for example those differing by whether $\zeta$ is fixed or free to vary, requires a distinct formulation. In contrast, \eqref{eq: extended continuation} is designed to support every possible choice by including among the monitor functions projections onto all variables that may or may not be designated as fixed during continuation. The decision to hold $\zeta$ fixed may thus be deferred to the moment of analysis, rather than imposed at the time of construction. Similar considerations apply to the problem considered in Section~\ref{sec: condition numbers}, where it is straightforward to construct an extended continuation problem with suitably chosen monitor functions that reduces to each of the four scenarios with appropriate selection of $\mathbb{I}$.

The constrained optimization examples in Section~\ref{sec: inflection points} actually point to a further extension to \eqref{eq: extended continuation} that recognizes the linearity and homogeneity of the modified adjoint conditions in the various Lagrange multipliers, $\eta$, and $\chi$. Additional study also of constrained optimization problems with inequality constraints~\cite{li2020optimization} inspires the definition of an augmented continuation problem $\mathbf{A}$ of the form
\begin{equation}
\label{eq:augmented continuation problem}
\left(\begin{array}{c}\Phi(u)\\\Psi(u)-\mu\\\Lambda^\ast(u)\lambda\\\Xi(u,\lambda,v)\\\Theta(u,\lambda,v)-\nu\end{array}\right)=0
\end{equation}
in terms of the \emph{zero functions} $\Phi$, \emph{monitor functions} $\Psi$, \emph{adjoint functions} $\Lambda^\ast$, \emph{complementary zero functions} $\Xi$, and \emph{complementary monitor functions} $\Theta$, as well as collections of \emph{continuation variables} $u$, \emph{continuation parameters} $\mu$, \emph{continuation multipliers} $\lambda$, \emph{complementary continuation variables} $v$, and \emph{complementary continuation parameters} $\nu$. The form~\eqref{eq:augmented continuation problem} again self-referentially contains both \eqref{eq: general continuation} and \eqref{eq: extended continuation} with $\mu$ and $\nu$ designated as variables that may be fixed or allowed to vary at the moment of analysis.

The augmented continuation problem in~\eqref{eq:augmented continuation problem} is a generalization of Eq.~(30) in~\cite{li2020optimization} for locating solutions to constrained optimization problems with equality and inequality constraints using continuation techniques:
\begin{equation}
\label{eq:Lipaper}
    \begin{pmatrix}
    \Phi(u)\\\Psi(u)-\mu\\\left(D\Phi(u)\right)^*\lambda+\left(D\Psi(u)\right)^*\eta+\left(D G(u)\right)^*\sigma\\\eta-\nu\\K(\sigma,-G(u))-\kappa
    \end{pmatrix}=0,
\end{equation}
where a finite subset of elements of $\Psi$ evaluate to the inequality function $G$. Here, the complementarity conditions of Karush-Kuhn-Tucker theory~\cite{ben1982unified} are expressed in terms of complementarity functions $K$ that must vanish at extrema. We obtain \eqref{eq:Lipaper} as a special case of \eqref{eq:augmented continuation problem} by defining $\Lambda^
\ast$ in terms of the collection of adjoint operators $\left(D\Phi(u)\right)^*$, $\left(D\Psi(u)\right)^*$, and $\left(D G(u)\right)^*$, and by designating the collection of Lagrange multipliers $(\lambda,\eta,\sigma)$ as the corresponding vector of continuation multipliers. Linearity follows from the additive decomposition of the constraint Lagrangian into terms coupling individual constraints with the corresponding adjoint variables. In this notation, $(\nu,\kappa)$ are complementary continuation parameters and the last two rows of \eqref{eq:Lipaper} define the corresponding complementary monitor functions $\Theta$. In practice, we often substitute relaxed complementarity functions that are smooth everywhere (the complementarity functions used in~\cite{li2020optimization} are nonsmooth at origin). Such relaxed functions are then parameterized by additional complementary continuation variables that, in turn, may be associated with complementary continuation parameters in order to consider variations that stiffen the constraint. We do not consider inequality constraints in this paper, but will have use for both $\Xi$ and $\Theta$ in later sections.

Suitably discretized, the augmented continuation problem $\mathbf{A}$ in \eqref{eq:augmented continuation problem} is the most general type of continuation problem supported by the most recent release of the \textsc{coco} platform~\citep{COCO}. Here, $u\in\mathbb{R}^{n_u}$, $\lambda\in\mathbb{R}^{n_\lambda}$, $v\in\mathbb{R}^{n_v}$, $\mu\in\mathbb{R}^{n_\Psi}$ and $\nu\in\mathbb{R}^{n_\Theta}$, while $\Phi:\mathbb{R}^{n_u}\rightarrow\mathbb{R}^{n_\Phi}$, $\Psi:\mathbb{R}^{n_u}\rightarrow\mathbb{R}^{n_\Psi}$, $\Lambda:\mathbb{R}^{n_u}\rightarrow\mathbb{R}^{n_\lambda\times n_\Lambda}$, $\Xi:\mathbb{R}^{n_u}\times\mathbb{R}^{n_\lambda}\times\mathbb{R}^{n_v}\rightarrow\mathbb{R}^{n_\Xi}$, and $\Theta:\mathbb{R}^{n_u}\times\mathbb{R}^{n_\lambda}\times\mathbb{R}^{n_v}\rightarrow\mathbb{R}^{n_\Theta}$, and $\Lambda^\ast=\Lambda^\text{T}$. We obtain a restricted continuation problem $\mathbf{A}\big|_{\mathbb{I}_\mu,\mathbb{I}_\nu}$ by designating subsets $\{\mu_i,\,i\in\mathbb{I}_\mu\}$ and $\{\nu_i,\,i\in\mathbb{I}_\nu\}$ as fixed. The resulting restricted continuation problem then has nominal dimensional deficit equal to $n_u+n_\lambda+n_v-n_\Phi-n_\Lambda-n_\Xi-\left|\mathbb{I}_\mu\right|-\left|\mathbb{I}_\nu\right|$.

While there may be some merit in the level of abstraction of the augmented continuation problem purely from an organizational viewpoint, it truly comes into its own when coupled with a systematic paradigm of problem construction. This is one of the features of the \textsc{coco} software platform. The reader may refer to Chapter 3 of~\cite{dankowicz2013recipes} for an earlier discussion that applies to the extended continuation problem~\eqref{eq: extended continuation}. %In the case of constrained optimization, this is motivated by the linearity of the constraint Lagrangian~\cite{li2017staged,li2020optimization}. 

It is a truism that a given (finite-dimensional) augmented continuation problem $\mathbf{A}$ may be interpreted as the largest element of a chain
\begin{equation}
\label{eq:chain}
    \emptyset=\mathbf{A}_0\subsetpluseq\mathbf{A}_1\subsetpluseq\cdots\subsetpluseq\mathbf{A}_N=\mathbf{A}
\end{equation}
of augmented continuation problems $\mathbf{A}_i$, where $\tilde{\mathbf{A}}\subsetpluseq\hat{\mathbf{A}}$ if 
\begin{gather}
    \tilde{n}_{u/\lambda/v/\Phi/\Psi/\Lambda/\Xi/\Theta}\le \hat{n}_{u/\lambda/v/\Phi/\Psi/\Lambda/\Xi/\Theta},\\
    \hat{\Phi}_{\left(1:\tilde{n}_\Phi\right)}(u)=\tilde{\Phi}\left(u_{\left(1:\tilde{n}_u\right)}\right),\,\hat{\Psi}_{\left(1:\tilde{n}_\Psi\right)}(u)=\tilde{\Psi}\left(u_{\left(1:\tilde{n}_u\right)}\right),\\
    \hat{\Lambda}_{\left(1:\tilde{n}_\lambda,1:\tilde{n}_\Lambda\right)}(u)=\tilde{\Lambda}\left(u_{\left(1:\tilde{n}_u\right)}\right),\,
    %\hat{\Lambda}_{\left(\tilde{n}_\lambda+1:\hat{n}_\lambda,1:\tilde{n}_\Lambda\right)}(u)=0\label{eq:LambdaHat},
    \hat{\Lambda}_{\left(1:\tilde{n}_\lambda,\tilde{n}_\Lambda+1:\hat{n}_\Lambda\right)}(u)=0\\
    \hat{\Xi}_{\left(1:\tilde{n}_\Xi\right)}(u,\lambda,v)=\tilde{\Xi}\left(u_{\left(1:\tilde{n}_u\right)},\lambda_{\left(1:\tilde{n}_\lambda\right)},v_{\left(1:\tilde{n}_v\right)}\right),\\
    \hat{\Theta}_{\left(1:\tilde{n}_\Theta\right)}(u,\lambda,v)=\tilde{\Theta}\left(u_{\left(1:\tilde{n}_u\right)},\lambda_{\left(1:\tilde{n}_\lambda\right)},v_{\left(1:\tilde{n}_v\right)}\right),
\end{gather}
and
\begin{equation}
    \hat{\mu}_{\left(1:\tilde{n}_\Psi\right)}=\tilde{\mu},\,\hat{\nu}_{\left(1:\tilde{n}_\Xi\right)}=\tilde{\nu},
\end{equation}
and where $\emptyset$ denotes an \emph{empty} continuation problem with $n_u=n_\lambda=n_v=n_\Phi=n_\Psi=n_\Lambda=n_\Xi=n_\Theta=0$. The chain in \eqref{eq:chain} represents a sequential embedding of partial realizations of $\mathbf{A}$ into successively larger problems with additional unknowns and additional constraints. Since $\emptyset\subsetpluseq\mathbf{A}$ for any $\mathbf{A}$, we obtain a nontrivial decomposition of $\mathbf{A}$ in the form of \eqref{eq:chain} when at least one of the partial realizations is nonempty and distinct from $\mathbf{A}$. Given an augmented continuation problem $\mathbf{A}$ with $n_\Phi+n_\Psi+n_\Lambda+n_\Xi+n_\Theta>1$, it is always possible to find a nontrivial decomposition \eqref{eq:chain} for some equivalent augmented continuation problem obtained by reordering the elements of $\Phi$, $\Psi$, $\Lambda$, $\Xi$, $\Theta$, $u$, $\lambda$, $v$, $\mu$, and $\nu$.

Given a chain decomposition \eqref{eq:chain}, there exists, for each $i$, four ordered index sets 
\begin{gather}
    \{n_{u,i-1}+1,\ldots,n_{u,i}\}\subseteq\mathbb{K}_{u,i}\subseteq\{1,\ldots,n_{u,i}\},\\
    \{n_{\lambda,i-1}+1,\ldots,n_{\lambda,i}\}\subseteq\mathbb{K}_{\lambda,i}\subseteq\{1,\ldots,n_{\lambda,i}\},\\
    \{n_{\Lambda,i-1}+1,\ldots,n_{\Lambda,i}\}\subseteq\mathbb{K}_{\Lambda,i}\subseteq\{1,\ldots,n_{\Lambda,i}\},\\
    \{n_{v,i-1}+1,\ldots,n_{v,i}\}\subseteq\mathbb{K}_{v,i}\subseteq\{1,\ldots,n_{v,i}\},
\end{gather}
such that
\begin{align}
    \Phi_{\left(n_{\Phi,i-1}+1:n_{\Phi,i}\right)}(u)&=\mathfrak{phi}^{(i)}\left(u_{\mathbb{K}_{u,i}}\right),\\
    \Psi_{\left(n_{\Psi,i-1}+:n_{\Psi,i}\right)}(u)&=\mathfrak{psi}^{(i)}\left(u_{\mathbb{K}_{u,i}}\right),\\
    %\Lambda_{\left(\mathbb{K}_{\lambda,i},n_{\Lambda,i-1}+1:n_{\Lambda,i}\right)}(u)&=\mathfrak{lambda}^{(i)}\left(u_{\mathbb{K}_{u,i}}\right),\label{eq:lambdaAti}\\
    \Lambda_{\left(n_{\lambda,i-1}+1:n_{\lambda,i},\mathbb{K}_{\Lambda,i}\right)}(u)&=\mathfrak{lambda}^{(i)}\left(u_{\mathbb{K}_{u,i}}\right),\label{eq:lambdaAti}\\
    \Xi_{\left(n_{\Xi,i-1}+1:n_{\Xi,i}\right)}(u,\lambda,v)&=\mathfrak{xi}^{(i)}\left(u_{\mathbb{K}_{u,i}},\lambda_{\mathbb{K}_{\lambda,i}},v_{\mathbb{K}_{v,i}}\right),\\
    \Theta_{\left(n_{\Theta,i-1}+1:n_{\Theta,i}\right)}(u,\lambda,v)&=\mathfrak{theta}^{(i)}\left(u_{\mathbb{K}_{u,i}},\lambda_{\mathbb{K}_{\lambda,i}},v_{\mathbb{K}_{v,i}}\right),
\end{align}
and
\begin{equation}
\label{eq:LambdaOldNewUzero}
    %\Lambda_{\left(\{1,\ldots,n_{\lambda,i}\}\setminus\mathbb{K}_{\lambda,i},n_{\Lambda,i-1}+1:n_{\Lambda,i}\right)}(u)=0
    \Lambda_{\left(n_{\lambda,i-1}+1:n_{\lambda,i},\{1,\ldots,n_{\Lambda,i}\}\setminus\mathbb{K}_{\Lambda,i}\right)}(u)=0
\end{equation}
for some functions $\mathfrak{phi}^{(i)}:\mathbb{R}^{|\mathbb{K}_{u,i}|}\rightarrow\mathbb{R}^{n_{\Phi,i}-n_{\Phi,i-1}}$, $\mathfrak{psi}^{(i)}:\mathbb{R}^{|\mathbb{K}_{u,i}|}\rightarrow\mathbb{R}^{n_{\Psi,i}-n_{\Psi,i-1}}$, $\mathfrak{lambda}^{(i)}:\mathbb{R}^{|\mathbb{K}_{u,i}|}\rightarrow\mathbb{R}^{(n_{\lambda,i}-n_{\lambda,i-1})\times|\mathbb{K}_{\Lambda,i}|}$, $\mathfrak{xi}^{(i)}:\mathbb{R}^{|\mathbb{K}_{u,i}|}\times\mathbb{R}^{|\mathbb{K}_{\lambda,i}|}\times\mathbb{R}^{|\mathbb{K}_{v,i}|}\rightarrow\mathbb{R}^{n_{\Xi,i}-n_{\Xi,i-1}}$, and $\mathfrak{theta}^{(i)}:\mathbb{R}^{|\mathbb{K}_{u,i}|}\times\mathbb{R}^{|\mathbb{K}_{\lambda,i}|}\times\mathbb{R}^{|\mathbb{K}_{v,i}|}\rightarrow\mathbb{R}^{n_{\Theta,i}-n_{\Theta,i-1}}$. We refer to these functions as representations of the corresponding left-hand sides and to $\mathbb{K}_{u,i}$, $\mathbb{K}_{\lambda,i}$, $\mathbb{K}_{\Lambda,i}$, and $\mathbb{K}_{v,i}$ as the corresponding dependency index sets.

We now arrive at a paradigm of decomposition of an augmented continuation problem $\mathbf{A}$ through a sequence of partial realizations $\mathbf{A}_i$ constructed sequentially in terms of the representations $\mathfrak{phi}^{(i)}$, $\mathfrak{psi}^{(i)}$, $\mathfrak{lambda}^{(i)}$, $\mathfrak{xi}^{(i)}$, and $\mathfrak{theta}^{(i)}$ and the dependency index sets $\mathbb{K}_{u,i}$, $\mathbb{K}_{\lambda,i}$, $\mathbb{K}_{\Lambda,i}$, and $\mathbb{K}_{v,i}$. Since we must associate an initial solution guess $(u_0,\lambda_0,v_0)$ to $\mathbf{A}$, we may construct the dependency index sets $\mathbb{K}_{u,i}$, $\mathbb{K}_{\lambda,i}$, and $\mathbb{K}_{v,i}$ in terms of the index sets 
\begin{gather}
    \mathbb{K}^\mathrm{o}_{u,i}=\mathbb{K}_{u,i}\setminus \{n_{u,i-1}+1,\ldots,n_{u,i-1}+(n_{u,i}-n_{u,i-1})\},\\
    \mathbb{K}^\mathrm{o}_{\lambda,i}=\mathbb{K}_{\lambda,i}\setminus \{n_{\lambda,i-1}+1,\ldots,n_{\lambda,i-1}+(n_{\lambda,i}-n_{\lambda,i-1})\},\\
    \mathbb{K}^\mathrm{o}_{v,i}=\mathbb{K}_{v,i}\setminus \{n_{v,i-1}+1,\ldots,n_{v,i-1}+(n_{v,i}-n_{v,i-1})\},
\end{gather}
and the cardinalities $|u_{0,n_{u,i-1}+1:n_{u,i}}|=n_{u,i}-n_{u,i-1}$, $|\lambda_{0,n_{\lambda,i-1}+1:n_{\lambda,i}}|=n_{\lambda,i}-n_{\lambda,i-1}$, and $|v_{0,n_{v,i-1}+1:n_{v,i}}|=n_{v,i}-n_{v,i-1}$, respectively. Similarly, we obtain the index set $\mathbb{K}_{\Lambda,i}$ from the index set
\begin{equation}
    \mathbb{K}^\mathrm{o}_{\Lambda,i}=\mathbb{K}_{\Lambda,i}\setminus \{n_{\Lambda,i-1}+1,\ldots,n_{\Lambda,i-1}+(n_{\Lambda,i}-n_{\Lambda,i-1})\}
\end{equation}
and the difference between the number of columns of $\mathfrak{lambda}^{(i)}(\cdot)$ and $\mathbb{K}^\mathrm{o}_{\Lambda,i}$, since this must equal $n_{\Lambda,i}-n_{\Lambda,i-1}$.

Rather than considering the decomposition of an existing augmented continuation problem, we may consider its staged construction through the successive application of a sequence of operators on the space of augmented continuation problems. Given an augmented continuation problem $\mathbf{A}$ with initial solution guess $(u_0,\lambda_0,v_0)$ we construct an augmented continuation problem $\hat{\mathbf{A}}$ with initial solution guess $(\hat{u}_0,\hat{\lambda}_0,\hat{v}_0)$ by the application of the operator
\begin{equation}
\label{eq:cocooperator}
    \hat{}:=\big[\mathfrak{phi},\mathfrak{psi},\mathfrak{lambda},\mathfrak{xi},\mathfrak{theta},\mathbb{K}^\mathrm{o}_u,\mathbb{K}^\mathrm{o}_\lambda,\mathbb{K}^\mathrm{o}_\Lambda,\mathbb{K}^\mathrm{o}_v,u^\mathrm{n}_0,\lambda^\mathrm{n}_0,v^\mathrm{n}_0\big]
\end{equation}
in terms of the index sets $\mathbb{K}^\mathrm{o}_u,\mathbb{K}^\mathrm{o}_\lambda,\mathbb{K}^\mathrm{o}_\Lambda,\mathbb{K}^\mathrm{o}_v$, vectors $u^\mathrm{n}_0\in\mathbb{R}^{k_u}$, $\lambda^\mathrm{n}_0\in\mathbb{R}^{k_\lambda}$, and $v^\mathrm{n}_0\in\mathbb{R}^{k_v}$, functions $\mathfrak{phi}:\mathbb{R}^{|\mathbb{K}^\mathrm{o}_u|+k_u}\rightarrow\mathbb{R}^{k_\Phi}$, $\mathfrak{psi}:\mathbb{R}^{|\mathbb{K}^\mathrm{o}_u|+k_u}\rightarrow\mathbb{R}^{k_\Psi}$, $\mathfrak{lambda}:\mathbb{R}^{|\mathbb{K}^\mathrm{o}_u|+k_u}\rightarrow\mathbb{R}^{k_\lambda\times (|\mathbb{K}^{o}_\Lambda|+k_\Lambda)}$, $\mathfrak{xi}:\mathbb{R}^{|\mathbb{K}^\mathrm{o}_u|+k_u}\times\mathbb{R}^{|\mathbb{K}^\mathrm{o}_\lambda|+k_\lambda}\times\mathbb{R}^{|\mathbb{K}^\mathrm{o}_v|+k_v}\rightarrow\mathbb{R}^{k_\Xi}$ and $\mathfrak{theta}:\mathbb{R}^{|\mathbb{K}^\mathrm{o}_u|+k_u}\times\mathbb{R}^{|\mathbb{K}^\mathrm{o}_\lambda|+k_\lambda}\times\mathbb{R}^{|\mathbb{K}^\mathrm{o}_v|+k_v}\rightarrow\mathbb{R}^{k_\Theta}$,
such that $\hat{u}_0=\begin{pmatrix}u_0,u^\mathrm{n}_0\end{pmatrix}$, $\hat{\lambda}_0=\begin{pmatrix}\lambda_0,\lambda^\mathrm{n}_0\end{pmatrix}$, and $\hat{v}_0=\begin{pmatrix}v_0,v^\mathrm{n}_0\end{pmatrix}$, $\mathbb{K}_u=\mathbb{K}^\mathrm{o}_u\cup\{n_u+1,\ldots,n_u+k_u\}$, $\mathbb{K}_\lambda=\mathbb{K}^\mathrm{o}_\lambda\cup\{n_\lambda+1,\ldots,n_\lambda+k_\lambda\}$,  $\mathbb{K}_\Lambda=\mathbb{K}^\mathrm{o}_\Lambda\cup\{n_\Lambda+1,\ldots,n_\Lambda+k_\Lambda\}$, and $\mathbb{K}_v=\mathbb{K}^\mathrm{o}_v\cup\{n_v+1,\ldots,n_v+k_v\}$, 
\begin{align}
    \hat{\Phi}:\hat{u}&\mapsto\begin{pmatrix}\Phi\left(\hat{u}_{\left(1:n_u\right)}\right)\\\mathfrak{phi}\left(\hat{u}_{\mathbb{K}_u}\right)\end{pmatrix},\\
    \hat{\Psi}:\hat{u}&\mapsto\begin{pmatrix}\Psi\left(\hat{u}_{\left(1:n_u\right)}\right)\\\mathfrak{psi}\left(\hat{u}_{\mathbb{K}_u}\right)\end{pmatrix},\\
    \hat{\Xi}:\left(\hat{u},\hat{\lambda},\hat{v}\right)&\mapsto\begin{pmatrix}\Xi\left(\hat{u}_{\left(1:n_u\right)},\hat{\lambda}_{\left(1:n_\lambda\right)},\hat{v}_{\left(1:n_v\right)}\right)\\\mathfrak{xi}\left(\hat{u}_{\mathbb{K}_u},\hat{\lambda}_{\mathbb{K}_\lambda},\hat{v}_{\mathbb{K}_v}\right)\end{pmatrix},\\
    \hat{\Theta}:\left(\hat{u},\hat{\lambda},\hat{v}\right)&\mapsto\begin{pmatrix}\Theta\left(\hat{u}_{\left(1:n_u\right)},\hat{\lambda}_{\left(1:n_\lambda\right)},\hat{v}_{\left(1:n_v\right)}\right)\\\mathfrak{theta}\left(\hat{u}_{\mathbb{K}_u},\hat{\lambda}_{\mathbb{K}_\lambda},\hat{v}_{\mathbb{K}_v}\right)\end{pmatrix},
\end{align}
and
\begin{align}
    \hat{\Lambda}_{\left(1:n_\lambda,1:n_\Lambda\right)}:\hat{u}&\mapsto\Lambda\left(\hat{u}_{\left(1:n_u\right)}\right),\\
    \hat{\Lambda}_{\left(1:n_\lambda,n_\Lambda+1:n_\Lambda+k_\Lambda\right)}:\hat{u}&\mapsto 0,\\
    \hat{\Lambda}_{\left(n_\lambda+1:n_\lambda+k_\lambda,\mathbb{K}_\Lambda\right)}:\hat{u}&\mapsto\mathfrak{lambda}\left(\hat{u}_{\mathbb{K}_u}\right),\\
    \hat{\Lambda}_{\left(n_\lambda+1:n_\lambda+k_\lambda,\{1,\ldots,n_\Lambda+k_\Lambda\}\setminus\mathbb{K}_\Lambda\right)}:\hat{u}&\mapsto 0.
\end{align}

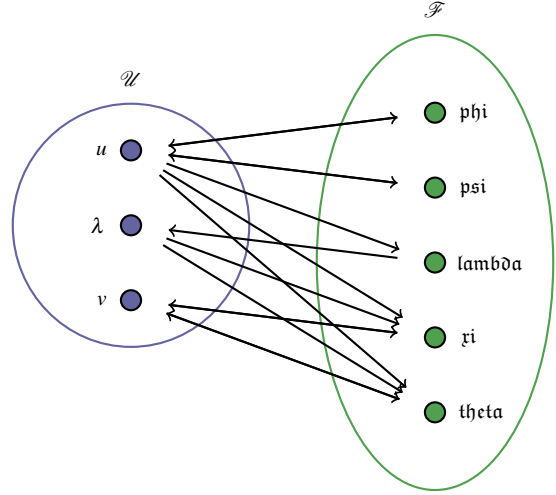
\begin{figure}[ht]
\centering
\begin{tikzpicture}[thick,
  every node/.style={draw,circle},
  fsnode/.style={fill=myblue},
  ssnode/.style={fill=mygreen},
  every fit/.style={ellipse,draw,inner sep=0pt,text width=2.2cm},
  ->,shorten >= 10pt,shorten <= 10pt
]

% the vertices of U
\begin{scope}[start chain=going below,node distance=7mm]
\node[fsnode,on chain,label=left:$u$] (f1) {};
\node[fsnode,on chain,label=left:$\lambda$] (f2) {};
\node[fsnode,on chain,label=left:$v$] (f3) {};
\end{scope}

% the vertices of V
\begin{scope}[xshift=4cm,yshift=0.5cm,start chain=going below,node distance=7mm]
\node[ssnode,on chain,label=right:$\mathfrak{phi}$] (s1) {};
\node[ssnode,on chain,label=right:$\mathfrak{psi}$] (s2) {};
\node[ssnode,on chain,label=right:$\mathfrak{lambda}$] (s3) {};
\node[ssnode,on chain,label=right:$\mathfrak{xi}$] (s4) {};
\node[ssnode,on chain,label=right:$\mathfrak{theta}$] (s5) {};
\end{scope}

% the set U
\node [myblue,fit=(f1) (f3),label=above:$\mathcal{U}$] {};
% the set V
\node [mygreen,fit=(s1) (s5),label=above:$\mathcal{F}$] {};

% the edges
\draw (f1) -- (s1);
\draw (f1) -- (s2);
\draw (f1) -- (s3);
\draw (f1) -- (s4);
\draw (f1) -- (s5);
\draw (f2) -- (s4);
\draw (f2) -- (s5);
\draw (f3) -- (s4);
\draw (f3) -- (s5);
\draw (s1) -- (f1);
\draw (s2) -- (f1);
\draw (s3) -- (f2);
\draw (s4) -- (f3);
\draw (s5) -- (f3);
\end{tikzpicture}
\caption{A directed bipartite graph illustration of the core constructors in \eqref{eq: coco_add_func1}-\eqref{eq: coco_add_adjt}. Here $\mathcal{U}=\{u,\lambda,v\}$ is a set of variables and $\mathcal{F}=\{\mathfrak{phi},\mathfrak{psi},\mathfrak{lambda},\mathfrak{xi},\mathfrak{theta}\}$ is a set of functions. A directed edge from node $A$ in $\mathcal{U}$ to node $B$ in $\mathcal{F}$ indicates that a variable of type $A$ is an argument of a function of type $B$. A directed edge from node $C$ in $\mathcal{F}$ to node $D$ in $\mathcal{U}$ indicates that a new variables of type $D$ can be introduced with the construction of a function of type $C$.}
\label{fig:core-constructors}
\end{figure}

In \textsc{coco}, an operator of the form \eqref{eq:cocooperator} is called a \textit{constructor}. Its core constructors correspond to the special operators
\begin{gather}
    \big[\mathfrak{phi},\emptyset,\emptyset,\emptyset,\emptyset,\mathbb{K}^\mathrm{o}_u,\emptyset,\emptyset,\emptyset,u^\mathrm{n}_0,\emptyset,\emptyset\big],\label{eq: coco_add_func1}\\
    \big[\emptyset,\mathfrak{psi},\emptyset,\emptyset,\emptyset,\mathbb{K}^\mathrm{o}_u,\emptyset,\emptyset,\emptyset,u^\mathrm{n}_0,\emptyset,\emptyset\big],\label{eq: coco_add_func2}\\
    \big[\emptyset,\emptyset,\emptyset,\mathfrak{xi},\emptyset,\mathbb{K}^\mathrm{o}_u,\mathbb{K}^\mathrm{o}_\lambda,\emptyset,\mathbb{K}^\mathrm{o}_v,\emptyset,\emptyset,v^\mathrm{n}_0\big],\label{eq: coco_add_comp1}\\
    \big[\emptyset,\emptyset,\emptyset,\emptyset,\mathfrak{theta},\mathbb{K}^\mathrm{o}_u,\mathbb{K}^\mathrm{o}_\lambda,\emptyset,\mathbb{K}^\mathrm{o}_v,\emptyset,\emptyset,v^\mathrm{n}_0\big],\label{eq: coco_add_comp2}
\end{gather}
and
\begin{equation}
    \big[\emptyset,\emptyset,\mathfrak{lambda},\emptyset,\emptyset,\mathbb{K}^\mathrm{o}_u,\emptyset,\mathbb{K}^\mathrm{o}_\Lambda,\emptyset,\emptyset,\lambda^\mathrm{n}_0,\emptyset\big],\label{eq: coco_add_adjt}
\end{equation}
where, in the last case, $\mathbb{K}^\mathrm{o}_u$ equals $\mathbb{K}_u$ for a previous call to one of the first two core constructors. A bipartite graph illustration of these core constructors and their variable dependence is presented in Fig.~\ref{fig:core-constructors}. Each call to a core constructor is associated with a unique \textit{function identifier} allowing subsequent stages of construction, for example, to reference its index sets. Each definition of a (complementary) monitor function is also associated with unique labels for the corresponding (complementary) continuation parameters, allowing each to be fixed or free to vary during the subsequent continuation analysis. Composition of calls to these core constructors defines the space of operators of the form \eqref{eq:cocooperator} that may be realized in \textsc{coco}.

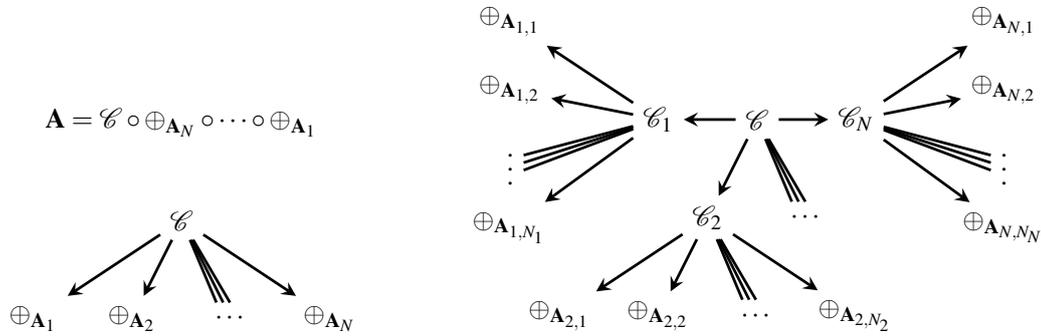
\begin{figure}[ht]
\centering
\scalebox{1.3}{
\begin{tikzpicture}
% define nodes
\node (aug) {$\mathbf{A}=\mathcal{C}\circ\oplus_{\mathbf{A}_N}\circ\cdots\circ\oplus_{\mathbf{A}_1}$};
\node (close) [below of=aug]{$\mathcal{C}$};
\node (a1) [below of=close,xshift=-1.5cm] {$\oplus_{\mathbf{A}_1}$};
\node (a2) [right of=a1] {$\oplus_{\mathbf{A}_2}$};
\node (a34) [right of=a2] {$\cdots$};
\node (an) [right of=a34] {$\oplus_{\mathbf{A}_N}$};
% define edges
\draw [arrow] (close) -- (a1);
\draw [arrow] (close) -- (a2);
\draw [thick] (close) -- (a34);
\draw [thick] ([xshift=1.5em]close) -- (a34);
\draw [thick] ([xshift=-1.5em]close) -- (a34);
\draw [arrow] (close) -- (an);
\end{tikzpicture}}
\qquad\qquad % <----------------- SPACE BETWEEN PICTURES
\scalebox{1.3}{
\begin{tikzpicture}
\node (c) {$\mathcal{C}$};
\node (c2) [below of=c,xshift=-0.5cm] {$\mathcal{C}_2$};
\node (c34) [right of=c2] {$\cdots$};
\node (cn) [right of=c] {$\mathcal{C}_N$};
\node (c1) [left of=c] {$\mathcal{C}_1$};
\draw [arrow] (c) -- (c1);
\draw [arrow] (c) -- (c2);
\draw [arrow] (c) -- (cn);
\draw [thick] (c) -- (c34);
\draw [thick] ([xshift=1.5em]c) -- (c34);
\draw [thick] ([xshift=-1.5em]c) -- (c34);
% subproblem for c1
\node (a11) [left of=c1,xshift=-0.5cm,yshift=1.0cm] {$\oplus_{\mathbf{A}_{1,1}}$};
\node (a12) [below of=a11,yshift=0.3cm] {$\oplus_{\mathbf{A}_{1,2}}$};
\node (a134) [below of=a12,yshift=0.3cm] {$\vdots$};
\node (a1n) [below of=a134,yshift=0.3cm] {$\oplus_{\mathbf{A}_{1,N_1}}$};
\draw [arrow] (c1) -- (a11);
\draw [arrow] (c1) -- (a12);
\draw [thick] (c1) -- (a134);
\draw [thick] ([yshift=-3.0em]c1) -- (a134);
\draw [thick] ([yshift=-6.0em]c1) -- (a134);
\draw [arrow] (c1) -- (a1n);
% subproblem for c2
\node (a21) [below of=c2,xshift=-1.5cm] {$\oplus_{\mathbf{A}_{2,1}}$};
\node (a22) [right of=a21] {$\oplus_{\mathbf{A}_{2,2}}$};
\node (a234) [right of=a22] {$\cdots$};
\node (a2n) [right of=a234] {$\oplus_{\mathbf{A}_{2,N_2}}$};
\draw [arrow] (c2) -- (a21);
\draw [arrow] (c2) -- (a22);
\draw [thick] (c2) -- (a234);
\draw [thick] ([xshift=1.5em]c2) -- (a234);
\draw [thick] ([xshift=-1.5em]c2) -- (a234);
\draw [arrow] (c2) -- (a2n);
% subproblem for cn
\node (an1) [right of=cn,xshift=0.5cm,yshift=1.0cm] {$\oplus_{\mathbf{A}_{N,1}}$};
\node (an2) [below of=an1,yshift=0.3cm] {$\oplus_{\mathbf{A}_{N,2}}$};
\node (an34) [below of=an2,yshift=0.3cm] {$\vdots$};
\node (ann) [below of=an34,yshift=0.3cm] {$\oplus_{\mathbf{A}_{N,N_N}}$};
\draw [arrow] (cn) -- (an1);
\draw [arrow] (cn) -- (an2);
\draw [thick] (cn) -- (an34);
\draw [thick] ([yshift=-3.0em]cn) -- (an34);
\draw [thick] ([yshift=-6.0em]cn) -- (an34);
\draw [arrow] (cn) -- (ann);
\end{tikzpicture}}
\caption{(left) A simple tree representation of the construction of the augmented continuation problem $\mathbf{A}$ in terms of a canonical sum of uncoupled problems coupled through the imposition of gluing conditions. (right) A recursive generalization.}
\label{fig:tree-representations}
\end{figure}

In the special case that $\mathbb{K}^\mathrm{o}_u=\mathbb{K}^\mathrm{o}_\lambda=\mathbb{K}^\mathrm{o}_\Lambda=\mathbb{K}^\mathrm{o}_v=\emptyset$, the augmented continuation problem $\hat{\mathbf{A}}$ obtained by application of the operator $\hat{}\,$ in \eqref{eq:cocooperator} can be defined as the \emph{canonical sum} of two \emph{uncoupled} augmented continuation problems $\mathbf{A}$ and $\tilde{\mathbf{A}}$, such that  $\oplus_{\tilde{\mathbf{A}}}(\mathbf{A})=\mathbf{A}\oplus\tilde{\mathbf{A}}:=\hat{\mathbf{A}}$. 
An arbitrary augmented continuation problem $\mathbf{A}$ may be constructed as the canonical sum of a sequence of uncoupled augmented continuation problems $\{\mathbf{A}_i\}_{i=1}^N$, glued together by the application of an operator $\mathcal{C}$:
\begin{equation}
    \mathbf{A}=\mathcal{C}\circ\oplus_{\mathbf{A}_N}\circ\cdots\circ\oplus_{\mathbf{A}_1},
\end{equation}
and represented graphically in the left panel of Fig.~\ref{fig:tree-representations}. Such a formulation is especially convenient in problems where the  individual operators $\oplus_{\mathbf{A}_i}$ may be sampled from a smaller set of operators, for example when modeling multi-segment boundary-value problems, where the  $\oplus_{\mathbf{A}_i}$ represent contributions associated with individual segments and $\mathcal{C}$ imposes the corresponding boundary conditions, as well as gluing conditions on the problem parameters. This paradigm of construction is naturally nested and recursive, as suggested in the right panel of Fig.~\ref{fig:tree-representations}.

The particular choice of core constructors in \textsc{coco} is not accidental and obviously reflects the unique position of the continuation multipliers $\lambda$ and complementary continuation variables $v$ in the problem hierarchy. This is best appreciated through examples.

\subsection{Data Assimilation}
\label{sec: data assimilation}
We consider in this section an augmented continuation problem obtained naturally from the optimization of an objective functional in the presence of delay differential constraints adapted from~\cite{TraversoMagri2019}. In contrast to this reference, we emphasize below the form of the resultant necessary conditions and describe a solution strategy similar to that presented in Section~\ref{sec: lessons and inspirations}.

\subsubsection{Problem formulation}
\label{sec: data assimilation problem formulation}

%https://link.springer.com/chapter/10.1007%2F978-3-030-22747-0_12

\begin{figure}[ht]
\centering
\subfloat{\includegraphics[width=0.45\columnwidth,trim={0 {0.0\textwidth} 0 0.0\textwidth},clip]{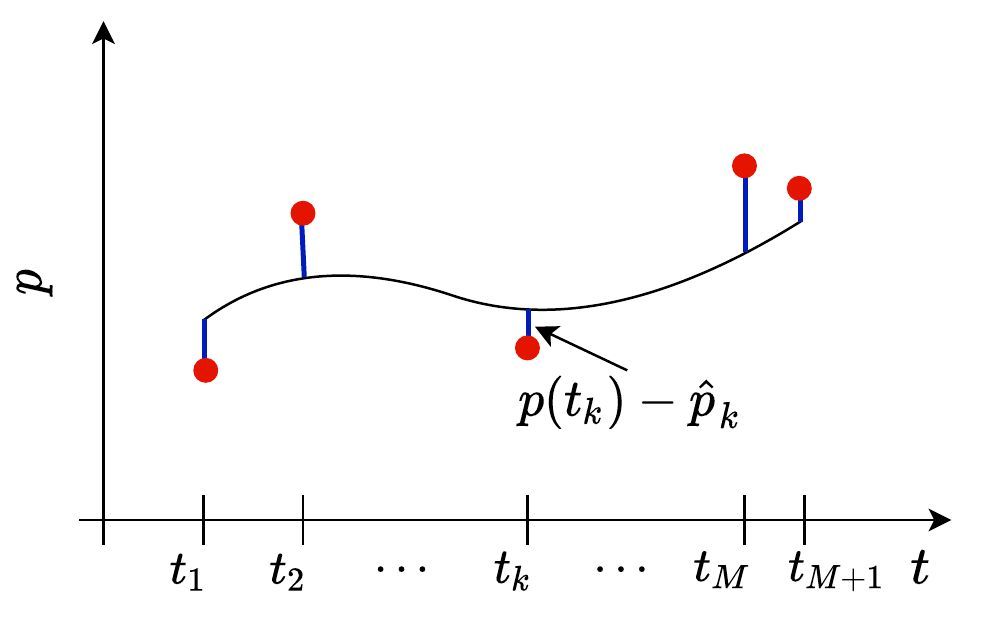}}
\caption{Illustration of data assimilation problem in Section~\ref{sec: data assimilation} seeking the optimal selection of $u(0)$ and $p(0)$ to minimize a weighted quadratic sum in the deviations $p(t_k)-\hat{p}_k$.}
\label{fig: data assimilation}
\end{figure}

From \cite{TraversoMagri2019} we obtain the data assimilation problem of finding initial values $u(0),p(0)$ for two functions $u,p:[0,T]\mapsto\mathbb{R}^n$ that minimize the cost functional (cf.\ Fig.~\ref{fig: data assimilation})
\begin{equation}
    J:=\sum_{k=1}^{M+1}w_k\left|p\big(t_k\big)-\hat{p}_k\right|^2
\end{equation}
in terms of the given sequence of observations $\hat{p}_k\in\mathbb{R}^n$,  non-negative weight vector $w\in\mathbb{R}^{M+1}$, and time sequence $0=t_1<\cdots<t_{M+1}=T$ under the differential constraints
\begin{equation}
\label{eq:diffconstraints}
    \dot{u}_j=-j\pi p_j,\,\dot{p}_j=-\zeta_jp_j+j\pi u_j+q_j
\end{equation}
for $t\in(0,\alpha)\cup(\alpha,T)$, continuity conditions
\begin{equation}
    \lim_{t\rightarrow\alpha-}u(t)=\lim_{t\rightarrow\alpha+}u(t),\,\lim_{t\rightarrow\alpha-}p(t)=\lim_{t\rightarrow\alpha+}p(t),
\end{equation}
and coupling constraints
\begin{equation}
\label{eq:coupconstraints}
    q(t)=\bigg\{\begin{array}{cc}0, & t\in(0,\alpha),\\2\beta F\big(\gamma^\text{T}u(t-\alpha)\big), & t\in(\alpha,T)\end{array}
\end{equation}
for $\beta,\gamma\in\mathbb{R}^n$ and a continuously differentiable function $F:\mathbb{R}\mapsto\mathbb{R}$. We treat this problem using standard techniques from the calculus of variations.

We anticipate discontinuities in the derivative of $p$ at $t=\alpha$ and in the Lagrange multipliers associated with the differential constraints at $t=t_i$ for $i=2,\ldots,M$. For simplicity, assume that $\alpha=t_I$ for some $2\le I\le M$. For notational convenience, let $T_k:=t_{k+1}-t_k$ and
\begin{equation}
    \hat{t}_{i,k}:=(t_i-t_k+t_I)/T_k
\end{equation}
for $k\ge I$. We then replace \eqref{eq:diffconstraints}-\eqref{eq:coupconstraints} with a sequential multi-point boundary-value problem for the functions $u^{(k)},p^{(k)}:[0,1]\rightarrow\mathbb{R}^n$ for $k=1,\ldots,M$ given by
\begin{align}
\label{eq:de}
    \dot{u}_j^{(k)}=T_k\left(-j\pi p_j^{(k)}\right),\,\dot{p}_j^{(k)}=T_k\left(-\zeta_jp_j^{(k)}+j\pi u_j^{(k)}+q_j^{(k)}\right)
\end{align}
for $t\in(0,1)$ and \begin{equation}
\label{eq:ic}
    u^{(k)}(0)=u^{(k-1)}(1),\,p^{(k)}(0)=p^{(k-1)}(1)
\end{equation}
for $k>1$, with
\begin{equation}
\label{eq:cp}
    q^{(k)}(t)=2\beta F\left(\gamma^\text{T}u^{(l)}\left(\frac{T_k}{T_l}\left(t-\hat{t}_{l,k}\right)\right)\right)
\end{equation}
for $k\ge I$,
\begin{equation}
    l\in \mathbb{I}_k:=\{1\le i\le M\,\big|\,[0,1]\cap(\hat{t}_{i,k}, \hat{t}_{i+1,k})\ne\emptyset\}
\end{equation}
and $t\in\left[\max\left(0,\hat{t}_{l,k}\right),\min\left(1,\hat{t}_{l+1,k}\right)\right]$, and $q^{(k)}(t)=0$ otherwise. In this notation,
\begin{equation}
    J=w_1\left|p^{(1)}(0)-\hat{p}_1\right|^2+\sum_{k=2}^{M+1}w_k\left|p^{(k-1)}(1)-\hat{p}_k\right|^2.
\end{equation}
We seek an optimal choice for $u^{(1)}(0)$ and $p^{(1)}(0)$ that corresponds to an extremum of $J$ along the corresponding constraint manifold.

\subsubsection{Adjoint conditions}
\label{sec: data assimilation adjoint conditions}

To locate such an extremum, consider the Lagrangian (which differs from~\cite{TraversoMagri2019} in the purposeful introduction of the auxiliary variable $\mu$)
\begin{equation}
    \mu+\eta\left(J-\mu\right)+L_\mathrm{de}+L_\mathrm{ic}+L_\mathrm{cp},
\end{equation}
where
\begin{align}
    L_\mathrm{de}&:=\sum_{k=1}^M\sum_{j=1}^n\int_0^1\kappa_j^{(k)}\left(\dot{u}_j^{(k)}+T_kj\pi p_j^{(k)}\right)\,\mathrm{d}t+\sum_{k=1}^M\sum_{j=1}^n\int_0^1\lambda_j^{(k)}\left(\dot{p}_j^{(k)}+T_k\zeta_jp_j^{(k)}-T_kj\pi u_j^{(k)}-T_kq_j^{(k)}\right)\,\mathrm{d}t
\end{align}
in terms of the Lagrange multiplier functions $\kappa^{(k)},\lambda^{(k)}:[0,1]\mapsto\mathbb{R}^n$,
\begin{align}
    L_\mathrm{ic}&:=\nu^{(1)\text{T}}\left(u^{(1)}(0)-u_0\right)+\sum_{k=2}^M\nu^{(k)\text{T}}\left(u^{(k)}(0)-u^{(k-1)}(1)\right)+\omega^{(1)\text{T}}\left(p^{(1)}(0)-p_0\right)+\sum_{k=2}^M\omega^{(k)\text{T}}\left(p^{(k)}(0)-p^{(k-1)}(1)\right)
\end{align}
in terms of the Lagrange multipliers $\nu^{(k)},\omega^{(k)}\in\mathbb{R}^n$, and
\begin{align}
    L_\mathrm{cp}&:=\sum_{k< I}\int_0^1\mu^{(k)\text{T}}q^{(k)}\,\mathrm{d}t+\sum_{k\ge I}\sum_{l\in L_k}\int_{\max\left(0,\hat{t}_{l,k}\right)}^{\min\left(1,\hat{t}_{l+1,k}\right)}\mu^{(k)\text{T}}\left(q^{(k)}-2\beta F\left(\gamma^Tu^{(l)}\left(\frac{T_k}{T_l}\left(t-\hat{t}_{l,k}\right)\right)\right)\right)\,\mathrm{d}t
\end{align}
in terms of the Lagrange multiplier function $\mu^{(k)}:[0,1]\mapsto\mathbb{R}^n$. Here, the Lagrange multiplier $\eta$ imposes the relationship between $J$ and the auxiliary variable $\mu$, while the auxiliary variables $u_0$ and $p_0$ are introduced to track $u^{(1)}(0)$ and $p^{(1)}(0)$. We assume below that $\kappa^{(k)}$ and $\lambda^{(k)}$ are continuous and piecewise differentiable, and that $\mu^{(k)}$ is continuous.

For further notational convenience, let
\begin{equation}
    \bar{\mathbb{I}}_k:=\{1\le i\le M\,\big|\,[0,1]\cap(\hat{t}_{k,i},\hat{t}_{k+1,i})\ne\emptyset\}.
\end{equation}
Independent variations of the constraint Lagrangian with respect to the components of $u^{(k)}(\cdot)$, $p^{(k)}(\cdot)$, $q^{(k)}(\cdot)$, $u_0$, $p_0$, and $\mu$ then yields the adjoint necessary conditions for an extremum given by
\begin{equation}
\label{eq:adjfirst}
    -\dot{\kappa}_j^{(k)}-T_kj\pi \lambda_j^{(k)}-2\gamma_j\mu^{(l)\text{T}}\left(\frac{T_k}{T_l}t+\hat{t}_{k,l}\right)\beta F'\left(\gamma^\text{T}u^{(k)}\right)=0
\end{equation}
for $t\in\left(\max\left(0,-\hat{t}_{k,l}T_l/T_k\right),\min\left(1,(1-\hat{t}_{k,l})T_l/T_k\right)\right)$ for some $l\in\bar{\mathbb{I}}_k$,
\begin{equation}
\label{eq:adjsecond}
    -\dot{\kappa}_j^{(k)}-T_kj\pi \lambda_j^{(k)}=0
\end{equation}
for $t\notin\left[\max\left(0,-\hat{t}_{k,l}T_l/T_k\right),\min\left(1,(1-\hat{t}_{k,l})T_l/T_k\right)\right]$ for any $l\in\bar{\mathbb{I}}_k$,
\begin{equation}
\label{eq:adjthird}
    -\dot{\lambda}_j^{(k)}+T_k\zeta_j\lambda_j^{(k)}+T_kj\pi\kappa_j^{(k)}=0,
\end{equation}
for $t\in(0,1)$,
\begin{equation}
\label{eq:adjfourth}
    \mu^{(k)}-T_k\lambda^{(k)}=0
\end{equation}
for $t\in[0,1]$,
\begin{gather}
\label{eq:adjfifth}
    \kappa^{(k-1)}(1)-\nu^{(k)}=0,\,\lambda^{(k-1)}(1)-\omega^{(k)}+2\eta w_k\left(p^{(k-1)}(1)-\hat{p}_k\right)=0,\\
    \label{eq:adjsixth}
    -\kappa^{(k)}(0)+\nu^{(k)}=0,\,-\lambda^{(k)}(0)+\omega^{(k)}=0
\end{gather}
for $k=2,\ldots,M$,
\begin{gather}
\label{eq:adjseventh}
    \kappa^{(M)}(1)=0,\,\lambda^{(M)}(1)+2\eta w_{M+1}\left(p^{(M)}(1)-\hat{p}_{M+1}\right)=0,\\
    -\kappa^{(1)}(0)+\nu^{(1)}=0,\,-\lambda^{(1)}(0)+\omega^{(1)}+2\eta w_1\left(p^{(1)}(0)-\hat{p}_1\right)=0,\label{eq:adjlast}
\end{gather} $\nu^{(1)}=\omega^{(1)}=0$, and $1-\eta=0$. As was the case in a previous section, these conditions are linear in the Lagrange multipliers and, apart from the final condition on $\eta$, homogeneous.

\subsubsection{Problem construction}
\label{sec: data assimilation problem construction}
We obtain an augmented continuation problem $\mathbf{A}$ of the form in \eqref{eq:augmented continuation problem}
by associating 
\begin{itemize}
    \item $\Phi$ with the multi-point boundary-value problem in \eqref{eq:de}-\eqref{eq:cp} in terms of the continuation variables $u^{(k)}$, $p^{(k)}$, and $q^{(k)}$;
    \item $\Psi$ with the vector $\begin{pmatrix}J,u^{(1)}(0),p^{(1)}(0)\end{pmatrix}$ and corresponding continuation parameters $\mu$, $u_0$, and $p_0$; and 
    \item $\Lambda^\ast$ with the linear operator in \eqref{eq:adjfirst}-\eqref{eq:adjlast} acting on the continuation multipliers $\kappa^{(k)}$, $\lambda^{(k)}$, $\mu^{(k)}$, $\nu^{(k)}$, $\omega^{(k)}$, and $\eta$. 
\end{itemize}
This problem has dimensional deficit $2n+1$ which reduces to $0$ once a solution is found with $\nu^{(1)}=\omega^{(1)}=0$ and $\eta=1$.

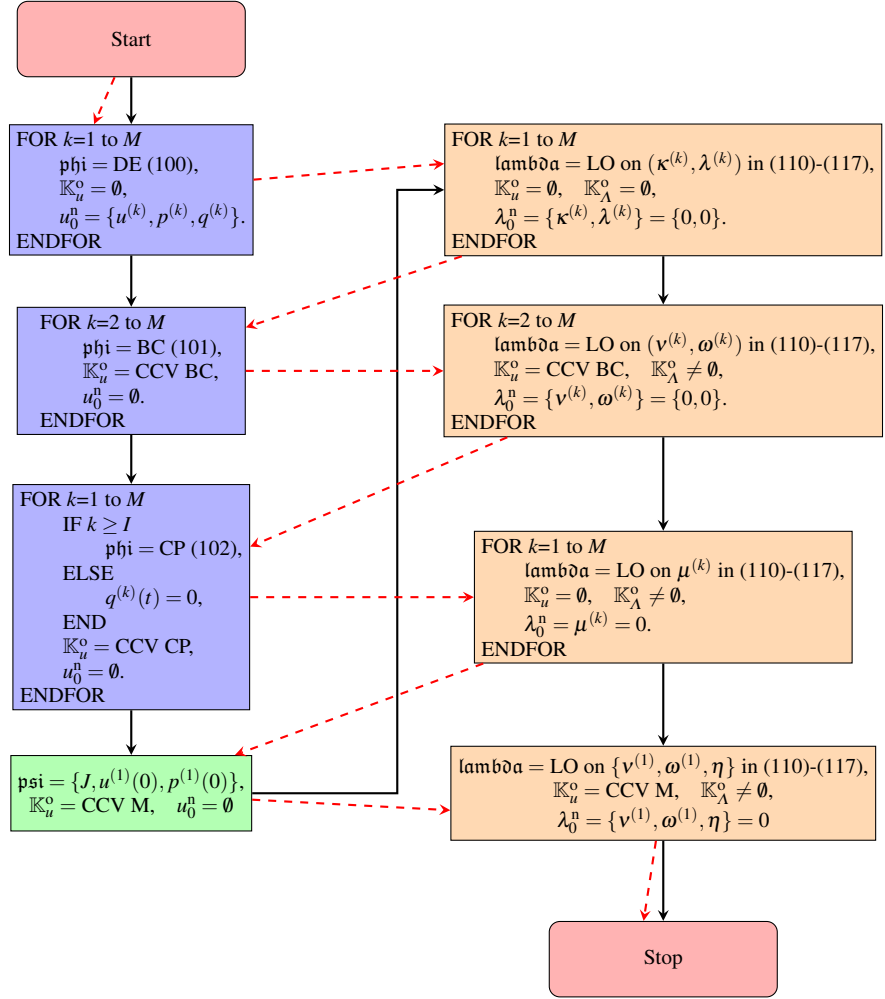
\begin{figure}[ht]
\centering
\begin{tikzpicture}
[node distance=2.0cm]
\node (start) [startstop] {Start};
\node (diff-eqs) [process-phi,below of=start,align=left] {FOR $k$=1 to $M$\\
\quad\quad $\mathfrak{phi}=\mathrm{DE}$~\eqref{eq:de},\\
\quad\quad $\mathbb{K}^\mathrm{o}_u=\emptyset$,\\
\quad\quad$u^\mathrm{n}_0=\{u^{(k)}, p^{(k)}, q^{(k)}\}$.\\
ENDFOR};
\node (bc) [process-phi, below of=diff-eqs,yshift=-0.4cm,align=left] {FOR $k$=2 to $M$\\
\quad\quad $\mathfrak{phi}=\mathrm{BC}$~\eqref{eq:ic},\\
\quad\quad $\mathbb{K}^\mathrm{o}_u=$ CCV BC,\\
\quad\quad $u^\mathrm{n}_0=\emptyset$.\\
ENDFOR};
\node (cp) [process-phi, below of=bc,yshift=-1.0cm,align=left] {FOR $k$=1 to $M$\\
\quad\quad IF $k\geq I$\\
\qquad\qquad$\mathfrak{phi}=\mathrm{CP}$~\eqref{eq:cp},\\
\qquad ELSE\\
\qquad\qquad $q^{(k)}(t)=0$,\\
\qquad END\\
\quad\quad $\mathbb{K}^\mathrm{o}_u=$ CCV CP,\\
\quad\quad$u^\mathrm{n}_0=\emptyset$.\\
ENDFOR};
\node (monitors) [process-psi,yshift=-0.6cm,below of=cp,align=center]{$\mathfrak{psi}=\{J, u^{(1)}(0), p^{(1)}(0)\}$,\\
$\mathbb{K}^\mathrm{o}_u=$ CCV M,\quad $u^\mathrm{n}_0=\emptyset$};
% adjoint ones
\node (diff-eqs-adjt) [process-lambda,right of=diff-eqs,xshift=5.0cm,align=left] {FOR $k$=1 to $M$\\
\quad\quad $\mathfrak{lambda}=$ LO on $(\kappa^{(k)},\lambda^{(k)})$ in~\eqref{eq:adjfirst}-\eqref{eq:adjlast},\\
\quad\quad $\mathbb{K}^\mathrm{o}_u=\emptyset$,\quad $\mathbb{K}^\mathrm{o}_\Lambda=\emptyset$,\\
\quad\quad$\lambda^\mathrm{n}_0=\{\kappa^{(k)}, \lambda^{(k)}\}=\{0,0\}$.\\
ENDFOR};
\node (bc-adjt) [process-lambda, right of=bc,xshift=5cm,align=left] {FOR $k$=2 to $M$\\
\quad\quad $\mathfrak{lambda}=$ LO on $(\nu^{(k)},\omega^{(k)})$ in~\eqref{eq:adjfirst}-\eqref{eq:adjlast},\\
\quad\quad $\mathbb{K}^\mathrm{o}_u=$ CCV BC,\quad $\mathbb{K}^\mathrm{o}_\Lambda\neq\emptyset$,\\
\quad\quad$\lambda^\mathrm{n}_0=\{\nu^{(k)}, \omega^{(k)}\}=\{0,0\}$.\\
ENDFOR};
\node (cp-adjt) [process-lambda, right of=cp,xshift=5cm,align=left] {FOR $k$=1 to $M$\\
\quad\quad $\mathfrak{lambda}=$ LO on $\mu^{(k)}$ in~\eqref{eq:adjfirst}-\eqref{eq:adjlast},\\
\quad\quad $\mathbb{K}^\mathrm{o}_u=\emptyset$,\quad $\mathbb{K}^\mathrm{o}_\Lambda\neq\emptyset$,\\
\quad\quad$\lambda^\mathrm{n}_0=\mu^{(k)}=0$.\\
ENDFOR};
\node (monitors-adjt) [process-lambda,xshift=5.0cm,right of=monitors,align=center]{$\mathfrak{lambda}=$ LO on $\{\nu^{(1)}, \omega^{(1)}, \eta\}$ in~\eqref{eq:adjfirst}-\eqref{eq:adjlast},\\
$\mathbb{K}^\mathrm{o}_u=$ CCV M,\quad$\mathbb{K}^\mathrm{o}_\Lambda\neq\emptyset$,\\ $\lambda^\mathrm{n}_0=\{\nu^{(1)}, \omega^{(1)}, \eta\}=0$};
\node (stop) [startstop, below of=monitors-adjt,yshift=-0.2cm] {Stop};
% DRAW LINES
\draw [arrow] (start) -- (diff-eqs);
\draw [arrow] (diff-eqs) -- (bc);
\draw [arrow] (bc) -- (cp);
\draw [arrow] (cp) -- (monitors);
\draw [->,thick,>=stealth] (monitors) to[-|-] (diff-eqs-adjt);
\draw [arrow] (diff-eqs-adjt) -- (bc-adjt);
\draw [arrow] (bc-adjt) -- (cp-adjt);
\draw [arrow] (cp-adjt) -- (monitors-adjt);
\draw [arrow] (monitors-adjt) -- (stop);
% AN ALTERNATIVE
\draw [arrow,red,dashed] ([xshift=-4em]start) -- (diff-eqs);
\draw [arrow,red,dashed] ([yshift=3em]diff-eqs) -- (diff-eqs-adjt);
\draw [arrow,red,dashed] ([yshift=1em]diff-eqs-adjt) -- (bc);
\draw [arrow,red,dashed] (bc) -- (bc-adjt);
\draw [arrow,red,dashed] (bc-adjt) -- (cp);
\draw [arrow,red,dashed] (cp) -- (cp-adjt);
\draw [arrow,red,dashed] (cp-adjt) -- (monitors);
\draw [arrow,red,dashed] ([yshift=-2em]monitors) -- (monitors-adjt);
\draw [arrow,red,dashed] ([xshift=-4em]monitors-adjt) -- (stop);
\end{tikzpicture}
\caption{A flowchart depicting the construction of the augmented continuation problem $\mathbf{A}$ associated with the data assimilation problem in Section~\ref{sec: data assimilation}. Here rectangles filled with blue, green and orange colors represent core constructors associated with functions of the type $\mathfrak{phi}$, $\mathfrak{psi}$ and $\mathfrak{lambda}$, respectively. The workflow with black solid arrows illustrates the algorithm detailed in Section~\ref{sec: data assimilation problem construction}. In this workflow, adjoint contributions are constructed after the construction of \emph{all} zero and monitor functions. In contrast, in the alternate workflow represented by red dashed arrows, one constructs the adjoint contributions after the introduction of \emph{each} of the corresponding zero or monitor functions. Here, the abbreviations DE, BC, CP, CCV, M, and LO represent differential equations, boundary conditions, coupling conditions, corresponding continuation variables, monitor functions, and linear operators, respectively. In particular, $\mathbb{K}^\mathrm{o}_u=$ CCV BC/CP/M denote indexing the corresponding continuation variables for boundary conditions/coupling conditions/monitor functions from the ones defined when constructing the differential constraints.}
\label{fig:flowchart}
\end{figure}

After suitable discretization, we may construct $\mathbf{A}$  according to the following algorithm: 
\begin{itemize}
    \item[\textbf{Step 1:}] As $k$ increments from $1$ to $M$, repeatedly invoke the core constructor \eqref{eq: coco_add_func1} with $\mathfrak{phi}$ encoding the differential constraints \eqref{eq:de}, $\mathbb{K}^\mathrm{o}_u=\emptyset$, and $u^\mathrm{n}_0$ given by an initial solution guess for the continuation variables $u^{(k)}$, $p^{(k)}$, and $q^{(k)}$.
    \item[\textbf{Step 2:}] As $k$ increments from $2$ to $M$, repeatedly invoke the core constructor \eqref{eq: coco_add_func1} with $\mathfrak{phi}$ encoding the boundary conditions \eqref{eq:ic}, $\mathbb{K}^\mathrm{o}_u$ indexing the corresponding continuation variables from \textbf{Step 1}, and $u^\mathrm{n}_0=\emptyset$.
    \item[\textbf{Step 3:}] As $k$ increments from $1$ to $M$, repeatedly invoke the core constructor \eqref{eq: coco_add_func1} with $\mathfrak{phi}$ encoding the coupling \eqref{eq:cp} for $k\ge I$ and the condition $q^{(k)}(t)=0$ for $k<I$, $\mathbb{K}^\mathrm{o}_u$ indexing the corresponding continuation variables from \textbf{Step 1}, and $u^\mathrm{n}_0=\emptyset$.
    \item[\textbf{Step 4:}] Invoke the core constructor \eqref{eq: coco_add_func2} with $\mathfrak{psi}$ encoding the evaluation of $J$, $u^{(1)}(0)$, and $p^{(1)}(0)$, $\mathbb{K}^\mathrm{o}_u$ indexing the corresponding continuation variables from \textbf{Step 1}, and $u^\mathrm{n}_0=\emptyset$.
    \item[\textbf{Step 5:}] As $k$ increments from $1$ to $M$, repeatedly invoke the core constructor \eqref{eq: coco_add_adjt} with $\mathfrak{lambda}$ encoding the linear operators acting on $\kappa^{(k)}$ and $\lambda^{(k)}$ in the adjoint conditions \eqref{eq:adjfirst}-\eqref{eq:adjlast}, $\mathbb{K}^\mathrm{o}_u$ indexing the continuation variables introduced in the corresponding call in \textbf{Step 1}, $\mathbb{K}^\mathrm{o}_\Lambda=\emptyset$, and $\lambda^\mathrm{n}_0$ given by an initial solution guess for the continuation variables $\kappa^{(k)}$ and $\lambda^{(k)}$.
    \item[\textbf{Step 6:}] As $k$ increments from $2$ to $M$, repeatedly invoke the core constructor \eqref{eq: coco_add_adjt} with $\mathfrak{lambda}$ encoding the linear operators acting on $\nu^{(k)}$ and $\omega^{(k)}$ in the adjoint conditions \eqref{eq:adjfirst}-\eqref{eq:adjlast}, $\mathbb{K}^\mathrm{o}_u$ indexing the continuation variables associated with the corresponding call in \textbf{Step 2}, $\mathbb{K}^\mathrm{o}_\Lambda\ne\emptyset$, and $\lambda^\mathrm{n}_0$ given by an initial solution guess for the continuation variables $\nu^{(k)}$ and $\omega^{(k)}$.
    \item[\textbf{Step 7:}] As $k$ increments from $1$ to $M$, repeatedly invoke the core constructor \eqref{eq: coco_add_adjt} with $\mathfrak{lambda}$ encoding the linear operators acting on $\mu^{(k)}$ in the adjoint conditions \eqref{eq:adjfirst}-\eqref{eq:adjlast}, $\mathbb{K}^\mathrm{o}_u$ indexing the continuation variables associated with the corresponding call in \textbf{Step 3}, $\mathbb{K}^\mathrm{o}_\Lambda\ne\emptyset$, and $\lambda^\mathrm{n}_0$ given by an initial solution guess for the continuation variables $\mu^{(k)}$.
    \item[\textbf{Step 8:}] Invoke the core constructor \eqref{eq: coco_add_adjt} with $\mathfrak{lambda}$ encoding the linear operators acting on $\nu^{(1)}$, $\omega^{(1)}$, and $\eta$ in the adjoint conditions \eqref{eq:adjfirst}-\eqref{eq:adjlast}, $\mathbb{K}^\mathrm{o}_u$ indexing the continuation variables associated with the corresponding call in \textbf{Step 4}, $\mathbb{K}^\mathrm{o}_\Lambda\ne\emptyset$, and $\lambda^\mathrm{n}_0$ given by an initial solution guess for the continuation variables $\nu^{(1)}$, $\omega^{(1)}$, and $\eta$.
\end{itemize}
One advantage of this algorithm is that steps 5 through 8 can be implemented automatically~\cite{li2017staged,li2020optimization,ahsan2020optimization} from information provided in steps 1 through 4, rather than simply using the core constructor \eqref{eq: coco_add_func1} to implement a general continuation problem $\mathbf{P}$. A flowchart representation of this algorithm is presented in Fig.~\ref{fig:flowchart}. This figure also shows a resequenced algorithm for constructing the augmented continuation problem that interlaces construction of adjoint contributions immediately following the construction of the corresponding zero and monitor functions.

\subsubsection{Problem analysis}
\label{sec: data assimilation problem analysis}

Using a method of successive continuation (originally described in~\cite{kernevez1987optimization} with further developments in~\cite{li2017staged,li2020optimization,ahsan2020optimization}), we may reach the desired local extremum through a sequence of intermediate points at the intersection of the solution manifolds to different restricted continuation problem. To this end, invoke the core constructor \eqref{eq: coco_add_comp2} to append complementary monitor functions evaluating to $\eta$, $\nu^{(1)}$, and $\omega^{(1)}$ with corresponding complementary continuation parameters $\varphi_\eta$, $\varphi_\nu$, and $\varphi_\omega$. Here, $\mathbb{K}^\mathrm{o}_u=\mathbb{K}^\mathrm{o}_v=v^\mathrm{n}_0=\emptyset$ and $\mathbb{K}^\mathrm{o}_\lambda$ indexes the corresponding continuation multipliers. We construct the desired sequence of restricted continuation problems by fixing fewer than $2n+1$ (complementary) continuation parameters. 

For example, we obtain an augmented continuation problem with dimensional deficit equal to $1$ by fixing $p_0$, all but the first component of $u_0$, and the first component of $\varphi_\nu$. A local extremum in $\mu$ along a family of solutions to this problem with all vanishing Lagrange multipliers (such a family exists by homogeneity) then coincides with an intersection with a secondary family of solutions along which only the Lagrange multipliers vary. One point along this secondary family has $\eta=1$. The continuation problem obtained next by fixing $\varphi_\eta$ at $1$ and allowing, say, the second component of $u_0$ to vary is satisfied along a tertiary manifold through this point. If we locate a point on this manifold where the second component of $\varphi_\nu$ equals $0$, we may use this point to switch to a different restricted continuation problem with the first three components of $u_0$ allowed to vary and the first two components of $\varphi_\nu$ fixed. Along the corresponding solution manifold we look for a point where the third component of $\varphi_\nu$ equals $0$, and continue in the same fashion until a local extremum is reached. 

Alternatively, once the initial point with $\eta=1$ is reached, denote the corresponding value of $\nu^{(1)}$ by $\nu^{(1)\ast}$. We may now invoke the core constructor \eqref{eq: coco_add_comp1} to append complementary zero functions that evaluate to all but the first component of the combination
\begin{equation}
    \nu^{(1)}-(1-\chi)\nu^{(1)\ast}
\end{equation}
in terms of the complementary continuation variable $\chi\in\mathbb{R}$. Here, $\mathbb{K}^\mathrm{o}_u=\mathbb{K}^\mathrm{o}_v=\emptyset$, $\mathbb{K}^\mathrm{o}_\lambda$ indexes the continuation multipliers $\nu^{(1)}$, and $v^\mathrm{n}_0$ contains an initial solution guess for $\chi$. By again fixing $\varphi_\eta$ at $1$ and allowing all remaining components of $u_0$ to vary, we obtain a continuation problem with dimensional deficit equal to $1$ and may search along its solution manifold for a point with $\nu^{(1)}=0$. We drive $\omega^{(1)}$ to $0$ following similar principles.

\subsection{Phase response curves of periodic orbits}
\label{sec: phase response curves}
\subsubsection{Linear response theory for closed regular problems}
\label{sec: phase response curves linear response theory}
Instead of optimization, as in Section~\ref{sec: data assimilation}, we consider in this section the simpler case where the zero problem $\Phi(u)=0$ for $\Phi:\mathcal{U}_\Phi\to \mathcal{R}_\Phi$ has dimensional deficit equal to $0$ and is regular at some solution $\tilde u$ (i.e., such that the Frech\'{e}t derivative $D\Phi(\tilde u)$ is regular). We choose the case of a scalar-valued monitor function $\Psi:\mathcal{U}_\Phi\mapsto\mathbb{R}$ such that the continuation parameter $\mu$ given by $\Psi(u)-\mu=0$ is also scalar ($\Psi(u)$ is called the  \emph{observable}). Let $\mathcal{R}_\Phi^*$ denote the (dual) space of linear functionals on $\mathcal{R}_\Phi$. At extremal points $(\tilde u,\tilde\mu,\tilde \lambda,\tilde\eta)\in\mathcal{U}_\Phi\times\mathbb{R}\times\mathcal{R}_\Phi^*\times\mathbb{R}$ of the Lagrangian
\begin{align}\label{prc:lagrangian:general}
    L(u,\mu,\lambda,\eta)&=\mu+\eta(\Psi(u)-\mu)+\lambda \Phi(u)
\end{align}
the Lagrange multiplier $\tilde \lambda$ measures the linear sensitivity of $\Psi$ to changes in $\Phi$. Indeed, by considering vanishing variations of $L$, it follows that $(\tilde u,\tilde{\mu},\tilde \lambda,\tilde{\eta})$ must satisfy
\begin{equation}
    \Phi(\tilde u)=0,\,\Psi(\tilde{u})-\tilde{\mu}=0,\,\eta D\Psi(\tilde u)+\tilde{\lambda}D\Phi(\tilde u)=0
\end{equation}
and $1-\eta=0$, from which we obtain $\tilde\lambda=-D\Psi(\tilde u)\left(D\Phi(\tilde u)\right)^{-1}$. For all small perturbations $\delta\Phi\in\mathcal{R}_\Phi$ the perturbed zero problem $\Phi(u)=\delta\Phi$ has a locally unique solution $u=\tilde u+\delta u$, where $\delta\Phi=D\Phi(\tilde u)\delta u+O(\|\delta u\|^2)$. It follows that
\begin{align}\label{eq:prc:general:linear-response}
    \delta\Psi&=\Psi(u)-\Psi(\tilde u)=D\Psi(\tilde u)\delta u+O(\|\delta u\|^2)=-\tilde \lambda \delta\Phi+O(\|\delta\Phi\|^2)\mbox{}
\end{align}

In this section, we to apply this general observation to the derivation of phase response curves associated with limit cycles in ordinary and delay differential equations.

\subsubsection{Phase response curves as linear response}
\label{sec: phase response curves as linear response}
The construction in Section~\ref{sec: phase response curves linear response theory} can be applied to abstract autonomous periodic boundary-value problems and the observable $T$ (the unknown period) to obtain so-called phase response curves~\cite{ermentrout1996type,govaerts2006computation,izhikevich2007dynamical,langfield2020continuation}.

In the notation of this section, let
$\mathcal{U}_\Phi=C^1([0,1];\mathcal{U})\times\mathbb{R}$ be the space of continuation variables, where $\mathcal{U}$ is some Banach space, and let the zero problem take the form
\begin{align}
    \Phi(u)=(\Phi_\mathrm{DE},\Phi_\mathrm{BC},\Phi_\mathrm{PS})((x(\cdot),T))=(\dot x(\cdot)-T f(x(\cdot)),x(0)-x(1),h(x(0)))
\end{align}
corresponding to a periodic orbit $\mathbb{R}\ni t\mapsto x(t/T)\in\mathcal{U}$ of period $T$ of an autonomous vector field $f$ and phase determined by the Poincar{\'e} condition $h(x(0))=0$. We define the monitor function $\Psi$ as the projection onto the scalar component $T$ of $(x(\cdot),T)\in\mathcal{U}_\Phi$ such that $\Psi(x(\cdot),T)=T$. The Lagrangian for the linear response of $T$, given in \eqref{prc:lagrangian:general}, is then
\begin{align}
    &L(x(\cdot),T,\mu,\lambda_\mathrm{DE}(\cdot),\lambda_\mathrm{BC},\lambda_\mathrm{PS},\eta)=\mu+\eta(T-\mu)+
    \int_0^1\lambda_\mathrm{DE}(\tau)\left(\dot{x}(\tau)-Tf(x(\tau)\right)\d \tau+
    \lambda_\mathrm{BC}\left(x(0)-x(1)\right)+
    \lambda_\mathrm{PS}h(x(0))
\end{align}
defined on $\mathcal{U}_\Phi\times\mathbb{R}\times\left(C^0\left([0,1];\mathcal{U}\right)\right)^*\times\mathcal{U}^\ast\times\mathbb{R}\times\mathbb{R}$. In this case, vanishing variations of $L$ with respect to the Lagrange multipliers $\lambda_\mathrm{DE}(\cdot)$, $\lambda_\mathrm{BC}$, $\lambda_\mathrm{PS}$, and $\eta$ at an extremal point $(\tilde{x}(\cdot),\tilde{T},\tilde{\mu},\tilde{\lambda}_\mathrm{DE}(\cdot),\tilde{\lambda}_\mathrm{BC},\tilde{\lambda}_\mathrm{PS},\tilde{\eta})$ imply that
\begin{equation}
    \dot{\tilde{x}}(\tau)-\tilde{T}f(\tilde{x}(\tau))=0\text{ for }\tau\in(0,1),\,\tilde{x}(0)-\tilde{x}(1)=0,\,h(\tilde{x}(0))=0,\label{eq:odeprc}
\end{equation}
and $\tilde{T}-\tilde{\mu}=0$, i.e., that $t\mapsto\tilde{x}(t/\tilde{T})$ is a periodic solution with period $\tilde{\mu}=\tilde{T}$ of a dynamical system with autonomous vector field $f$ and with initial condition on the zero-level surface of the function $h$. 

Vanishing variations of $L$ with respect to $x(\cdot)$, $T$, and $\mu$ yield the necessary adjoint conditions
\begin{align}
    -\dot{\tilde{\lambda}}_\mathrm{DE}(\tau)-\tilde{\lambda}_\mathrm{DE}(\tau)\tilde{T}Df(\tilde{x}(\tau))&=0\text{ for }\tau\in(0,1),\label{eq:adjdiffodeprc}\\
    \tilde{\lambda}_\mathrm{BC}-\tilde{\lambda}_\mathrm{DE}(0)+\tilde{\lambda}_\mathrm{PS}Dh(\tilde{x}(0))&=0,\label{eq:adjbc0odeprc}\\
    -\tilde{\lambda}_\mathrm{BC}+\tilde{\lambda}_\mathrm{DE}(1)&=0,\label{eq:adjbc1odeprc}\\
    \eta-\int_0^1\tilde{\lambda}_\mathrm{DE}(\tau)f(\tilde{x}(\tau))\,\d \tau&=0,\label{eq:adjTodeprc}
\end{align}
and $1-\eta=0$. From \eqref{eq:odeprc} and \eqref{eq:adjdiffodeprc}, it follows that $\tilde{\lambda}_\mathrm{DE}(\tau)f(\tilde{x}(\tau))$ is constant, i.e., that
\begin{equation}
    \tilde{\lambda}_\mathrm{DE}(\tau)f(\tilde{x}(\tau))=1\text{ for all }\tau\in[0,1],
\end{equation}
where we used \eqref{eq:adjTodeprc} and the fact that $\eta=1$. Moreover, from \eqref{eq:adjbc0odeprc} and \eqref{eq:adjbc1odeprc}, we see that $\tilde{\lambda}_\mathrm{BC}=\tilde{\lambda}_\mathrm{DE}(1)$ and $\tilde{\lambda}_\mathrm{PS}Dh(\tilde{x}(0))=\tilde{\lambda}_\mathrm{DE}(0)-\tilde{\lambda}_\mathrm{DE}(1)$. By the periodicity of $\tilde{x}$, it follows that
\begin{equation}
    \tilde{\lambda}_\mathrm{PS}Dh(\tilde{x}(0))\tilde{T}f(\tilde{x}(0))=\tilde{\lambda}_\mathrm{DE}(0)\tilde{T}f(\tilde{x}(0))-\tilde{\lambda}_\mathrm{DE}(1)\tilde{T}f(\tilde{x}(1))=0,
\end{equation}
i.e., that $\tilde{\lambda}_{PS}=0$ provided that the periodic trajectory $\tilde{x}(\cdot)$ intersects $h=0$ transversally at $\tilde{x}(0)$. In this case, $t\mapsto\tilde \lambda_\mathrm{DE}(t/\tilde{T})$ is also periodic with period $\tilde{T}$.

The invertibility of the linearization of \eqref{eq:odeprc} is equivalent to a simple Floquet multiplier at $1$ for the corresponding periodic orbit. This invertibility implies the existence of a unique pair $(x(\cdot),T)\in\mathcal{U}_\Phi$ near $(\tilde{x}(\cdot),\tilde{T})$ for each pair of small perturbations $(\delta_\mathrm{BC},\delta_\mathrm{PS})$ such that
\begin{equation}
    \dot{x}(\tau)-Tf(x(\tau))=0\text{ for }\tau\in[0,1],\, x(0)-x(1)=\delta_\mathrm{BC},\,
    h(x(0))=\delta_\mathrm{PS}\mbox{.}
\end{equation}
The linear response formula \eqref{eq:prc:general:linear-response} from Section~\ref{sec: phase response curves linear response theory} then implies
\begin{equation}
    T-\tilde{T}=-\tilde{\lambda}_\mathrm{DE}(0)\delta_\mathrm{BC}+O\left(\|(\delta_\mathrm{BC},\delta_\mathrm{PS})\|^2\right),
\end{equation}
where we have used the fact that $\tilde{\lambda}_\mathrm{BC}=\tilde{\lambda}_\mathrm{DE}(1)=\tilde{\lambda}_\mathrm{DE}(0)$ and $\tilde{\lambda}_\mathrm{PS}=0$. In particular, since
\begin{equation}
    \delta_{BC}=x(0)-\tilde{x}(0)+\tilde{x}(1)-x(1),
\end{equation}
it follows that to first order in $\|x(0)-\tilde{x}(0)\|$ and $|T-\tilde{T}|$,
\begin{equation}
    \tilde{\lambda}_\mathrm{DE}(0)\delta_\mathrm{BC}=\tilde{T}-T.
\end{equation}

\begin{figure}[ht]
\centering
\subfloat{\includegraphics[width=0.6\columnwidth,trim={0 {0.0\textwidth} 0 0.0\textwidth},clip]{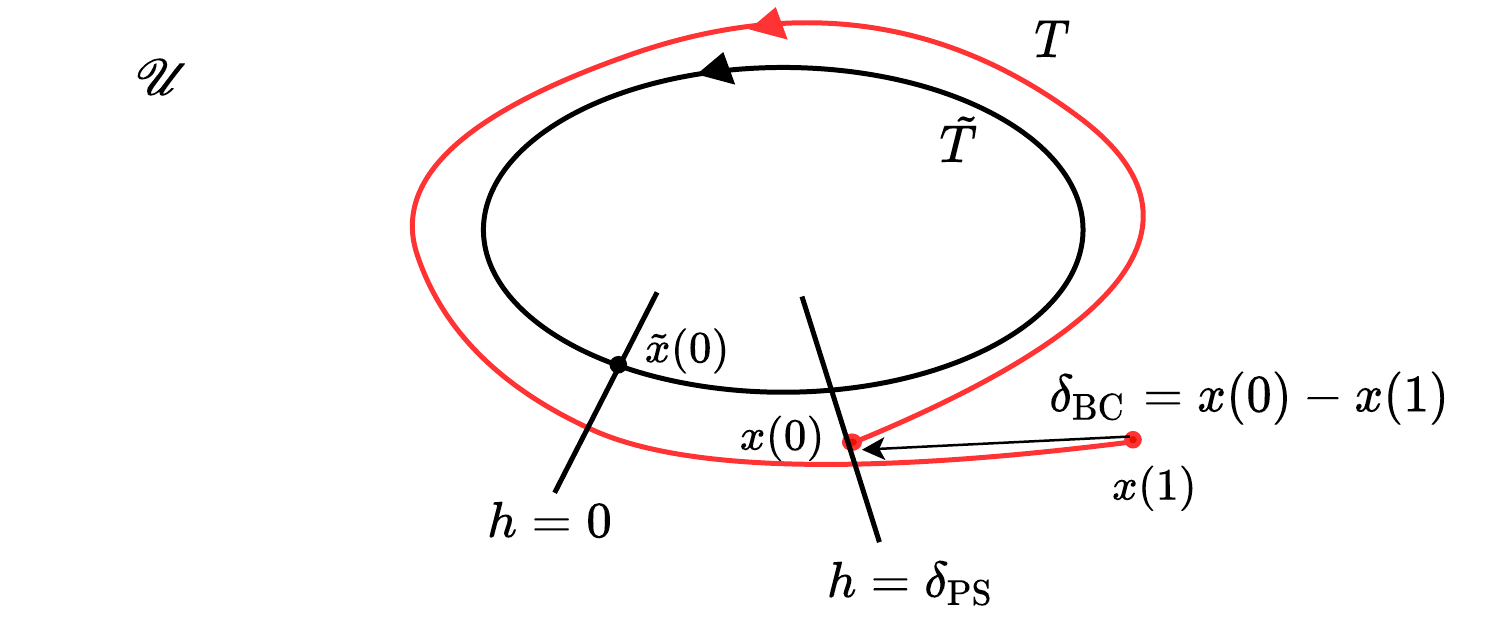}}
\caption{The periodic orbit $\mathbb{R}\ni t\mapsto\tilde{x}(t/\tilde{T})\in\mathcal{U}$ satisfies the Poincar\'{e} condition $h(x(0))=0$ and periodicity condition $x(0)-x(1)=0$. A violation by $\delta_{\text{PS}}$ of the Poincar\'{e} condition and by $\delta_{\text{BC}}$ of the periodicity condition, as shown in the figure, results in a change in the duration $T$ relative to the period $\tilde{T}$ of the periodic orbit by $-\tilde{\lambda}_{\text{DE}}(0)\delta_{\text{BC}}$ to first order in $\|\delta_{\text{BC}}\|$ and $\|\delta_{\text{PS}}\|$. For a stable limit cycle, $\tilde{\lambda}_{\text{DE}}(0)$ equals the Fr\'{e}chet derivative $D\varphi$ of the asymptotic phase evaluated at $\tilde{x}(0)$.}
\label{fig: prc}
\end{figure}

Consider the special case that the periodic function $\tilde{x}(t/\tilde{T})$ describes a linearly asymptotically stable limit cycle of the vector field $f(x)$. Then, there exists a unique map $\varphi$, called the \emph{asymptotic phase}~\cite{chicone2004asymptotic}, defined on the basin of attraction $\mathcal{B}$ of the limit cycle, such that $\varphi:\mathcal{B}\rightarrow[0,\tilde{T})$, $\varphi(\tilde{x}(0))=0$, and 
\begin{equation}
    \lim_{t\rightarrow\infty} x(t/T)-\tilde{x}\left((\varphi(x(0))+t)/\tilde{T}\right)=0
\end{equation}
for every solution $x(t)$ to $\dot{x}=Tf(x)$ with $x(0)$ in $\mathcal{B}$. The substitution $t\mapsto s+t$ into this limit identity shows that $\phi$ satisfies the
\begin{equation}
    \varphi(x(s/T))=\varphi(x(0))+s\text{ modulo } \tilde{T}
\end{equation}
for all $s$, which can also be used as a defining equation for $\phi$. In particular, for $s=T$, we obtain
\begin{equation}
    \varphi(x(1))-\varphi(x(0))=T\text{ modulo } \tilde{T}
\end{equation}
For $x(0)\approx\tilde{x}(0)$ and $T\approx\tilde{T}$, it follows that to first order in $\|x(0)-\tilde{x}(0)\|$ and $|T-\tilde{T}|$,
\begin{equation}
    D\varphi(\tilde{x}(0))\delta_\mathrm{BC}=\tilde{T}-T,
\end{equation}
i.e., that the Frech\'{e}t derivative $D\varphi(\tilde{x}(0))=\tilde{\lambda}_\mathrm{DE}(0)$. By considering an arbitrary $t\in[0,\tilde{T})$, we obtain the \emph{phase response functional} $D\varphi(\tilde{x}(t/\tilde{T}))=\tilde{\lambda}_\mathrm{DE}(t/\tilde{T})$ for the period-$\tilde{T}$ orbit $\tilde{x}(t/\tilde{T})$ of the vector field $f(x)$.

This phase response functional (or vector in the case of finite-dimensional $\mathcal{U}$) is reduced to the periodic \emph{phase response curve} $\prc(\tau)$ on $[0,1]$ for a particular perturbation $\delta_\mathrm{BC}\in{\cal U}$ by applying the functional to the perturbation at every time $\tau$:
\begin{align}
    \prc_{\delta_\mathrm{BC}}:[0,1]\ni\tau\mapsto \tilde\lambda_\mathrm{DE}(\tau)\delta_\mathrm{BC}\in\R.
\end{align}
This measures the first-order shift in the asymptotic phase due to a perturbation of the state $x(\tau)$ at time $\tau$ by $\delta_\mathrm{BC}$ (ignoring terms of order $\|(\delta_\mathrm{BC},\delta_\mathrm{PS})\|^2$).

\subsubsection{Delay differential equations}
\label{sec: phase response curves delay-differential equations}

The results in the previous section apply also to periodic orbits in delay differential equations (DDEs), including the special case of a single discrete delay $\alpha$, given by
\begin{equation}
    \dot{z}(t)=f(z(t),z(t-\alpha)).
    \label{eq:ddedisc}
\end{equation}
The defining conditions for $\tilde \lambda_{DE}(\cdot)$, $\tilde{\lambda}_{BC}$, and $\tilde{\lambda}_{PS}$ may be obtained from the general conditions \eqref{eq:adjdiffodeprc}--\eqref{eq:adjTodeprc} by writing the delay differential equation in the form
\begin{equation}
    \dot{z}(t)=f(z(t),\zeta(-\alpha,t)),\,\zeta_{,t}(s,t)=\zeta_{,s}(s,t),\,\zeta(0,t)=z(t)
\end{equation}
in which the delayed term is obtained from the solution of an advective boundary-value problem~\cite{novivcenko2012phase}. Suitable choices of the Banach space $\mathcal{U}$ and of the action of the Lagrange multiplier $\lambda_\mathrm{DE}$ yields the corresponding adjoint boundary-value problem after solving the advective problem and its adjoint along characteristics. In this section, we go a different route. We apply the linear response approximation \eqref{eq:prc:general:linear-response} from Section~\ref{sec: phase response curves as linear response} directly to the following form of \eqref{eq:ddedisc}:
\begin{align}\label{prc:dde:coco}
    \dot{x}(\tau)&=Tf(x(\tau),y(\tau)),\,t\in(0,1),\\
    y(\tau)&=\begin{cases}
    x(\tau+1-\alpha/T),& \tau\in[0,\alpha/T],\\
    x(\tau-\alpha/T),& \tau\in(\alpha/T,1].
    \end{cases}
\end{align}
where $x(\tau)=z(T\tau)$ and $y(\tau)=z(T\tau-\alpha)$. For this coupled system, we consider vanishing variations of the Lagrangian
\begin{align}
    &L(x(\cdot),y(\cdot),T,\mu,\lambda_\mathrm{DE}(\cdot),\lambda_\mathrm{CP}(\cdot),\lambda_\mathrm{BC},\lambda_\mathrm{PS},\eta) = \mu+\eta(T-\mu)\nonumber\\
    &\quad+\int_0^1 \lambda_\mathrm{DE}^\text{T}(t)\left(\dot{x}(t)-Tf(x(t),y(t))\right)\d t+\int_{0}^{\alpha/T}\lambda_\mathrm{CP}^\text{T}(t)\left(y(t)-x(t+1-\alpha/T)\right)\d t\nonumber\\
    &\quad+\int_{\alpha/T}^1\lambda_\mathrm{CP}^\text{T}(t)\left(y(t)-x(t-\alpha/T)\right)\d t+\lambda_\mathrm{BC}^\text{T}\left(x(0)-x(1)\right)+\lambda_\mathrm{PS}h(x(0))
\end{align}
with respect to its arguments at an extremal point. We assume that $x,\lambda_\mathrm{DE}\in C^1([0,1];\mathbb{R}^n)$, while $y,\lambda_\mathrm{CP}\in C^0([0,1];\mathbb{R}^n)$. As in the previous section, variations with respect to the Lagrange multipliers $\lambda_\mathrm{DE}(\cdot)$, $\lambda_\mathrm{CP}(\cdot)$, $\lambda_\mathrm{BC}$, $\lambda_\mathrm{PS}$, and $\eta$ yield the boundary-value problem
\begin{equation}
    \dot{\tilde{x}}(t)-\tilde{T}f(\tilde{x}(t),\tilde{y}(t))=0\text{ for }t\in(0,1),\,\tilde{x}(0)-\tilde{x}(1)=0,\,h(\tilde{x}(0))=0,
    \label{eq:ddebvp}
\end{equation}
and $\tilde{T}-\tilde{\mu}=0$, where
\begin{align}
    \tilde{y}(t)-\tilde{x}(t+1-\alpha/\tilde{T})&=0\text{ for }t\in(0,\alpha/\tilde{T}),\\\tilde{y}(t)-\tilde{x}(t-\alpha/\tilde{T})&=0\text{ for }t\in(\alpha/\tilde{T},1).
    \label{eq:ddeqcp}
\end{align}
It follows that $\tilde{x}(t/\tilde{T})$ is a periodic solution with period $\tilde{\mu}=\tilde{T}$ of the delay differential equation \eqref{eq:ddedisc} and with initial condition on the zero-level surface of the function $h$.

Vanishing variations of $L$ with respect to $x(\cdot)$, $y(\cdot)$, $T$, and $\mu$, yields the necessary adjoint differential equations
\begin{equation}
    -\dot{\tilde{\lambda}}_\mathrm{DE}^\text{T}(t)-\tilde{\lambda}_\mathrm{DE}^\text{T}(t)\tilde{T}D_xf(\tilde{x}(t),\tilde{y}(t))-\tilde{\lambda}_\mathrm{CP}^\mathrm{T}(t+\alpha/\tilde{T})=0,
    \label{eq:adjdde1}
\end{equation}

for $t\in(0,1-\alpha/\tilde{T})$ and
\begin{equation}
    -\dot{\tilde{\lambda}}_\mathrm{DE}^\text{T}(t)-\tilde{\lambda}_\mathrm{DE}^\text{T}(t)\tilde{T}D_xf(\tilde{x}(t),\tilde{y}(t))-\tilde{\lambda}_\mathrm{CP}^\mathrm{T}(t+\alpha/\tilde{T}-1)=0,
    \label{eq:adjdde2}
\end{equation}
for $t\in(1-\alpha/\tilde{T},1)$, boundary conditions
\begin{equation}
    -\tilde{\lambda}_\mathrm{DE}^\text{T}(0)+\lambda_\mathrm{BC}^\text{T}+\lambda_\mathrm{PS} Dh(x(0))=0,\,\tilde{\lambda}_\mathrm{DE}^\text{T}(1)-\lambda_\mathrm{BC}^\text{T}=0,
    \label{eq:adjddebc}
\end{equation} 
coupling conditions
\begin{equation}
    -\tilde{\lambda}_\mathrm{DE}^\text{T}(t)\tilde{T}D_yf(\tilde{x}(t),\tilde{y}(t))+\tilde{\lambda}_\mathrm{CP}^\text{T}(t)=0,\,t\in(0,1)
    \label{eq:adjddecp}
\end{equation}
integral condition
\begin{align}
    \eta-\int_0^1\tilde{\lambda}_\mathrm{DE}^\text{T}(t)f(\tilde{x}(t),\tilde{y}(t))\,\d t&-\frac{\alpha}{\tilde{T}^2}\int_0^{\alpha/\tilde{T}}\tilde{\lambda}_\mathrm{CP}^\text{T}(t)\dot{\tilde{x}}(t+1-\alpha/\tilde{T})\,\d t-\frac{\alpha}{\tilde{T}^2}\int_{\alpha/\tilde{T}}^1\tilde{\lambda}_\mathrm{CP}^\text{T}(t)\dot{\tilde{x}}(t-\alpha/\tilde{T})\,\d t=0,
    \label{eq:adjddeint}
\end{align}
and $1-\eta=0$. 

As in the previous section, we show by differentiation and use of \eqref{eq:ddebvp}-\eqref{eq:ddeqcp}, \eqref{eq:adjdde1}-\eqref{eq:adjdde2}, and \eqref{eq:adjddecp} that the function
\begin{equation}
    \tilde{\lambda}_\mathrm{DE}^\text{T}(t)f(\tilde{x}(t),\tilde{y}(t))+\frac{1}{\tilde{T}}\int_t^{t+\alpha/\tilde{T}}\tilde{\lambda}_\mathrm{CP}^\text{T}(s)\dot{\tilde{x}}(s-\alpha/\tilde{T})\,\d s
    \label{eq:hfunc1}
\end{equation}
for $t\in[0,1-\alpha/\tilde{T})$ and
\begin{align}
   &\tilde{\lambda}_\mathrm{DE}^\text{T}(t)f(\tilde{x}(t),\tilde{y}(t))+\frac{1}{\tilde{T}}\int_t^{1}\tilde{\lambda}_\mathrm{CP}^\text{T}(s)\dot{\tilde{x}}(s-\alpha/\tilde{T})\,\d s+\frac{1}{\tilde{T}}\int_{0}^{t-1+\alpha/\tilde{T}}\tilde{\lambda}_\mathrm{CP}^\text{T}(s)\dot{\tilde{x}}(s-\alpha/\tilde{T})\,\d s
    \label{eq:hfunc2}
\end{align}
for $t\in[1-\alpha/\tilde{T},1]$ is continuous and constant, such that
\begin{equation}
    \tilde{\lambda}_\mathrm{DE}^\text{T}(1)f(\tilde{x}(1),\tilde{y}(1))=\tilde{\lambda}_\mathrm{DE}^\text{T}(0)f(\tilde{x}(0),\tilde{y}(0)).
\end{equation}
From \eqref{eq:adjddebc} it then follows that $\lambda_\mathrm{PS}=0$ provided that the periodic trajectory $\tilde{x}(t)$ intersects $h=0$ transversally at $\tilde{x}(0)$, since in this case $Dh(x(0))f(x(0),y(0))\ne 0$. In this case, $\tilde{\lambda}_\mathrm{DE}(t/\tilde{T})$ is also periodic with period $\tilde{T}$. By \eqref{eq:adjddecp} this also holds for the function $\tilde{\lambda}_\mathrm{CP}(t/\tilde{T})$. It follows that the constant function in \eqref{eq:hfunc1} and \eqref{eq:hfunc2} may be written in the form \eqref{eq:hfunc1} for all $t$. Integration of this function over $t\in[0,1]$ and changing the order of integration then yields
\begin{align}
    &\int_0^1\tilde{\lambda}_\mathrm{DE}^\text{T}(t)f(\tilde{x}(t),\tilde{y}(t))\,\d t+\frac{1}{\tilde{T}}\int_0^1\int_t^{t+\alpha/\tilde{T}}\tilde{\lambda}_\mathrm{CP}^\text{T}(s)\dot{\tilde{x}}(s-\alpha/\tilde{T})\,\d s\,\d t\nonumber\\
    &\quad=\int_0^1\tilde{\lambda}_\mathrm{DE}^\text{T}(t)f(\tilde{x}(t),\tilde{y}(t))\,\d t+\frac{\alpha}{\tilde{T}^2}\int_0^1\tilde{\lambda}_\mathrm{CP}^\text{T}(s)\dot{\tilde{x}}(s-\alpha/\tilde{T})\,\d s=1,
\end{align}
where we used \eqref{eq:adjddeint}, periodicity, and the fact that $\eta=1$. After substitution for $\dot{\tilde{x}}$ and of the integration variable, we obtain the normalization condition~\cite{novivcenko2012phase}
\begin{equation}
    \tilde{\lambda}_\mathrm{DE}^\text{T}(0)f(\tilde{x}(0),\tilde{y}(0))+\int_{-\alpha/\tilde{T}}^0\tilde{\lambda}_\mathrm{CP}^\text{T}(s+\alpha/\tilde{T})f(\tilde{x}(s),\tilde{y}(s))\,\d s=1.
\end{equation}

As in the previous section, the regularity of the periodic orbit implies the existence of a unique triplet $(x(\cdot),y(\cdot),T)$ near 
$(\tilde{x}(\cdot),\tilde{y}(t),\tilde{T})$ for each pair of small $\delta_\mathrm{BC}$ and $\delta_\mathrm{PS}$, such that
\begin{equation}
    \dot{x}(t)-Tf(x(t),y(t))=0\text{ for }t\in[0,1],x(0)-x(1)=\delta_\mathrm{BC},\,h(x(0))=\delta_\mathrm{PS},
\end{equation}
where
\begin{align}
    y(t)-x(t+1-\alpha/T)&=0\text{ for }t\in(0,\alpha/T),\\y(t)-x(t-\alpha/T)&=0\text{ for }t\in(\alpha/T,1).
\end{align}
From the analysis in Section~\ref{sec: phase response curves linear response theory}, we conclude that
\begin{equation}
    T-\tilde{T}=-\tilde{\lambda}_\mathrm{DE}^\text{T}(0)\delta_\mathrm{BC}+O\left(\|\delta_\mathrm{BC}\|^2\right).
\end{equation}
where we have used the fact that $\tilde{\lambda}_\mathrm{BC}=\tilde{\lambda}_\mathrm{DE}(1)=\tilde{\lambda}_\mathrm{DE}(0)$ and $\tilde{\lambda}_\mathrm{PS}=0$. For an asymptotically stable limit cycle, we may again associate $\tilde{\lambda}_\mathrm{DE}^\text{T}(0)$ with the Frech\'{e}t derivative $D\varphi(\tilde{x}(0))$ of the corresponding asymptotic phase~\cite{chicone2004asymptotic}.

\subsubsection{Problem construction and analysis}
\label{sec: phase response curves problem construction}
As in the previous section on optimization, we obtain an augmented continuation problem $\mathbf{A}$ of the form in \eqref{eq:augmented continuation problem} corresponding to the analysis of a periodic orbit with $\mathcal{U}=\mathbb{R}^n$
by associating
\begin{itemize}
    \item $\Phi$ with the boundary-value problem in \eqref{eq:odeprc} in terms of the continuation variables $x$ and $T$; 
    \item $\Psi$ with the scalar $T$ and corresponding continuation parameter $\mu$; and 
    \item $\Lambda^\ast$ with the linear operator in \eqref{eq:adjdiffodeprc}-\eqref{eq:adjTodeprc} acting on the continuation multipliers $\lambda_\mathrm{DE}$, $\lambda_\mathrm{BC}$, $\lambda_\mathrm{PS}$, and $\eta$.
\end{itemize}
This problem has dimensional deficit $1$ which reduces to $0$ once a solution is found with $\eta=1$.

After suitable discretization, we may construct $\mathbf{A}$ according to the following algorithm:
\begin{itemize}
    \item[\textbf{Step 1:}] Invoke the core constructor \eqref{eq: coco_add_func1} with $\mathfrak{phi}$ encoding the differential constraint in \eqref{eq:odeprc}, $\mathbb{K}^\mathrm{o}_u=\emptyset$, and $u^\mathrm{n}_0$ given by an initial solution guess for the continuation variables $x$ and $T$.
    \item[\textbf{Step 2:}] Invoke the core constructor \eqref{eq: coco_add_func1} with $\mathfrak{phi}$ encoding the periodic boundary conditions in \eqref{eq:odeprc}, $\mathbb{K}^\mathrm{o}_u$ indexing the corresponding continuation variables from \textbf{Step 1}, and $u^\mathrm{n}_0=\emptyset$.
    \item[\textbf{Step 3:}] Invoke the core constructor \eqref{eq: coco_add_func1} with $\mathfrak{phi}$ encoding the phase condition in \eqref{eq:odeprc}, $\mathbb{K}^\mathrm{o}_u$ indexing the corresponding continuation variables from \textbf{Step 1}, and $u^\mathrm{n}_0=\emptyset$.
    \item[\textbf{Step 4:}] Invoke the core constructor \eqref{eq: coco_add_func2} with $\mathfrak{psi}$ encoding the evaluation of $T$, $\mathbb{K}^\mathrm{o}_u$ indexing the corresponding continuation variable from \textbf{Step 1}, and $u^\mathrm{n}_0=\emptyset$.
    \item[\textbf{Step 5:}] Invoke the core constructor \eqref{eq: coco_add_adjt} with $\mathfrak{lambda}$ encoding the linear operators acting on $\lambda_{DE}$ in the adjoint conditions \eqref{eq:adjdiffodeprc}-\eqref{eq:adjTodeprc}, $\mathbb{K}^\mathrm{o}_u$ indexing the continuation variables introduced in the corresponding call in \textbf{Step 1}, $\mathbb{K}^\mathrm{o}_\Lambda=\emptyset$, and $\lambda^\mathrm{n}_0$ given by an initial solution guess for the continuation variables $\lambda_{DE}$.
    \item[\textbf{Step 6:}] Invoke the core constructor \eqref{eq: coco_add_adjt} with $\mathfrak{lambda}$ encoding the linear operators acting on $\lambda_{BC}$ in the adjoint conditions \eqref{eq:adjdiffodeprc}-\eqref{eq:adjTodeprc}, $\mathbb{K}^\mathrm{o}_u$ indexing the continuation variables associated with the corresponding call in \textbf{Step 2}, $\mathbb{K}^\mathrm{o}_\Lambda\ne\emptyset$, and $\lambda^\mathrm{n}_0$ given by an initial solution guess for the continuation variables $\lambda_{BC}$.
    \item[\textbf{Step 7:}] Invoke the core constructor \eqref{eq: coco_add_adjt} with $\mathfrak{lambda}$ encoding the linear operators acting on $\lambda_{PS}$ in the adjoint conditions \eqref{eq:adjdiffodeprc}-\eqref{eq:adjTodeprc}, $\mathbb{K}^\mathrm{o}_u$ indexing the continuation variables associated with the corresponding call in \textbf{Step 3}, $\mathbb{K}^\mathrm{o}_\Lambda\ne\emptyset$, and $\lambda^\mathrm{n}_0$ given by an initial solution guess for the continuation variables $\lambda_{PS}$.
    \item[\textbf{Step 8:}] Invoke the core constructor \eqref{eq: coco_add_adjt} with $\mathfrak{lambda}$ encoding the linear operators acting on $\eta$ in the adjoint conditions \eqref{eq:adjdiffodeprc}-\eqref{eq:adjTodeprc}, $\mathbb{K}^\mathrm{o}_u$ indexing the continuation variables associated with the corresponding call in \textbf{Step 4}, $\mathbb{K}^\mathrm{o}_\Lambda\ne\emptyset$, and $\lambda^\mathrm{n}_0$ given by an initial solution guess for the continuation variable $\eta$.
\end{itemize}
As suggested previously, steps 5 through 8 can be implemented automatically from information provided in steps 1 through 4, thereby reducing the task of construction to the definition of the vector field $f$. To solve for the corresponding phase response curve, we invoke the core constructor \eqref{eq: coco_add_comp1} to append a complementary zero function that evaluates to $1-\eta$. Here, $\mathbb{K}^\mathrm{o}_u=\mathbb{K}^\mathrm{o}_v=\emptyset$ and  $\mathbb{K}^\mathrm{o}_\lambda$ indexes the continuation multiplier $\eta$. 

For the analysis of a periodic orbit of the delay differential equation \eqref{eq:ddedisc}, we similarly obtain an augmented continuation problem $\mathbf{A}$ of the form in \eqref{eq:augmented continuation problem} by associating
\begin{itemize}
    \item $\Phi$ with the boundary-value problem in \eqref{eq:ddebvp}-\eqref{eq:ddeqcp} in terms of the continuation variables $x$, $y$, and $T$; \item $\Psi$ with the scalar $T$ and corresponding continuation parameter $\mu$; and \item $\Lambda^\ast$ with the linear operator in \eqref{eq:adjdde1}-\eqref{eq:adjddeint} acting on the continuation multipliers $\lambda_\mathrm{DE}$, $\lambda_\mathrm{CP}$, $\lambda_\mathrm{BC}$, $\lambda_\mathrm{PS}$, and $\eta$.
\end{itemize}
This problem again has dimensional deficit $1$ which reduces to $0$ once a solution is found with $\eta=1$. We leave a detailed description of the algorithm of construction to the reader.

\section{Toolbox construction}
\label{sec: Toolbox construction}

The examples on data assimilation and phase response curves in Sections~\ref{sec: data assimilation} and \ref{sec: phase response curves} have demonstrated how advanced examples of augmented continuation problems may be constructed through repeated calls to the core \textsc{coco} constructors in \eqref{eq: coco_add_func1}-\eqref{eq: coco_add_adjt}. Several observations follow from this discussion and inform our continued development in this section.

\subsection{Composite construction}
\label{sec: composite construction}
A \textsc{coco} constructor of the form \eqref{eq:cocooperator} with $\mathbb{K}^\mathrm{o}_u=\emptyset$, $\mathbb{K}^\mathrm{o}_\lambda=\emptyset$, $\mathbb{K}^\mathrm{o}_\Lambda=\emptyset$, and $\mathbb{K}^\mathrm{o}_v=\emptyset$ and with its remaining arguments defined according to an abstract paradigm is called a \emph{toolbox constructor}. We often use this terminology also to describe constructors whose index sets $\mathbb{K}^\mathrm{o}_u$, $\mathbb{K}^\mathrm{o}_\lambda$, $\mathbb{K}^\mathrm{o}_\Lambda$, and $\mathbb{K}^\mathrm{o}_v$ are associated with a preceding application of a toolbox constructor and whose remaining arguments are defined by the same abstract paradigm. A collection of toolbox constructors is called a (\textsc{coco}-compatible) toolbox. A paradigm for toolbox construction using the \textsc{coco} platform is described in Part II of textbook \cite{dankowicz2013recipes}.

As an example, a toolbox constructor may be designed to construct the zero problem and monitor functions associated with analyzing equilibria of an arbitrary autonomous vector field in terms of its problem parameters. In this case, the functions $\mathfrak{phi}$ and $\mathfrak{psi}$ may be defined to include calls to user-defined encodings of the vector field and its derivatives without assuming a particular state space dimension or a particular number of problem parameters. In the \textsc{coco}-compatible toolbox \mcode{'ep'}~(see EP-Tutorial.pdf in \cite{COCO}), for example, the toolbox constructors \mcode{ode_isol2ep} and \mcode{ode_ep2ep} accomplish this task using data provided by the user and data stored from a previous analysis, respectively. 

For constrained optimization along manifolds of equilibria, it is necessary to construct the contributions to the adjoint conditions corresponding to the equilibrium constraints. In the \mcode{'ep'} toolbox, this is accomplished with the \mcode{adjt_isol2ep} and \mcode{adjt_ep2ep} toolbox constructors, provided that the equilibrium constraint was constructed using \mcode{ode_isol2ep} and \mcode{ode_ep2ep}, respectively. Since the vector field, the state space dimension, and the number of problem parameters were provided already to the latter constructors, the adjoint constructors \mcode{adjt_isol2ep} and \mcode{adjt_ep2ep} require no further user data. This observation simplifies the calling syntax to these operators.

These same principles apply also to the constructors associated with the \textsc{coco}-compatible \mcode{'coll'}~(see COLL-Tutorial.pdf in \cite{COCO}) and \mcode{'po'}~(see PO-Tutorial.pdf in \cite{COCO}) toolboxes that are included with the \textsc{coco} release. For example, given a vector field, a sequence of time instants, a corresponding array of state-space vectors, and an array of numerical values for the problem parameters, the toolbox constructor \mcode{ode_isol2coll} constructs a discretization of the trajectory problem
\begin{equation}
    \dot{x}=f(t,x,p),\,t\in[T_0,T_0+T]
\end{equation}
in terms of the unknown values of $x(t)$ for a finite set of values of $t\in[T_0,T_0+T]$, unknown values of $p$, and unknown initial time $T_0$ and duration $T$. Given a preceding call to \mcode{ode_isol2coll}, the corresponding adjoint contributions may be appended to the augmented continuation problem using the \mcode{adjt_isol2coll} toolbox constructor. As before, no additional information about the problem is required in this call. Additional constructors associated with the \mcode{coll} toolbox encode multi-segment boundary-value problems and the corresponding contributions to the adjoint necessary conditions. The \mcode{'po'} toolbox encodes the special case of periodic boundary conditions, also for piecewise-defined vector fields and hybrid dynamical systems, according to the same fundamental paradigm. 

By definition, we arrive at a toolbox constructor by recognizing a universality among a class of individual problems and by encoding an abstract representation of this universality in suitable constructors. We see elements of such universality represented in the examples in Sections~\ref{sec: data assimilation} and \ref{sec: phase response curves} and proceed in the remainder of this section to derive the mathematical formalism of a corresponding toolbox.

\subsection{Delay graphs}
\label{sec: delay graphs}
The data assimilation example in Section~\ref{sec: data assimilation} may be abstractly represented according to a graph-theoretic framework~\cite{paul2000designing}. Specifically, let the term \emph{segment} here refer to a variable $x:[-\hat{\alpha},T]\rightarrow\mathbb{R}^n$ in terms of a \emph{maximal past} $\hat{\alpha}$ and \emph{duration} $T$ such that $\dot{x}=f(x(t),x(t-\alpha))$ for $t\in(0,T)$ in terms of a vector field $f$ and \emph{delay} $\alpha\le\hat{\alpha}$. For a collection of $M$ segments, we use a subscript $k=1,\ldots,M$ to identify individual segments. In the example in Section~\ref{sec: data assimilation}, $x_k(T_k\alpha)=(u^{(k)}(\alpha),p^{(k)}(\alpha))$ with $\alpha_k=0$ for $k<I$ and $\alpha_k=\alpha$ for $k\ge I$.

In this section, we say that the $i$-th and $j$-th segments are \emph{coupled}, in that order, if there exists a coupling matrix $B_{i,j}$ such that
\begin{equation}
\label{eq: coupledsegments}
    x_i(s)=B_{i,j}x_j(T_j+s)\text{ for all }s\in[-\hat{\alpha}_i,0]
\end{equation}
and, in particular, that $x_i(0)=B_{i,j}x_j(T_j)$.
Clearly, this is possible only if $-\hat{\alpha}_j\le T_j-\hat{\alpha}_i$. We say that the $j$-th segment is a \emph{predecessor} of the $i$-th segment. For a consistent definition, we require that there be at most one predecessor for each segment, but allow a segment to be a predecessor of multiple nodes. We naturally arrive at a directed graph with nodes representing segments and directed edges representing predecessor coupling. For a segment that is not a predecessor of any other segment, we may assume without loss of generality that $\hat{\alpha}_i=\alpha_i$. Similarly, for a segment that is not preceded by another segment, we may assume without loss of generality that $\hat{\alpha}_i=0$. It is not the case, however, that $\hat{\alpha}_j$ may be assumed to equal $0$ for a segment with $\alpha_j=0$, since the former is involved in coupling conditions of the form \eqref{eq: coupledsegments} with other segments.

\begin{figure}[ht]
\centering
\subfloat{\includegraphics[width=0.7\columnwidth,trim={0 {0.0\textwidth} 0 0.0\textwidth},clip]{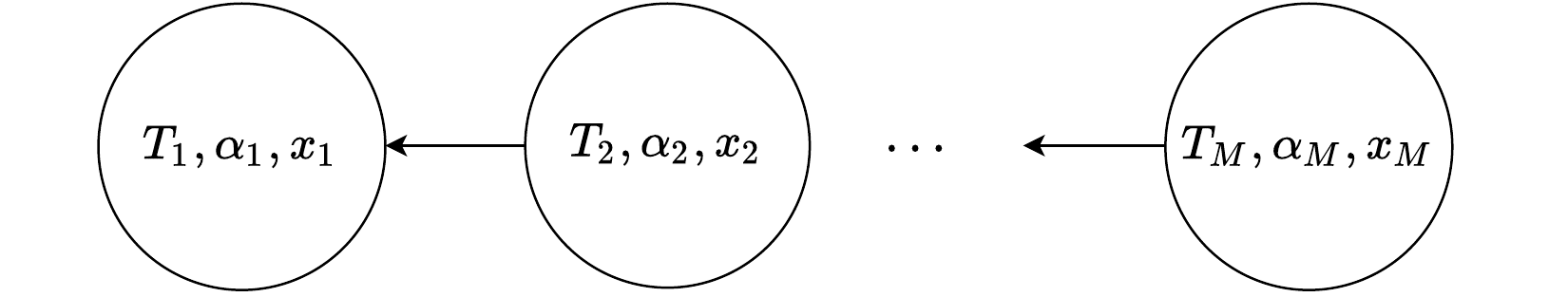}}
\caption{Directed graph representation of the data assimilation problem with $M$ segments. In the special case in the text, $M=5$, the segment lengths are $T_{1}=T_{3}=T_{4}=0.2T,T_{2}=0.3T,T_{5}=0.1T$ and the delays equal $\alpha_{1}=\alpha_{2}=0$, $\alpha_{3}=\alpha_{4}=\alpha_{5}=\alpha:=0.5T$.} 
\label{fig:networks1}
\end{figure}

The example in Section~\ref{sec: data assimilation} may be described in terms of the directed graph in Fig.~\ref{fig:networks1}, in which the $k$-th segment is the unique predecessor to the $k+1$-th segment with coupling matrix $I_n$. By definition
\begin{equation}
    \sum_{k=1}^{I-1} T_k=\alpha,\,\sum_{k=1}^M T_k=T.
\end{equation}
From this graph, we see that $x_{k}(s)=x_{k-1}(T_{k-1}+s)$ for all $s\in[-\hat{\alpha}_k,0]$. Provided that $T_{k-1}-\alpha_k\ge 0$, we may evaluate the delayed term $x_k(t-\alpha_k)$ for all $t\in(0,T_k)$ without further reference to the graph. If, instead, $T_{k-1}-\alpha_k<0$, we may use the fact that $x_{k-1}(s)=x_{k-2}(T_{k-2}+s)$ for all $s\in[-\hat{\alpha}_{k-1},0]$ to obtain $x_{k}(s)=x_{k-2}(T_{k-2}+T_{k-1}+s)$ for $s\in[-\hat{\alpha}_k,-T_{k-1}]$. Provided that $T_{k-2}+T_{k-1}-\alpha_k\ge 0$, we may evaluate the delayed term $x_k(t-\alpha_k)$ for all $t\in(0,T_k)$ without further reference to the graph. If not, then we proceed iteratively until the sum $\sum_{l=1}^{L}T_{k-l}$ equals or exceeds $\alpha_k$ for some $L<k$.

As an example, suppose that $M=5$, $\alpha=0.5T$, $T_1=T_3=T_4=0.2T$, $T_2=0.3T$, and $T_5=0.1T$. The algorithm in the previous paragraph shows that
\begin{align}
    x_5(s)&=x_2(T_2+T_3+T_4+s),\,s\in[-\alpha_5,T_5-\alpha_5],\\
    x_4(s)&=x_2(T_2+T_3+s),\,s\in[-\alpha_4,T_4-\alpha_4],\\
    x_3(s)&=x_1(T_1+T_2+s),\,s\in[-\alpha_3,T_3-\alpha_3],
\end{align}
and, consequently,
\begin{align}
    x_5(t-\alpha_5)&=x_2(T_2+T_3+T_4-\alpha_5+t),\,t\in[0,T_5],\\
    x_4(t-\alpha_4)&=x_2(T_2+T_3-\alpha_4+t),\,t\in[0,T_4],\\
    x_3(t-\alpha_3)&=x_1(T_1+T_2-\alpha_3+t),\,t\in[0,T_3].
\end{align}
It follows, for example, that
\begin{equation}
    u^{(5)}\left(t-\frac{\alpha_5}{T_5}\right)=u^{(2)}\left(\frac{T_2+T_3+T_4-\alpha_5}{T_2}+\frac{T_5}{T_2}t\right),\,t\in[0,1]
\end{equation}
in agreement with the general expression in \eqref{eq:cp}. The reader is encouraged to perform the corresponding calculations for other choices of the number of segments $M$, the time delay $\alpha$, and interval durations $T_k$.

From the general form of the adjoint contributions in Section~\ref{sec: data assimilation}, it may be correctly surmised that the directed graph in Fig.~\ref{fig:networks1} contains all the information required to construct these expressions for any number of segments $M$, time delay $\alpha$, and interval durations $T_k$. This is analogous to the possible application of the adjoint constructors \mcode{adjt_isol2coll} and \mcode{adjt_isol2po} in the \mcode{'coll'} and \mcode{'ep'} toolboxes, respectively, without requiring additional information than that provided to \mcode{ode_isol2coll} and \mcode{ode_isol2po}, respectively. 

We recognize in this construction of constraints and adjoint contributions a universal paradigm for multi-segment boundary-value problems with discrete delays. In the next section, we formulate an abstract toolbox template that is sufficiently flexible to handle a broad class of such problems. Several examples illustrate the reduction of the abstract framework to problems involving periodic and quasiperiodic orbits.

\subsection{A toolbox template for delay-coupled differential equations}
\label{sec: toolbox template}
\subsubsection{An abstract zero problem}
\label{sec: toolbox template zero problems}

Consider a collection of non-autonomous vector fields $f_i:\mathbb{R}\times\mathbb{R}^{n}\times\mathbb{R}^{n}\times\mathbb{R}^{q}\rightarrow\mathbb{R}^{n}$, for $i=1,\ldots,M$, governing the time histories of $M$ \emph{differential state variables} $x_i\in C^1\left([0,1],\mathbb{R}^{n}\right)$ according to the \emph{differential constraints}
\begin{equation}
    x_{i}^{\prime}(\tau)=T_{i}f_{i}\left(T_{0,i}+T_{i}\tau,x_{i}(\tau),y_{i}(\tau),p\right),\,\tau\in\left(0,1\right),\label{rep_eq1}
\end{equation}
in terms of the individual \emph{initial times} $T_{0,i}\in\mathbb{R}$, individual \emph{durations} $T_i\in\mathbb{R}$, and shared \emph{problem parameters} $p\in\mathbb{R}^q$. For a subset of indices $i$, we assume the imposition of \emph{boundary conditions}
\begin{equation}
\label{eq: abstractbcs}
    x_i(0)=\sum_{j=1}^{ M}B_{i,j}(p)\cdot x_j(1)
\end{equation}
in terms of (possibly zero) \emph{coupling matrices} $B_{i,j}(p)$.

The collection of \emph{algebraic state variables} $y_i\in C^0\left([0,1],\mathbb{R}^{n}\right)$ provide exogenous excitation to the system dynamics~\eqref{rep_eq1}. For indices $i$ such that $f_i$ depends explicitly on $y_i$, we obtain a closed abstract model by specifying $C_i$ \emph{coupling conditions} on a partition of $[0,1]$ into subintervals $\left[\gamma_{\mathrm{b},i,k},\gamma_{\mathrm{e},i,k}\right]$ with
\begin{equation}
    0=\gamma_{\mathrm{b},i,1}\le \gamma_{\mathrm{e},i,1}=\gamma_{\mathrm{b},i,2}\le\cdots\le\gamma_{\mathrm{e},i,C_i-1}=\gamma_{\mathrm{b},i,C_i}\le\gamma_{\mathrm{e},i,C_i}=1
    \label{rep_eq3}
\end{equation}
according to the expressions
\begin{equation}
    y_{i}(\tau)=\sum_{s=1}^{S_{i,k}}A_{i,k,s}(p)\cdot x_{j_{i,k,s}}\left( \frac{T_{i}}{T_{j_{i,k}}}\left(\tau - \varDelta_{i,k}\right) \right) ,\,\tau\in\left[\gamma_{\mathrm{b},i,k},\gamma_{\mathrm{e},i,k}\right]\label{rep_eq2}
\end{equation}
for $k=1,\ldots,C_{i}$ and $j_{i,k},j_{i,k,s}\in\{1,\ldots,M\}$, and where 
\begin{equation}
\label{eq_xi}
    0\le\xi_{\mathrm{b},i,k}:=\frac{T_i}{T_{j_{i,k}}}\left(\gamma_{\mathrm{b},i,k}-\varDelta_{i,k}\right)\le\xi_{\mathrm{e},i,k}:=\frac{T_i}{T_{j_{i,k}}}\left(\gamma_{\mathrm{e},i,k}-\varDelta_{i,k}\right)\le 1
\end{equation}
for all $i$ and $k$. In particular, the \emph{coupling delays} 
$\varDelta_{i,k}$ 
and interval boundaries  $\gamma_{\mathrm{b},i,k}$, and $\gamma_{\mathrm{e},i,k}$ are constrained such that
\begin{align}
    \xi_{\mathrm{e},i,k}&=1,\,k=1,\ldots,C_i-1,\label{dde:xi:e}\\
    \xi_{\mathrm{b},i,k}&=0,\,k=2,\ldots,C_i,\label{dde:xi:b}
\end{align}
and the \emph{composite coupling matrices} $A_{i,k,s}$ are chosen so that continuity of $y_i$ across $\tau=\gamma_{\mathrm{e},i,k}=\gamma_{\mathrm{b},i,k+1}$ for $k=1,\ldots,C_i-1$ is implied by the boundary conditions~\eqref{eq: abstractbcs}. It is clear that the functions $y_i$ are completely determined by the coupling conditions. Consequently, no additional constraints on the algebraic state variables can be added to the continuation problem. 

In practice, additional problem-specific algebraic constraints relate the quantities $T_{0,i}$, $T_i$, and $\gamma_{\mathrm{e},i,1}$ for $i=1,\ldots,M$ to each other and/or the problem parameters $p$. Per the additive principles of the constrained optimization paradigm, we omit consideration of such dependencies in the abstract framework, but make the relationships explicit in the context of particular examples.

\subsubsection{Examples}
\label{sec: toolbox template examples}
As a first example, suppose that $M=1$, $C_1=2$, and $S_{1,1}=S_{1,2}=j_{1,1,1}=j_{1,2,1}=j_{1,1}=j_{1,2}=1$, and let
\begin{equation}
    B_{1,1}=A_{1,1,1}=A_{1,2,1}=I_n,\,\gamma_{\mathrm{e},1,1}=\frac{\alpha}{T}
\end{equation}
with $T_1=T$ in terms of the problem parameters $\alpha$ and $T$, where $0\le\alpha\le T$. It follows that
\begin{equation}
    \gamma_{\mathrm{b},1,2}=\varDelta_{1,2}=\frac{\alpha}{T},\,\varDelta_{1,1}=\frac{\alpha}{T}-1.
\end{equation}
We obtain the equations
\begin{equation}
\begin{aligned}
 x'(\tau)&=Tf\left(T_{0}+T\tau,x(\tau),y(\tau),p\right),\,\tau\in\left(0,1\right),\\
 y(\tau)&=\begin{cases}
x\left(\tau+1-\alpha/T\right), & \tau\in\left(0,\alpha/T\right),\\
x\left(\tau-\alpha/T\right), & \tau\in\left(\alpha/T,1\right),
\end{cases}
\end{aligned}\label{eq:toolbox:example1}
\end{equation}
where we omitted the trivial subscript (cf.\ the governing boundary-value problem in the derivation of the phase response functional in Section~\ref{sec: phase response curves delay-differential equations} for an autonomous vector field). Here, the boundary condition $x(0)=x(1)$ implies continuity of $y(\tau)$ across $\tau=\alpha/T$, while evaluation of the coupling conditions at $\tau=0$ and $\tau=1$ shows that $y(0)=y(1)$. With $z(t):=x\left((t-T_0)/T\right)$ and $\tilde{z}(t):=y\left((t-T_0)/T\right)$, the coupling conditions may be condensed using the modulo operator to yield
\begin{equation}
    \tilde{z}(t)=z\left(t-\alpha\big|_{\mathrm{mod}[T_0,T_0+T)}\right),\,t\in[T_0,T_0
    +T),
\end{equation}
from which we conclude that the abstract problem corresponds to the existence of a continuously-differentiable, $T$-periodic solution $z(t)$ of the delay differential equation
\begin{equation}
    \dot{z}(t)=f(t,z(t),z(t-\alpha),p)
\end{equation}
provided that $f$ is periodic with period $T$ in its first argument. 

As a second example, suppose that $M=C_1=C_2=2$,  $S_{1,1}=S_{1,2}=S_{2,1}=S_{2,2}=1$, $j_{1,1,1}=j_{2,2,1}=j_{1,1}=j_{2,2}=2$, and $j_{1,2,1}=j_{2,1,1}=j_{1,2}=j_{2,1}=1$, and let
\begin{equation}
   B_{1,1}=B_{2,2}=0_n,\, B_{1,2}=B_{2,1}=A_{1,1,1}=A_{1,2,1}=A_{2,1,1}=A_{2,2,1}=I_n,\,\gamma_{\mathrm{e},1,1}=\frac{\alpha}{\beta},\,\gamma_{\mathrm{e},2,1}=\frac{\alpha}{T-\beta}
\end{equation}
with $T_1=\beta$, $T_2=T-\beta$ in terms of the problem parameters $\alpha$, $\beta$, and $T$, where $0\le\alpha<\beta<T-\alpha$. It follows that
\begin{equation}
    \gamma_{\mathrm{b},1,2}=\frac{\alpha}{\beta},\,\gamma_{\mathrm{b},2,2}=\frac{\alpha}{T-\beta},\,\varDelta_{1,1}=\frac{\alpha-T+\beta}{\beta},\,\varDelta_{1,2}=\frac{\alpha}{\beta},\,\varDelta_{2,1}=\frac{\alpha-\beta}{T-\beta},\,\varDelta_{2,2}=\frac{\alpha}{T-\beta}.
\end{equation}
Then, if $T_{0,1}=T_0$ and $T_{0,2}=T_0+\beta$, we obtain the equations
\begin{align}
x_{1}^{\prime}(\tau)&=\beta f_{1}\left(T_{0}+\beta\tau,x_{1}(\tau),y_{1}(\tau),p\right),\,\tau\in(0,1), \\
x_{2}^{\prime}(\tau) &=\left(T-\beta\right)f_{2}\left(T_{0}+\beta+\left(T-\beta\right)\tau,x_{2}(\tau),y_{2}(\tau),p\right),\,\tau\in(0,1),\\
y_{1}(\tau)&=\begin{cases}
x_{2}\left(\beta\tau/(T-\beta)+1-\alpha/(T-\beta)\right), & \tau\in\left(0,\alpha/\beta\right),\\
x_{1}\left(\tau-\alpha/\beta\right), & \tau\in\left(\alpha/\beta,1\right),
\end{cases} \label{per_hybrid1}\\
y_{2}(\tau)&=\begin{cases}
x_{1}\left((T-\beta)\tau/\beta+1-\alpha/\beta\right), & \tau\in\left(0,\alpha/(T-\beta)\right),\\
x_{2}\left(\tau-\alpha/(T-\beta)\right), & \tau\in\left(\alpha/(T-\beta),1\right).
\end{cases} \label{per_hybrid2}
\end{align}
Here, the boundary conditions $x_1(0)=x_2(1)$ and $x_1(1)=x_2(0)$ imply continuity of $y_1(\tau)$ across $\tau=\alpha/\beta$ and of $y_2(\tau)$ across $\tau=\alpha/(T-\beta)$, respectively, while evaluation of the coupling conditions at $\tau=0$ and $\tau=1$ shows that $y_1(0)=y_2(1)$ and $y_1(1)=y_2(0)$. With $z_1(t):=x_1\left((t-T_0)/\beta\right)$, $z_2(t):=x_2\left((t-T_0-\beta)/(T-\beta)\right)$, $\tilde{z}_1(t):=y_1\left((t-T_0)/\beta\right)$, and $\tilde{z}_2(t):=y_2\left((t-T_0-\beta)/(T-\beta)\right)$, the coupling conditions may now be condensed using the modulo operator to yield
\begin{align}
\tilde{z}_{1}(t)&=\begin{cases}
z_{2}\left(t-\alpha\big|_{\mathrm{mod}[T_0,T_0+T]}\right), & t\in\left[T_0,T_0+\alpha\right],\\
z_{1}\left(t-\alpha\big|_{\mathrm{mod}[T_0,T_0+T]}\right), & t\in\left[T_0+\alpha,T_0+\beta\right],
\end{cases}\\
\tilde{z}_{2}(t)&=\begin{cases}
z_{1}\left(t-\alpha\big|_{\mathrm{mod}[T_0,T_0+T]}\right), & t\in\left[T_0+\beta,T_0+\alpha+\beta\right],\\
z_{2}\left(t-\alpha\big|_{\mathrm{mod}[T_0,T_0+T]}\right), & t\in\left[T_0+\alpha+\beta,T_0+T\right],
\end{cases}
\end{align}
from which we conclude that the abstract problem corresponds to the existence of a continuous, $T$-periodic,  piecewise differentiable function $z(t)$ of the system of delay differential equations
\begin{align}
    \dot{z}(t)&=f_1(t,z(t),z(t-\alpha),p),\,t\big|_{\mathrm{mod}[T_0,T_0+T)}\in(T_0,T_0+\beta),\\
    \dot{z}(t)&=f_2(t,z(t),z(t-\alpha),p),\,t\big|_{\mathrm{mod}[T_0,T_0+T)}\in(T_0+\beta,T_0+T)
\end{align}
provided that $f_1$ and $f_2$ are periodic with period $T$ in their first argument.

As a final example inspired by the analysis of quasiperiodic invariant tori~\cite{ahsan2020optimization}, suppose that $M$ equals an odd integer, and that $C_i=2$, $S_{i,1}=M$, $S_{i,2}=1$, $j_{i,1,1}=1,\ldots,j_{i,1,M}=j_{i,1}=M$, and $j_{i,2,1}=j_{i,2}=i$ for $i=1,\ldots,M$. For each $i$, let
\begin{equation}
    B_{i,j}=A_{i,1,j}=A_{i,j},\,j=1,\ldots,M
\end{equation}
and
\begin{equation}
        A_{i,2,1}=I_n,\,\gamma_{\mathrm{e},i,1}=\frac{\alpha}{T}
\end{equation}
with $T_1=\cdots=T_{M}=T$ in terms of the problem parameters $\alpha$ and $T$, where $0\le\alpha\le T$. It follows that
\begin{equation}
\gamma_{\mathrm{b},i,2}=\frac{\alpha}{T},\,\varDelta_{i,1}=\frac{\alpha}{T}-1,\,\varDelta_{i,2}=\frac{\alpha}{T}.    
\end{equation}
Then, if $T_{0,1}=\cdots=T_{0,M}=T_0$ and $f_1=\cdots=f_{M}=f$, we obtain
\begin{equation}
\begin{gathered}
x_i'(\tau)=Tf\left(T_{0}+T\tau,x_{i}(\tau),y_{i}(\tau),p\right),\,\tau\in\left(0,1\right),\label{eq:quasi1}\\
y_{i}(\tau)=\begin{cases}
\sum_{s=1}^{M}A_{i,s}
\cdot x_{s}\left(\tau+1-\alpha/T\right), & \tau\in\left(0,\alpha/T\right),\\
x_{i}\left(\tau-\alpha/T\right), & \tau\in\left(\alpha/T,1\right). 
\end{cases}  
\end{gathered}
\end{equation}
Here, the boundary conditions
\begin{equation}
\label{sample:eq2}
x_i(0)=\sum_{j=1}^M A_{i,j}\cdot x_j(1), 
\end{equation}
imply continuity of $y_i(\tau)$ across $\tau=\alpha/T$, while evaluation of the coupling conditions at $\tau=0$ and $\tau=1$ shows that
\begin{equation}
    y_i(0)=\sum_{j=1}^M A_{i,j}\cdot y_j(1).
\end{equation}

We specialize to the case with $A_{i,s}$ given by the $(i,s)$-th $n\times n$ block in the $Mn\times Mn$ matrix $(F\otimes I_n)^{-1}\left((RF)\otimes I_n\right)$, where $F$ denotes the symmetric square matrix whose $(i,j)$-th entry equals $e^{-2\pi\mathrm{j}(i-1)(j-1)/M}$
and $R$ denotes the diagonal matrix whose diagonal elements equal $1,e^{-2\pi\rho},\ldots,e^{-2\pi \lfloor M/2\rfloor\rho},e^{2\pi \lfloor M/2\rfloor\rho},\ldots,e^{2\pi\rho}$ in terms of the \emph{rotation number} $\rho$. Then, if $x_i(\tau)=x\left(\varphi_i,\tau\right)$ and $y_i(\tau)=y\left(\varphi_i,\tau\right)$ in terms of some functions
\begin{equation}
    x\left(\varphi,\tau\right):=\sum_{m=-\lfloor M/2\rfloor}^{\lfloor M/2\rfloor} c_{x,m}(\tau)e^{\mathrm{j}m\varphi},\,
    y\left(\varphi,\tau\right):=\sum_{m=-\lfloor M/2\rfloor}^{\lfloor M/2\rfloor} c_{y,m}(\tau)e^{\mathrm{j}m\varphi},
\end{equation}
and $\varphi_i=2\pi(i-1)/M$, it follows from the coupling conditions that the discrete Fourier transform of $y(\varphi,\tau)$ sampled at $\varphi=\varphi_i$ equals the discrete Fourier transform of $x(\varphi-2\pi\rho,\tau+1-\alpha/T)$ sampled at $\varphi=\varphi_i$ when $\tau\in[0,\alpha/T]$ and, consequently, that
\begin{equation}
\label{eq_coupling_quasi}
    y(\varphi,\tau)=\begin{cases}x(\varphi-2\pi\rho,\tau+1-\alpha/T), & \tau\in(0,\alpha/T),\\
    x(\varphi,\tau-\alpha/T),&\tau\in(\alpha/T,1).
    \end{cases}
\end{equation}
In this case, continuity of $y(\varphi,\tau)$ across $\tau=\alpha/T$ requires that $x(\varphi,0)=x(\varphi-2\pi\rho,1)$, while evaluation at $\tau=0$ and $\tau=1$ shows that $y(\varphi,0)=y(\varphi-2\pi\rho,1)$.
Suppose now, additionally, that
\begin{equation}
\label{eq: pde_quasi}
    x_{,\tau}(\varphi,\tau)=Tf(T_0+T\tau,x(\varphi,\tau),y(\varphi,\tau),p),\,\tau\in(0,1)
\end{equation}
for all $\varphi\in\mathbb{S}$ rather than only at $\varphi=\varphi_i$. Then, with
\begin{equation}
    \tilde{z}(\theta_1(t),\theta_2(t)):=x\left(\theta_1(t)-\rho\theta_2(t),\frac{\theta_2(t)}{2\pi}\right),\,\theta_1(t)=\varphi+2\pi\rho\frac{t-T_0}{T},\,\theta_2(t)=2\pi\frac{t-T_0}{T},
\end{equation}
we obtain
\begin{equation}\label{eq:qp:torus}
    \frac{d}{dt}\tilde{z}(\theta_1(t),\theta_2(t))=\frac{1}{T}x_{,\tau}\left(\theta_1(t)-\rho\theta_2(t),\frac{\theta_2(t)}{2\pi}\right)=f\left(t,z(\theta_1(t),\theta_2(t)),z\left(\theta_1(t)-\frac{2\pi\rho\alpha}{T},\theta_2(t)-\frac{2\pi\alpha}{T}\right),p\right)
\end{equation}
and, consequently, that the function $z(t)=\tilde{z}(\theta_1(t),\theta_2(t))$  describes a continuously-differentiable, quasiperiodic solution with angular frequencies $2\pi/T$ and $2\pi\rho/T$ of the delay differential equation
\begin{equation}\label{eq:qp:dde}
    \dot{z}(t)=f(t,z(t),z(t-\alpha),p)
\end{equation}
provided that $f$ is periodic with period $T$ in its first argument and $\rho$ is irrational. 

\begin{figure}[ht]
\centering
\subfloat{\includegraphics[width=0.45\columnwidth,trim={0 {0.0\textwidth} 0 0.0\textwidth},clip]{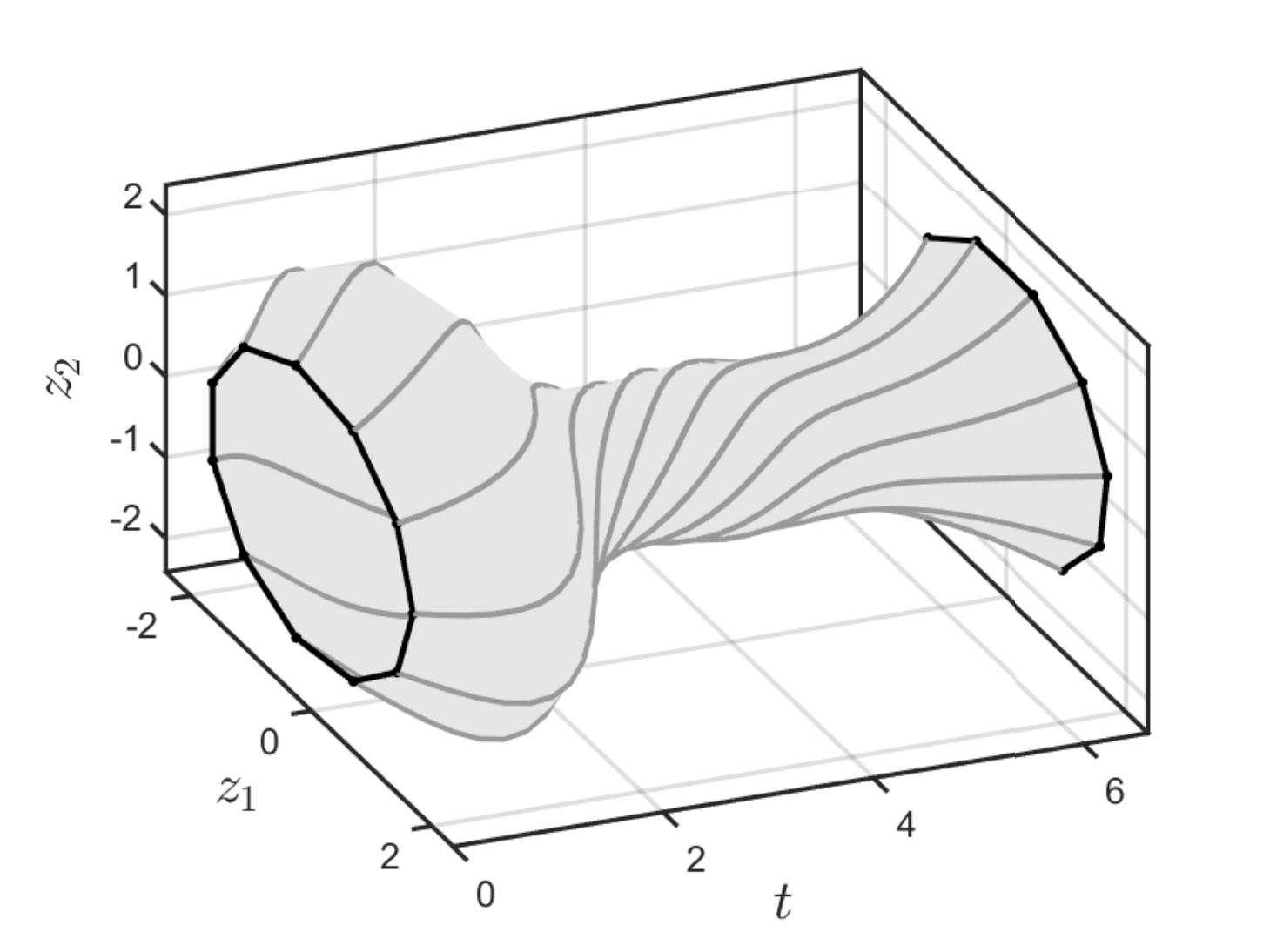}}
\caption{A discretized representation of a quasiperiodic invariant torus for a two-dimensional non-autonomous equation of the form \eqref{eq:qp:dde} using the multi-segment formalism described in (\ref{eq:quasi1})-(\ref{sample:eq2}) with $M=11$ and $T=2\pi$.}
\label{fig:quasi_approx}
\end{figure}

An example of the use of (\ref{eq:quasi1})-(\ref{sample:eq2}) to approximate a quasiperiodic invariant torus for a two-dimensional equation of the form \eqref{eq:qp:dde} (taken from~\cite{ahsan2020optimization}) using a finite collection of trajectory segments is shown in Fig.~\ref{fig:quasi_approx}. The boundary conditions~\eqref{sample:eq2} correspond to a discretized representation of a relative rotation by $2\pi\rho$ between the intersections of the torus with the $t=0$ and $t=T$ (here, $T=2\pi$) surfaces, as is also the case for problems without delay~\cite{dankowicz2013recipes}. In the presence of delay, the history for each segment is obtained by Fourier interpolation over the family of trajectory segments, rotated by $2\pi\rho$ and shifted by $T$. 

\subsubsection{Adjoint conditions}
\label{sec: toolbox template adjoint conditions}
We proceed to consider the problem of optimizing a scalar valued function of the differential state variables $x_i$, algebraic state variables $y_i$, initial times $T_{0,i}$, durations $T_i$, coupling delays $\varDelta_{i,k}$, interval limits $\gamma_{\mathrm{b},i,k}$ and $\gamma_{\mathrm{e},i,k}$, and problem parameters $p$, subject to the proposed differential constraints~\eqref{rep_eq1}, boundary conditions~\eqref{eq: abstractbcs}, mesh conditions~\eqref{rep_eq3}, coupling conditions~\eqref{rep_eq2}, and algebraic conditions \eqref{dde:xi:e}-\eqref{dde:xi:b}. In this section, we derive the corresponding contributions to the necessary adjoint conditions for stationary points. Several examples of such optimization problems may be found in~\cite{ahsan2020optimization} and are revisited here in the context of the abstract framework.

By analogy with the $i$-th differential constraint \eqref{rep_eq1} and boundary conditions \eqref{eq: abstractbcs} (when applicable), consider the partial Lagrangian 
\begin{equation}
\int_{0}^{1}\lambda^{\text{T}}(\tau)\cdot \left(x'(\tau)-Tf\left(T_{0}+T\tau,x(\tau),y(\tau),p\right)\right)\,\mathrm{d}\tau +\zeta^\text{T}\cdot\left(x(0)-\sum_{j=1}^MB_{i,j}(p)x_j(1)\right) \label{eq: lagr1} 
\end{equation}
in terms of the Lagrange multipliers $\lambda:\mathbb{R}\rightarrow\mathbb{R}^n$ and $\zeta\in\mathbb{R}^n$. We obtain additive contributions to the necessary adjoint conditions by considering independent variations with respect to $x(\cdot)$, $x_j(1)$, $y(\cdot)$, $T_{0}$, $T$, and $p$, followed by identification of the coefficients of $\delta x(\cdot)$, $\delta x_j(1)$, $\delta y(\cdot)$, $\delta T_{0}$, $\delta T$, and $\delta p$, respectively. For example, using integration by parts, we obtain the contributions
\begin{equation}
\label{eq:diffadjx}
    -\lambda^{'\text{T}}(\tau)-T\lambda^\text{T}(\tau)\cdot f_{,x}\left(T_{0}+T\tau,x(\tau),y(\tau),p\right),
\end{equation}
corresponding to variations $\delta x(\tau)$ for $\tau\in(0,1)$, and 
\begin{equation}
\label{eq:diffadjx01}
\zeta^\text{T}-\lambda^\text{T}(0),\quad\lambda^\text{T}(1)
\end{equation}
for variations $\delta x(0)$ and $\delta x(1)$, respectively. Similarly, contributions corresponding to $\delta y(\tau)$ for $\tau\in[0,1]$, $\delta T_0$, $\delta T$, $\delta p$ and $\delta x_j(1)$ are given by
\begin{align}
    &-T\lambda^\text{T}(\tau)\cdot f_{,y}\left(T_{0}+T\tau,x(\tau),y(\tau),p\right),\label{eq:diffadjy}\\
    &-T\int_0^1\lambda^\text{T}(\tau)\cdot f_{,t}\left(T_{0}+T\tau,x(\tau),y(\tau),p\right)\,\mathrm{d}\tau,\label{eq:diffadjT0}\\
    &-\int_0^1\lambda^\text{T}(\tau)\cdot \left(\tau Tf_{,t}\left(T_{0}+T\tau,x(\tau),y(\tau),p\right)+f\left(T_{0}+T\tau,x(\tau),y(\tau),p\right)\right)\,\mathrm{d}\tau,\label{eq:diffadjT}\\
    &-T\int_0^1 \lambda^\text{T}(\tau)\cdot f_{,p}\left(T_{0}+T\tau,x(\tau),y(\tau),p\right)\,\mathrm{d}\tau-\zeta^\text{T}\cdot\sum_{j=1}^M\left(B_{i,j}(p)x_j(1)\right)_{,p},\label{eq:diffadjp}\\
    & -\zeta^\text{T}\cdot B_{i,j}(p),\label{eq:diffadjxj1}
\end{align}
respectively.

Next, by analogy with the mesh conditions in \eqref{rep_eq3} and $C_i$ coupling conditions in \eqref{rep_eq2}, consider the partial Lagrangian
\begin{align}
    \sum_{k=1}^C\int_{\gamma_{\mathrm{b},k}}^{\gamma_{\mathrm{e},k}}\mu^{\text{T}}(\tau)\cdot \left(y(\tau)-\sum_{s=1}^{S_k} A_{k,s}(p)\cdot x_{j_{k,s}}\left( \frac{T}{T_{j_{k}}}\left(\tau - \varDelta_{k}\right) \right)\right)\,\mathrm{d}\tau\label{eq: lagr2} +\eta_1\gamma_{\mathrm{b},1}+\sum_{k=1}^{C-1}\eta_{k+1}\left(\gamma_{\mathrm{b},k+1}-\gamma_{\mathrm{e},k }\right)+\eta_{C+1}\left(1-\gamma_{\mathrm{e},C}\right)
\end{align}
in terms of the Lagrange multipliers $\mu:\mathbb{R}\rightarrow\mathbb{R}^n$ and $\eta_k\in\mathbb{R}$ for $k=1,\ldots,C+1$. We obtain additive contributions to the necessary adjoint conditions by considering independent variations with respect to $y(\cdot)$, $p$, and $T$, as well as $T_{j_{k}}$, $\varDelta_{k}$, $x_{j_{k,s}}(\cdot)$, $\gamma_{\mathrm{b},k}$, and $\gamma_{\mathrm{e},k}$ for $s=1,\ldots,S_k$ and $k=1,\ldots,C$, followed by identification of the coefficients of $\delta y(\cdot)$, $\delta p$, and $\delta T$, as well as $\delta T_{j_{k}}$, $\delta \varDelta_{k}$, $\delta x_{j_{k,s}}(\cdot)$, $\delta \gamma_{\mathrm{b},k}$, and $\delta \gamma_{\mathrm{e},k}$, for $s=1,\ldots,S_k$ and $k=1,\ldots,C$, respectively. In order, these equal
\begin{align}
    &\mu^{\text{T}}(\tau),\,\tau\in\left[0,1\right],\label{eq:adjtcoupy}\\
    &-\sum_{k=1}^C\int_{\gamma_{\mathrm{b},k}}^{\gamma_{\mathrm{e},k}}\mu^{\text{T}}(\tau)\cdot \sum_{s=1}^{S_k} \left(A_{k,s}(p)\cdot x_{j_{k,s}}\left( \frac{T}{T_{j_{k}}}\left(\tau - \varDelta_{k}\right) \right)\right)_{,p}\,\mathrm{d}\tau,\label{eq:adjtcoupp}\\
        &-\sum_{k=1}^C\frac{1}{T_{j_{k}}}\int_{\gamma_{\mathrm{b},k}}^{\gamma_{\mathrm{e},k}}(\tau-\varDelta_{k})\mu^{\text{T}}(\tau)\cdot \sum_{s=1}^{S_k}A_{k,s}(p)\cdot x_{j_{k,s}}^{\prime}\left(\frac{T}{T_{j_{k}}}\left(\tau-\varDelta_{k}\right)\right)\,\mathrm{d}\tau, \label{eq:adjtcoupT}\\
        &\frac{T}{T_{j_{k}}^2}\int_{\gamma_{\mathrm{b},k}}^{\gamma_{\mathrm{e},k}}(\tau-\varDelta_{k})\mu^{\text{T}}(\tau)\cdot \sum_{s=1}^{S_k}A_{k,s}(p)\cdot x_{j_{k,s}}^{\prime}\left(\frac{T}{T_{j_{k}}}\left(\tau -\varDelta_{k}\right)\right)\,\mathrm{d}\tau,\,k=1,\ldots,C,\label{eq:adjtcoupTj}\\
    &\frac{T}{T_{j_{k}}}\int_{\gamma_{\mathrm{b},k}}^{\gamma_{\mathrm{e},k}}\mu^{\text{T}}(\tau)\cdot \sum_{s=1}^{S_k}A_{k,s}(p)\cdot x_{j_{k,s}}^{\prime}\left(\frac{T}{T_{j_{k}}}\left(\tau -\varDelta_{k}\right)\right)\,\mathrm{d}\tau,\,k=1,\ldots,C,\label{eq:adjtcoupvarDelta}\\
    &-\frac{T_{j_{k}}}{T} \mu^{\text{T}}\left(\frac{T_{j_{k}}}{T}\tau+\varDelta_{k}\right)\cdot A_{k,s}(p),\,\tau\in\left(\xi_{\mathrm{b},k},\xi_{\mathrm{e},k}\right),\,s=1,\ldots,S_k,\,k=1,\ldots,C\label{eq:adjtcoupxjks}
\end{align}
and, by the assumed continuity of $y(\tau)$, the sequences
\begin{equation}
\label{eq:adjt-mesh}
    \eta_1,\ldots,\eta_C,\quad\text{and}\quad -\eta_2,\ldots,-\eta_{C+1},
\end{equation}
respectively. When several subscripts $j_{k}$ evaluate to the same integer, the corresponding contributions may be added to each other to obtain the adjoint contributions associated with a particular differential state variable or duration. Consistent with the smoothness assumptions on $x$ and $y$, we assume that $\lambda$ and $\mu$ are continuous, piecewise-differentiable and continuous, respectively.

Finally, by analogy with the additional conditions \eqref{dde:xi:e} and \eqref{dde:xi:b} on the quantities $\xi_{\mathrm{b},k}$ and $\xi_{\mathrm{e},k}$, consider the partial Lagrangian
\begin{align}
\sum_{k=1}^{C-1}\chi_{\mathrm{e},k}\left(\frac{T}{T_{j_{k}}}\left(\gamma_{\mathrm{e},k}-\varDelta_{k}\right)-1\right)+\sum_{k=2}^{C}\chi_{\mathrm{b},k}\left(\frac{T}{T_{j_{k}}}\left(\gamma_{\mathrm{b},k}-\varDelta_{k}\right)\right),
\end{align}
in terms of the Lagrange multipliers $\chi_{\mathrm{e},k}$ for  $k=1,\ldots,C-1$ and $\chi_{\mathrm{b},k}$ for  $k=2,\ldots,C$. We obtain additive contributions to the necessary adjoint conditions by considering independent variations with respect to $T$, $T_{j_{k}}$ and $\varDelta_{k}$ for $k=1,\ldots,C$, $\gamma_{\mathrm{e},k}$ for $k=1,\ldots,C-1$, and $\gamma_{\mathrm{b},k}$ for $k=2,\ldots,C$. For example, identification of the coefficient of $\delta T$ yields the contribution
\begin{gather}
\label{eq:adjt-alg-T}
\sum_{k=1}^{C-1}\chi_{\mathrm{e},k}\frac{\gamma_{\mathrm{e},k}-\varDelta_k}{T_{j_{k}}}+\sum_{k=2}^{C}\chi_{\mathrm{b},k}\frac{\gamma_{\mathrm{b},k}-\varDelta_k}{T_{j_{k}}}.
\end{gather}
We find the contributions
\begin{align}
-\chi_{\mathrm{e},1}\frac{T}{T_{j_{1}}^2}\left(\gamma_{\mathrm{e},1}-\varDelta_k\right),\,-\chi_{\mathrm{e},1}\frac{T}{T_{j_1}}
\end{align}
by considering coefficients of $\delta T_{j_{1}}$ and $\delta\varDelta_{1}$,
\begin{equation}
    -\chi_{\mathrm{b},C}\frac{T}{T_{j_{C}}^2}\left(\gamma_{\mathrm{b},C}-\varDelta_C\right),\,-\chi_{\mathrm{b},C}\frac{T}{T_{j_C}}
\end{equation}
by considering coefficients of $\delta T_{j_{C}}$ and $\delta\varDelta_{C}$, and
\begin{align}
&-\left(\chi_{\mathrm{e},k}\gamma_{\mathrm{e},k}+\chi_{\mathrm{b},k}\gamma_{\mathrm{b},k}\right)\frac{T}{T_{j_{k}}^2},\,-\left(\chi_{\mathrm{e},k}+\chi_{\mathrm{b},k}\right)\frac{T}{T_{j_k}}
\end{align}
by considering coefficients of $\delta T_{j_{k}}$ and $\delta\varDelta_{k}$ for $k=2,\ldots,C-1$, respectively. As before contributions may be added to each other when several subscripts $j_{k}$ evaluate to the same integer.
Finally, identification of the coefficients of $\delta\gamma_{\mathrm{e},k}$ for $k=1,\ldots,C-1$ and $\delta\gamma_{\mathrm{b},k}$ for $k=2,\ldots,C$
yields the contributions
\begin{equation}
    \chi_{\mathrm{e},k}\frac{T}{T_{j_{k}}}
\end{equation}
and
\begin{equation}
\label{eq:adjt-alg-gammab}
   \chi_{\mathrm{b},k}\frac{T}{T_{j_{k}}},
\end{equation}
respectively.
%%%%%%%%%%
% Examples
%%%%%%%%%%

\subsubsection{Examples, continued}
\label{sec: toolbox template examples continued}
For the first example in Section~\ref{sec: toolbox template examples}, we obtain the adjoint contributions
\begin{align}
    &-\lambda^{'\text{T}}(\tau)-T\lambda^\text{T}(\tau)\cdot f_{,x}\left(T_0+T\tau,x(\tau),y(\tau),p\right)-\mu^\text{T}\left(\tau+\frac{\alpha}{T}\right),\,\tau\in\left(0,1-\frac{\alpha}{T}\right),\\
    &-\lambda^{'\text{T}}(\tau)-T\lambda^\text{T}(\tau)\cdot f_{,x}\left(T_0+T\tau,x(\tau),y(\tau),p\right)-\mu^\text{T}\left(\tau+\frac{\alpha}{T}-1\right),\,\tau\in\left(1-\frac{\alpha}{T},1\right),
\end{align}
corresponding to variations $\delta x(\tau)$ for $\tau\in(0,1)$, and
\begin{equation}
    -T\lambda^\text{T}(\tau)\cdot f_{,y}(T_0+T\tau,x(\tau),y(\tau),p)+\mu^\text{T}(\tau),\,\tau\in[0,1]
\end{equation}
corresponding to variations $\delta y(\tau)$ for $\tau\in[0,1]$. All three of these must equal zero in the absence of any additional constraints involving $x(\tau)$ on $(0,1)$ or $y(\tau)$ on $[0,1]$. In this case,
\begin{equation}
\lim_{\tau\uparrow 1-\frac{\alpha}{T}}\lambda^{'\text{T}}\left(\tau\right)-\lim_{\tau\downarrow 1-\frac{\alpha}{T}}\lambda^{'\text{T}}\left(\tau\right)=\mu^\text{T}(0)-\mu^\text{T}(1)=T\lambda^\text{T}(0)\cdot f_{,y}(T_0,x(0),y(0),p)-T\lambda^\text{T}(1)\cdot f_{,y}(T_0+T,x(1),y(1),p).
\end{equation}
Since by continuity $y(0)=y(1)$ and $x(0)=x(1)$, the right-hand side equals
\begin{equation}
    T\left(\lambda^\text{T}(0)-\lambda^\text{T}(1)\right)\cdot f_{,y}(T_0,x(0),y(0),p)
\end{equation}
in the case that $f$ is periodic in its first argument with period $T$. In general, this is nonzero and a discontinuity in the derivative of $\lambda$ occurs at $\tau=1-\alpha/T$.

For the second example in Section~\ref{sec: toolbox template examples}, we obtain the adjoint contributions
\begin{align}
    &-\lambda_1^{'\text{T}}(\tau)-\beta \lambda_1^\text{T}(\tau)\cdot f_{1,x}\left(T_0+\beta\tau,x_1(\tau),y_1(\tau),p\right)-\mu_1^\text{T}\left(\tau+\frac{\alpha}{\beta}\right),\,\tau\in\left(0,1-\frac{\alpha}{\beta}\right),\\    &-\lambda_1^{'\text{T}}(\tau)-\beta \lambda_1^\text{T}(\tau)\cdot f_{1,x}\left(T_0+\beta\tau,x_1(\tau),y_1(\tau),p\right)-\frac{\beta}{T-\beta}\mu_2^\text{T}\left(\frac{\beta}{T-\beta}\tau+\frac{\alpha-\beta}{T-\beta}\right),\,\tau\in\left(1-\frac{\alpha}{\beta},1\right),\\
    &-\lambda_2^{'\text{T}}(\tau)-(T-\beta) \lambda_2^\text{T}(\tau)\cdot f_{2,x}\left(T_0+\beta+(T-\beta)\tau,x_2(\tau),y_2(\tau),p\right)-\mu_2^\text{T}\left(\tau+\frac{\alpha}{T-\beta}\right),\,\tau\in\left(0,1-\frac{\alpha}{T-\beta}\right),\\
    &-\lambda_2^{'\text{T}}(\tau)-(T-\beta) \lambda_2^\text{T}(\tau)\cdot f_{2,x}\left(T_0+\beta+(T-\beta)\tau,x_2(\tau),y_2(\tau),p\right)-\frac{T-\beta}{\beta}\mu_1^\text{T}\left(\frac{T-\beta}{\beta}\tau+\frac{\alpha-T+\beta}{\beta}\right),\,\tau\in\left(1-\frac{\alpha}{T-\beta},1\right),
\end{align}
corresponding to variations $\delta x_1(\tau)$ and $\delta x_2(\tau)$ for $\tau\in(0,1)$, and
\begin{align}
    -\beta\lambda_1^\text{T}(\tau)\cdot f_{1,y}(T_0+\beta\tau,x_1(\tau),y_1(\tau),p)+\mu_1^\text{T}(\tau),&\,\tau\in[0,1],\\
    -(T-\beta)\lambda_2^\text{T}(\tau)\cdot f_{2,y}(T_0+\beta+(T-\beta)\tau,x_2(\tau),y_2(\tau),p)+\mu_2^\text{T}(\tau),&\,\tau\in[0,1],
\end{align}
corresponding to variations $\delta y_1(\tau)$ and $\delta y_2(\tau)$ for $\tau\in[0,1]$. All six of these must equal zero in the absence of any additional constraints involving $x_1(\tau)$ or $x_2(\tau)$ on $(0,1)$ or $y_1(\tau)$ or $y_2(\tau)$ on $[0,1]$.

We encourage the reader to apply the general expressions in the previous section to derive the adjoint contributions for the third example in Section~\ref{sec: toolbox template examples} and to compare the resulting necessary conditions with the adjoint equations derived in~\cite{ahsan2020optimization} from the
partial differential equation~\eqref{eq: pde_quasi}, coupling conditions~\eqref{eq_coupling_quasi}, and boundary conditions $x(\varphi,0)=x(\varphi-2\pi\rho,1)$.

%%%%%%%%%%%%%%%%
% Discretization
%%%%%%%%%%%%%%%%
\subsection{Discretization}
\label{sec: discretization} 

We proceed to describe a natural discretization of the governing zero problem and adjoint contributions, first as an identity in spaces of piecewise polynomials and then in terms of the resulting large systems of algebraic equations with unknowns in $\R^{\hat{N}}$ with $\hat{N}\gg1$. Given such a discretization, we may formulate Jacobians with respect to the set of discrete unknowns, but omit their explicit expressions in this text.

\subsubsection{Abstract formulation of collocation discretization}
\label{sec: discretization abstract} 

Given a partition $\tau_{\mathrm{pt}}$ with $0=\tau_{\mathrm{pt},1}<\ldots<\tau_{\mathrm{pt},N+1}=1$ of the interval $[0,1]$, two spaces of, potentially discontinuous, piecewise-polynomial functions are relevant to our discussion, namely
\begin{align}
{\cal P}_{m\phantom{,0}}&:=\left\{x:[0,1]\to\R^n:x\vert_{[\tau_{\mathrm{pt},j},\tau_{\mathrm{pt},j+1}]}\mbox{\ is polynomial of degree $m$ for $j=1,\ldots,N$}\right\},
\end{align}
and the subspace ${\cal P}_{m,0}:=\left\{x\in{\cal P}_m: x(0)=0\right\}$. In particular, ${\cal P}_{0,0}$ is the space of piecewise-constant functions that equal $0$ at $t=0$. Elements of ${\cal P}_m$ are permitted to be discontinuous and multivalued
on the interior partition points $\{\tau_{\mathrm{pt},j}\}_{j=2}^N$ %% $\{t_{\mathrm{pt},j}\}_{j=2}^N$ %% 
and we use the notation $x(\tau_{\mathrm{pt},j}^\pm)$ to distinguish between left- and right-sided limits.  The spaces depend on the partition $\tau_\mathrm{pt}$ (of length $N$), the degree $m$ and
space dimension $n$. Specifically, $\dim{\cal P}_m=nN(m+1)$, while $\dim{\cal P}_{m,0}=nN(m+1)-n$ (we use the space ${\cal P}_{0,0}$ of piecewise constant functions with $x(0)=0$, which has $\dim{\cal P}_{0,0}=n(N-1)$). 
We observe that differentiation maps ${\cal P}_m$ into ${\cal P}_{m-1}$.

Our proposed discretization (consistent with the approach in the \textsc{coco}-compatible \mcode{'coll'} toolbox) is expressed in terms of four projections into the spaces
${\cal P}_m$ and ${\cal P}_{0,0}$, using interpolation at either the collection of Gauss-Legendre points of degree $m-1$ on each subinterval, $\{\tau_{\mathrm{cn},j}\}_{j=1}^{Nm}$, a mesh of $N(m+1)$ base points, $\{\tau_{\mathrm{bp},j}\}_{j=1}^{N(m+1)}$, or the interior partition points $\{\tau_{\mathrm{pt},j}\}_{j=2}^N$:
\begin{align}
\begin{aligned}
  \pcn&:(x:[0,1]\to\R^n)\mapsto \tilde{x}\in{\cal P}_{m-1}
  \mbox{\ with $\tilde{x}(\tau_{\mathrm{cn},j})=x(\tau_{\mathrm{cn},j})$ for all $j=1,\ldots,Nm$.}\\
  \pbp&:(x:[0,1]\to\R^n)\mapsto \tilde{x}\in{\cal P}_{m\phantom{+1}}
  \mbox{\ with $\tilde{x}(\tau_{\mathrm{bp},j})=x(\tau_{\mathrm{bp},j})$ for all $j=1,\ldots,N(m+1)$,}\\
  \pcont^\pm&:(x:[0,1]\to\R^n)\mapsto \tilde{x}\in{\cal P}_{0,0}
  \mbox{\ \ \ with $\tilde{x}(\tau_{\mathrm{pt},j}^+)=x(\tau_{\mathrm{pt},j}^\pm)$ for all $j=2,\ldots,N$.}\\
\end{aligned}
\label{disc:proj}
\end{align}
Note the use of the right limit $\tilde{x}(\tau_{\mathrm{pt},j}^+)$ for the result $\tilde{x}$ for both projections $\pcont^\pm$ (which will enforce continuity of the solution in \eqref{disc:eq3:poly} below). All projections depend on the partition $\tau_\mathrm{pt}$, collection of base points $\tau_\mathrm{bp}$, and space dimension $n$ (without indicating these dependencies as subscripts). They can be applied to functions that are continuous on each partition interval and have well-defined left and right limits.

For the $i$-th segment (omitting the $i$ subscript), let $S:=\sum_{k=1}^{C} S_{k}$ and $J_{k}:=(\gamma_{\mathrm{b},k},\gamma_{\mathrm{e},k})$. Recall the differential constraint \eqref{rep_eq1}
\begin{equation}
\label{disc:eq1:inf}
    x'(\tau)=Tf(T_{0}+T\tau,x(\tau),y(\tau),p),\,\tau\in(0,1)
\end{equation}
and coupling conditions \eqref{rep_eq2}, which are of the form
\begin{equation}
    y(\tau)=\sum_{\ell=1}^{S}a_\ell(\tau)z_\ell\left(b_\ell\tau-c_\ell\right),\,\tau\in[0,1]\mbox{,}\label{disc:eq2:inf}
\end{equation}
where
\begin{equation}
    a_\ell(\tau):=A_{k(\ell),s(\ell)}(p)\mathds{1}_{\textstyle J_{k(\ell)}}(\tau),\,z_\ell(\tau):=x_{\textstyle j_{k(\ell),s(\ell)}}(\tau),\,b_\ell:=\frac{T}{T_{j_{k}}},\,
    c_\ell:=\frac{T}{T_{j_{k}}}\Delta_{k}\label{disc:eq:abcdef}
\end{equation}
and
\begin{equation}
    k(\ell):=\min\left\{\nu:\sum_{j=1}^\nu S_j\geq\ell\right\},\,
    s(\ell):=\ell-\sum_{j=1}^{k-1}S_j\mbox{,}
\end{equation}
in terms of the unknown differential and algebraic state variables $x(\cdot)$ and $y(\cdot)$.

We included the indicator function $\mathds{1}_{\textstyle J_{k(\ell)}}(\tau)$ in the definition of $a_\ell(\tau)$ to make explicit that the functions $x_{\textstyle j_{k(\ell),s(\ell)}}$ are only evaluated on the subinterval $J_{k(\ell)}$. The
unknowns $x$ and $y$ and $z_\ell$ are continuous (or more regular) functions on the interval $[0,1]$. Our chosen method of discretization looks for functions $\tilde{x},\tilde{y}\in{\cal P}_m$, coupled to $\tilde{z}_\ell\in{\cal P}_m$ from possibly other segments, that satisfy the finite-dimensional constraints
\begin{align}
  \tilde{x}'&=\pcn Tf\left(T_0+(\cdot)T,\tilde{x}(\cdot),\tilde{y}(\cdot),p\right)&&\mbox{discretized ODE in ${\cal P}_{m-1}$, dimension $nNm$}\label{disc:eq1:poly}\\
  \tilde{y}&=\pbp\sum_{\ell=1}^S a_\ell(\cdot)\tilde{z}_\ell((\cdot)b_\ell-c_\ell)
  &&\mbox{discretized algebraic constraint in ${\cal P}_m$, dimension $nN(m+1)$}\label{disc:eq2:poly}\\
  \pcont^-\,\tilde x&=\pcont^+\,\tilde x
    &&\mbox{zero gaps for $\tilde x$ in ${\cal P}_{0,0}$, dimension $n(N-1)$.}\label{disc:eq3:poly}
\end{align}
Equation \eqref{disc:eq1:poly} is an identity between discontinuous piecewise polynomials of degree $m-1$, while \eqref{disc:eq2:poly} is an identity between discontinuous piecewise polynomials of degree $m$, and \eqref{disc:eq3:poly} is an identity between piecewise constant functions that are $0$ in $t=0$. The equations \eqref{disc:eq1:poly}--\eqref{disc:eq3:poly} have a dimensional deficit $n$ (considering $\tilde x$ and $\tilde y$ as the variables). Evaluation of $f$, the time shift and scaling $\tilde z_\ell\mapsto \tilde z_\ell((\cdot)b_\ell-c_\ell)$ and multiplication by the piecewise continuous function $a_\ell(\cdot)$ are exact such that an approximation is only performed when the respective projections $\pcn$ and $\pbp$ are applied in \eqref{disc:eq1:poly} and \eqref{disc:eq2:poly}. 

In lieu of the Lagrange multipliers, we seek functions $\tilde\lambda,\tilde\mu\in{\cal P}_m$. The equations resulting from vanishing variations with respect to $x$ are projected by $\pcn$, giving piecewise polynomial identities in ${\cal P}_{m-1}$, while equations from vanishing variations with respect to $y$  are projected by $\pbp$, giving piecewise polynomial identities in ${\cal P}_m$. Thus, continuity has to be enforced only for $\tilde\lambda$ (by imposing $\pcont^-\,\tilde\lambda=\pcont^+\,\tilde\lambda$). Integrals over subintervals $J$ of $[0,1]$ occurring in finite-dimensional adjoint contributions, such as \eqref{eq:adjtcoupTj} and \eqref{eq:adjtcoupvarDelta}, are approximated using $\pcn$ on the expression truncated by the indicator function. Specifically, for arbitrary $g\in C([0,1];\R^j)$, the integral $\int_{J} g(\tau)\d \tau$ is approximated as $\int_0^1\pcn\left[\mathds{1}_{J}(\cdot)g(\cdot)\right](\tau)\d \tau$. 

The polynomial identities \eqref{disc:eq1:poly}-- \eqref{disc:eq3:poly} and the corresponding contributions to the adjoint conditions have to be evaluated for all segments $i=1,\ldots,M$. The identities in the spaces of piecewise polynomials, ${\cal P}_m$, ${\cal P}_{m-1}$ and ${\cal P}_{0,0}$, are reduced to algebraic equations in $\R^{nN(m+1)}$, $\R^{nNm}$ and $\R^{n(N-1)}$ by evaluating them on the meshes $\tau_\mathrm{bp}$, $\tau_\mathrm{cn}$ and the interior partition points $\{\tau_{\mathrm{pt},j}^+\}_{j=2}^N$, respectively.

Before proceeding to consider a detailed implementation of the abstract discretization scheme presented in this section, we conclude with a comment on convergence analysis. Indeed, it is notable that rigorous convergence analysis---showing that solutions of \eqref{disc:eq1:poly}--\eqref{disc:eq3:poly} converge to solutions of \eqref{disc:eq1:inf} and \eqref{disc:eq2:inf} under appropriate regularity assumptions---is an open problem. Even for autonomous single-segment periodic boundary-value problems with delay, a complete convergence proof has only been presented recently~\cite{andoSIAM2020}. The difficulty with this analysis is the non-differentiability of the discretized nonlinear system with respect to many unknowns away from the solution (for example, the period $T$, the delay $\alpha$ and the quantities $b_\ell,c_\ell$ in \eqref{disc:eq:abcdef}). Similar concerns apply to the convergence analysis for the Lagrange multipliers. In our framework, we use discretizations of the adjoints of the infinite-dimensional problem \eqref{disc:eq1:inf}-\eqref{disc:eq2:inf}, not the adjoints of the discretized problem \eqref{disc:eq1:poly}--\eqref{disc:eq3:poly} (which are different and possibly not well defined away from the solution manifold; the same concern applies to Jacobians of all equations).

\subsubsection{Implementation as large systems of algebraic equations}
\label{sec: discretization implementation}
We now describe a \textsc{coco}-implementable form of the zero problems and contributions to adjoint conditions derived in Sections~\ref{sec: toolbox template zero problems} and \ref{sec: toolbox template adjoint conditions}, and their discretization as polynomials described in Section~\ref{sec: discretization abstract}. As in the abstract discussion above, we consider an arbitrary segment $i$, but omit the index $i$ in the description below.  Where possible, we rely on vectorized notation to suppress a jungle of indices. For example, we use the $\mathfrak{vec}$ operator to convert its argument into a one-dimensional array of scalars. If $A$ is an array of possible values for the argument of a function $f$, then $f(A)$ is an array of the same size as $A$ of values of $f$ applied to each element of $A$.

Consistent with the abstract discussion and following~\cite{dankowicz2013recipes}, let $N$ and $m$ be two positive integers and define the uniform partition $\tau_{\mathrm{pt},j}=(j-1)/N$, for $j=1,\ldots,N+1$, and time sequence
\begin{equation}
    \tau_{\mathrm{bp},(m+1)(j-1)+k}=\tau_{\mathrm{pt},j}+\frac{k-1}{Nm},\,j=1,\ldots,N,\,k=1,\ldots,m+1.
\end{equation}
In particular, $\tau_{\mathrm{bp},1}=\tau_{\mathrm{pt},1}=0$, $\tau_{\mathrm{bp},N(m+1)}=\tau_{\mathrm{pt},N}+1/N=\tau_{\mathrm{pt},N+1}=1$, and
\begin{equation}
\label{eq:meshcont}
    \tau_{\mathrm{bp},(m+1)j+1}=\tau_{\mathrm{pt},j+1}=\tau_{\mathrm{bp},(m+1)(j-1)+m+1}=\tau_{\mathrm{bp},(m+1)j},\,j=1,\ldots,N-1.
\end{equation}
We represent the discretized differential state variable $\tilde x(\cdot)$ as the one-dimensional array $x_{\mathrm{bp}}:=\tilde x\left(\tau_{\mathrm{bp}}\right)$ of $N(m+1)$  vectors in $\mathbb{R}^n$ and proceed, similarly, to represent the discretized algebraic state variable $y(\cdot)$ and Lagrange multipliers $\lambda(\cdot)$ and $\mu(\cdot)$ as the one-dimensional arrays $y_{\mathrm{bp}}:=\tilde y\left(\tau_{\mathrm{bp}}\right)$, $\lambda_{\mathrm{bp}}:=\tilde\lambda\left(\tau_{\mathrm{bp}}\right)$, and $\mu_{\mathrm{bp}}:=\tilde \mu\left(\tau_{\mathrm{bp}}\right)$, respectively, of $N(m+1)$ vectors in $\mathbb{R}^n$ each. In this representation the continuity equation \eqref{disc:eq3:poly} (and its correspondent for $\tilde\lambda$) takes the form of $2(N-1)n$ continuity conditions
\begin{equation}
    x_{\mathrm{bp}}\cdot \mathcal{C}=\lambda_{\mathrm{bp}}\cdot \mathcal{C}=0\mbox{.}
\end{equation}
Since the base points include the partition points according to \eqref{eq:meshcont}, the $j$-th column of the $N(m+1)\times(N-1)$ matrix $\mathcal{C}$ equals $\mathbf{e}_{(m+1)j+1}-\mathbf{e}_{(m+1)j}$ in the standard basis of $\mathbb{R}^{N(m+1)}$. Equivalently, in vectorized form, we write $Q\cdot \mathfrak{vec}\left(x_{\mathrm{bp}}\right)=Q\cdot \mathfrak{vec}\left(\lambda_{\mathrm{bp}}\right)=0$ with the $(N-1)n\times N(m+1)n$ matrix $Q=\mathcal{C}^\text{T}\otimes I_n$. Since the discretized coupling conditions \eqref{disc:eq2:poly} and their adjoint contributions with respect to variations $\delta y(\cdot)$ are imposed in ${\cal P}_m$, there is no need to impose explicit continuity conditions for $y_{\mathrm{bp}}$ and $\mu_{\mathrm{bp}}$.

It is convenient to define an index function that maps times $\tau$ to the corresponding subinterval indices for the partition of ${\cal P}_m$. This is useful, for example, when applying the time scale and -shift operator $\tilde x(\cdot)\mapsto \tilde x((\cdot)b-c)$ on piecewise polynomials in ${\cal P}_m$ in the discretized coupling conditions in \eqref{disc:eq2:poly}. To this end, given a partition $\sigma$ of $[0,1]$ into $P$ intervals, let $\tau\mapsto\iota\left(\tau;\sigma\right)$ denote the linear interpolation of the pairing $\sigma\rightarrow\{1,\ldots,P+1\}$ at $\tau$. In terms of the floor function $\lfloor\cdot\rfloor$,
\begin{equation}
    \tau\in\left[\sigma_{p},\sigma_{p+1}\right)\Rightarrow \pi(\tau;\sigma):=\left\lfloor\iota\left(\tau;\sigma\right)\right\rfloor=p
\end{equation} for any $p=1,\ldots,P$. Then, with $\pi(\tau;\sigma)=1$ for $\tau<0$ and $\pi(\tau;\sigma)=P$ for $\tau\ge 1$, the piecewise-constant, non-decreasing index function $\pi(\cdot;\sigma)$ maps $(-\infty,\infty)$ to $\{1,\ldots,P\}$. Consider, for example, the sequence $\{\mathcal{L}_l\}_{l=1}^{m+1}$ of $m$-th degree Lagrange polynomials defined on the uniform partition of $[-1,1]$, such that
\begin{equation}
    \mathcal{L}_l\left(-1+2\frac{k-1}{m}\right)=\delta_{l,k},\,k=1,\ldots,m+1.
\end{equation}
Then, since
\begin{equation}
    \lim_{\tau\rightarrow\tau_{\mathrm{pt},j}-}\mathcal{L}_l\left(2N\tau+1-2\pi(\tau;\tau_{\mathrm{pt}})\right)=\delta_{l,m+1},\,\lim_{\tau\rightarrow\tau_{\mathrm{pt},j}+}\mathcal{L}_l\left(2N\tau+1-2\pi(\tau;\tau_{\mathrm{pt}})\right)=\delta_{l,1},
\end{equation}
the continuous, piecewise-polynomial interpolants given by
\begin{align}
\label{eq:interpolants}
    \tau\mapsto\sum_{l=1}^{m+1}\mathcal{L}_l\left(2N\tau+1-2\pi(\tau;\tau_{\mathrm{pt}})\right)x_{\mathrm{bp},(m+1)(\pi(\tau;\tau_{\mathrm{pt}})-1)+l},\,
    \tau\mapsto\sum_{l=1}^{m+1}\mathcal{L}_l\left(2N\tau+1-2\pi(\tau;\tau_{\mathrm{pt}})\right)y_{\mathrm{bp},(m+1)(\pi(\tau;\tau_{\mathrm{pt}})-1)+l},\\
    \tau\mapsto\sum_{l=1}^{m+1}\mathcal{L}_l\left(2N\tau+1-2\pi(\tau;\tau_{\mathrm{pt}})\right)\lambda_{\mathrm{bp},(m+1)(\pi(\tau;\tau_{\mathrm{pt}})-1)+l},\,
    \tau\mapsto\sum_{l=1}^{m+1}\mathcal{L}_l\left(2N\tau+1-2\pi(\tau;\tau_{\mathrm{pt}})\right)\mu_{\mathrm{bp},(m+1)(\pi(\tau;\tau_{\mathrm{pt}})-1)+l}
\end{align}
allow us to evaluate the piecewise polynomials $\tilde x(\tau)$, $\tilde y(\tau)$, $\tilde \lambda(\tau)$, and $\tilde \mu(\tau)$ and, as appropriate, their derivatives at arbitrary $\tau$ in $[0,1]$ in terms of linear combinations of the elements of $x_{\mathrm{bp}}$, $y_{\mathrm{bp}}$, $\lambda_{\mathrm{bp}}$, and $\mu_{\mathrm{bp}}$, respectively. As a special case, let $z$ denote the one-dimensional array of $m$-th order Gauss-Legendre quadrature nodes on the interval $[-1,1]$ in increasing order, and define the time sequence
\begin{equation}
\label{eq:tau_cn}
    \tau_{\mathrm{cn},m(j-1)+k}=\tau_{\mathrm{pt},j}+\frac{1+z_k}{2N},\,j=1,\ldots,N,\,k=1,\ldots,m,
\end{equation}
such that $2N\tau_{\mathrm{cn}}-2\pi\left(\tau_{\mathrm{cn}};\tau_{\mathrm{pt}}\right)=J_{N,1}\otimes (z-1)$. Then,
\begin{gather}
    \tilde{x}\left(\tau_{\mathrm{cn}}\right)= x_{\mathrm{cn}}:=x_{\mathrm{bp}}\cdot \mathcal{L}_{\mathrm{cn}},\,\tilde{x}'\left(\tau_{\mathrm{cn}}\right)= x'_{\mathrm{cn}}:=2Nx_{\mathrm{bp}}\cdot\mathcal{L}'_{\mathrm{cn}},\,\tilde{y}\left(\tau_{\mathrm{cn}}\right)= y_{\mathrm{cn}}:=y_{\mathrm{bp}}\cdot \mathcal{L}_{\mathrm{cn}},\\
    \tilde{\lambda}\left(\tau_{\mathrm{cn}}\right)= \lambda_{\mathrm{cn}}:=\lambda_{\mathrm{bp}}\cdot \mathcal{L}_{\mathrm{cn}},\,\tilde{\lambda}'\left(\tau_{\mathrm{cn}}\right)= \lambda'_{\mathrm{cn}}:=2N\lambda_{\mathrm{bp}}\cdot\mathcal{L}'_{\mathrm{cn}},\,\tilde{\mu}\left(\tau_{\mathrm{cn}}\right)= \mu_{\mathrm{cn}}:=\mu_{\mathrm{bp}}\cdot \mathcal{L}_{\mathrm{cn}},
\end{gather}
where
\begin{equation}
    \mathcal{L}_{\mathrm{cn},(m+1)(a-1)+c,m(b-1)+d}=\delta_{a,b}\mathcal{L}_c(z_d),\,\mathcal{L}'_{\mathrm{cn},(m+1)(a-1)+c,m(d-1)+d}=\delta_{a,b}\mathcal{L}'_c(z_d)
\end{equation}
for $a,b=1,\ldots,N$, $c=1,\ldots,m+1$, and $d=1,\ldots,m$. Equivalently, we write $\mathfrak{vec}\left(x_{\mathrm{cn}}\right)=W\cdot \mathfrak{vec}\left(x_{\mathrm{bp}}\right)$, $\mathfrak{vec}\left(x_{\mathrm{cn}}'\right)=2NW'\cdot \mathfrak{vec}\left(x_{\mathrm{bp}}\right)$,  $\mathfrak{vec}\left(y_{\mathrm{cn}}\right):=W\cdot \mathfrak{vec}\left(y_{\mathrm{bp}}\right)$, $\mathfrak{vec}\left(\lambda_{\mathrm{cn}}\right)=W\cdot \mathfrak{vec}\left(\lambda_{\mathrm{bp}}\right)$, $\mathfrak{vec}\left(\lambda_{\mathrm{cn}}'\right)=2NW'\cdot \mathfrak{vec}\left(\lambda_{\mathrm{bp}}\right)$, and $\mathfrak{vec}\left(\mu_{\mathrm{cn}}\right)=W\cdot \mathfrak{vec}\left(\mu_{\mathrm{bp}}\right)$ in terms of the $Nmn\times N(m+1)n$ matrices $W= \mathcal{L}^\text{T}_{\mathrm{cn}}\otimes I_n$ and $W'= \mathcal{L}'^{\text{T}}_{\mathrm{cn}}\otimes I_n$. 

We use evaluation at $\tau_{\mathrm{cn}}$ to represent the discretization \eqref{disc:eq1:poly} of  differential constraint \eqref{rep_eq1} as the vectorized algebraic constraints
\begin{equation}
    0=\mathfrak{vec}\left(x'_{\mathrm{cn}}-Tf_{\mathrm{cn}}\right)=2NW'\cdot\mathfrak{vec}\left(x_{\mathrm{bp}}\right)-T\mathfrak{vec}\left(f_{\mathrm{cn}}\right),
\end{equation}
where $f_{\mathrm{cn}}:=f(T_0+T\tau_{\mathrm{cn}},x_{\mathrm{cn}},y_{\mathrm{cn}},p)$ is a one-dimensional array of $Nm$ vectors in $\mathbb{R}^n$. By a similar use of notation and evaluation at $\tau_{\mathrm{cn}}$, we represent the   discretizations of the corresponding contributions in \eqref{eq:diffadjx} to the necessary adjoint conditions associated with variations $\delta x(\tau)$ for $\tau\in(0,1)$ as the vectorized algebraic expression
\begin{equation}
    -\mathfrak{vec}\left(\lambda'_{\mathrm{cn}}\right)-T \mathfrak{diag}\left(f_{x,\mathrm{cn}}\right)^\text{T}\cdot \mathfrak{vec}\left(\lambda_{\mathrm{cn}}\right)=\left(-2NW'-T\mathfrak{diag}\left(f_{x,\mathrm{cn}}\right)^\text{T}\cdot W\right)\cdot \mathfrak{vec}\left(\lambda_{\mathrm{bp}}\right),
\end{equation}
where $f_{x,\mathrm{cn}}:=f_{,x}(T_0+T\tau_{\mathrm{cn}},x_{\mathrm{cn}},y_{\mathrm{cn}},p)$ is a one-dimensional array of $Nm$ matrices in $\mathbb{R}^{n\times n}$ and
$\mathfrak{diag}\left(f_{x,\mathrm{cn}}\right)$ is an $Nmn\times Nmn$ block-diagonal matrix with the elements of $f_{x,\mathrm{cn}}$ along the diagonal. The discretized contributions from the Lagrangian~\eqref{eq: lagr1} associated with variations $\delta x(0)$ and $\delta x(1)$ are similarly given by $ \zeta-\begin{pmatrix}I_{n} & 0 & \cdots & 0\end{pmatrix}\cdot\mathfrak{vec}\left(\lambda_{\mathrm{bp}}\right)$ and $\begin{pmatrix}0 & \cdots & 0 & I_{n} \end{pmatrix}\cdot\mathfrak{vec}\left(\lambda_{\mathrm{bp}}\right)$, while that corresponding to $\delta x_j(1)$ equals $-\left(B_{i,j}(p)\right)^\text{T}\cdot \zeta$. By replacing the $x$ subscript with $y$ and $\mathrm{cn}$ with $\mathrm{bp}$, evaluation at $\tau_{\mathrm{bp}}$ yields the discretization of the contributions in \eqref{eq:diffadjy} to the necessary adjoint condition associated with variations $\delta y(\tau)$ for $\tau\in[0,1]$ in terms of the vectorized algebraic expression
\begin{gather}
    -T\mathfrak{diag}\left(f_{y,\mathrm{bp}}\right)^\text{T}\cdot \mathfrak{vec}\left(\lambda_{\mathrm{bp}}\right).
\end{gather}
Finally, Gaussian quadrature on the partition $\tau_{\mathrm{pt}}$ with collocation nodes at $\tau_{\mathrm{cn}}$ and with $t$ or $p$ subscripts in place of $x$ or $y$ yields the discretization of the contributions to the necessary adjoint conditions in \eqref{eq:diffadjT0}, \eqref{eq:diffadjT}, and \eqref{eq:diffadjp} associated with variations $\delta T_0$, $\delta T$, and $\delta p$,  given in order by the vectorized algebraic expressions
\begin{gather}
    -\frac{T}{2N}\,\mathfrak{vec}\left(f_{t,\mathrm{cn}}\right)^\text{T}\cdot\Omega_{\mathrm{cn}}\cdot W\cdot\mathfrak{vec}\left(\lambda_{\mathrm{bp}}\right)\\
    -\frac{1}{2N}\left(\mathfrak{vec}\left(f_{\mathrm{cn}}\right)+T\left(\mathfrak{diag}\left(\tau_{\mathrm{cn}}\right)\otimes I_n\right)\cdot \mathfrak{vec}\left(f_{t,\mathrm{cn}}\right)\right)^\text{T}\cdot\Omega_{\mathrm{cn}}\cdot W\cdot\mathfrak{vec}\left(\lambda_{\mathrm{bp}}\right),\\
    -\frac{T}{2N}\mathfrak{transp}\left(f_{p,\mathrm{cn}}\right)^\text{T}\cdot\Omega_{\mathrm{cn}}\cdot W\cdot\mathfrak{vec}\left(\lambda_{\mathrm{bp}}\right)-\sum_{j=1}^{M}\left(\mathfrak{vec}\left(B_{i,j}(p)\right)_{,p}\right)^\text{T}\cdot \left(\begin{pmatrix}0 & \cdots & 0 & I_n\end{pmatrix}\cdot \mathfrak{vec}\left(x_{j,\mathrm{bp}}\right)\otimes I_n\right)\cdot\zeta
\end{gather}
in terms of the $Nmn\times Nmn$ matrix $\Omega_{\mathrm{cn}}=I_N\otimes\mathfrak{diag}(\omega)\otimes I_n$, the one-dimensional array $\omega$ of $m$-th order quadrature weights associated with the nodes $z_1,\ldots,z_m$, and the $Nmn\times q$ block-vertical matrix $\mathfrak{transp}\left(f_{p,\mathrm{cn}}\right)$ with the elements of $f_{p,\mathrm{cn}}$ stacked vertically.

For the discretization \eqref{disc:eq2:poly} of the coupling conditions in \eqref{rep_eq2} and the corresponding adjoint contributions in \eqref{eq:adjtcoupp}-\eqref{eq:adjtcoupxjks}, it is necessary to consider interpolation at time instants defined by the arguments of the  differential state variables $x_{j_{k,s}}$ and Lagrange multiplier $\mu$ for $\tau\in\tau_{\mathrm{bp}}$ or $\tau_{\mathrm{cn}}$. In contrast to interpolation at $\tau_{\mathrm{cn}}$, the corresponding interpolation matrices by necessity depend on the durations $T$ and $T_{j_{k}}$, coupling delays $\varDelta_{k}$, and mesh limits $\gamma_{\mathrm{b},k}$ and $\gamma_{\mathrm{e},k}$. For example, let $k_{\mathrm{bp}}=\pi(\tau_{\mathrm{bp}};\{0,\gamma_{\mathrm{b},2}\ldots,\gamma_{\mathrm{b},C},1\})$ as a non-decreasing sequence of coupling interval indices associated with the elements of the $\tau_{\mathrm{bp}}$ array. For each subarray of successive elements $\tau_{\mathrm{bp}}^k$ that share an interval index $k\in k_{\mathrm{bp}}$, associate the shifted time instants
\begin{equation}
    \tau_{\mathrm{bp}\downarrow \mathrm{sh}}^{k}:=\frac{T}{T_{j_{k}}}\left(\tau_{\mathrm{bp}}^k-J_{|\tau_{\mathrm{bp}}^k|,1}\otimes\varDelta_{k}\right)
\end{equation}
with the non-decreasing sequence $j^{k}_{\mathrm{bp}\downarrow \mathrm{sh}}=\pi\left(\tau_{\mathrm{bp}\downarrow \mathrm{sh}}^{k},\tau_{\mathrm{pt}}\right)$ of interval indices $j$. For each such $k$, the coupling conditions \eqref{rep_eq2} discretized at $\tau_{\mathrm{bp}}^k$ then take the form
\begin{equation}
    y_{\mathrm{bp}}^k -\sum_{s=1}^{S_k}A_{k,s}(p)\cdot x_{j_{k,s},\mathrm{bp}}\cdot\mathcal{L}_{\mathrm{bp}\downarrow \mathrm{sh}}^{k}=0,
\end{equation}
where $y_{\mathrm{bp}}^k$ denotes the corresponding elements of $y_{\mathrm{bp}}$ and
\begin{equation}
    \mathcal{L}^{k}_{\mathrm{bp}\downarrow \mathrm{sh},(m+1)(a-1)+c,b}=\delta_{a,j^{k}_{\mathrm{bp}\downarrow  sh,b}}\mathcal{L}_c\left(2N\tau^{k}_{\mathrm{bp}\downarrow \mathrm{sh},b}+1-2j^{k}_{\mathrm{bp}\downarrow  sh,b}\right),
\end{equation}
for $a=1,\ldots,N$, $c=1,\ldots,m+1$, and $b=1,\ldots,|\tau_{\mathrm{bp}}^k|$. Equivalently, in vectorized form,
\begin{equation}
    \mathfrak{vec}\left(y_{\mathrm{bp}}^k\right)-\sum_{s=1}^{S_k}\left(\left(\mathcal{L}^{k}_{\mathrm{bp}\downarrow \mathrm{sh}}\right)^\text{T}\otimes A_{k,s}(p)\right)\cdot\mathfrak{vec}\left(x_{j_{k,s},\mathrm{bp}}\right)=0.
\end{equation}

For the adjoint contributions in \eqref{eq:adjtcoupxjks} associated with variations $\delta x_{j_{k,s}}(\cdot)$ and discretized at $\tau_{\mathrm{cn}}$, we obtain nonzero contributions only on the subset $\tau_{\mathrm{cn}}^{k}\subseteq\tau_{\mathrm{cn}}$ of time instances in $\left(\xi_{\mathrm{b},k},\xi_{\mathrm{e},k}\right)$, corresponding to the non-decreasing sequence $j^{k}_{\mathrm{cn}\uparrow \mathrm{sh}}=\pi\left(\tau^{k}_{\mathrm{cn}\uparrow \mathrm{sh}},\tau_{\mathrm{pt}}\right)$ of interval indices for $\tau^{k}_{\mathrm{cn}\uparrow \mathrm{sh}}:=T_{j_{k}}\tau^{k}_{\mathrm{cn}}/T+\varDelta_{k}$. The corresponding vectorized expression is now given by
\begin{equation}
    -\frac{T_{j_{k}}}{T}\left(\left(\mathcal{L}_{\mathrm{cn}\uparrow \mathrm{sh}}^{k}\right)^\text{T}\otimes A_{k,s}^\text{T}(p)\right)\cdot\mathfrak{vec}\left(\mu_{\mathrm{bp}}\right),
\end{equation}
where
\begin{equation}
    \mathcal{L}_{\mathrm{cn}\uparrow \mathrm{sh},(m+1)(a-1)+c,b}^{k}=\delta_{a,j^{k}_{\mathrm{cn}\uparrow  sh,b}}\mathcal{L}_c\left(2N\tau^{k}_{\mathrm{cn}\uparrow \mathrm{sh},b}+1-2j^{k}_{\mathrm{cn}\uparrow  sh,b}\right),
\end{equation}
for $a=1,\ldots,N$, $c=1,\ldots,m+1$, and $b=1,\ldots,|\tau^{k}_{\mathrm{cn}\uparrow \mathrm{sh}}|$.

For the adjoint contributions in \eqref{eq:adjtcoupy} associated with variations $\delta y(\cdot)$ and discretized at $\tau_{\mathrm{bp}}$, the vectorization is simply given by $\mathfrak{vec}\left(\mu_{\mathrm{bp}}\right)$. In contrast, for the adjoint contributions in \eqref{eq:adjtcoupp}-\eqref{eq:adjtcoupvarDelta} associated with variations $\delta p$, $\delta T$, $\delta T_{j_k}$, and $\delta\varDelta_{k}$  and discretized using quadrature on the partition $\tau_{\mathrm{pt}}$ with collocation nodes at $\tau_{\mathrm{cn}}$, let $k_{\mathrm{cn}}=\pi\left(\tau_{\mathrm{cn}};\{0,\gamma_{\mathrm{b},2},\ldots,\gamma_{\mathrm{b},C},1\}\right)$ be a non-decreasing sequence of coupling interval indices associated with the elements of the $\tau_{\mathrm{cn}}$ array. For each subarray of successive elements $\tau_{\mathrm{cn}}^k$ that share an interval index $k\in k_{\mathrm{cn}}$, let $\Omega^k_{\mathrm{cn}}$ denote the corresponding subset of $\Omega_{\mathrm{cn}}$ and associate the shifted time instants
\begin{equation}
    \tau_{\mathrm{cn}\downarrow \mathrm{sh}}^{k}:=\frac{T}{T_{j_{k}}}\left(\tau_{\mathrm{cn}}^k-J_{|\tau_{\mathrm{cn}}^k|,1}\otimes\varDelta_{k}\right)
\end{equation}
with the non-decreasing sequence $j^{k}_{\mathrm{cn}\downarrow \mathrm{sh}}=\pi\left(\tau_{\mathrm{cn}\downarrow \mathrm{sh}}^{k},\tau_{\mathrm{pt}}\right)$ of interval indices $j$. The sought vectorized contributions are then given by
\begin{gather}
    -\sum_{k=1}^C\left(\sum_{s=1}^{S_k}\left(\mathfrak{vec}\left(A_{k,s}(p)\right)_{,p}\right)^\text{T}\cdot\left(x_{j_{k,s},\mathrm{bp}}\cdot\mathcal{L}^{k}_{\mathrm{cn}\downarrow \mathrm{sh}}\otimes I_n\right)\right)\cdot \Omega_{\mathrm{cn}}^k\cdot\mathfrak{vec}\left(\mu_{\mathrm{cn}}^k\right),\\
    -\frac{1}{T}\sum_{k=1}^C\left(\sum_{s=1}^{S_k}\mathfrak{vec}\left( x_{j_{k,s},\mathrm{bp}}\right)^\text{T}\cdot\left(\mathcal{L}'^{k}_{\mathrm{cn}\downarrow \mathrm{sh}}\otimes A_{k,s}^\text{T}(p)\right)\right)\cdot\left(\mathfrak{diag}(\tau_{\mathrm{cn}\downarrow \mathrm{sh}}^{k})\otimes I_n\right)\cdot\Omega^k_{\mathrm{cn}}\cdot\mathfrak{vec}\left(\mu_{\mathrm{cn}}^k\right),\\
     \frac{1}{T_{j_{k}}}\left(\sum_{s=1}^{S_k}\mathfrak{vec}\left( x_{j_{k,s},\mathrm{bp}}\right)^\text{T}\cdot\left(\mathcal{L}'^{k}_{\mathrm{cn}\downarrow \mathrm{sh}}\otimes A_{k,s}^\text{T}(p)\right)\right)\cdot\left(\mathfrak{diag}(\tau_{\mathrm{cn}\downarrow \mathrm{sh}}^{k})\otimes I_n\right)\cdot\Omega^k_{\mathrm{cn}}\cdot\mathfrak{vec}\left(\mu_{\mathrm{cn}}^k\right),\\
    \frac{T}{T_{j_k}}\left(\sum_{s=1}^{S_k}\mathfrak{vec}\left( x_{j_{k,s},\mathrm{bp}}\right)^\text{T}\cdot\left(\mathcal{L}'^{k}_{\mathrm{cn}\downarrow \mathrm{sh}}\otimes A_{k,s}^\text{T}(p)\right)\right)\cdot\Omega^k_{\mathrm{cn}}\cdot\mathfrak{vec}\left(\mu_{\mathrm{cn}}^k\right),
\end{gather}
where $\mu_{\mathrm{cn}}^k$ denote the corresponding elements of $\mu_{\mathrm{cn}}$ and
\begin{equation}
    \mathcal{L}'^{k}_{\mathrm{cn}\downarrow \mathrm{sh},(m+1)(a-1)+c,b}=\delta_{a,j^{k}_{\mathrm{cn}\downarrow  sh,b}}\mathcal{L}'_c\left(2N\tau^{k}_{\mathrm{cn}\downarrow \mathrm{sh},b}+1-2j^{k}_{\mathrm{cn}\downarrow  sh,b}\right),
\end{equation}
for $a=1,\ldots,N$, $c=1,\ldots,m+1$, and $b=1,\ldots,|\tau^{k}_{\mathrm{cn}\uparrow \mathrm{sh}}|$.

%%%%%%%%%%%%%%%%%%%%%
% Dimensional Deficit
%%%%%%%%%%%%%%%%%%%%%
\subsection{Dimensional Deficit}
\label{sec: toolbox template dimensional deficit}
The dimensional deficit of the zero problem in Section~\ref{sec: toolbox template zero problems} equals $M(n+2)+q-Kn+G$, where $K$ is the total number of boundary conditions~\eqref{eq: abstractbcs} and $G$ is the number of segments with $f_i$ explicitly dependent on $y_i$. With the imposition of the adjoint conditions, the dimensional deficit is reduced by $M(n+2)+q-Kn+G$ for a net value of $0$.

The discretization in the previous section is consistent with these counts, since i) there the $(M+G)N(m+1)n$ unknown components of $x_{\mathrm{bp}}$ and $y_{\mathrm{bp}}$ are constrained by $M(N-1)n$ continuity conditions, $MNmn$ discretized differential conditions, and $GN(m+1)n$ discretized coupling conditions, and ii) the $(M+G)N(m+1)n$ unknown components of $\lambda_{\mathrm{bp}}$ and $\mu_{\mathrm{bp}}$ are constrained by $M(N-1)n$ continuity conditions, $MNmn$ discretized adjoint differential conditions, $GN(m+1)n$ discretized adjoint coupling conditions, and $2Mn$ adjoint boundary conditions.

In the first example in Section~\ref{sec: toolbox template examples}, we supplement with the algebraic constraints $T_{0,1}=T_0$, $T_1=T$, and $\gamma_{\mathrm{e},1,1}=\alpha/T$ resulting in a composite zero problem with dimensional deficit $q$ provided that $T,T_0,\alpha\in p$. Similarly, in the second example, we supplement with the algebraic constraints $\gamma_{\mathrm{e},1,1}=\alpha/\beta$, $\gamma_{\mathrm{e},2,1}=\alpha/(T-\beta)$, $T_1=\beta$, $T_2=T-\beta$, $T_{0,1}=T_0$, and $T_{0,2}=T_0+\beta$ resulting in a composite zero problem with dimensional deficit $q$ provided that $\alpha,\beta,T,T_0\in p$. Finally, in the last example in this section, we supplement with the algebraic constraints $\gamma_{\mathrm{e},i,1}=\alpha/T$, $T_i=T$, and $T_{0,i}=T_0$ for $i=1,\ldots,M$ resulting in a composite zero problem with dimensional deficit $q$ provided that $\alpha,T,T_0\in p$. 

\subsection{Toolbox construction}
\label{sec: toolbox template toolbox construction}

We may pattern the development of a \textsc{coco}-compatible toolbox on the abstract template introduced in Section~\ref{sec: toolbox template zero problems} and the corresonding contributions to adjoint conditions in Section~\ref{sec: toolbox template adjoint conditions}. A zero problem constructor naturally decomposes into repeated calls to a toolbox constructor for the differential constraint~\eqref{rep_eq1} (the \mcode{ode_isol2coll} toolbox constructor already accomplishes this for vector fields that do not depend on an algebraic state variable), followed by a constructor for the boundary conditions~\eqref{eq: abstractbcs}, followed or interspersed by repeated calls to a constructor for the mesh conditions~\eqref{rep_eq3}, coupling conditions~\eqref{rep_eq2},  and algebraic conditions~\eqref{dde:xi:e}-\eqref{dde:xi:b}. Assuming, as is typically the case, that the constructor for the differential constraint assumes independent problem parameters in each constructor call, it is necessary to introduce additional algebraic constraints (a.k.a.~\emph{gluing conditions}) to ensure that the problem parameters are shared across all segments. Similar considerations may also apply to the construction of the coupling conditions.

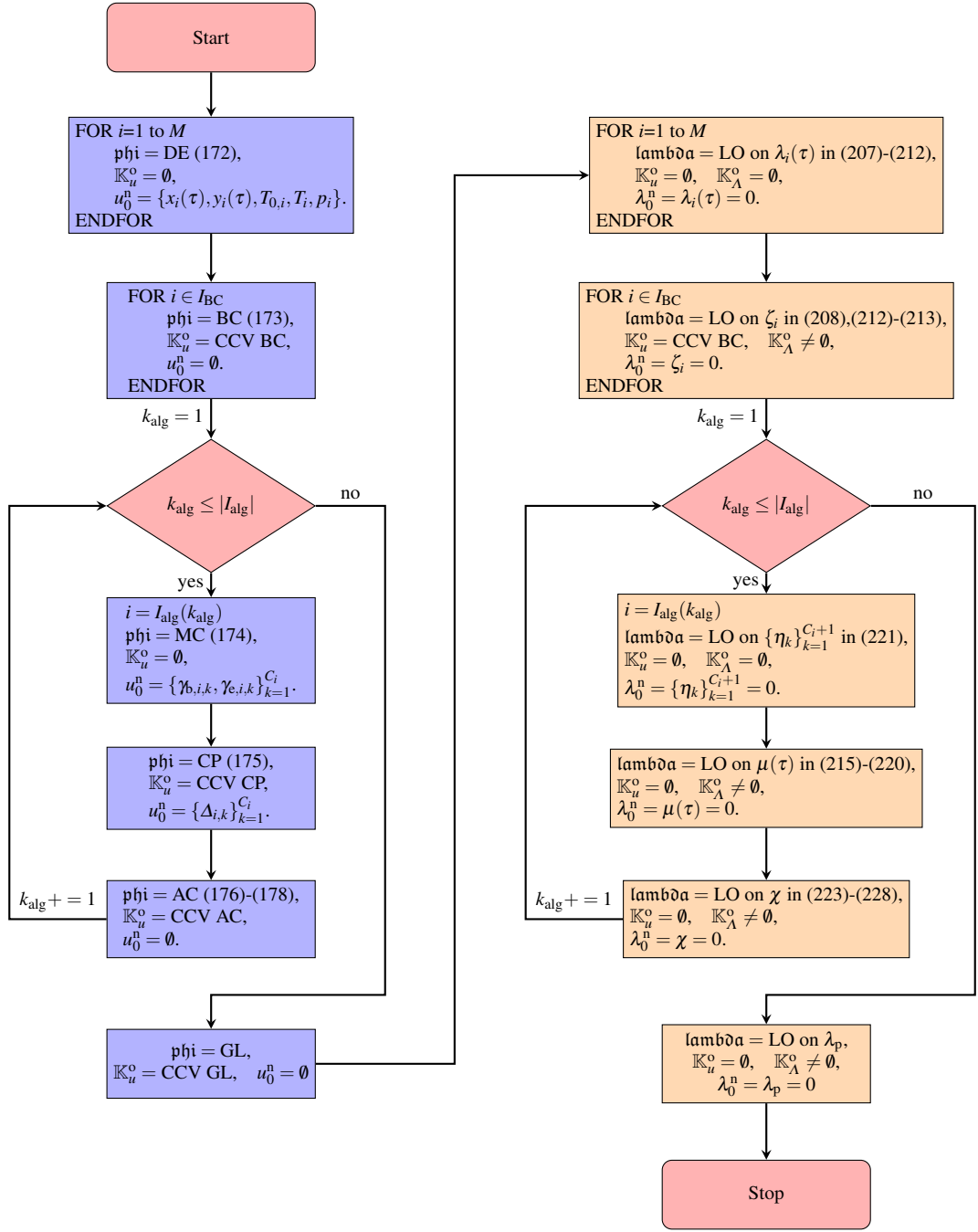
\begin{figure}
\centering
\begin{tikzpicture}
[node distance=2.0cm]
\node (start) [startstop] {Start};
\node (diff-eqs) [process-phi,below of=start,align=left] {FOR $i$=1 to $M$\\
\quad\quad $\mathfrak{phi}=\mathrm{DE}$~\eqref{rep_eq1},\\
\quad\quad $\mathbb{K}^\mathrm{o}_u=\emptyset$,\\
\quad\quad$u^\mathrm{n}_0=\{x_i(\tau), y_i(\tau),T_{0,i},T_i, p_i\}$.\\
ENDFOR};
\node (bc) [process-phi, below of=diff-eqs,yshift=-0.4cm,align=left] {FOR $i\in {I}_\mathrm{BC}$\\
\quad\quad $\mathfrak{phi}=\mathrm{BC}$~\eqref{eq: abstractbcs},\\
\quad\quad $\mathbb{K}^\mathrm{o}_u=$ CCV BC,\\
\quad\quad $u^\mathrm{n}_0=\emptyset$.\\
ENDFOR};
\node(check-alg) [decision, below of=bc, yshift=-0.4cm] {$k_\mathrm{alg}\leq|I_\mathrm{alg}|$};
\node (mc) [process-phi, below of=check-alg,yshift=-0.1cm,align=left] {$i=I_\mathrm{alg}(k_\mathrm{alg})$\\
$\mathfrak{phi}=\mathrm{MC}$~\eqref{rep_eq3},\\
$\mathbb{K}^\mathrm{o}_u=\emptyset$,\\
$u^\mathrm{n}_0=\{\gamma_{\mathrm{b},i,k},\gamma_{\mathrm{e},i,k}\}_{k=1}^{C_i}$.};
\node (cp) [process-phi, below of=mc,yshift=-0.0cm,align=left] {
$\mathfrak{phi}=\mathrm{CP}$~\eqref{rep_eq2},\\
$\mathbb{K}^\mathrm{o}_u=$ CCV CP,\\
$u^\mathrm{n}_0=\{\varDelta_{i,k}\}_{k=1}^{C_i}$.};
%\varDelta_{i,k}\right) \right) ,\,\tau\in\left[\gamma_{\mathrm{b},i,k},\gamma_{\mathrm{e},i,k}\right]
\node (ac) [process-phi, below of=cp,yshift=0.1cm,align=left] {
$\mathfrak{phi}=\mathrm{AC}$~\eqref{eq_xi}-\eqref{dde:xi:b},\\
$\mathbb{K}^\mathrm{o}_u=$ CCV AC,\\
$u^\mathrm{n}_0=\emptyset$.};
\node (glue) [process-phi,yshift=-0.1cm,below of=ac,align=center]{$\mathfrak{phi}=\mathrm{GL}$,\\
$\mathbb{K}^\mathrm{o}_u=$ CCV GL,\quad $u^\mathrm{n}_0=\emptyset$};
% adjoint ones
\node (diff-eqs-adjt) [process-lambda,right of=diff-eqs,xshift=6.0cm,align=left] {FOR $i$=1 to $M$\\
\quad\quad $\mathfrak{lambda}=$ LO on $\lambda_i(\tau)$ in~\eqref{eq:diffadjx}-\eqref{eq:diffadjp},\\
\quad\quad $\mathbb{K}^\mathrm{o}_u=\emptyset$,\quad $\mathbb{K}^\mathrm{o}_\Lambda=\emptyset$,\\
\quad\quad$\lambda^\mathrm{n}_0=\lambda_i(\tau)=0$.\\
ENDFOR};
\node (bc-adjt) [process-lambda, right of=bc,xshift=6cm,align=left] {FOR $i\in I_\mathrm{BC}$\\
\quad\quad $\mathfrak{lambda}=$ LO on $\zeta_i$ in~\eqref{eq:diffadjx01},\eqref{eq:diffadjp}-\eqref{eq:diffadjxj1},\\
\quad\quad $\mathbb{K}^\mathrm{o}_u=$ CCV BC,\quad $\mathbb{K}^\mathrm{o}_\Lambda\neq\emptyset$,\\
\quad\quad$\lambda^\mathrm{n}_0=\zeta_i=0$.\\
ENDFOR};
\node(check-alg-adjt) [decision, below of=bc-adjt, yshift=-0.4cm] {$k_\mathrm{alg}\leq|I_\mathrm{alg}|$};
\node (mc-adjt) [process-lambda, right of=mc,xshift=6cm,align=left] {$i=I_\mathrm{alg}(k_\mathrm{alg})$\\
$\mathfrak{lambda}=$ LO on $\{\eta_k\}_{k=1}^{C_i+1}$ in~\eqref{eq:adjt-mesh},\\
$\mathbb{K}^\mathrm{o}_u=\emptyset$,\quad $\mathbb{K}^\mathrm{o}_\Lambda=\emptyset$,\\
$\lambda^\mathrm{n}_0=\{\eta_k\}_{k=1}^{C_i+1}=0$.};
\node (cp-adjt) [process-lambda, right of=cp,xshift=6cm,align=left] {
$\mathfrak{lambda}=$ LO on $\mu(\tau)$ in~\eqref{eq:adjtcoupy}-\eqref{eq:adjtcoupxjks},\\
$\mathbb{K}^\mathrm{o}_u=\emptyset$,\quad $\mathbb{K}^\mathrm{o}_\Lambda\neq\emptyset$,\\
$\lambda^\mathrm{n}_0=\mu(\tau)=0$.};
\node (ac-adjt) [process-lambda, right of=ac,xshift=6cm,align=left] {
$\mathfrak{lambda}=$ LO on $\chi$ in~\eqref{eq:adjt-alg-T}-\eqref{eq:adjt-alg-gammab},\\
$\mathbb{K}^\mathrm{o}_u=\emptyset$,\quad $\mathbb{K}^\mathrm{o}_\Lambda\neq\emptyset$,\\
$\lambda^\mathrm{n}_0=\chi=0$.};
\node (glue-adjt) [process-lambda,xshift=6.0cm,right of=glue,align=center]{$\mathfrak{lambda}=$ LO on $\lambda_\mathrm{p}$,\\
$\mathbb{K}^\mathrm{o}_u=\emptyset$,\quad$\mathbb{K}^\mathrm{o}_\Lambda\neq\emptyset$,\\ $\lambda^\mathrm{n}_0=\lambda_\mathrm{p}=0$};
\node (stop) [startstop, below of=glue-adjt,yshift=0.1cm] {Stop};
% DRAW LINES
\draw [arrow] (start) -- (diff-eqs);
\draw [arrow] (diff-eqs) -- (bc);
\draw [arrow] (bc) -- node[anchor=east]{$k_\mathrm{alg}=1$}(check-alg);
\draw [arrow] (check-alg) -- node[anchor=east] {yes} (mc);
\draw [arrow] (mc) -- (cp);
\draw [arrow] (cp) -- (ac);
\draw [->,thick,>=stealth]
    (ac.west) -- node[anchor=south] {$k_\mathrm{alg}+=1$}
    ++(-1.4cm,0) |-
    (check-alg.west)
    ;
\draw [->,thick,>=stealth]
    (check-alg.east) -- node[anchor=south] {no}
    ++ (1.0cm,0) -- ++(0,-7.1cm) -| (glue.north);
\draw [arrow] (diff-eqs-adjt) -- (bc-adjt);
\draw [arrow] (bc-adjt) -- node[anchor=east]{$k_\mathrm{alg}=1$}(check-alg-adjt);
\draw [arrow] (check-alg-adjt) -- node[anchor=east] {yes} (mc-adjt);
\draw [arrow] (mc-adjt) -- (cp-adjt);
\draw [arrow] (cp-adjt) -- (ac-adjt);
\draw [->,thick,>=stealth]
    (ac-adjt.west) -- node[anchor=south] {$k_\mathrm{alg}+=1$}
    ++(-1.4cm,0) |-
    (check-alg-adjt.west)
    ;
\draw [->,thick,>=stealth]
    (check-alg-adjt.east) -- node[anchor=south] {no}
    ++ (1.5cm,0) -- ++(0,-7.1cm) -| (glue-adjt.north);
\draw [arrow] (glue-adjt) -- (stop);
\draw [->,thick,>=stealth]
    (glue.east) -- ++ (2.0cm,0) |- (diff-eqs-adjt.west);
%\draw [->,thick,>=stealth] (ac) to[-|-] (check-alg);
\end{tikzpicture}
\caption{A flowchart depicting the construction of the abstract zero problem in Section~\ref{sec: toolbox template zero problems} and corresponding contributions to adjoint conditions in Section~\ref{sec: toolbox template adjoint conditions}. Here rectangles filled with blue and orange colors represent core constructors associated with functions of the type $\mathfrak{phi}$ and $\mathfrak{lambda}$, respectively. The abbreviations DE, BC, MC, CP, CCV, AC, GL, and LO represent differential equations, boundary conditions, mesh conditions, coupling conditions, corresponding continuation variables, algebraic conditions, glue conditions, and linear operators, respectively. In particular, $\mathbb{K}^\mathrm{o}_u=$ CCV BC/CP/AC/GL denote indexing the corresponding continuation variables for boundary conditions/coupling conditions/algebraic conditions/glue conditions from the ones defined when constructing the differential constraints. $I_\mathrm{BC}\subset\{1,\cdots,M\}$ gives the set of differential state variables involving boundary conditions. $I_\mathrm{alg}\subset\{1,\cdots,M\}$ gives the collection of algebraic state variables involving coupling conditions. $k_\mathrm{alg}+=1$ should be interpreted as $k_\mathrm{alg}=k_\mathrm{alg}+1$. Note that the indices $i$ in Section~\ref{sec: toolbox template adjoint conditions} has been omitted. Such indices can be added properly to the adjoints derived in the section.}
\label{fig:flowchart-toolbox}
\end{figure}

A flow chart similar in character to that in Section~\ref{sec: data assimilation problem construction} corresponding here to the construction of an augmented continuation problem per the abstract toolbox template is shown in Fig.~\ref{fig:flowchart-toolbox}. Since the $i$-th coupling condition depends on segments $x_{j_{i,k,s}}$, for $s=1,\ldots,S_k$ and $k=1,\ldots,C_i$, it must be constructed after the corresponding differential constraints have been appended to the continuation problem. It is not necessary, however, to wait until all differential constraints have been introduced. Similar considerations apply to the constructors for the contributions to the adjoint conditions. These may be invoked only after the entire zero problem has been constructed, or at opportune moments following the construction of the corresponding zero functions.

In the next section, we proceed to illustrate the application of such a \textsc{coco}-compatible toolbox through several numerical examples. These demonstrate the versatility of the tool and the opportunity to use such a toolbox, and the \textsc{coco} construction paradigm described in Section~\ref{sec: The coco formalism}, to build sophisticated special-purpose toolboxes, dedicated to particular classes of delay differential equations, say.

%%%%%%%%%%%%%%%%%%%%
% Numerical Examples
%%%%%%%%%%%%%%%%%%%%
\section{Numerical examples}
\label{sec: Numerical examples}
\subsection{Generalizations of the abstract framework}
\label{sec: numerical examples generalizations}

Before we consider numerical examples illustrating the ability of the proposed toolbox to perform continuation and constrained optimization for boundary-value problems with delay, we consider two possible generalizations that allow one to handle initial-value problems and unknown exogenous driving, as well as problems involving multiple discrete delays.  

Consider, for example, the initial-value problem
\begin{align}
    \dot{z}\left(t\right)&=f\left(t,z(t),z\left(t-\alpha\right),p\right),\,t\in\left(T_{0},T_{0}+T\right),\\
    z\left(t\right)&=g\left(t+\alpha-T_0,p\right),\,t\in\left[T_{0}-\alpha,T_{0}\right]
\end{align}
for some known function $g(s,p)$, for $s\in[0,\alpha]$. The substitution $x(\tau)=z(T_0+T\tau)$ then yields
\begin{equation}
\begin{aligned}
 x'(\tau)&=Tf\left(T_0+T\tau,x(\tau),y(\tau),p\right),\,\tau\in\left(0,1\right),\\
y\left(\tau\right) &=\begin{cases}
g\left(T\tau, p \right), & \tau\in\left(0,\alpha/T\right),\\
x\left(\tau-\alpha/T\right), & \tau\in\left(\alpha/T,1\right).
\end{cases}
\end{aligned}
\end{equation}
Inspired by this example, we consider coupling conditions of the form
\begin{equation}
\label{abstract_eq2}
  y_{i}\left(\tau\right)=g_{i,k}\left(T_{i}\tau-\varDelta_{i,k,1},p\right),\tau\in\left(\gamma_{\mathrm{b},i,k},\gamma_{\mathrm{e},i,k}\right)  
\end{equation}
for some $g_{i,k}:\mathbb{R}\times\mathbb{R}^{q}\rightarrow \mathbb{R}^{n}$ as an alternative to the form given in \eqref{rep_eq2}. We leave it to the reader to derive the associated form of the contributions to the adjoint necessary conditions.

A further generalization of the framework introduced in Section~\ref{sec: Toolbox construction} is support for vector fields of the form $f:\mathbb{R}\times\mathbb{R}^n\times\mathbb{R}^{n\times n_\delta}\times\mathbb{R}^q\rightarrow\mathbb{R}^n$ with coupling conditions of the form \eqref{rep_eq2} for each column of the $n\times n_\delta$ matrix of algebraic state variables. This may be used to analyze problems with multiple discrete time delays. In combination with the implementation of coupling conditions of the form \eqref{abstract_eq2}, it allows for analysis of problems with or without delay and with unknown exogenous driving terms to be determined through the optimization of an objective function, as in the theory of optimal control. 

Both of these techniques are illustrated through the sequence of examples below.

\subsection{Connecting orbits}
\label{sec: numerical examples connecting orbits}
The dynamical system
\begin{equation}
\label{eq1:connecting}
\dot{z}(t)=f(z(t),z(t-\alpha),p):=\begin{pmatrix}
z_{2}(t)\\
z_{1}(t)-z_{1}(t)z_{1}\left(t-\alpha\right)+p_2 z_{2}(t)+p_1 z_{1}(t)z_{2}(t)
\end{pmatrix}
\end{equation}
admits equilibrium solutions at $(z_1,z_2)=(0,0)$ and $(1,0)$ for all values of the parameters $p_1$ and $p_2$ and delay $\alpha$. The equilibrium at the origin is always a saddle, while that at $(1,0)$ undergoes a Hopf bifurcation with angular frequency $\omega$ provided that
\begin{equation}
    \cos\alpha\omega=\omega^2,\,\sin\alpha\omega=-(p_1+p_2)\omega.
\end{equation}
For $\alpha=0.8255$ and $p_1=0.5$ (cf.\ \cite{samaey2002numerical}), the bifurcation occurs for $p_2=p_{2,\text{HB}}\approx-1.257$ with $\omega\approx0.868$. A supercritical family of limit cycles of approximate limiting form $(z_1(t)-1,z_2(t))\sim\sqrt{p_2-p_{2,\text{HB}}}(\omega^{-1}\sin\omega t,\cos\omega t)$ grows out of this point as $p_2$ increases past the critical value $p_{2,\text{HB}}$.

We obtain a continuation problem for single-segment periodic orbits of the form developed in the first example in Section~\ref{sec: toolbox template examples} by straightforward substitution. Since the vector field is autonomous, we let $T_0=0$ without loss of generality and impose the condition $x_2(0)=0$ to fix the solution phase. We use interpolation according to~\eqref{eq:interpolants} to impose the condition $x_{2}(\tau_\text{cr})=0$ in terms of the additional continuation variable $\tau_\text{cr}$. Finally, we introduce monitor functions evaluating to $p_1$, $p_2$, $\alpha$, $T$, and $\tau_\text{cr}$ and denote the corresponding continuation parameters by $\mu_{p_1}$, $\mu_{p_2}$, $\mu_\alpha$, $\mu_T$ and $\mu_{\tau_\text{cr}}$. An initial solution guess is then given by $(x_1(\tau),x_2(\tau))=(1,0)+0.01(0.868^{-1}\sin(2\pi\tau-\pi/2),\cos(2\pi\tau-\pi/2))$, $T=2\pi/0.868$, $p_1=0.5$, $p_2=-1.257$, $\alpha=0.8255$, and $\tau_\text{cr}=0.5$. 

\begin{figure}[ht]
\centering
\subfloat[]{\includegraphics[width=0.49\columnwidth,trim={0 {0.0\textwidth} 0 0.0\textwidth},clip]{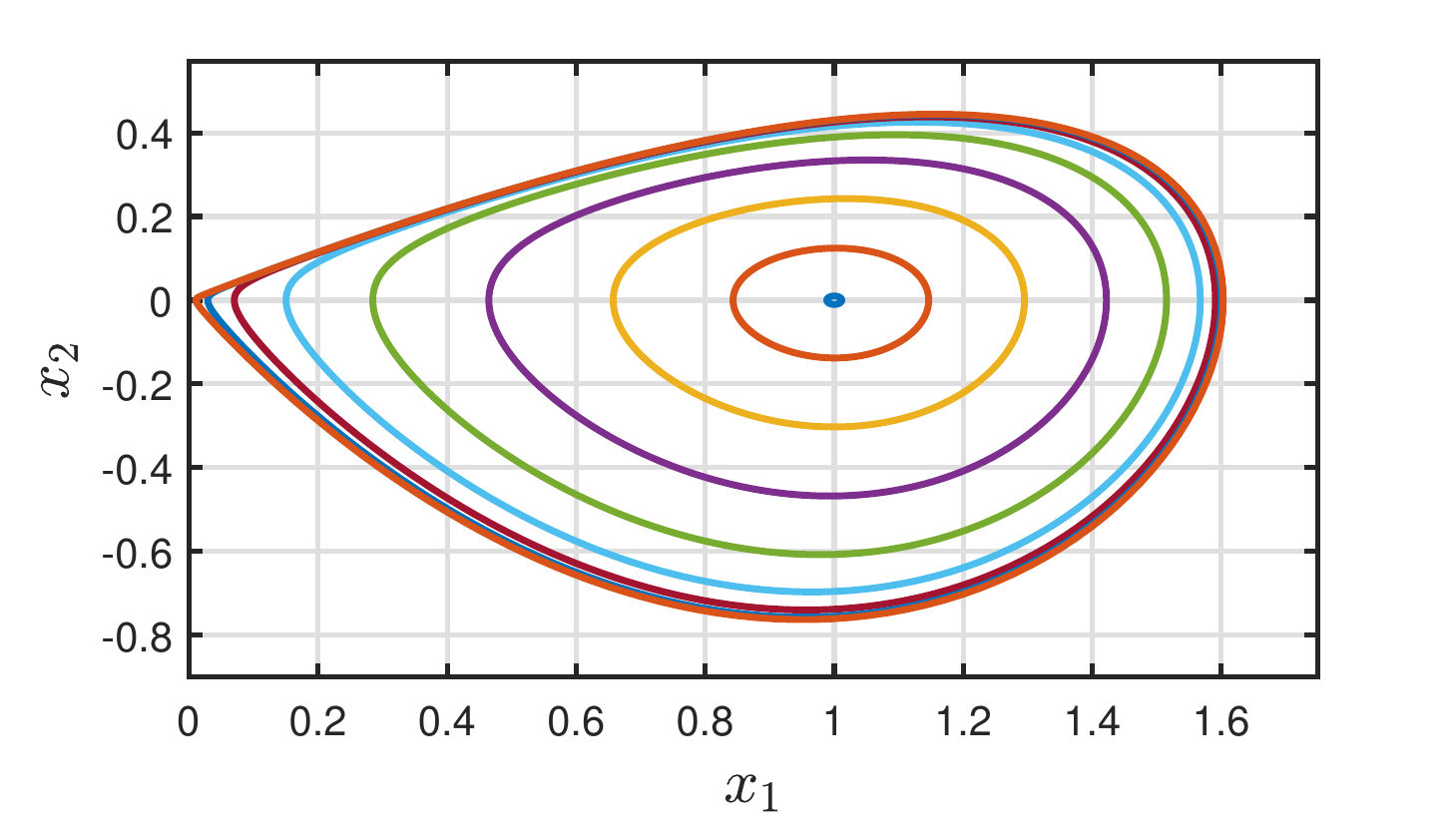}}
\subfloat[]{\includegraphics[width=0.49\columnwidth,trim={0 {0.0\textwidth} 0 0.0\textwidth},clip]{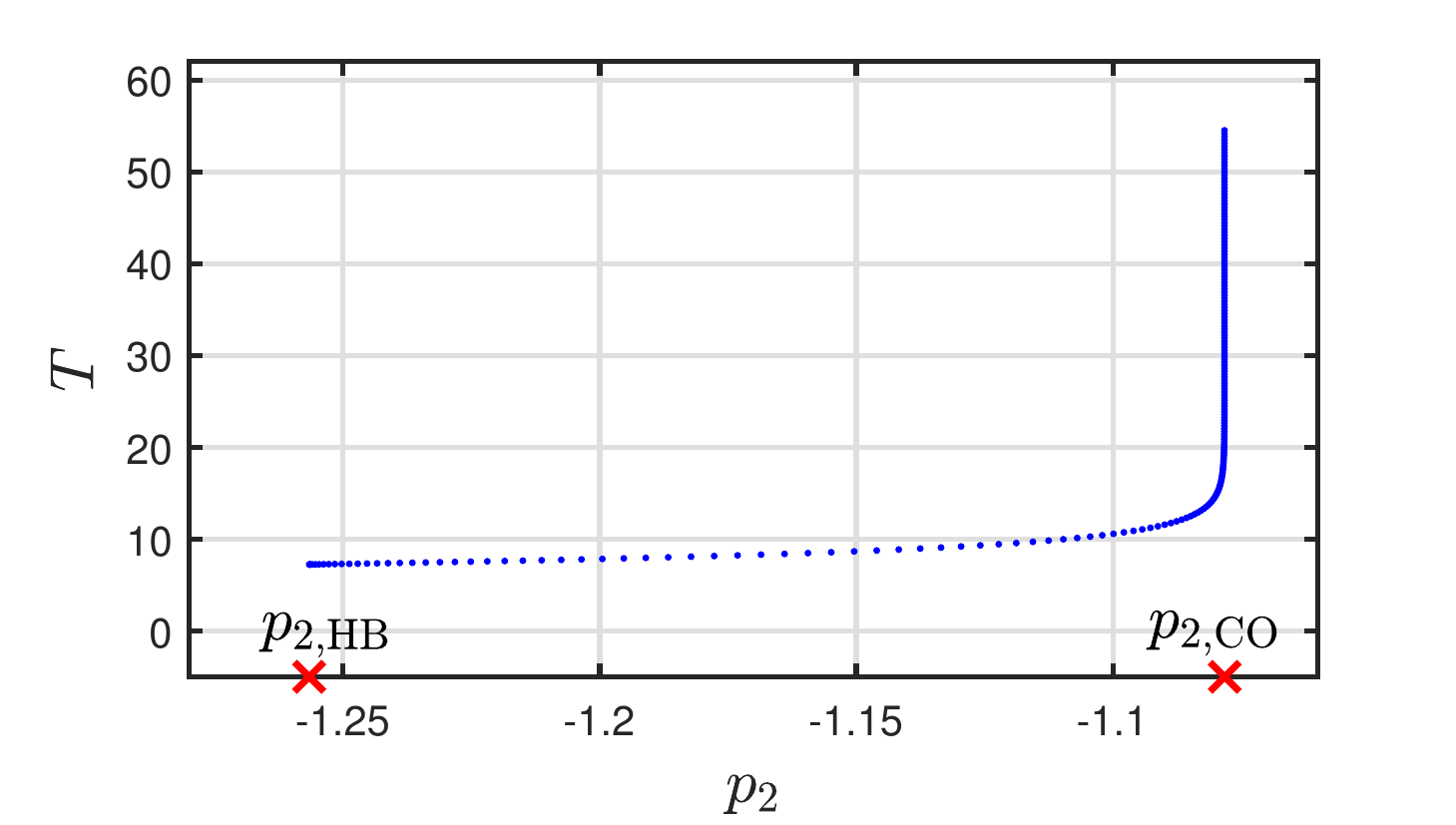}}
\caption{(left) A family of periodic orbits of the dynamical system in~(\ref{eq1:connecting}) born from a Hopf bifurcation at $p_2=p_{2,\text{HB}}$ with fixed $p_1=0.5$ and $\alpha=0.8255$ and varying $T$ and $p_2$. (right) The corresponding period $T$ shows unbounded growth as a function of $p_2$ as a homoclinic bifurcation at $p_2\approx p_{2,\text{CO}}$ is approached.} 
\label{fig2:connecting}
\end{figure}

With $N=100$, one-dimensional continuation with $\mu_{p_1}$ and $\mu_\alpha$ fixed and $\mu_{p_2}$, $\mu_T$, and $\mu_{\tau_\text{cr}}$ free to vary yields the one-dimensional family of limit cycles sampled in the left panel of Fig.~\ref{fig2:connecting}. Notably, as seen in the right panel, the period $T$ increases without bound as $p_2$ approaches a number close to $-1.0782$, suggestive of the existence of a homoclinic connecting orbit based at the saddle equilibrium at the origin, as also evident in the left panel. We may consider the periodic orbit with $T=20$ as a first-order approximation to such a connecting orbit.

We proceed to construct a continuation problem for locating and tracking an approximate family of such connecting orbits under variations in the problem parameters and delay by replacing the single-segment periodic orbit coupling conditions with the equations
\begin{align}
    y(\tau)&=\begin{cases}
    \epsilon ve^{\ell T\left( \tau-\alpha/T\right)}, & \tau\in\left(0,\alpha/T\right),\\
    x\left(\tau-\alpha/T\right), & \tau\in\left(\alpha/T,1\right),
    \end{cases}
    \label{eq2:connecting1}
\end{align}
and the periodic boundary condition with $x(0)=\epsilon v$ in terms of an unknown scalar $\epsilon$ and a normalized right eigenvector $v$ of the Jacobian
\begin{equation}
    A=\begin{pmatrix}0 & 1\\1 & p_2\end{pmatrix}
\end{equation}
of the  linearization at the origin corresponding to the unstable eigenvalue $\ell$. We impose the additional boundary condition
\begin{equation}
    w^\text{T}x\left(1\right)=0
\end{equation}
in terms of a normalized left eigenvector $w$ of $A$ in order to ensure that the end point $x(1)$ lies in the stable eigenspace of the equilibrium at the origin. We again let $T_0=0$ and use interpolation according to~\eqref{eq:interpolants} to impose the condition $x_{2}(\tau_\text{cr})=0$ in terms of the additional continuation variable $\tau_\text{cr}$. We introduce monitor functions evaluating to $p_1$, $p_2$, $\alpha$, $T$, $\tau_\text{cr}$, and $\epsilon$ and denote the corresponding continuation parameters by $\mu_{p_1}$, $\mu_{p_2}$, $\mu_\alpha$, $\mu_T$, $\mu_{\tau_\text{cr}}$, and $\mu_\epsilon$.

\begin{figure}[ht]
\centering
\subfloat[]{\includegraphics[width=0.49\columnwidth,trim={0 {0.0\textwidth} 0 0.0\textwidth},clip]{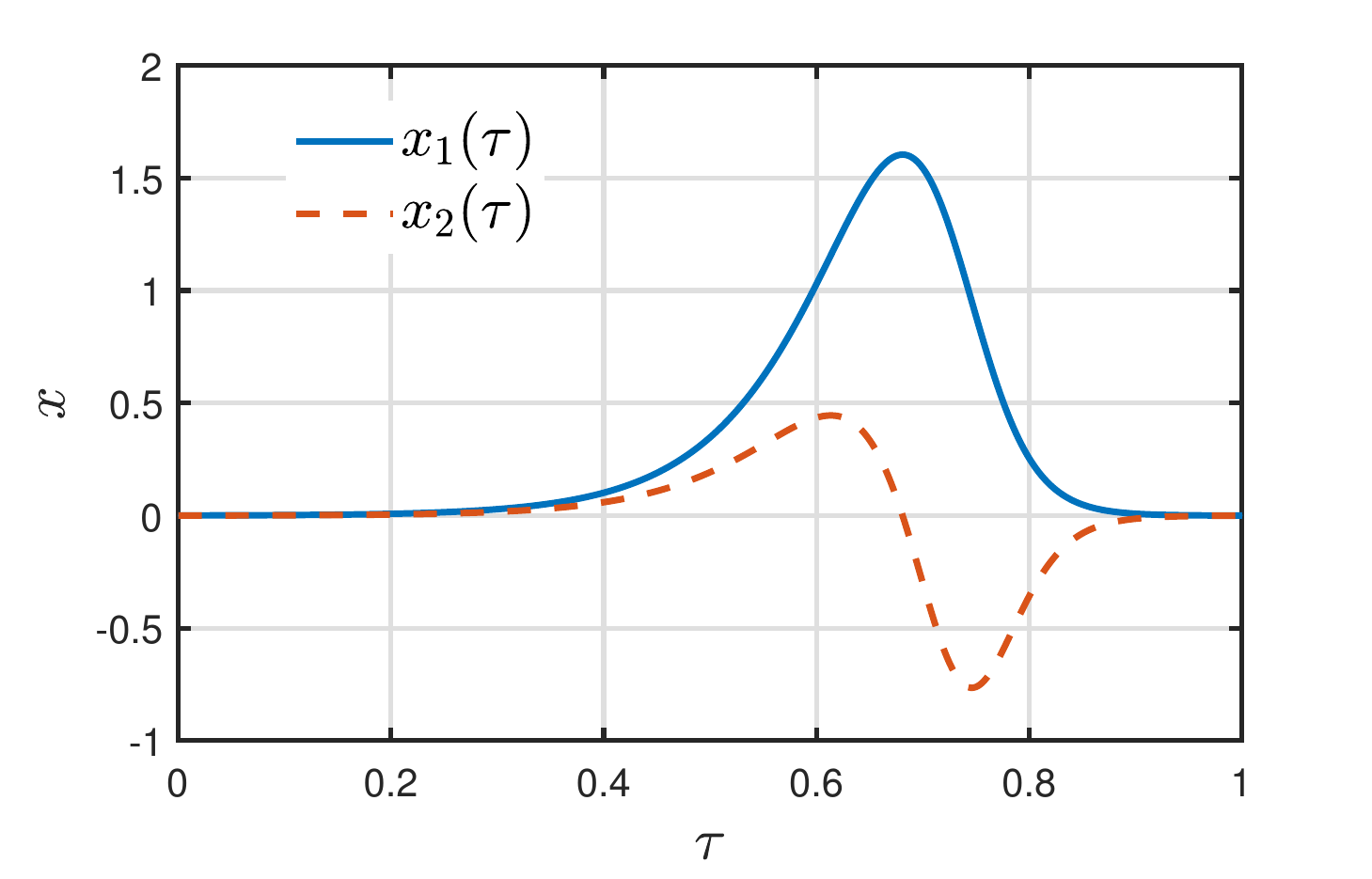}}
\subfloat[]{\includegraphics[width=0.49\columnwidth,trim={0 {0.0\textwidth} 0 0.0\textwidth},clip]{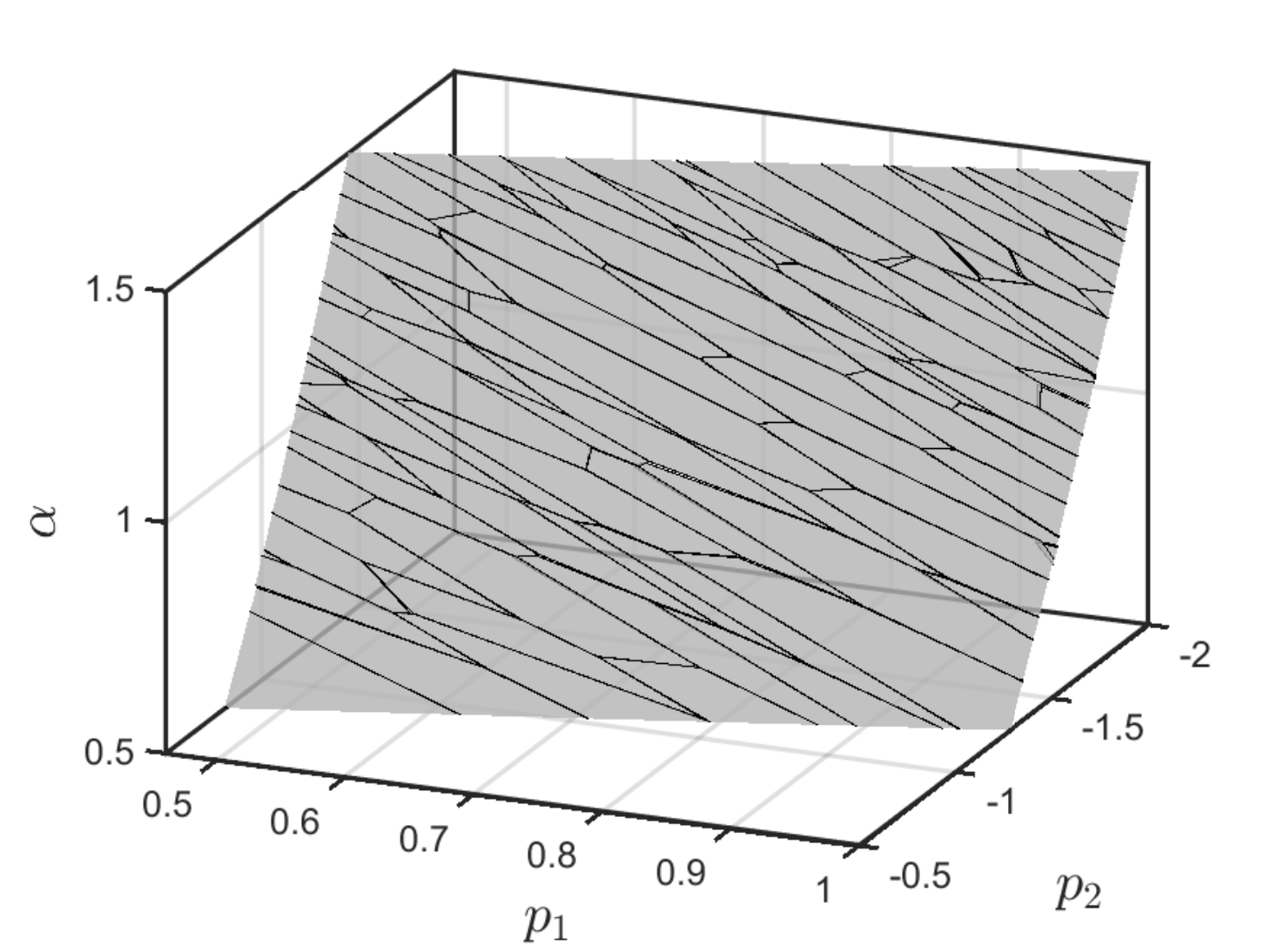}}
\caption{(a) Approximate connecting orbit of the delay differential equation~(\ref{eq1:connecting}) obtained using the proposed framework for  $\alpha=0.8255$, $p_{1}=0.5$, $p_2\approx -1.0782$, $T=20$, $\tau_\text{cr}\approx0.6756$, $\epsilon\approx 1.3907\times 10^{-3}$. (b) Homoclinic bifurcation surface obtained using two-dimensional continuation with \textsc{coco}.}
\label{fig1:connecting}
\end{figure}

For fixed $\mu_T$, $\mu_{p_{1}}$, $\mu_{p_{2}}$, $\mu_\alpha$, $
\mu_{\tau_{\mathrm{cr}}}$, and $\mu_\epsilon$ we obtain a continuation problem with dimensional deficit $-2$. We proceed to first allow $\mu_{p_2}$ and $\mu_\epsilon$ to vary and obtain a unique solution by applying a nonlinear solver (\mcode{atlas\_0d} in \textsc{coco}) with initial solution guess given by the periodic orbit with $T=20$, the corresponding values of $p_1$, $p_2$, $\alpha$, and $\tau_\text{cr}$, and $\epsilon=10^{-3}$. Next, we fix $\mu_\epsilon$, and allow additionally $\mu_{p_1}$, $\mu_\alpha$, and $\mu_T$ to vary in order to obtain a two-dimensional solution manifold and seed \mcode{atlas\_kd} with the unique solution found in the previous step. The left panel of Fig.~\ref{fig1:connecting} shows the approximate homoclinic obtained with $p_1=0.5$ and $\alpha=0.8255$, while the right panel shows the corresponding two-dimensional homoclinic bifurcation surface under simultaneous variations in $p_1$, $p_2$, and $\alpha$.

\subsection{Phase response curves}
\label{sec: numerical examples phase response curves}
\begin{figure}[ht]
\centering
\subfloat{\includegraphics[width=0.49\columnwidth,trim={0 {0.0\textwidth} 0 0.0\textwidth},clip]{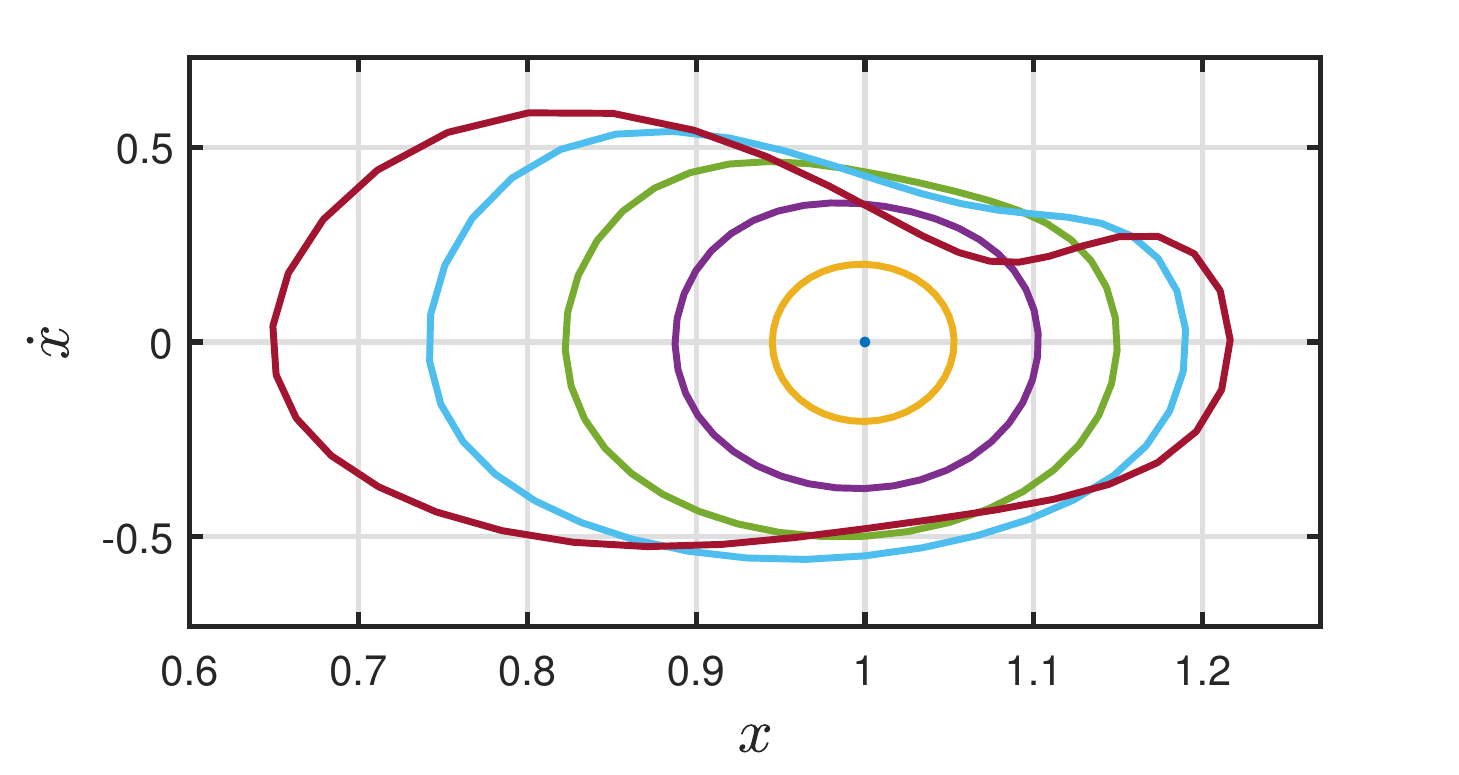}}
\caption{Sample orbits from a family of limit cycles of the Mackey-Glass system~\eqref{mackey} with $a=2$ and $b=10$ emanating from a supercritical Hopf bifurcation under variations in the delay $\alpha$ past the critical value $\alpha_{\text{HB}}\approx 0.4708$.
} 
\label{fig:prc0}
\end{figure}

We use the methodology in Section~\ref{sec: phase response curves linear response theory} to compute the phase response curve of a stable limit cycle of the Mackey-Glass equation~\cite{glass2010mackey} 
\begin{align}
\label{mackey}
    \dot{z}(t)=\frac{az\left(t-\alpha\right)}{1+z^{b}\left(t-\alpha\right)}-z(t),
\end{align}
an oft-used model for the dynamics of physiological systems such as respiratory dynamics~\cite{berezansky2012mackey} and the production of white blood cells~\cite{berezansky2006mackey}. Specifically, a family of such limit cycles is born from the equilibrium at $z=(a-1)^{1/b}$ for $a=2$ and $b=10$ when $\alpha$ increases through a Hopf bifurcation at $\alpha_\text{HB}\approx0.4708$ as shown in Fig.~\ref{fig:prc0}. We proceed to analyze the phase response curve for the orbit obtained with $\alpha=0.7$.

With the decomposition in Section~\ref{sec: toolbox template examples}, we obtain
\begin{equation}
\begin{aligned}
 x'(\tau)&=T\left(\frac{ay\left(\tau\right)}{1+y^{b}\left(\tau\right)}-x(\tau)\right),\,\tau\in\left(0,1\right),\\
 y(\tau)&=\begin{cases}
x\left(\tau+1-\alpha/T\right), & \tau\in\left(0,\alpha/T\right),\\
x\left(\tau-\alpha/T\right), & \tau\in\left(\alpha/T,1\right).
\end{cases}
\end{aligned}
\end{equation}
The toolbox developed in Section~\ref{sec: Toolbox construction} may be applied out of the box to construct the corresponding zero problems and adjoint contributions provided that we append the algebraic conditions $\gamma_{\mathrm{e},1,1}=\alpha/T$, $T_0=0$, and the periodicity condition $x(0)=x(1)$. Since the vector field is autonomous, we append the phase condition $x(0)=1$ to remove the invariance to time shifts.

We proceed to append monitor functions that evaluate to $T$, $\alpha$, $a$, and $b$ with corresponding continuation parameter $\mu_T$, $\mu_\alpha$, $\mu_a$, and $\mu_b$, respectively. We denote the continuation multipliers associated with the corresponding adjoint contributions by $\eta_T$, $\eta_\alpha$, $\eta_a$, and $\eta_b$. Then, the phase response curve is obtained from the continuation multiplier $\lambda_{\text{DE}}(\cdot)$ associated with the differential constraint at a solution to the augmented continuation problem with $\eta_T=1$. To locate such a solution, we append a complementary zero function that evaluates to $\eta_T-1$.

\begin{figure}[ht]
\centering
\subfloat{\includegraphics[width=0.49\columnwidth,trim={0 {0.0\textwidth} 0 0.0\textwidth},clip]{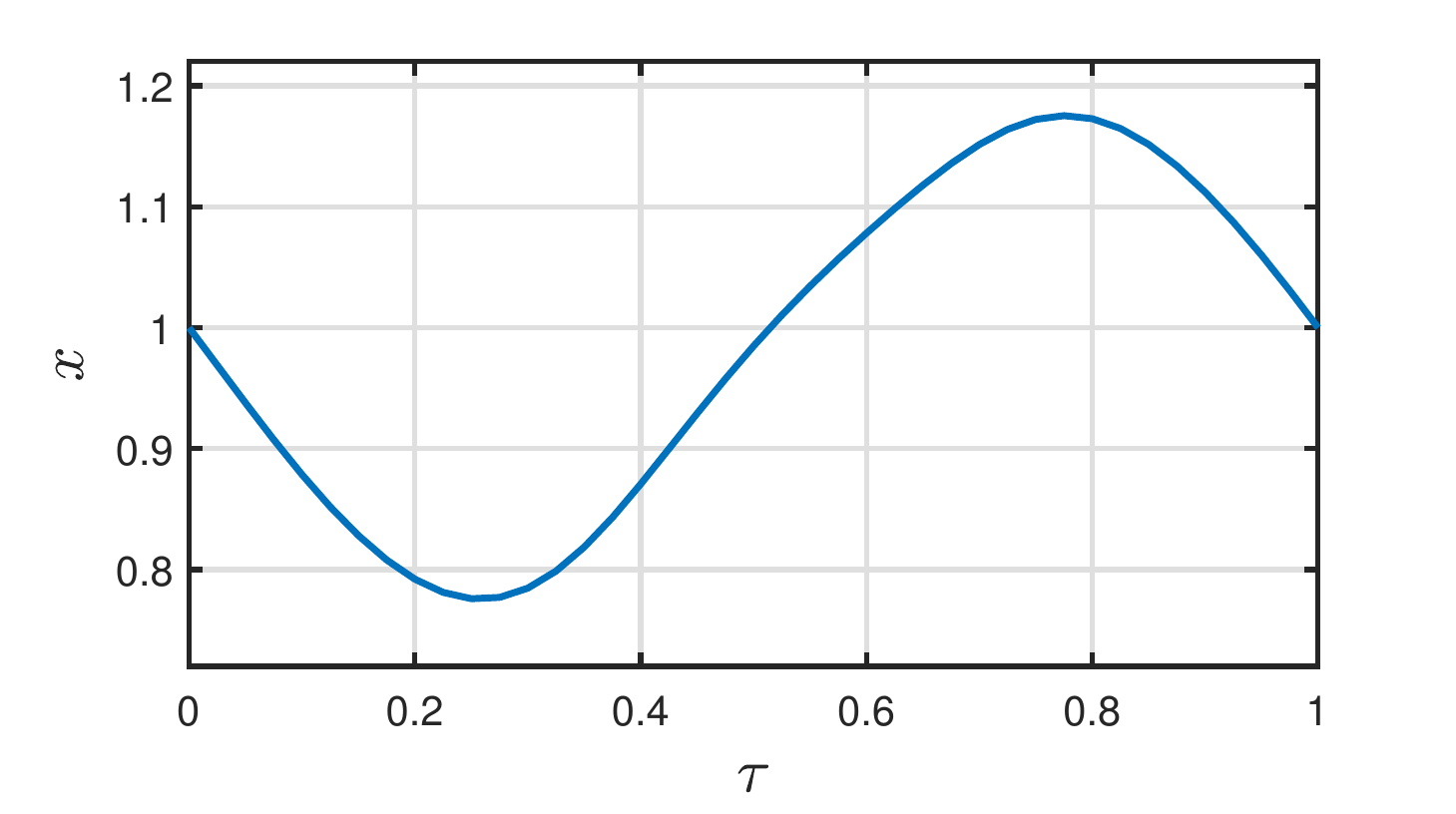}}
\subfloat{\includegraphics[width=0.49\columnwidth,trim={0 {0.0\textwidth} 0 0.0\textwidth},clip]{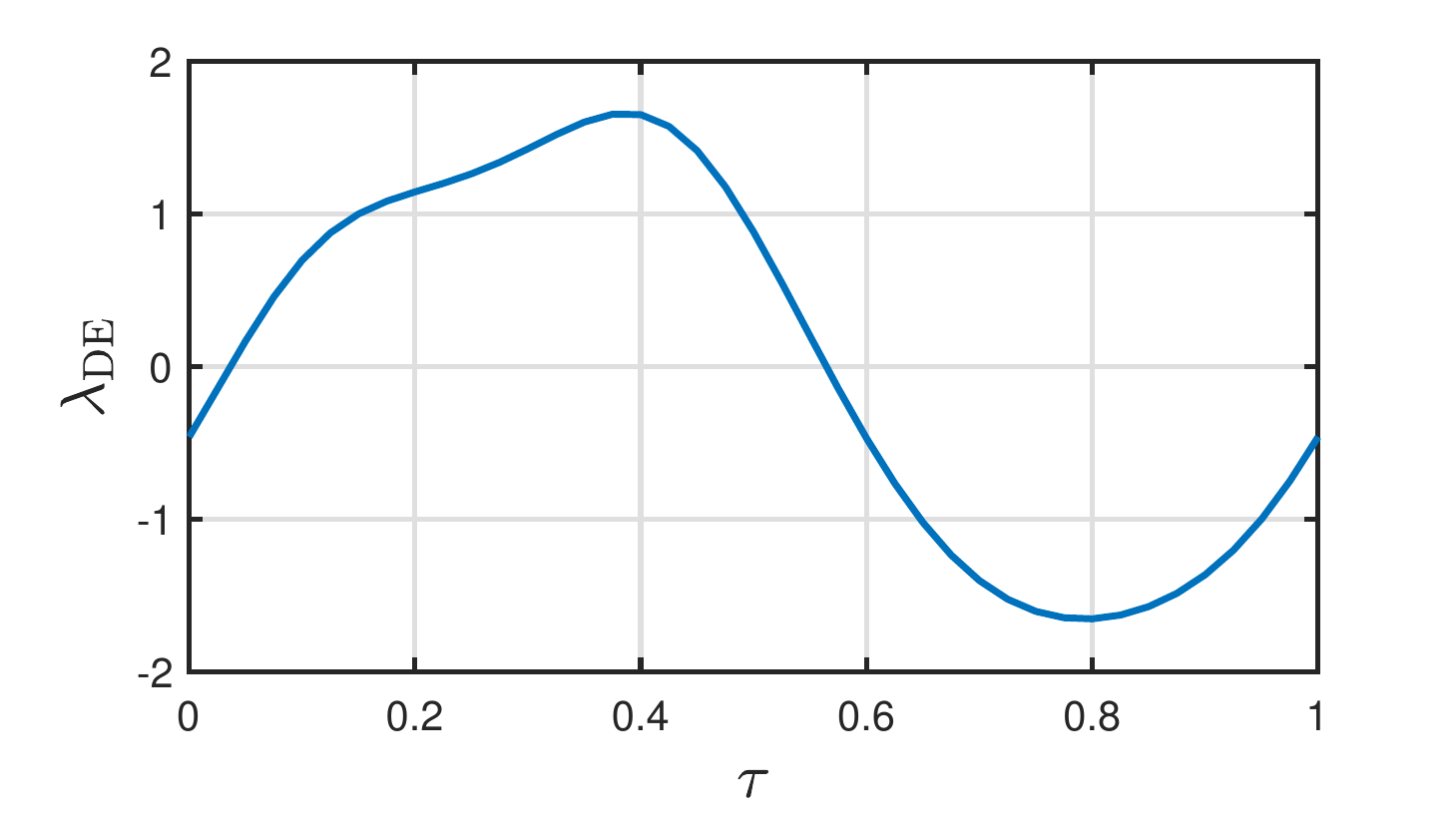}}
\caption{Time histories for the differential variable $x(\cdot)$ and corresponding phase response curve $\lambda_\text{DE}(\cdot)$ for a limit cycle of the delay differential equation \eqref{mackey} with $a=2$, $b=10$, $\alpha=0.7$, and $T\approx 2.2958$.} 
\label{fig:prc1}
\end{figure}

\begin{figure}[ht]
\centering
\subfloat[]{\includegraphics[width=0.49\columnwidth,trim={0 {0.0\textwidth} 0 0.0\textwidth},clip]{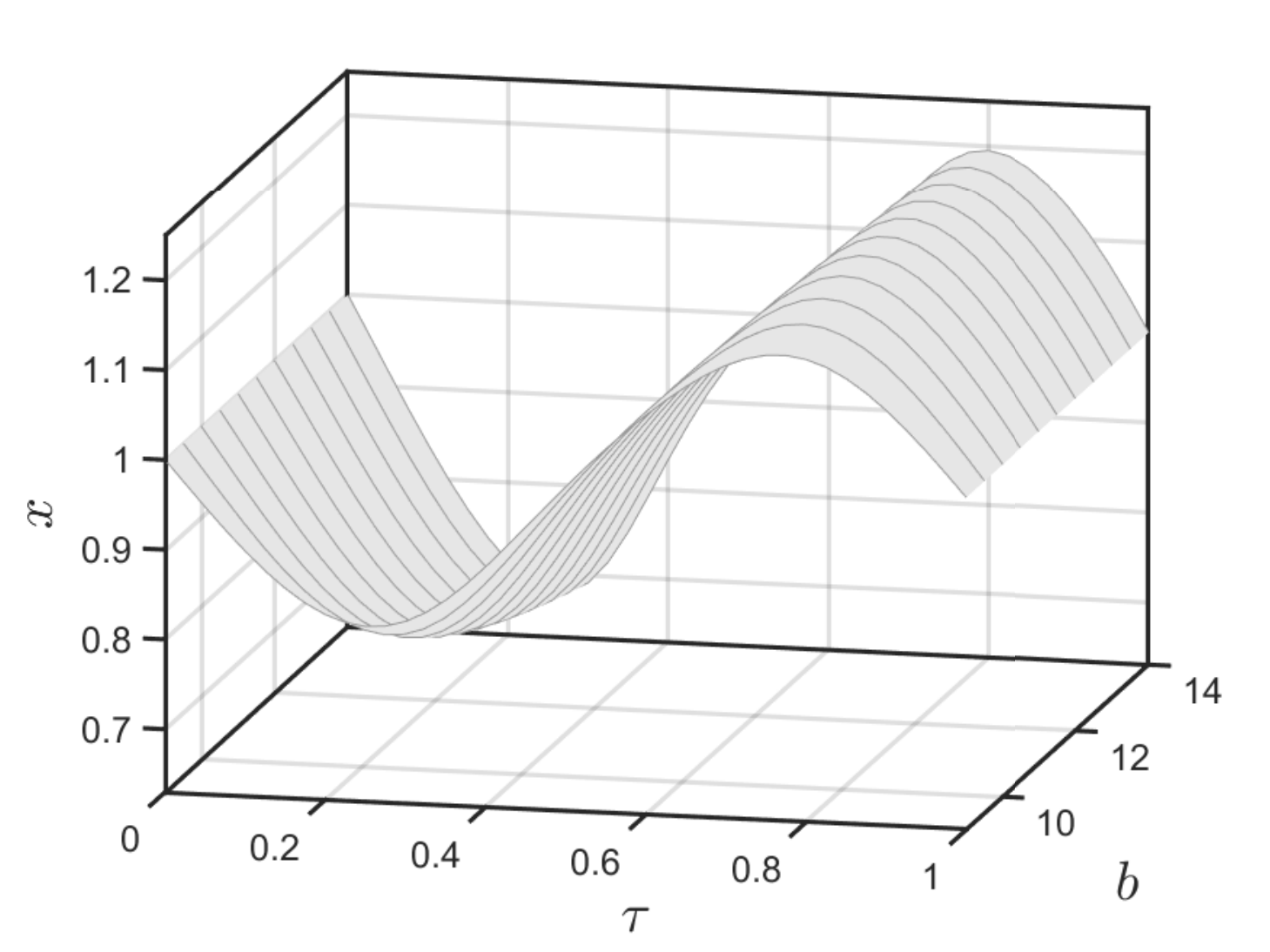}}
\subfloat[]{\includegraphics[width=0.49\columnwidth,trim={0 {0.0\textwidth} 0 0.0\textwidth},clip]{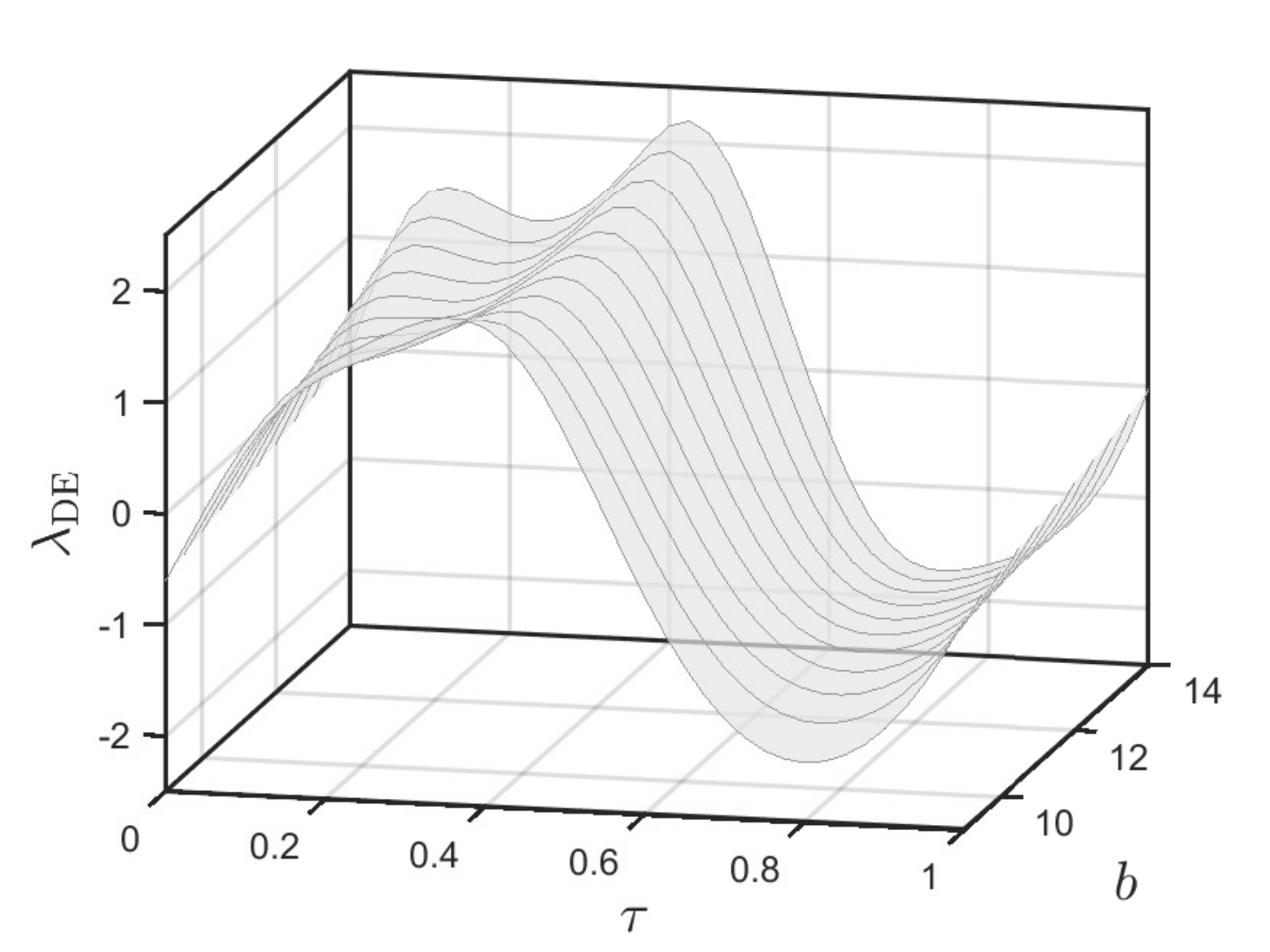}}
\caption{A one-parameter family of time histories for (a) the differential state variable $x(\cdot)$ and (b) corresponding phase response curve $\lambda_\text{DE}(\cdot)$ of the delay differential equation~\eqref{mackey} with $a=2$ and $\alpha=0.7$ under variations in $b$.} 
\label{fig:prc2}
\end{figure}

With $\mu_T$, $\mu_\alpha$, $\mu_a$, and $\mu_b$ fixed, we obtain a reduced continuation problem with dimensional deficit $-1$. With the Lagrange multipliers initially set to $0$, we allow $\mu_T$ to vary and obtain the graphs of $x(\cdot)$ and $\lambda_\text{DE}(\cdot)$ shown in Fig.~\ref{fig:prc1}. These agree with phase-shifted versions of the curves obtained in~\cite{novivcenko2012phase} using backward integration. A family of phase response curves computed using one-dimensional continuation under simultaneous variations of $\mu_T$ and $\mu_b$ is shown in Fig.~\ref{fig:prc2}.

\subsection{Sensitivity analysis of a system with delay}
The regularizing effect of multi-parameter continuation benefits all \textsc{coco} toolboxes. We demonstrate this for the toolbox for delay-coupled systems by applying it to the harmonically forced Duffing oscillator with delayed feedback
\begin{equation}
\label{duffing}
    \ddot{z}+2\zeta\dot{z}+z+z^{3}=\gamma z\left(t-\alpha\right) + a\cos\omega t+b\sin\omega t
\end{equation}
with $\zeta\ll 1$ and feedback gain $\gamma$. Such models with feedback delay arise frequently in experiments when internal actuator dynamics are taken into account 
\cite{sieber2008tracking,wallace2005adaptive}. For this system, as done in Section~\ref{sec: condition numbers}, we analyze the condition numbers associated with different continuation problems along a single frequency response curve. Even though the oscillator with delayed feedback in principle has infinitely many degrees of freedom, the linearization of the problem in the zero solution at zero forcing ($a=b=0$, $\zeta=5\times 10^{-3}$, $\gamma=-0.01$, $z(t)=\dot z(t)=0$) has only a single pair of dominant eigenvalues $\lambda_{1,\pm}\approx -7.9\times 10^{-4}\pm \mathrm{i}$ near the imaginary axis (the eigenvalue pair with the next largest real part is $\lambda_{2,\pm}\approx -9.22\pm3.94\mathrm{i}$). Thus, for small $\gamma$ we expect a situation qualitatively similar to the Duffing oscillator without delayed feedback shown in Fig.~\ref{fig:sens:linear}.

We again follow the decomposition introduced in Section~\ref{sec: toolbox template examples} for a single-segment periodic orbit problem with $T_0=0$ and $T=2\pi/\omega$. 
%the periodic orbit example to obtain
%\begin{align}\label{delay:duffing:de}
%\left(\begin{array}{c}
%x_{1}^{\prime}\\
%x_{2}^{\prime}
%\end{array}\right) &=T\left(\begin{array}{c}
%x_{2}\\
%-2\zeta x_{2}-x_{1}-x_{1}^{3}+ky_{1}+a\cos\omega\left(T_{0}+T\tau\right)+b\sin\omega\left(T_{0}+T\tau\right)
%\end{array}\right),\\
%y_{l}&=\begin{cases}
%x_{l}\left(\tau+1-\alpha/T\right), & \tau\in\left(0,\alpha/T\right),\\
%x_{l}\left(\tau-\alpha/T\right), & \tau\in\left(\alpha/T,1\right),
%\end{cases}\mbox{\quad with $l=1,2$.}\label{delay:duffing:alg}
%\end{align}
%The toolbox developed in Section~\ref{sec: Toolbox construction} may then be applied out of the box for \eqref{delay:duffing:de}--\eqref{delay:duffing:alg} to construct the corresponding zero problems. The algebraic conditions coupling internal interval boundaries, coupling delays and periods are in this example
%\begin{align}\label{delay:duffing:cond}
%     \gamma_{\mathrm{e},1,1}^{(1)}&=\alpha/T\mbox{,}&
%     \varDelta_{1,\ell}&=\alpha/T-j_\ell\mbox{,}& T=2\pi/\omega\mbox{.}
%\end{align} with $j_\ell=\ell\vert_{\operatorname{mod}(2)}$ for $\ell=1,2$. We look for periodic $x$, such that the boundary condition is $x(0)=x(1)$, and we start at initial time $T_{0}=0$.
In order to compute the sensitivity of the Jacobian with respect to the different problem parameters, we  proceed to append monitor functions evaluating to $\omega$, $a$, and $b$ and denote the corresponding continuation parameters by $\mu_{\omega}$, $\mu_{a}$, and $\mu_{b}$. Throughout our analysis, we fix $\zeta=5\times10^{-3}$, $\gamma=-0.01$, and $\alpha=1$. The initial guess for continuation is constructed by first simulating the system dynamics using the \texttt{dde23} solver in \textsc{matlab} with $\omega=0.5$, $a=1.5\zeta$, and $b=0$. A periodic orbit approximation of the differential state variable is then obtained by allowing the system to attain steady state and extracting a segment of length $2\pi/\omega$ from the terminal end of the solution profile. For the algebraic state variable, corresponding to the delayed state, we use linear interpolation to obtain the desired initial guess.

\begin{figure}[ht]
\centering
\subfloat{\includegraphics[width=0.49\columnwidth,trim={0 {0.0\textwidth} 0 0.0\textwidth},clip]{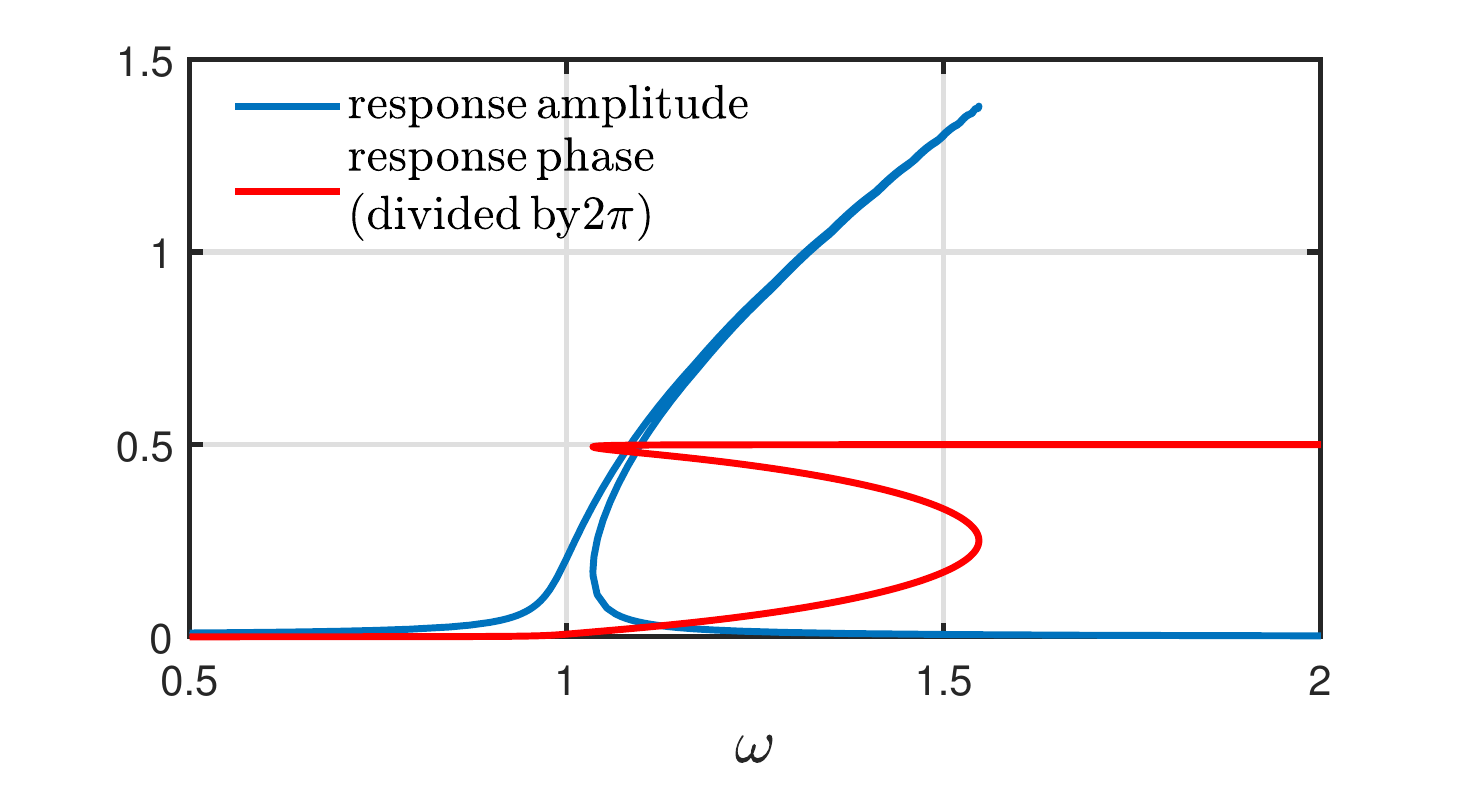}}
\subfloat{\includegraphics[width=0.49\columnwidth,trim={0 {0.0\textwidth} 0 0.0\textwidth},clip]{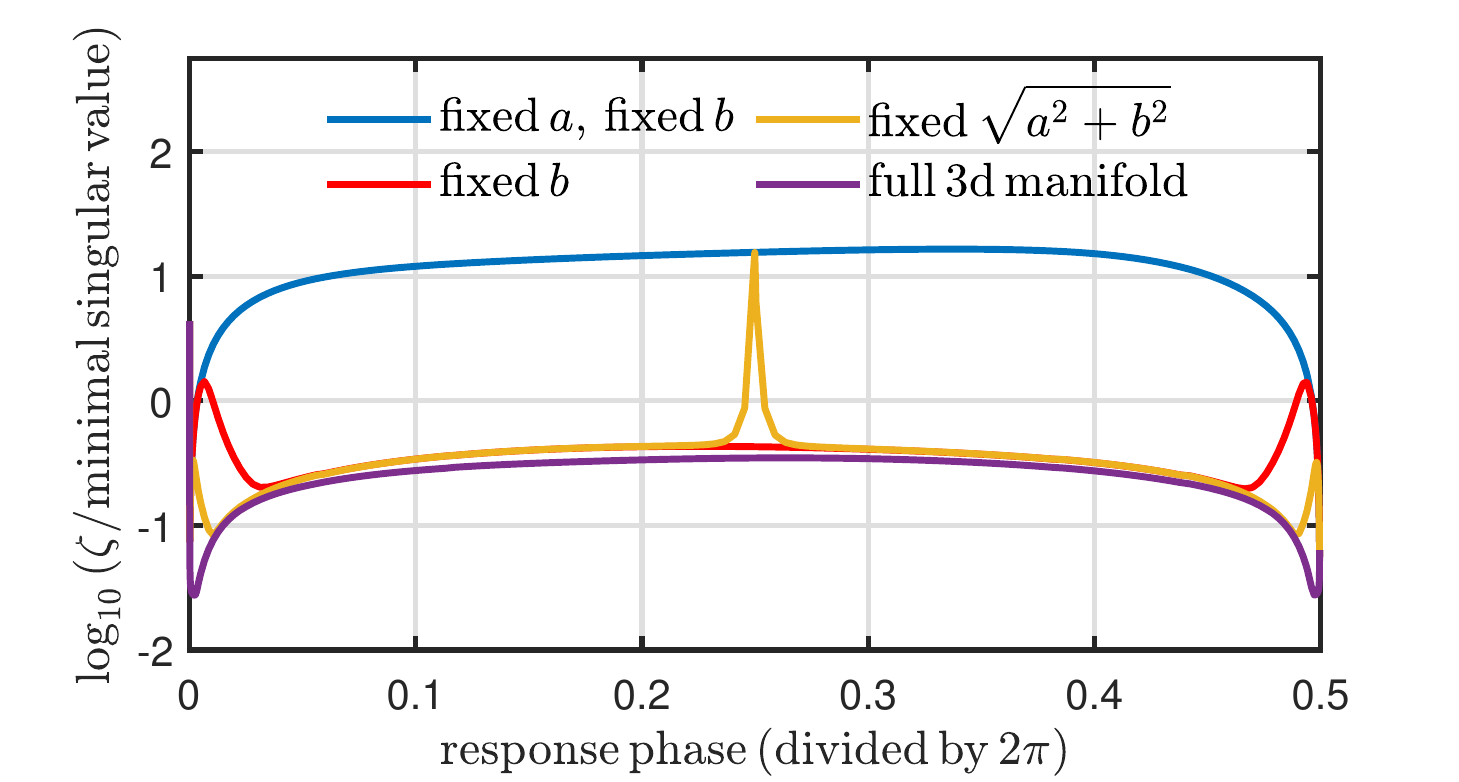}}
\caption{(left panel) Frequency response curve of the harmonically forced delayed duffing oscillator~\eqref{duffing} with $\zeta=5\times 10^{-3}$, $\gamma=-0.01$, $a=1.5\zeta$, $b=0$, and $\alpha=1$. Here, the response amplitude denotes the maximum displacement along the orbit and the response phase is the displacement phase measured relative to the phase of the harmonic forcing. (right panel) Inverse minimal singular values under four different continuation scenarios obtained using the Jacobian of the discretization of the corresponding periodic boundary-value problem.} 
\label{fig:duffing1}
\end{figure}

With $\mu_{\omega}$, $\mu_{a}$, and $\mu_{b}$ fixed, the problem has a dimensional deficit of $0$. We obtain the amplitude and phase response curves shown in the left panel of Fig.~\ref{fig:duffing1} by allowing $\mu_{\omega}$ to vary. Since $\gamma$ is small, these resemble the corresponding shapes for the Duffing oscillator without delay (the dashed curves in Fig.~\ref{fig:sens:linear}). In the right panel, we plot the norm of the inverse of the Jacobians of the problem discretization (scaled by $\zeta^{-1}$) under four different continuation scenarios obtained by fixing both $\mu_{a}$ and $\mu_{b}$, fixing either $\mu_a$ or $\mu_b$, or fixing neither $\mu_a$ nor $\mu_b$. Here, the case with $\mu_a$ fixed is equivalent for the purpose of computing condition numbers to holding the forcing amplitude $\sqrt{a^2+b^2}$ fixed when evaluating the Jacobian at $b=0$ (as was done in Section~\ref{sec: condition numbers}). We observe that the qualitative statements derived for the linear oscillator in Section~\ref{sec: condition numbers} still hold: the norm of the inverse of the Jacobian is of order $\zeta^{-1} N$ across the entire resonance peak when both $\mu_{a}$ and $\mu_{b}$ are fixed. (The factor $N$ is the number of collocation intervals and is here equal to 10.) When we allow $\mu_{a}$ to vary while $\mu_b$ is fixed, the norm of the inverse is independent of $\zeta$ everywhere except for increases near the base of the resonance peak. When $\mu_{b}$ is allowed to vary while $\mu_a$ is fixed, the norm of the inverse (when evaluated at $b=0$) has a near-singularity $\sim \zeta^{-1}N$ near the tip of the resonance peak. When $\mu_{a}$ and $\mu_{b}$ are both free to vary, the norm of the inverse is small throughout the frequency range. This analysis therefore reinforces the observation from Section~\ref{sec: condition numbers} on the regularizing benefits of multi-dimensional continuation.

\begin{figure}[ht]
\centering
\subfloat{\includegraphics[width=0.49\columnwidth,trim={0 {0.0\textwidth} 0 0.0\textwidth},clip]{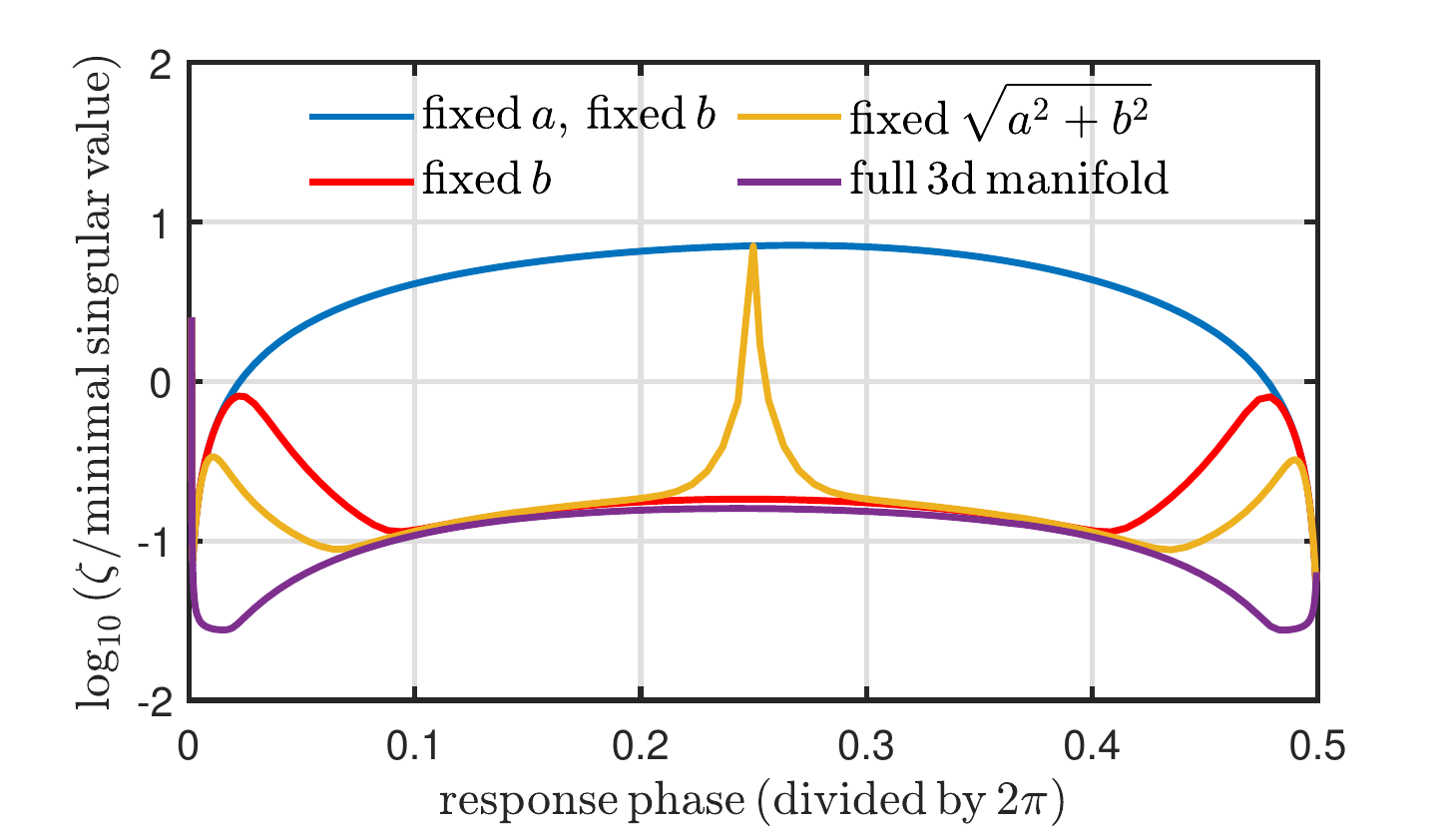}}
\subfloat{\includegraphics[width=0.49\columnwidth,trim={0 {0.0\textwidth} 0 0.0\textwidth},clip]{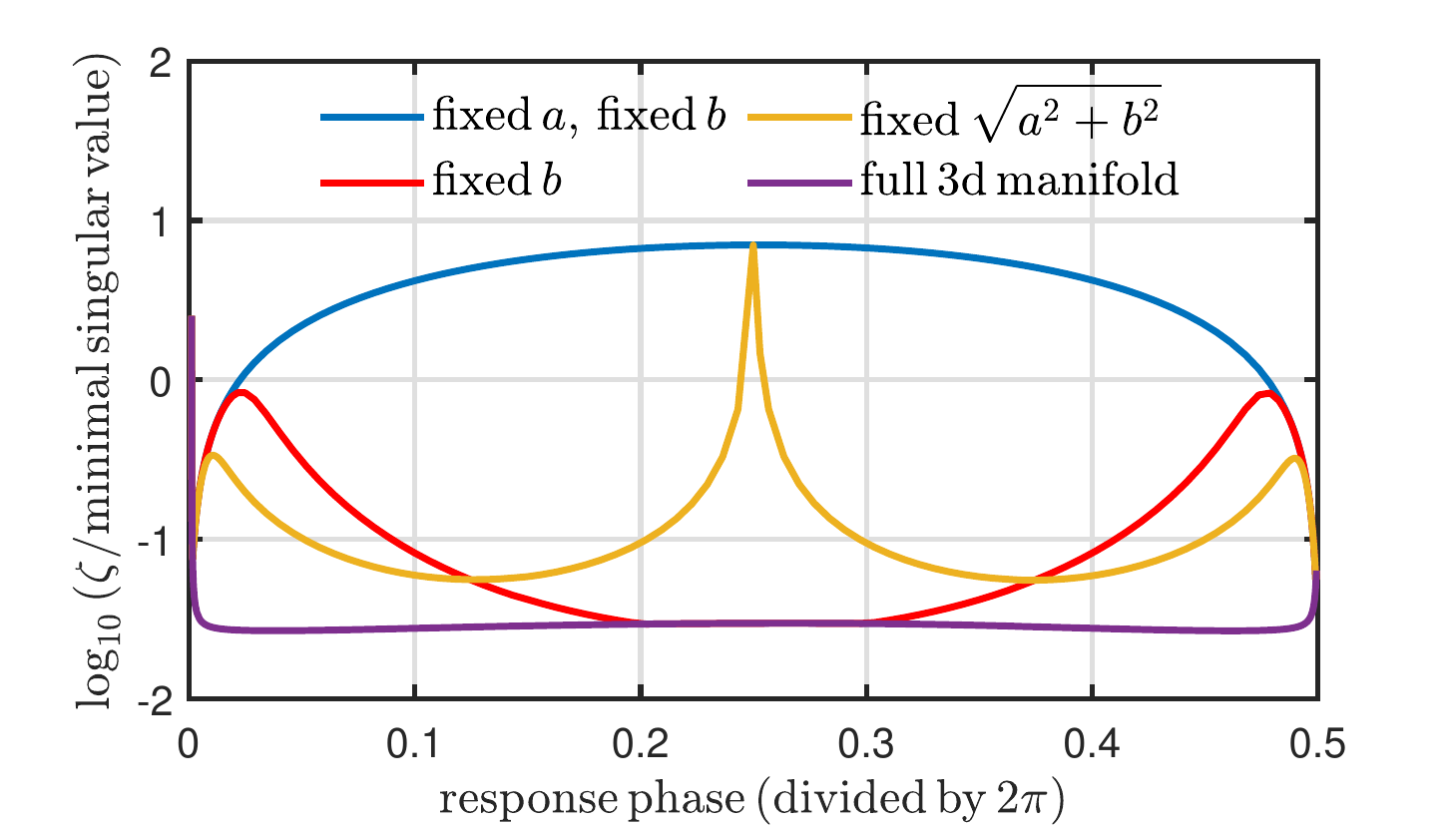}}
\caption{(left panel) Inverse minimal singular values under four different continuation scenarios obtained using the Jacobian of the discretization of the periodic boundary-value problem for the harmonically forced delayed duffing oscillator~\eqref{duffing} with $\zeta=5\times 10^{-3}$, $\gamma=0$, $a=1.5\zeta$, $b=0$, and $\alpha=1$. (right panel) Same analysis as in the left panel, but after division by $N$ of the entries of the Jacobian corresponding to the sensitivity of the algebraic state variables with respect to the coupling delays.} 
\label{fig:duffing2}
\end{figure}

Although we obtain qualitative agreement with the trends for the norms from Fig.~\ref{fig:sens:linear} in Section~\ref{sec: condition numbers}, here we observe pronounced upward bulging of the lower three curves across the resonance peak (for nearly all phases between $0$ and $\pi$). As shown in the left panel of Fig.~\ref{fig:duffing2}, this effect occurs even if we set the feedback gain $\gamma$ to $0$, such that the differential constraint becomes independent of the algebraic state variable. An inspection of the Jacobian matrix reveals that this effect is caused by the derivatives $x'(\tau-\varDelta_k)$ of the coupling conditions $y(\tau)-x(\tau-\varDelta_k)=0$ with respect to the coupling delays $\varDelta_{k}$. The latter are implicitly coupled to the variable period $T$ via the algebraic conditions on $\xi_{\mathrm{e},1}$ and $\xi_{\mathrm{b},2}$ in \eqref{dde:xi:e}-\eqref{dde:xi:b}. After discretization, the solution for the algebraic state variable $y(\tau)$ then exhibits a sensitivity of order $N$ (the number of collocation intervals) with respect to $\varDelta_k$ if $x'(\tau)$ is of order $1$ (as is the case across the resonance peak). When we manipulate the Jacobian before computing the norm of its inverse, first dividing the corresponding terms by $N$, the bulging disappears as shown in the right panel of Fig.~\ref{fig:duffing2} for the case when $\gamma=0$. 
%, leading to the observed bulging. %When we ignore the corresponding column of the Jacobian, the bulging disappears.

\subsection{Optimal control problems}
\label{sec: numerical examples optimal control}
From~\cite{smithoptimal}, we obtain the problem of choosing a control input $u(t)\in\mathbb{R}$ that minimizes the objective functional
\begin{equation}
    J=\int_{0}^{2}\left(z^{2}+u^{2}\right)\text{d}t
    \label{optimal_dde_eq1}
\end{equation}
subject to the initial-value problem
\begin{equation}
\label{optimal_dde_eq2}
    \begin{gathered}
 \dot{z}=tz+ z\left(t-1\right)+u\left(t\right),\,t\in(0,2),\\
 z\left(t\right)=1,\,t\in[-1,0].
    \end{gathered}
\end{equation}
We parameterize the sought optimal control input in terms of a truncated expansion of normalized Chebyshev polynomials of the first kind $T_{c,j}$ defined on the interval $[-1,1]$. Specifically, with $x(\tau):=z(2\tau)$, we obtain
\begin{align}
 x'(\tau)&=4\tau x(\tau)+2y^{(1)}(\tau)+2y^{(2)}(\tau),\,\tau\in\left(0,1\right),\\
y^{(1)}\left(\tau\right) &=\begin{cases}
1, & \tau\in\left(0,1/2\right),\\
x\left(\tau-1/2\right), & \tau\in\left(1/2,1\right),
\end{cases}\\
y^{(2)}\left(\tau\right) &=\sum_{j=1}^q p_j T_{c,j}(2\tau-1),\,\tau\in\left(0,1\right),
\end{align}
and boundary condition $x(0)=1$ in terms of the \emph{a priori} unknown coefficients $p_j$. With the generalization described in Section~\ref{sec: numerical examples generalizations} for matrix-valued algebraic state variables, the toolbox developed in Section~\ref{sec: Toolbox construction} may be applied out of the box to construct the corresponding zero problems and adjoint contributions provided that we append the algebraic conditions $\gamma^{(1)}_{\mathrm{e},1,1}=1/2$, $T_0=0$, $T=2$.

In order to search for optimal choices of the expansion coefficients, we append monitor functions that evaluate to $J$ and $\{p_j\}_{j=1}^q$, as well as the corresponding contributions to the adjoint conditions. We denote the continuation parameters associated with $J$ and $\{p_j\}_{j=1}^q$ by $\mu_J$ and $\{\mu_{p_j}\}_{j=1}^q$ and the corresponding continuation multipliers by $\eta_J$ and $\{\eta_{p_j}\}_{j=1}^q$. Next, we introduce additional complementary monitor functions that evaluate to $\eta_J$ and $\{\eta_{p_j}\}_{j=1}^q$ and denote the corresponding complementary continuation parameters by $\nu_J$ and $\{\nu_{p_j}\}_{j=1}^q$. At the sought extremum, $\nu_J=1$ and $\nu_{p_j}=0$ for $j=1,\ldots,q$. We construct an initial solution guess for the discretization of the differential state variables by integrating the initial-value problem \eqref{optimal_dde_eq2} using the \textit{dde23} solver in \textsc{matlab}. We use linear interpolation to construct an initial solution guess for the discretization of $y^{(1)}(\cdot)$ and initially let $p_j=0$ for $j=1,\ldots,q$. Finally, the Lagrange multipliers are all initialized with zero values.

\begin{figure}[h]
\centering
\subfloat{\includegraphics[width=0.45\columnwidth,trim={0 {0.0\textwidth} 0 0.0\textwidth},clip]{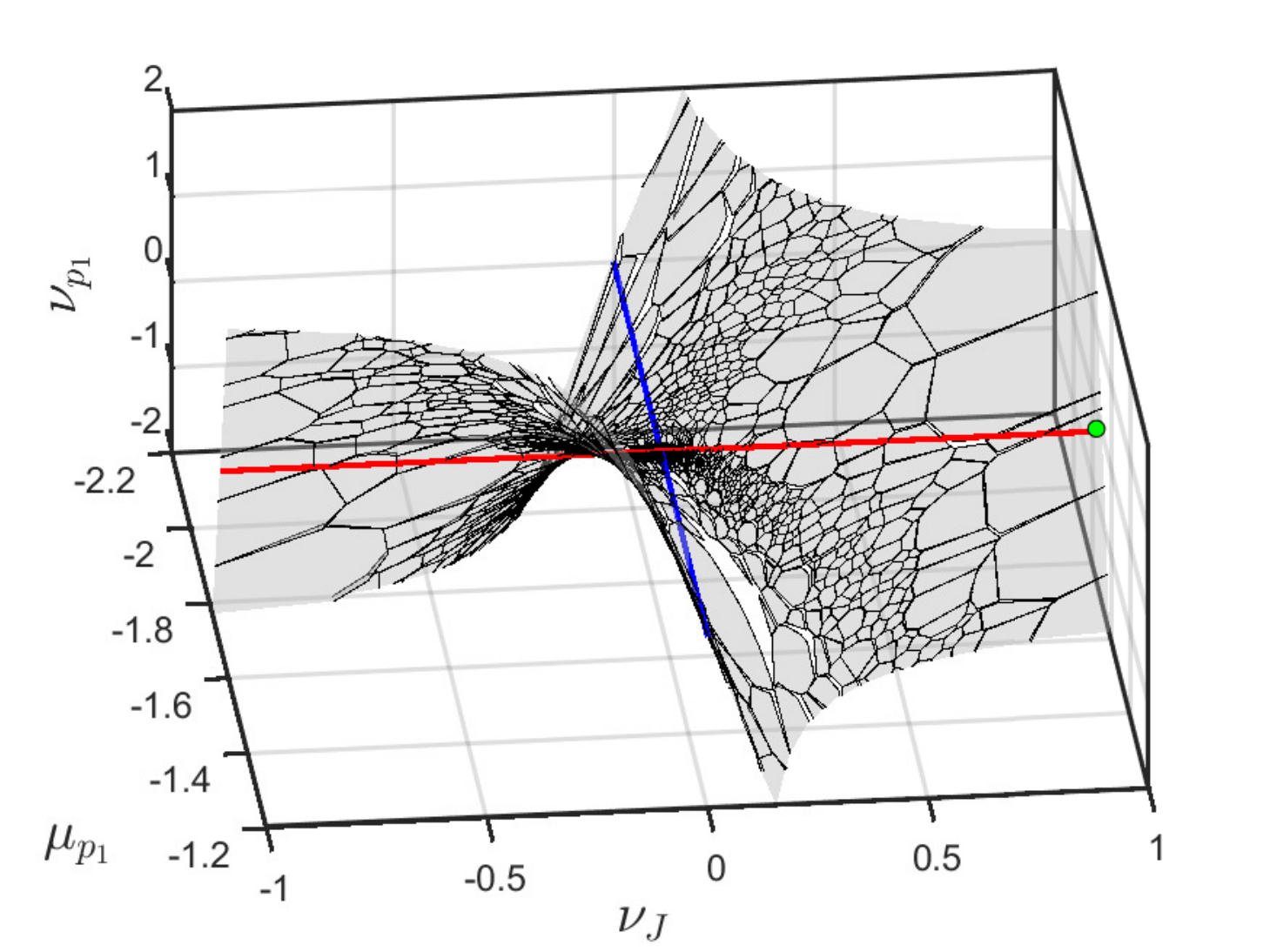}}
\caption{Projection of a two-dimensional solution manifold for the optimal control problem (\ref{optimal_dde_eq1}--\ref{optimal_dde_eq2}) and the corresponding adjoint equations obtained using the \mcode{atlas\_kd} atlas algorithm by allowing $\mu_J$, $\nu_J$, $\mu_{p_1}$, and $\{\nu_{p_j}\}_{j=1}^q$ to vary while holding $\{\mu_{p_j}\}_{j=2}^q$ fixed. As in Section~\ref{sec: inflection points}, the zero-level curves of $\nu_{p_1}$ on this manifold are two straight lines with $\nu_J=0$ (blue) and $\mu_{p_1}$ (red) constant, respectively, that intersect at a stationary point of $\mu_J$ along the first curve. The solution with $\nu_J=1$ and $\nu_{p_1}=0$ (green circle) can be located (to within desired tolerance) by continuation along the first of these straight lines, followed by branch switching and continuation along the second of these lines. Alternatively, it may be approximated by the solution point on the intersection of the two-dimensional manifold with the $\nu_J=1$ coordinate plane (located within desired tolerance) with $\nu_{p_{1}}$ closest to zero.
}
\label{fig1:optimal}
\end{figure}

With $\mu_J$, $\{\mu_{p_j}\}_{j=1}^q$, and $\{\nu_{p_j}\}_{j=1}^q$ fixed, the dimensional deficit equals $-q-1$. We obtain the two-dimensional solution manifold in Fig.~\ref{fig1:optimal} by allowing $\mu_J$, $\nu_J$, $\mu_{p_1}$, and $\{\nu_{p_j}\}_{j=1}^q$ to vary while holding $\{\mu_{p_j}\}_{j=2}^q$ fixed. We select the point with $\nu_{J}=1$ and $\nu_{p_{1}}$ closest to $0$ as initial solution guess for a second stage of continuation obtained by fixing $\nu_{J}$ and $\nu_{p_{1}}$ to $1$ and $0$, respectively, and allowing $\mu_{p_2}$ to vary. If we locate a point with $\nu_{p_2}=0$ along the corresponding one-dimensional manifold, then we may continue from this point along a new one-dimensional manifold obtained by fixing $\nu_{p_{2}}$ at $0$ and allowing $\mu_{p_{3}}$ to vary. Continuing in this fashion, we locate the sought extremum.

\begin{figure}[h!]
\centering
\subfloat{\includegraphics[width=0.49\columnwidth,trim={0 {0.0\textwidth} 0 0.0\textwidth},clip]{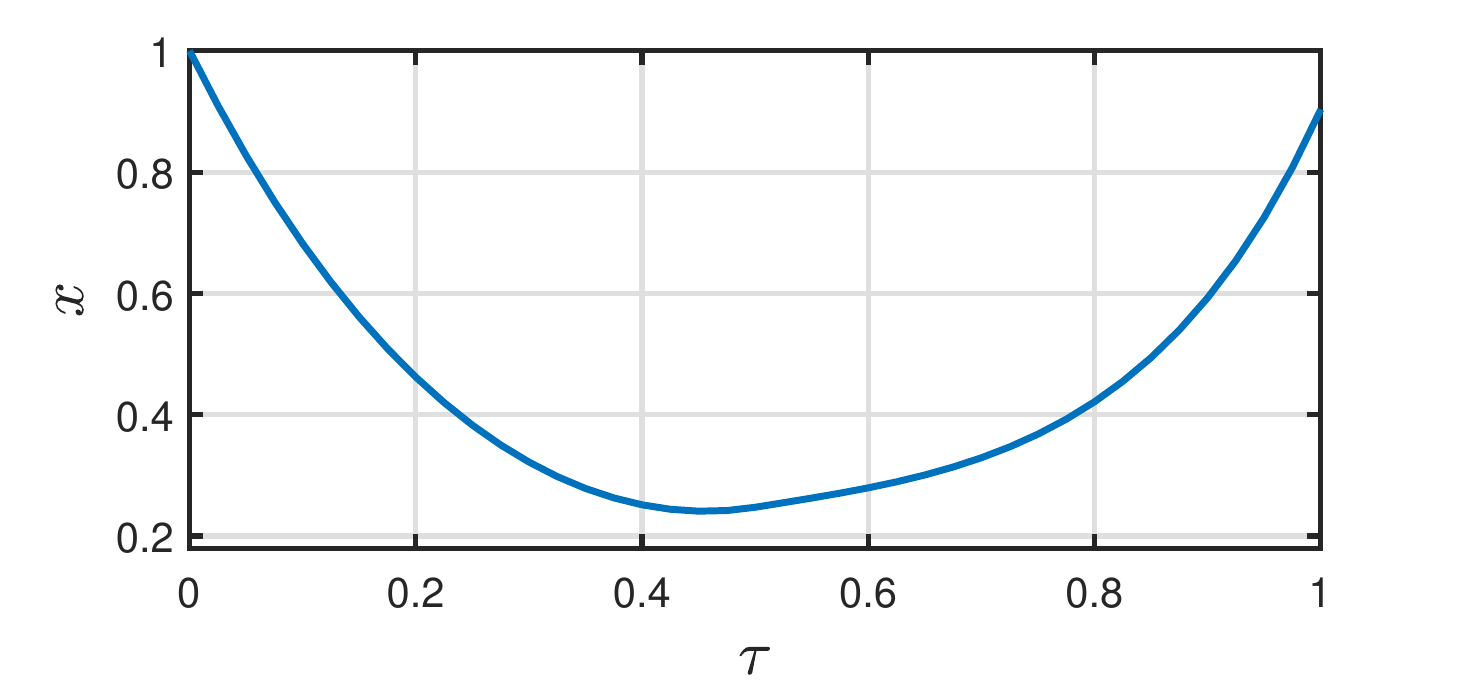}}
\subfloat{\includegraphics[width=0.49\columnwidth,trim={0 {0.0\textwidth} 0 0.0\textwidth},clip]{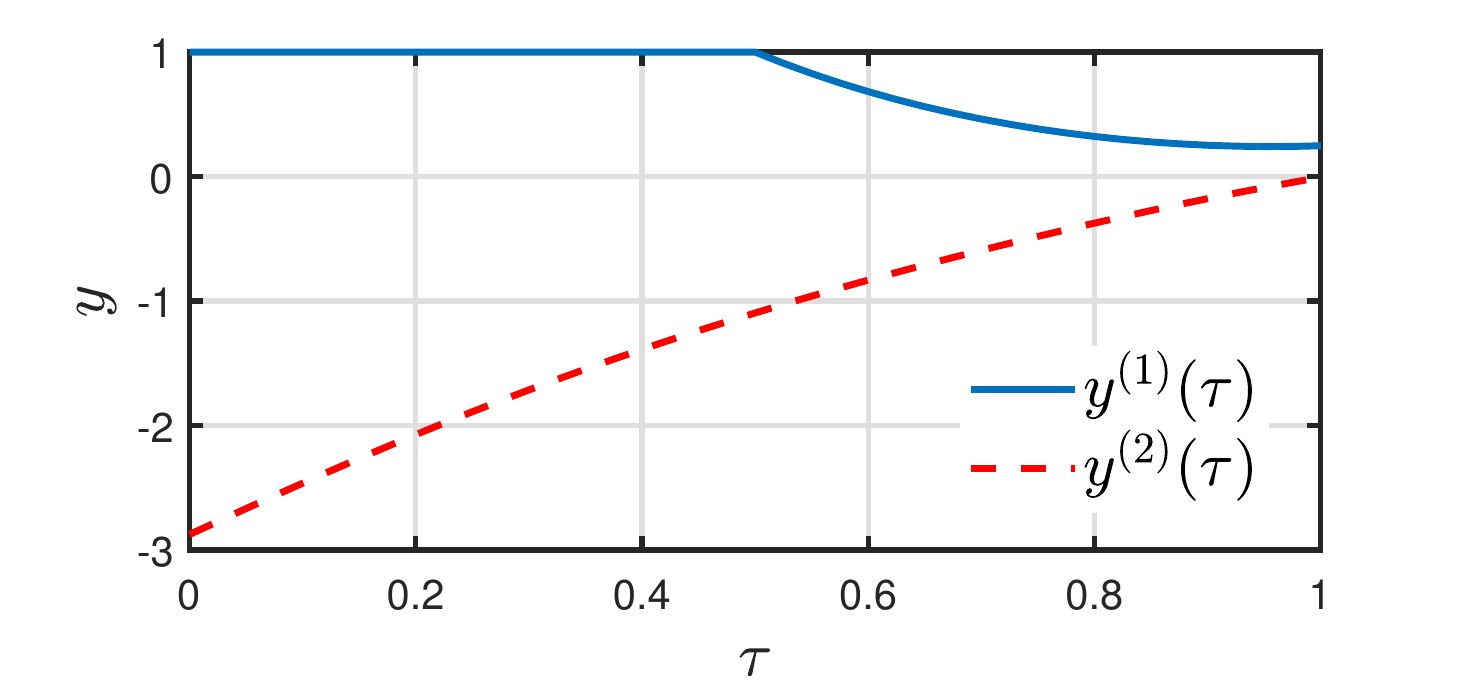}}\\
\subfloat{\includegraphics[width=0.49\columnwidth,trim={0 {0.0\textwidth} 0 0.0\textwidth},clip]{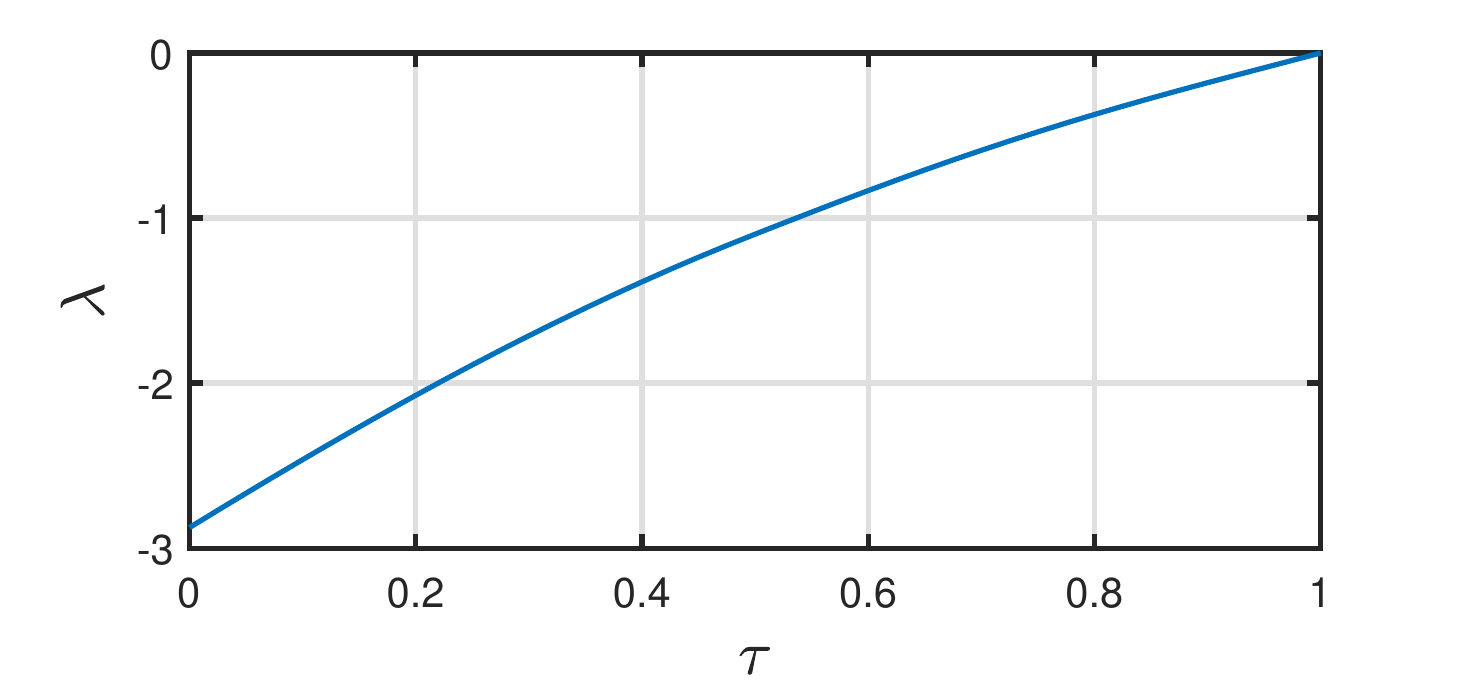}}
\subfloat{\includegraphics[width=0.49\columnwidth,trim={0 {0.0\textwidth} 0 0.0\textwidth},clip]{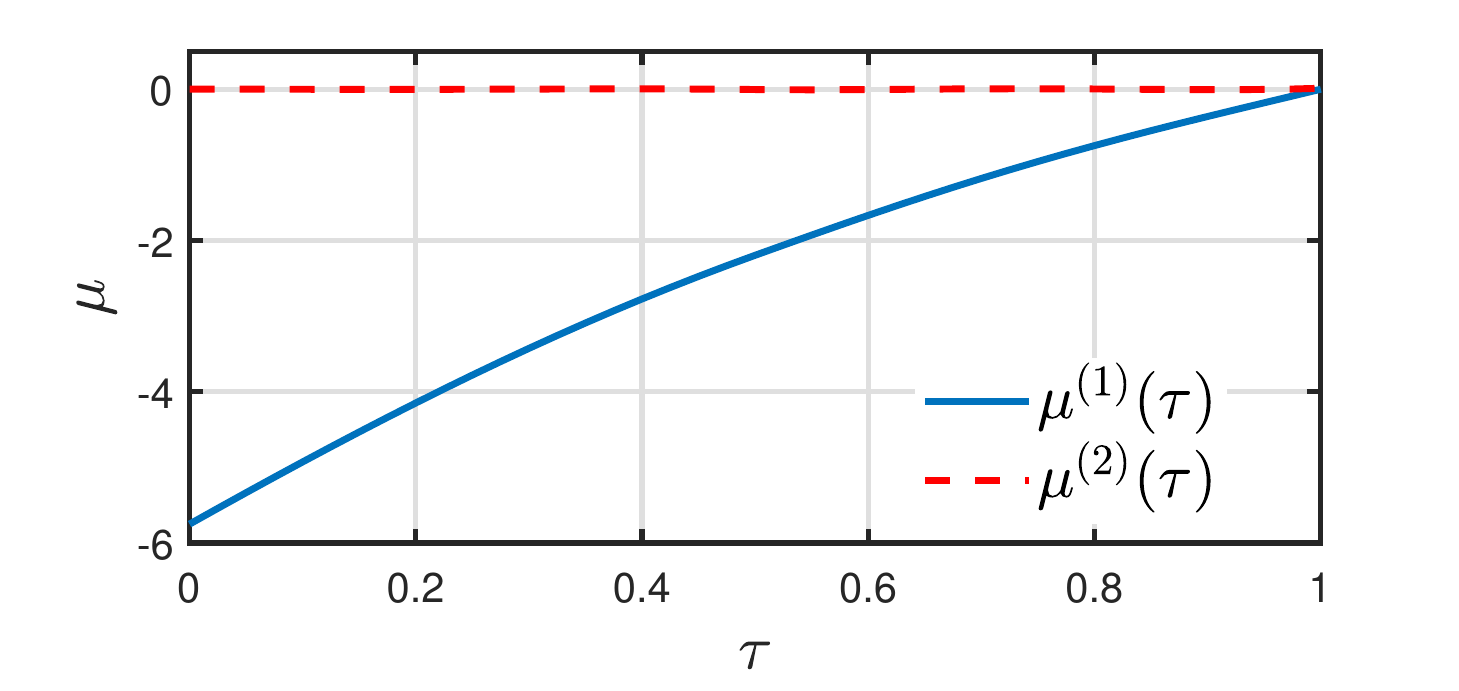}}
\caption{State variables and Lagrange multipliers at a local minimum for the optimal control problem (\ref{optimal_dde_eq1}--\ref{optimal_dde_eq2}).} 
\label{fig2:optimal}
\end{figure}

Figure~\ref{fig2:optimal} shows the locally optimal solution obtained using this methodology with $q=8$ (here, $N=10$ and $m=4$). At this solution, $J=J_8\approx 4.797$, is in excellent agreement with the minimum reported in ~\cite{smithoptimal} using stochastic optimization. Table~\ref{tbl1:optimal} demonstrates the rapid convergence of the locally optimal value $J=J_q$ with the expansion order $q$ anticipated from the use of Chebyshev polynomials and the smoothness of the optimal solution.

\begin{table}[h]
    \caption{The error $\left|\left(J_{q}-J_8\right)/J_8\right|$ decays rapidly with increasing truncation order $q$ for the optimal control problem~(\ref{optimal_dde_eq1})-(\ref{optimal_dde_eq2}).}
    \centering
    {\begin{tabular}{cccccccc}%
        \hline
        \hline
        $q$ & $1$ & $2$ & $3$ & $4$ & $5$ & $6$ & $7$  \\  	
        \hline
				Error & $4.3\times 10^{-1}$ & $4.8\times 10^{-3}$ & $9.7\times 10^{-6}$ & $2.6\times 10^{-5}$ & $2.7\times 10^{-5}$ & $2.8\times 10^{-5}$ & $8\times 10^{-6}$ \\
        \hline
        \hline
    \end{tabular}}
    \label{tbl1:optimal}
\end{table}

\subsection{Quasiperiodic orbits}
\label{sec: numerical examples quasiperiodic orbits}
As a final example, studied previously in~\cite{ahsan2020optimization}, we consider the problem of locating stationary values of $\omega$ along a family of quasiperiodic invariant tori of the delay differential equations
\begin{equation}
\label{quasi_ex}
\dot{z}(t)=f(t,z(t),z(t-\alpha),p):=\left(\begin{array}{c}
-\omega z_{2}(t)+z_{1}\left(t-\alpha\right)\left(1+r(t)\left(\cos2\pi t/T-1\right)\right)\\
\omega z_{1}(t)+z_{2}\left(t-\alpha\right)\left(1+r(t)\left(\cos2\pi t/T-1\right)\right)
\end{array}\right) 
\end{equation}
with $r=\sqrt{z_{1}^2+z_{2}^2}$. By applying the analysis of the third example in Section~\ref{sec: toolbox template examples} in reverse, we obtain the delay-coupled multi-segment boundary-value problem in \eqref{eq:quasi1}-\eqref{sample:eq2}.
The toolbox developed in Section~\ref{sec: Toolbox construction} may be applied out of the box to construct the corresponding zero problem and adjoint contributions provided that we append the algebraic conditions $T_{0,i}=0$, $T_i=T$, and $\gamma^{(1)}_{\mathrm{e},i,1}=\alpha/T$ for $i=1,\ldots,M$. The additional phase condition $x_{1,2}(0)=1$ restricts attention to a unique family of trajectory segments discretizing a quasiperiodic invariant torus.

In order to search for extremal values of $\omega$, we proceed to append monitor functions that evaluate to $\omega,T$ and $\alpha$, respectively, and denote the corresponding continuation parameters by $\mu_{\omega}$, $\mu_{T}$, and $\mu_{\alpha}$. We denote the continuation multipliers associated with the corresponding adjoint contributions by $\eta_{\omega},\eta_{T}$ and $\eta_{\alpha}$. Next, we introduce complementary monitor functions that evaluate to $\eta_{\omega},\eta_{T},\eta_{\alpha}$ and denote the corresponding complementary continuation parameters by $\nu_{\omega},\nu_{T}$ and $\nu_{\alpha}$. At a stationary point of $\omega$, we must have $\nu_{\omega}=1,\nu_{T}=\nu_{\alpha}=0$.

For the special case with $\alpha=0$ (this case is considered in the tutorial documentation for the \textsc{coco}-compatible \mcode{coll} toolbox), a one-parameter family of quasiperiodic orbits covering an invariant torus for the delay differential equation~\eqref{quasi_ex} is given by
\begin{equation}
z(T_0+T\tau;\varphi)=x(\varphi,\tau)=\left(r\left(\tau\right)\cos\left(2\pi\varrho\tau+\varphi\right),r\left(\tau\right)\sin\left(2\pi\varrho\tau+\varphi\right)\right),\,\varphi\in\mathbb{S}
\end{equation}
where
\begin{equation}
r\left(\tau\right)=\frac{1+\Omega^{2}}{1+\Omega^{2}-\cos 2\pi\tau-\Omega\sin 2\pi\tau},
\end{equation}
$\Omega=2\pi/T$, and $\varrho$ is an irrational number. With the substitution $x_i(\tau)=x(2\pi(i-1)/M,\tau)$ we obtain an initial solution guess for the $i$-th trajectory segment of the corresponding boundary-value problem. Here, we initialize continuation with $\omega=1$ and $\Omega=1.5$ and set the rotation number $\rho$ to $0.6618$. We use the initial guess for the problem with $\alpha=0$ along with continuation in $\alpha$ to obtain the differential state variables for nonzero $\alpha$.

\begin{figure}[ht!]
\centering
\subfloat[]{\includegraphics[width=0.45\columnwidth,trim={0 {0.0\textwidth} 0 0.0\textwidth},clip]{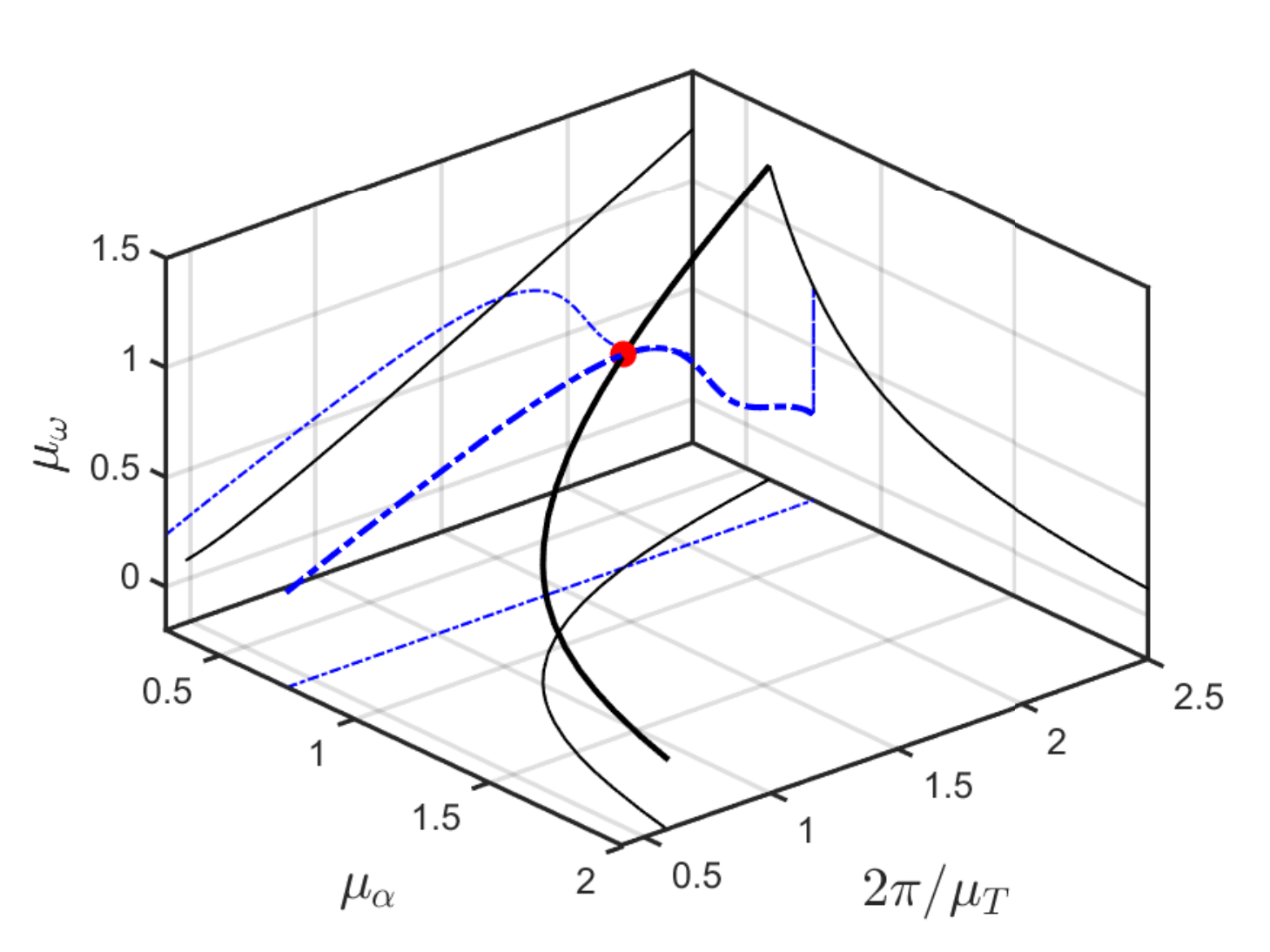}}\\
\subfloat[]{\includegraphics[width=0.45\columnwidth,trim={0 {0.0\textwidth} 0 0.0\textwidth},clip]{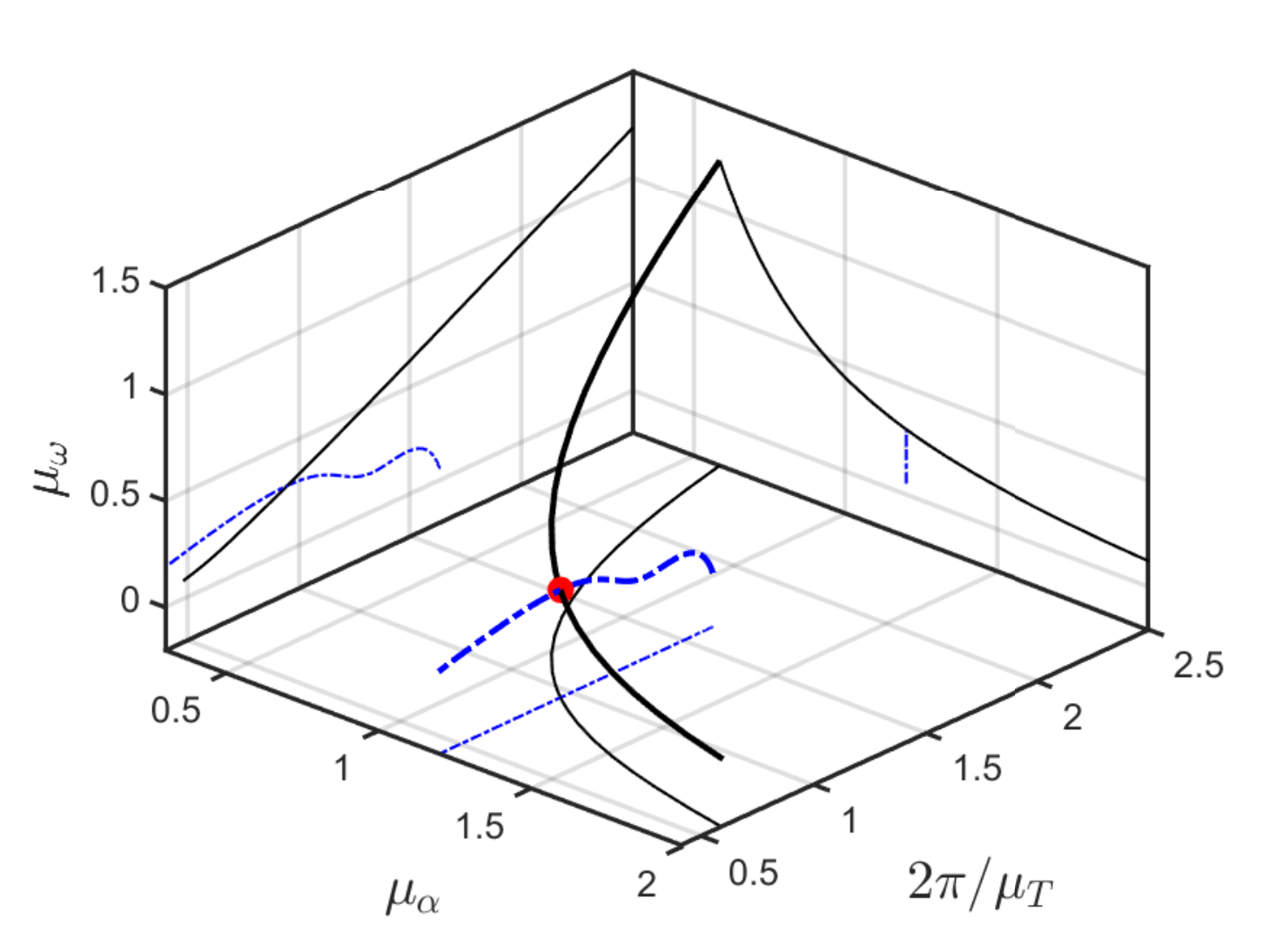}}
\subfloat[]{\includegraphics[width=0.45\columnwidth,trim={0 {0.0\textwidth} 0 0.0\textwidth},clip]{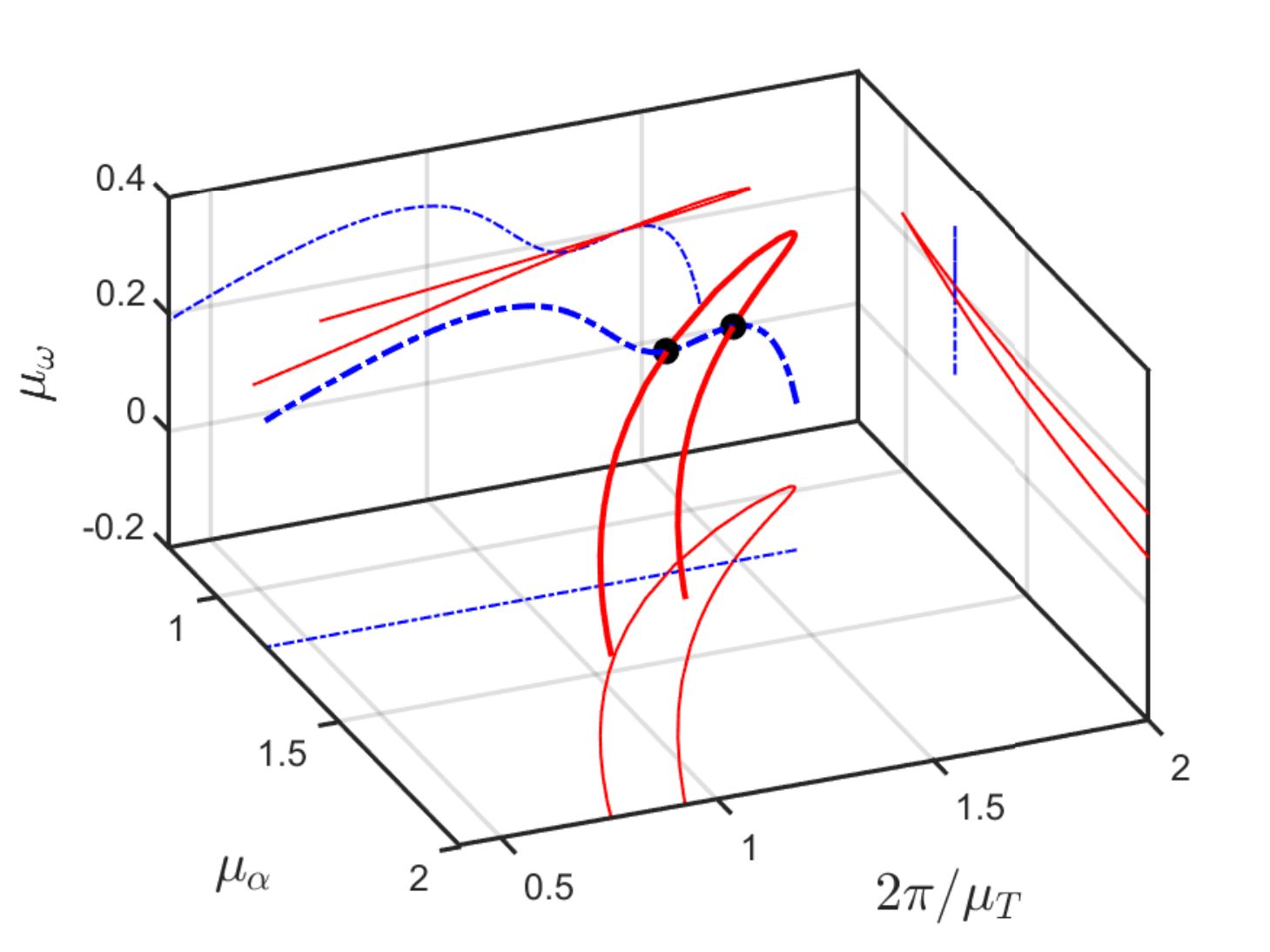}}
\caption{Curves (solid) of local stationary points of $\mu_\omega$ with respect to variations in $\mu_{T}$ along approximate families of quasiperiodic invariant tori for the delay differential equation \eqref{quasi_ex} obtained by first continuing along a one-dimensional manifold with trivial Lagrange multipliers and fixed $\mu_\alpha$ (dashed), then switching at a local stationary point to a branch with varying Lagrange multipliers, and finally fixing $\nu_\omega$ at $1$ and allowing $\mu_{T}$ to vary. (a) Initial continuation with $\mu_{\alpha}=0.75$ and branch switching at a unique local maximum. (b) Initial continuation with $\mu_{\alpha}=1.2$ and branch switching at one of the two local maxima. The final manifold coincides with that obtained in (a). (c) Initial continuation with $\mu_{\alpha}=1.2$ and branch switching at either of the local minimum or other local maximum. The final manifolds obtained in these two cases coincide.}
\label{fig:quasi1}
\end{figure}

With $\mu_{\omega},\mu_{T},\mu_{\alpha},\nu_{\omega},\nu_{T},\nu_{\alpha}$ fixed, the problem has a dimensional deficit of $-3$. We obtain a dimensional deficit of $1$ by allowing  $\mu_{\omega},\mu_{T},\nu_{\omega},\nu_{\alpha}$ to vary. The dashed-dotted line in Fig.~\ref{fig:quasi1}(a) shows the corresponding solution manifold with trivial Lagrange multipliers obtained with $\mu_{\alpha}=0.75$. Through the local maximum in $\omega$ (denoted by the red sphere) at $\left(\mu_{\omega}^{\star},\mu_{T}^{\star},\mu_{\alpha}^{\star}\right) \approx \left(0.773,3.601,0.75\right)$ runs a secondary branch of solutions along which $\nu_{\omega}$ (and the other Lagrange multipliers) vary. We locate a unique point on this branch with $\nu_{\omega}=1$. Continuation along the one-dimensional solution manifold through this point obtained by fixing $\nu_{\omega}$ at $1$ and allowing $\mu_\alpha$ to vary yields a family of local maxima of $\omega$ under variations in $T$ (solid black curve in Fig.~\ref{fig:quasi1}(a)). Within the chosen computational domain, we do not find a point along this manifold where $\nu_{\alpha}=0$. 

The successive continuation approach used here may yield different families of stationary points depending on the initial solution guess and the order in which initially-fixed continuation parameters are allowed to vary. For example, if we repeat the above analysis with an initial value of $\mu_{\alpha}=1.2$, then the one-dimensional solution manifold with trivial Lagrange multipliers has two local maxima and one local minimum in $\mu_\omega$ within the computational domain, as depicted by the dashed-dotted line in Fig.~\ref{fig:quasi1}(b,c). The red sphere in Fig.~\ref{fig:quasi1}(b) denotes one such local maximum at $\left(\mu_{\omega}^{\star},\mu_{T}^{\star},\mu_{\alpha}^{\star}\right)\approx \left(0.305,6.62,1.2\right)$. We may again switch to a secondary branch through this point in order to locate a point with $\nu_\omega=1$ and then continue along a one-dimensional manifold with $\nu_\omega$ fixed at $1$ and varying $\mu_\alpha$ corresponding to a family of stationary points of $\mu_\omega$ under variations in $\mu_{T}$. The latter manifold is found to equal that obtained in the previous paragraph. The identical methodology applied to the other stationary points on the original solution branch located at $\left(\mu_{\omega}^{\star},\mu_{T}^{\star},\mu_{\alpha}^{\star}\right)\approx \left(0.202,4.229,1.2\right)$ and $\left(0.181,4.724,1.2\right)$ (black spheres in Fig.~\ref{fig:quasi1}(c)) yields a single curve of stationary points  (solid red curve in Fig.~\ref{fig:quasi1}(c)) of $\mu_\omega$ under variations in $\mu_{T}$. Although $\mu_\omega$ achieves a local maximum along this manifold, $\nu_{\alpha}$ does not equal $0$ at any point along this curve, and we again conclude that there does not exist a stationary point of $\mu_\omega$ with respect to variations in both $\mu_{T}$ and $\mu_\alpha$ within the chosen computational domain. 

\begin{figure}[ht]
\centering
\subfloat{\includegraphics[width=0.45\columnwidth,trim={0 {0.0\textwidth} 0 0.0\textwidth},clip]{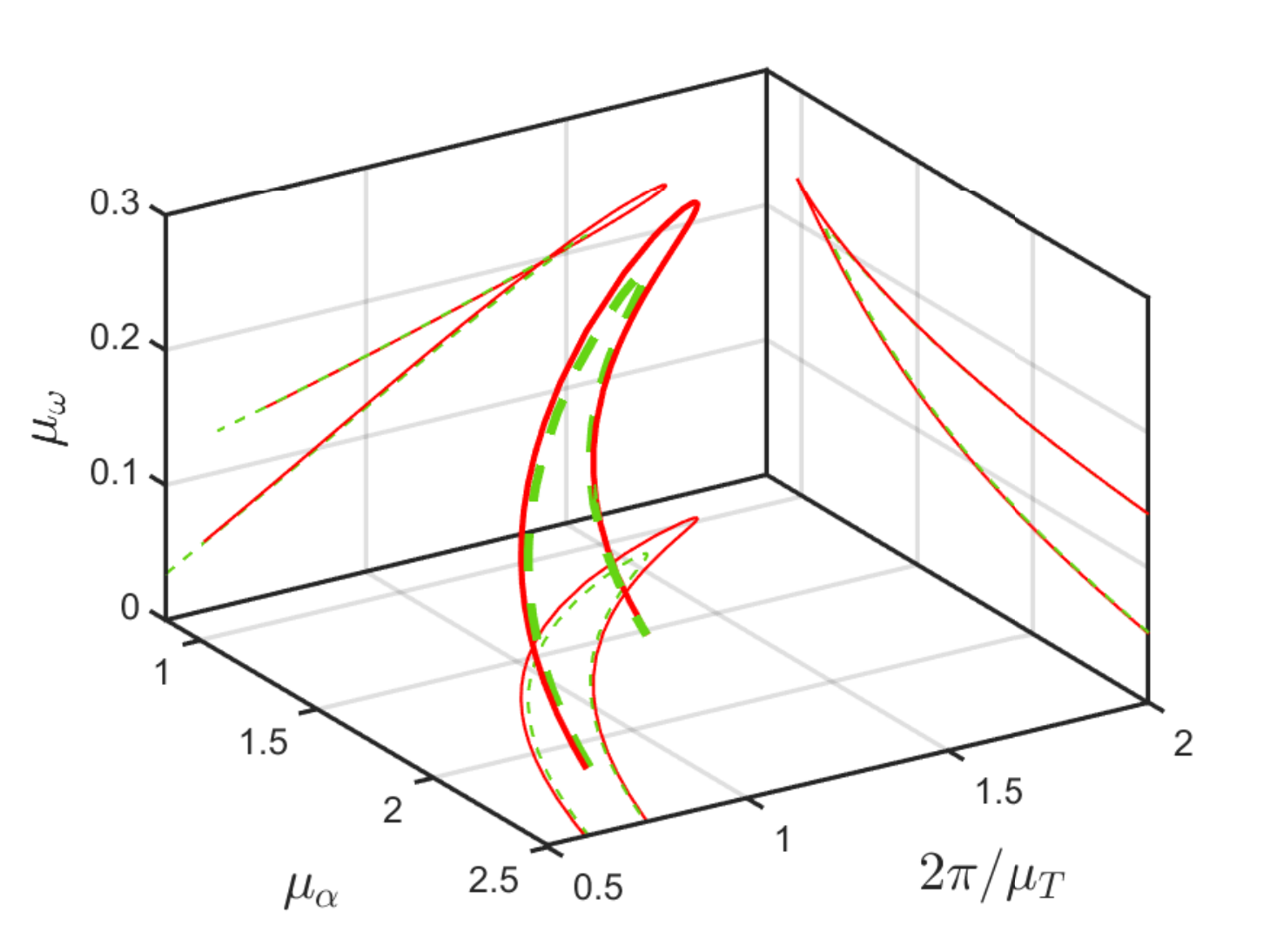}}
\caption{Curves of stationary points of $\mu_{\omega}$ with respect to $\mu_{T}$ (solid) and $\mu_{\alpha}$ (dashed) along approximate families of quasiperiodic invariant tori for the delay differential equation~\eqref{quasi_ex}. The two curves never intersect within the chosen computational domain.} 
\label{fig:quasi2}
\end{figure}

\begin{figure}[ht]
\centering
\subfloat[View 1]{\includegraphics[width=0.49\columnwidth,trim={0 {0.0\textwidth} 0 0.0\textwidth},clip]{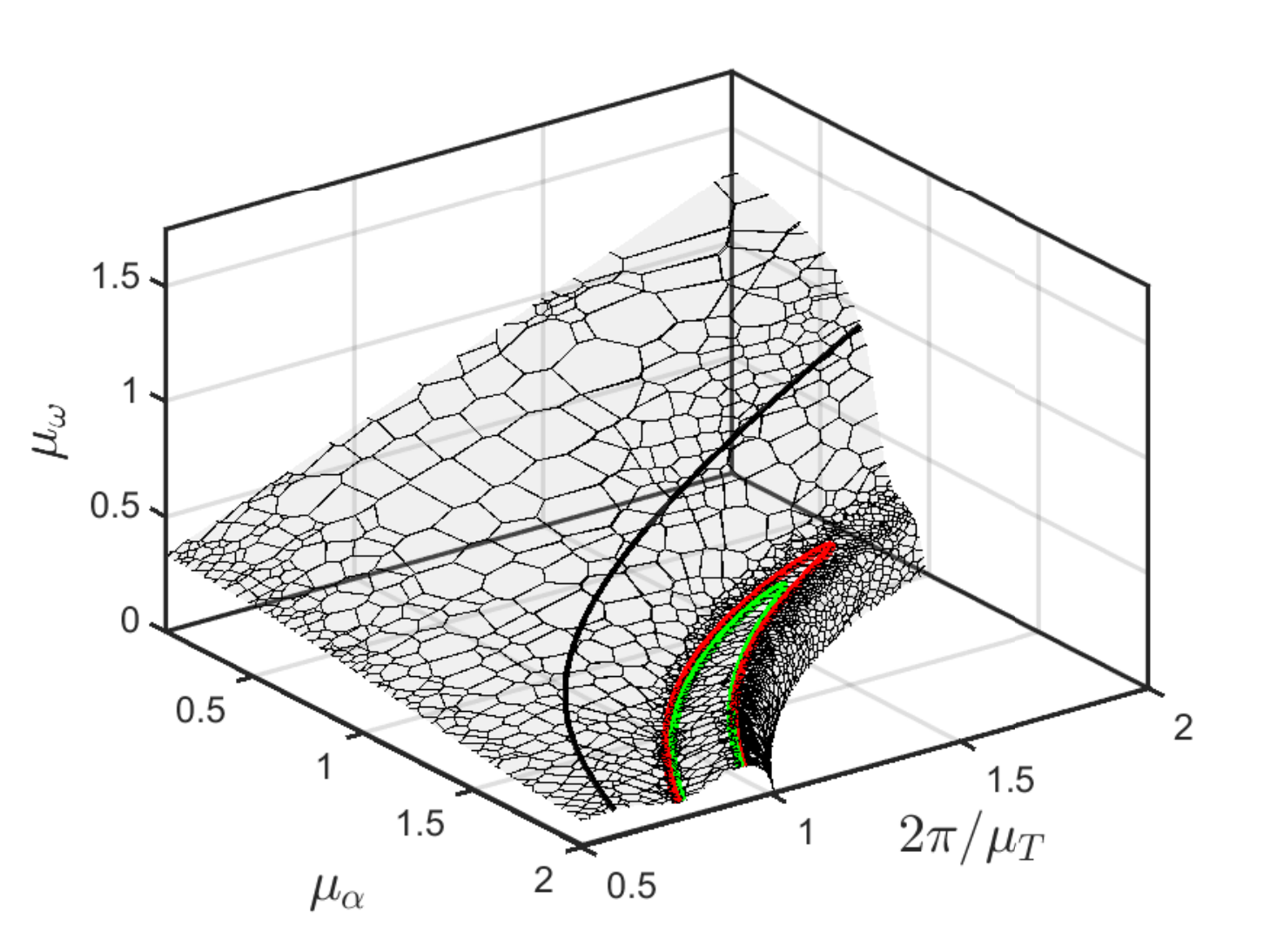}}
\subfloat[View 2]{\includegraphics[width=0.49\columnwidth,trim={0 {0.0\textwidth} 0 0.0\textwidth},clip]{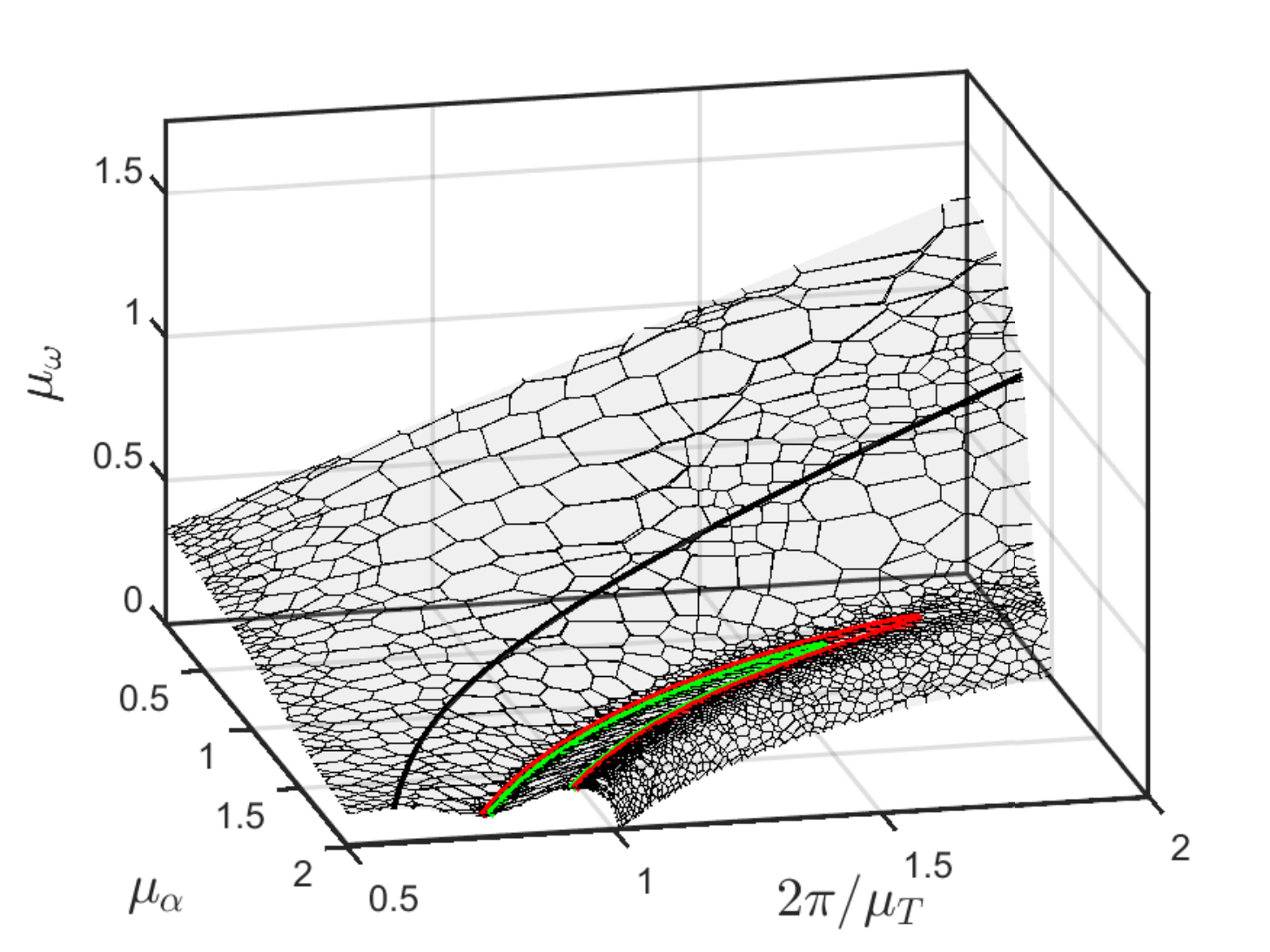}}
\caption{An approximate family of quasiperiodic invariant tori for the delay differential equation~\eqref{quasi_ex} obtained using two-dimensional continuation in \textsc{coco}. Highlighted curves consist of stationary points of $\mu_{\omega}$ with respect to $\mu_{T}$ (red and black) and stationary points of $\mu_{\omega}$ with respect to $\mu_{\alpha}$ (green) located using the corresponding augmented continuation problem.} 
\label{fig:quasi3}
\end{figure}

The same conclusion is obtained by holding $\mu_{T}$ fixed initially while allowing $\mu_\alpha$ to vary along a one-dimensional manifold with trivial Lagrange multipliers. For example, with $\mu_{T}$ fixed at $8.37$, we locate a local minimum and a local maximum in the value of $\mu_\omega$ at $\mu_\alpha=2.04$ and $2.49$. Through each of these points runs a secondary branch along which $\nu_\omega$ varies. The unique point with $\nu_\omega=1$ on each of these secondary branches may be used to continue along a one-dimensional solution manifold with $\nu_\omega$ fixed at $1$ and $\mu_{T}$ allowed to vary. Such a manifold corresponds to a curve of stationary points of $\mu_\omega$ with respect to $\mu_\alpha$. As shown in Fig.~\ref{fig:quasi2}, we obtain a single such curve which does not intersect the curve of stationary points of $\mu_\omega$ with respect to $\mu_{T}$ at any point in the computational domain. These observations are consistent with the visualization in Fig.~\ref{fig:quasi3} of the two-dimensional solution manifold obtained by allowing $\mu_\omega$, $\mu_\alpha$, and $\mu_{T}$ to vary.

\section{Further development}
\label{sec: Further development}

We conclude this paper with a discussion of opportunities for further development along the lines described in the previous sections, some of which are ongoing work. These pertain to the derivation of a compact data structure that uniquely defines a zero problem of the form considered in Section~\ref{sec: toolbox template zero problems} and the corresponding adjoint contributions, as well as modifications to support adaptive discretization.

\subsection{Construction from delay graphs} 
\label{sec: further developments graphs}
The collection of differential and algebraic constraints provided in Section~\ref{sec: toolbox template zero problems} are still one step away from the form in which a user typically formulates a multi-segment boundary-value problem with delay(s). The discussion in Section~\ref{sec: delay graphs} provides a template for how to automate this step by encapsulating the data required for a toolbox to construct the associated boundary-value problem and adjoint contributions using a suitable graph representation. In this section, we propose a general theory that is compatible with the abstract toolbox template and apply this to the examples in Section~\ref{sec: toolbox template examples}.

As before, associate with each segment
\begin{compactitem}
\item a duration $T>0$: we develop all notation for unscaled intervals $[0,T]$ and assume that the involved differential equations are autonomous for simpler notation; the scaling can be performed in a separate final step;
\item a delay $\alpha\geq0$: we only consider a single delay per segment for simpler notation;
\item a variable $x:[-\hat{\alpha},T]\to\R^{n}$ with
  $\hat{\alpha}\geq\alpha$ that is involved in the temporal coupling between segments;
\item directed links to its \emph{predecessors}; associated to a link
  from node $i$ to node $j$ is a non-zero coupling
  matrix $B_{ij}$, such that
  \begin{align}\label{network:coupling}
  \sum_{j}B_{ij}x_j(T_j+s)&=x_i(s)\mbox{ for all } s\in[-\hat{\alpha}_i,0],
\end{align}
where, in contrast to Section~\ref{sec: delay graphs}, we sum over all predecessors (allowing for more than one). The graph in Fig.~\ref{fig:networks1} in Section~\ref{sec: delay graphs} shows that the term predecessor refers to \emph{direct} predecessors (so, a predecessor of a predecessor of segment (node) $i$ is not automatically also a predecessor of node $i$.)
\end{compactitem}
A node without predecessors can have $\hat{\alpha}=0$, but in general $\hat{\alpha}$ must be larger than $\alpha$, even if $\alpha=0$ for a particular segment. For example, a node $i$ with $\alpha_i=0$ and duration $T_i=1$ may be predecessor of a node $k$ with $\alpha_k=2$. Then $\hat{\alpha}_i$ has to be at least equal to $1$. Notably, no information about $\hat{\alpha}$ will be required during construction of the associated coupling conditions. The graph representation immediately implies the \emph{boundary condition}
\begin{equation}
\label{eq:graphbcs}
    \sum_{j}B_{ij}x_j(T_j)=x_i(0)
\end{equation}
obtained by letting $s=0$ in \eqref{network:coupling}.

\begin{figure}[ht]
\centering
\includegraphics[width=0.9\columnwidth,trim={0 {0.0\textwidth} 0 0.0\textwidth},clip]{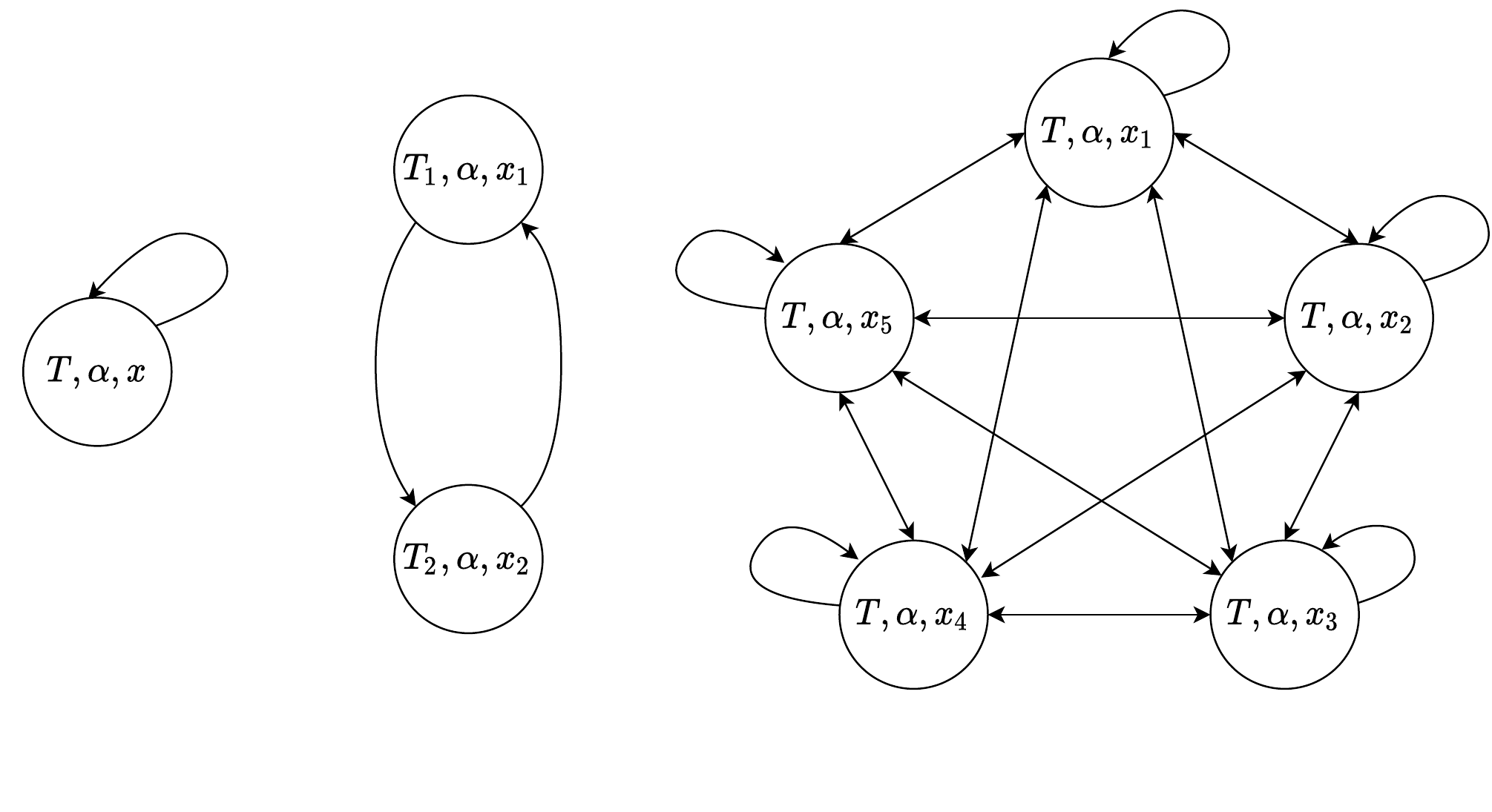}
\vspace{-0.4in}
\caption{Graph representations for delay-coupled boundary-value problems representing (left) a single segment periodic orbit (middle) a two-segment periodic orbit and (right) a quasiperiodic invariant torus approximated with 5 segments.} 
\label{fig:networks2}
\end{figure}

Consider, for example, the single-node graph shown in the left panel of Fig.~\ref{fig:networks2}. This encapsulates a delay differential equation $\dot{x}(t)=f(x(t),x(t-\alpha))$ for $t\in(0,T)$ and algebraic condition $x(s)=x(T+s)$ for $s\in[-\hat{\alpha},0]$ corresponding to the search for a periodic solution of period $T$. We apply the predecessor coupling as many times as necessary to ensure evaluation of $x(\tau)$ only for $\tau\in(0,T)$ in the coupling conditions and obtain
\begin{align}
    \dot{x}(t)&=f(x(t),y(t)),\,t\in(0,T),\\
    y(t)&=\begin{cases}\vspace{0.1cm}
    x\left(T+t-\alpha\big|_{\mathrm{mod}[0,T]}\right), & t\in\left(0,\alpha\big|_{\mathrm{mod}[0,T]}\right),\\x\left(t-\alpha\big|_{\mathrm{mod}[0,T]}\right), & t\in\left(\alpha\big|_{\mathrm{mod}[0,T]},T\right).
\end{cases} 
\end{align}
Finally, \eqref{eq:graphbcs} implies the periodic boundary condition $x(T)=x(0)$.

As a second example, consider the two-node graph shown in the middle panel of Fig.~\ref{fig:networks2}. This encapsulates the system of delay differential equations
\begin{align}
    \dot{x}_1(t)&=f_1(x_1(t),x_1(t-\alpha)),\,t\in(0,T_1),\\
    \dot{x}_2(t)&=f_2(x_2(t),x_2(t-\alpha)),\,t\in(0,T_2)
\end{align}
and algebraic conditions $x_1(s)=x_2(T_2+s)$ for $s\in[-\hat{\alpha}_1,0]$ and $x_2(s)=x_1(T_1+s)$ for $s\in[-\hat{\alpha}_2,0]$ corresponding to the search for a periodic solution of period $T_1+T_2$ for a piecewise-defined vector field with delay $\alpha$. We again apply the predecessor coupling as many times as necessary to ensure evaluation of $x_1(\tau)$ and $x_2(\tau)$ only for $\tau\in(0,T_1)$ and $\tau\in(0,T_2)$, respectively, in the coupling conditions. For example, if $\alpha<T_1,T_2$, and with $y_1(t)=x_1(t-\alpha)$ and $y_2(t)=x_2(t-\alpha)$, we obtain the coupling conditions
\begin{align}
    y_1(t)&=\begin{cases}
    x_2\left(T_2+t-\alpha\right), & t\in\left(0,\alpha\right),\\x_1\left(t-\alpha\right), & t\in\left(\alpha,T_1\right),
    \end{cases}\\
    y_2(t)&=\begin{cases}
    x_1\left(T_1+t-\alpha\right), & t\in\left(0,\alpha\right),\\x_2\left(t-\alpha\right), & t\in\left(\alpha,T_2\right).
    \end{cases}
\end{align}
If, instead, $T_1<\alpha<T_2$, we obtain
\begin{align}
    y_1(t)&=
    x_2\left(T_2+t-\alpha\right),\, t\in\left(0,T_1\right),
    \\
    y_2(t)&=\begin{cases}
    x_2\left(T_1+T_2+t-\alpha\right), & t\in\left(0,\alpha-T_1\right),\\
    x_1\left(T_1+t-\alpha\right), & t\in\left(\alpha-T_1,\alpha\right),\\
    x_2(t-\alpha), & t\in\left(\alpha,T_2\right).
    \end{cases}
\end{align}
In either case, \eqref{eq:graphbcs} implies the boundary conditions $x_1(T)=x_2(0)$ and $x_2(T)=x_1(0)$.

For a general construction, we associate with each node one or several finite paths through our graph consisting of predecessors to this node and their predecessors. Specifically, a sequence $\kappa=(\kappa_1,\ldots,\kappa_\ell)$ with $\kappa_1=i$ is a \emph{history} for the $i$-th segment if $\kappa_{k+1}$ is a predecessor to $\kappa_k$ for $k=1,\ldots,\ell-1$ and
\begin{equation}
    \sum_{j=2}^{\ell-1}T_{\kappa_j}<\alpha_i\le\sum_{j=2}^\ell T_{\kappa_j}.
\end{equation}
In the second example above, the sequence $(1,2)$ is a history for segment $1$ when $\alpha<T_2$ and the sequence $(2,1)$ is a history for segment $2$ when $\alpha<T_1$. In contrast, when $T_1<\alpha<T_1+T_2$, the sequence $(2,1,2)$ is a history for segment $2$. 

For a given history, there exists a smallest index $\iota$ such that
\begin{equation}
    \alpha_i<T_i+\sum_{j=2}^\iota T_{\kappa_j}.
\end{equation}
In particular, when $\alpha_i<T_i$, $\iota=1$. In contrast, for the history $(2,1,2)$ in the second example above with $T_1<\alpha<T_1+T_2$, $\iota=2$ and $\kappa_\iota=1$. All indices between $\iota$ and $\ell$ result in $\ell-\iota$ internal boundaries at
\begin{equation}
    \alpha_i-\sum_{j=2}^{\ell-k} T_{\kappa_j},\,k=1,\ldots,\ell-\iota
\end{equation}
and a partition of $[0,T_i]$ into $\ell-\iota+1$ subintervals separated by these boundaries. For example, if $\ell-\iota=1$, then $[0,T_i]$ is decomposed into the subintervals $[0,\alpha_i]$ and $[\alpha_i,T_i]$. Given these definitions, the $k$-th coupling condition for the $i$-th segment associated with the history $\kappa$ is given by
\begin{equation}
    \label{eq: generalcoupling}
    y_i(t)=A_k(\kappa)x_{\kappa_{\ell-k+1}}\left(t-\alpha_i+\sum_{j=2}^{\ell-k+1}T_{\kappa_j}\right),\,t\in\left[\max\left(0,\alpha_i-\sum_{j=2}^{\ell-k+1}T_{\kappa_j}\right),\min\left(T_i,\alpha_i-\sum_{j=2}^{\ell-k}T_{\kappa_j}\right)\right],
\end{equation}
where 
\begin{equation}
    A_k(\kappa)=B_{\kappa_1\kappa_2}\cdots B_{\kappa_{\ell-k}\kappa_{\ell-k+1}}
\end{equation}
and an empty matrix product is interpreted as an identity matrix. From the boundary condition
\begin{equation}
    B_{\kappa_{l-k}\kappa_{l-k+1}}x_{\kappa_{l-k+1}}(T_{\kappa_{l-k+1}})=x_{\kappa_{l-k}}(0),
\end{equation}
we obtain
\begin{equation}
   A_k(\kappa)x_{\kappa_{\ell-k+1}}\left(T_{\kappa_{\ell-k+1}}\right)=A_{k+1}(\kappa)x_{\kappa_{\ell-k}}\left(0\right),
\end{equation}
which, in turn implies continuity of $y_i(t)$ at the $k$-th internal boundary.

In a general graph, a segment may be associated with two or more distinct histories. In the simplest case, any two such histories $\kappa$ and $\kappa'$ of a segment $i$ correspond to identical sequences of durations $(T_{i,1},\ldots,T_{i,\ell})$. In this special case, the partition of $[0,T_i]$ into subintervals is a property of the segment. We sum over the set of all histories, $K_i$, of the $i$-th segment to obtain a composite coupling condition
\begin{equation}
\label{eq: compcoup}
    y_i(t)=\sum_{\kappa\in K_i}A_k(\kappa)x_{\kappa_{\ell-k+1}}\left(t-\alpha_i+\sum_{j=2}^{\ell-k+1}T_{\kappa_j}\right)
\end{equation}
on the $k$-th subinterval. After appropriate time rescaling, this form matches the general coupling condition in Section~\ref{sec: toolbox template zero problems} and ensures that continuity of the algebraic state variables follows from the boundary conditions on the differential state variables. We leave it to the reader to show that the graph in the right-most panel of Fig.~\ref{fig:networks2} is an example of the special case for which \eqref{eq: compcoup} applies and that the corresponding coupling conditions match the third example in Section~\ref{sec: toolbox template examples}.

In problems where $\alpha_i$ and/or $T_i$ vary during continuation, the indices $\ell$ and $\iota$ may change discretely at critical junctures necessitating a switch between different sets of coupling conditions. The constructive methodology introduced in this section may be deployed to yield a set of equations that remain valid on both sides of such junctures. For example, in the case of the single segment in the left panel of Fig.~\ref{fig:networks2}, the predecessor coupling relationship $x(s)=x(T+s)$ yields
\begin{equation}
    x(t-\alpha)=\begin{cases}
    x(t-\alpha), & t\in(\alpha,T)\\
    x(T+t-\alpha), & t\in(\alpha-T,\alpha)\\
    x(2T+t-\alpha), & t\in(\alpha-2T,\alpha-T)\\
    \qquad\vdots
    \end{cases}
\end{equation}
For $\alpha<T$, the third condition could be omitted, since the left-hand side is never evaluated outside $[0,T]$. If we, nevertheless, retain this condition in our formulation, we need to omit the redundant imposition of continuity across $t=\alpha-T$, since this implies that $x(0)=x(1)$, something that already follows from continuity at $t=\alpha$. Similarly, for $\alpha>T$, the first condition could be omitted, since an inverted interval is assumed to be empty. Again nothing prevents us from retaining this condition also for this case provided that we omit imposing continuity across $t=\alpha$. With proper treatment of continuity, retaining all three conditions allows for variations of $\alpha/T$ across $1$. By adding the next condition in the sequence, we include the possibility of variations of $\alpha/T$ across $2$, and so on.

We leave it to the reader to derive the appropriate generalizations for each of the boundary-value problems considered in this paper.

\subsection{Adaptive discretization}
\label{sec: adaptive discretization}

When a(n augmented) continuation problem is defined on an infinite-dimensional Banach space $\mathcal{U}_\Phi$, it may be appropriate to change discretization (or \emph{remesh} the problem) during continuation, e.g., in order to stay within pre-imposed bounds on the discretization errors (see Part V of~\cite{dankowicz2013recipes} for an extensive discussion of such adaptive meshing). In \textsc{coco}, a continuation problem is said to be \emph{adaptive} if
\begin{compactitem}
\item it is accompanied by instructions for switching between different discretizations without changing the dimensional deficit, and
\item all monitor functions are defined independently of the problem discretization and then discretized accordingly.
\end{compactitem}
During continuation, \textsc{coco} will remesh an adaptive continuation problem at some frequency defined by the corresponding atlas algorithm and according to an algorithm particular to the discretization scheme. Since the monitor functions must be defined independently of the discretization, they span the coordinate axes of an invariant, finite-dimensional projection of $\mathcal{U}_\Phi$ which may serve to visualize an arbitrary solution manifold. Indeed, as long as a sufficient number of independent monitor functions are included with the continuation problem, continuation may proceed along such a solution manifold in terms of a geometry defined in the projected space, independently of any adaptive changes to the mesh. This is the solution implemented in the \mcode{atlas\_kd} atlas algorithm~\cite{dankowicz2020multidimensional}.

For the abstract toolbox template presented in Section~\ref{sec: toolbox template zero problems}, the corresponding discretization in Section~\ref{sec: discretization implementation} is uniquely determined by the order $N$ and polynomial degree $m$, since the mesh points $\tau_{\mathrm{pt},j}=(j-1)/N$ were assumed to be evenly distributed over the interval $[0,1]$ (even though this was not required by the abstract form of problem discretization discussed in Section~\ref{sec: discretization abstract}). A simple form of adaptation would allow discrete changes to $N$ and/or $m$ during continuation, in order to accommodate variations of an estimated discretization error. Since such changes would inevitably change the relationship between individual base points and the corresponding time instants, it would be inappropriate to define a monitor function that evaluated, e.g., to  $x_{\mathrm{bp},j}$ for $j\in\{2,\ldots,N(m+1)-1\}$. In contrast, a monitor function that evaluated to the value of $x(\cdot)$ at a particular fixed time or the integral of $x(\cdot)$ over the interval $[0,1]$ would be defined (if not computed) independently of the particular mesh, since the piecewise polynomial $\tilde x(\cdot)$ is a continuous function at every point of the solution manifold.

In a more sophisticated form of adaptation, not only could the number $N$ of mesh intervals (or, less commonly, the polynomial degree $m$) vary during continuation, but one would also allow for non-uniform mesh intervals with unevenly spaced time meshes $\{\tau_{\mathrm{pt},j}\}_{j=1}^{N+1}$.  The \textsc{coco} toolbox \mcode{coll} implements a mesh-selection strategy that chooses the order $N$ and the mesh points $\tau_{\mathrm{pt},j}$ such that they equidistribute an estimated (positive) density $e(\tau)$ of a given error measure according to
\begin{equation}
\label{eq:adaptive:mesh}
    N=\frac{\langle e\rangle^{(m+1)/m}}{\mathtt{tol}^{1/m}},\,\int_0^{\tau_{\mathrm{pt},j}}e(\tau)\d \tau=\frac{j-1}{N}\langle e\rangle,\,j=2,\ldots,N
\end{equation}
where 
\begin{equation}
    \langle e\rangle=\int_0^1e(\tau)\d \tau
\end{equation}
and \texttt{tol} is a user-defined tolerance. For a delay-coupled system of differential constraints, a similar strategy would need to be concerned about possible loss of orders of differentiability of the exact solution at certain breakpoints even in the presence of smooth problem coefficients. Such breakpoints occur, for example, in initial-value problems with delay, including the data assimilation problem from Section~\ref{sec: data assimilation problem construction} or the optimal control problem in Section~\ref{sec: numerical examples optimal control}. Similar breakpoints would be expected in the adjoint variables for periodic delay-coupled boundary-value problems when the corresponding objective functional is not invariant with respect to time shifts. The reduced regularity of the solution at these breakpoints may lead to poles in the estimates for the error measure density $e(\tau)$ and, consequently, to inefficient placement of mesh points or reduced accuracy. For initial-value problems, the interaction of mesh selection and breakpoints has been discussed extensively \cite{guglielmi2001implementing,shampine2009numerical}.

%%%%%%%%%%%%%
% Conclusions
%%%%%%%%%%%%%
\section{Conclusions}
\label{sec: Conclusions}

The staged approach to problem construction supported by \textsc{coco} permits the user to build up nonlinear problems gradually by adding new variables and systems of equations and coupling them flexibly to variables defined previously, at each step increasing or decreasing the dimensional deficit (nominally the dimension of the corresponding solution manifold). This is the natural way of thinking about problem construction if algorithms for multi-dimensional continuation are at one's disposal. The initial examples of the paper showed how bifurcations or function extrema are embedded within higher-dimensional solution manifolds, and how such higher-dimensional solution manifolds are computable even in cases when low precision and poor condition numbers obstruct classical one-parameter continuation.

Our paper then described in detail the abstract staged construction formalism in the full generality currently supported by \textsc{coco}. A major innovation since its original realization in \cite{dankowicz2013recipes} is that the formalism now supports the simultaneous gradual build-up of adjoint information and also includes a new layer that permits construction of complementarity conditions associated with design optimization in the presence of inequality constraints (as partially described in~\cite{li2020optimization}). The utility of such staged construction with automatic accumulation of adjoints was illustrated using two detailed examples. The first example, a data assimilation problem, was an optimization problem with multiple delay-coupled time segments. The second, a phase response analysis of periodic orbits, was formulated as a linear sensitivity analysis of the orbital duration with respect to perturbations in the boundary conditions. 

Both examples showed that it is, in principle, possible to perform staged construction of a boundary-value problem associated with multiple segments, coupled to each other by discrete time delays, while automatically accumulating adjoints. The underlying structure turned out to be a network of delay-coupled systems of ordinary differential equations, linked by algebraic coupling constraints. A general representation in terms of delay graphs inspired the formulation of an abstract toolbox for delay-coupled problems, where each building block (a differential constraint and a set of algebraic coupling conditions) is sufficiently general but also simple enough to implement its adjoint at the toolbox level. The paper went on to formulate the discretized version of this abstract network of equations, first in terms of abstract projections, then with a detailed vectorized description of the resulting algebraic equations. The generality of the toolbox was demonstrated in the context of several numerical examples of coupled systems with delay as they arise for connecting orbits, optimal control problems, and quasiperiodic invariant tori. The final section commented on what is missing before the toolbox is ``ready for production'': automated decomposition of the delay graph (which is the form in which a user presents the problem) into the building blocks of the toolbox, as well as a means of error control through adaptive meshing and theoretical convergence analysis.

Beyond such improvements to the proposed toolbox, several opportunities for further work follow from the treatment in this paper. Among them are generalizations of the phase-response analysis to other normally hyperbolic invariant manifolds with a natural definition of an asymptotic phase, including periodic orbits in piecewise-smooth dynamical systems~\cite{PhysRevE.95.012212} and quasiperiodic invariant tori~\cite{5654185}. Indeed, we anticipate that our derivation of the corresponding adjoint boundary-value problems using sensitivity analysis, while carried out here only for a single-segment periodic orbit problem, should carry over without significant modification or additional overhead to multi-segment periodic orbit problems with discrete delays. It would follow that such problems could be analyzed using existing \textsc{coco} toolboxes without the need for special-purpose solutions.

We illustrated the theoretical use of multi-dimensional continuation for regularizing nearly singular problems in the presence of low precision numerics as would be expected in data-driven applications, e.g., experiments using control-based continuation. For such a methodology to work well in practice, we anticipate the need for a more purposeful design of the \mcode{atlas\_kd} algorithm to allow, for example, continuation along strips of higher-dimensional manifolds. This would permit the benefits of regularization without excessive (and costly) excursions into additional dimensions. For problems with underlying continuous symmetries (such as the rotational symmetry for the linear harmonic oscillator), it would be beneficial to develop appropriate modifications to \mcode{atlas\_kd} to again retain the benefits of higher-dimensional continuation without incurring its full cost.

Several classes of problems involving delay are not covered by the template toolbox developed here. These include problems with state- or time-dependent delays, as well as those with distributed delays. Even for discrete delays, we have assumed an explicit form of the differential constraints with similarly explicit algebraic coupling conditions. In contrast, the defining problem in \textsc{dde-biftool} admits delay differential equations with a nontrivial (and possibly singular) coefficient matrix on the left-hand side, thereby enabling analysis of problems with nontrivial algebraic coupling conditions. Our general approach to recognizing universality and encoding such universality in the \textsc{coco} framework, including with attention to the automated construction of adjoints, should inform such further development.

\section{Declarations}
\subsection{Funding}
JS gratefully acknowledges support by EPSRC Fellowship EP/N023544/1 and EPSRC grant EP/V04687X/1.

\subsection{Conflicts of interest}
The authors declare that they have no conflict of interest.

\subsection{Data availability}
The data used to generate the numerical results included in this paper are available from the corresponding author on reasonable request.

\subsection{Code availability}
Code used to generate the numerical results included in this paper is available from the corresponding author on reasonable request.

\subsection{Authors' contributions}
All authors contributed equally to the conception of this paper and the formal analysis. Zaid Ahsan and Mingwu Li were leads on the software development with Harry Dankowicz and Jan Sieber acting in supporting roles. Zaid Ahsan and Mingwu Li were leads on the development of graphical illustrations with Jan Sieber acting in a supporting role. All authors contributed equally to the writing, review, and editing of the text.

\bibliographystyle{spmpsci}

\begin{thebibliography}{100}
\providecommand{\url}[1]{{#1}}
\providecommand{\urlprefix}{URL }
\expandafter\ifx\csname urlstyle\endcsname\relax
  \providecommand{\doi}[1]{DOI~\discretionary{}{}{}#1}\else
  \providecommand{\doi}{DOI~\discretionary{}{}{}\begingroup
  \urlstyle{rm}\Url}\fi

\bibitem{abbas2011parametric}
Abbas, L.K., Rui, X., Marzocca, P., Abdalla, M., De~Breuker, R.: A parametric
  study on supersonic/hypersonic flutter behavior of aero-thermo-elastic
  geometrically imperfect curved skin panel.
\newblock Acta Mechanica \textbf{222}(1), 41--57 (2011)

\bibitem{acharya2020non}
Acharya, V., Lieuwen, T.: Non-monotonic flame response behaviors in
  harmonically forced flames.
\newblock Proceedings of the Combustion Institute  (2020)

\bibitem{ahsan2020optimization}
Ahsan, Z., Dankowicz, H., Sieber, J.: Optimization along families of periodic
  and quasiperiodic orbits in dynamical systems with delay.
\newblock Nonlinear Dynamics \textbf{99}(1), 837--854 (2020)

\bibitem{allgower2003introduction}
Allgower, E.L., Georg, K.: Introduction to numerical continuation methods.
\newblock SIAM (2003)

\bibitem{amandio2014stochastic}
Amandio, L., Marta, A., Afonso, F., Vale, J., Suleman, A., Araujo, A.:
  Stochastic optimization in aircraft design.
\newblock In: Engineering Optimization, pp. 267--272. CRC Press (2014)

\bibitem{andoSIAM2020}
And{\`o}, A., Breda, D.: Convergence analysis of collocation methods for
  computing periodic solutions of retarded functional differential equations.
\newblock SIAM Journal on Numerical Analysis \textbf{58}(5), 3010--3039 (2020)

\bibitem{arnold1972lectures}
Arnold, V.I.: Lectures on bifurcations in versal families.
\newblock In: Vladimir I. Arnold-Collected Works, pp. 271--340. Springer (1972)

\bibitem{back1992dstool}
Back, A., Guckenheimer, J., Myers, M., Wicklin, F., Worfolk, P.: Ds{T}ool:
  Computer assisted exploration of dynamical systems.
\newblock Notices Amer. Math. Soc \textbf{39}(4), 303--309 (1992)

\bibitem{barton2012control}
Barton, D., Mann, B., Burrow, S.: Control-based continuation for investigating
  nonlinear experiments.
\newblock Journal of Vibration and Control \textbf{18}(4), 509--520 (2012)

\bibitem{barton2009stability}
Barton, D.A.: Stability calculations for piecewise-smooth delay equations.
\newblock International Journal of Bifurcation and Chaos \textbf{19}(02),
  639--650 (2009)

\bibitem{barton2017control}
Barton, D.A.: Control-based continuation: bifurcation and stability analysis
  for physical experiments.
\newblock Mechanical Systems and Signal Processing \textbf{84}, 54--64 (2017)

\bibitem{bartoszewski2011solving}
Bartoszewski, Z.: Solving boundary value problems for delay differential
  equations by a fixed-point method.
\newblock Journal of Computational and Applied Mathematics \textbf{236}(6),
  1576--1590 (2011)

\bibitem{ben1982unified}
Ben-Tal, A., Zowe, J.: A unified theory of first and second order conditions
  for extremum problems in topological vector spaces.
\newblock In: Optimality and stability in mathematical programming, pp. 39--76.
  Springer (1982)

\bibitem{berezansky2006mackey}
Berezansky, L., Braverman, E.: Mackey-{G}lass equation with variable
  coefficients.
\newblock Computers \& Mathematics with Applications \textbf{51}(1), 1--16
  (2006)

\bibitem{berezansky2012mackey}
Berezansky, L., Braverman, E., Idels, L.: The {M}ackey--{G}lass model of
  respiratory dynamics: review and new results.
\newblock Nonlinear Analysis: Theory, Methods \& Applications \textbf{75}(16),
  6034--6052 (2012)

\bibitem{beyn1990numerical}
Beyn, W.J.: The numerical computation of connecting orbits in dynamical
  systems.
\newblock IMA Journal of Numerical Analysis \textbf{10}(3), 379--405 (1990)

\bibitem{blyth2020}
Blyth, M., Renson, L., Marucci, L.: Tutorial of numerical continuation and
  bifurcation theory for systems and synthetic biology.
\newblock \url{https://https://arxiv.org/pdf/2008.05226.pdf}.
\newblock Accessed: 2021-03-24

\bibitem{boender1982stochastic}
Boender, C.G.E., Kan, A.R., Timmer, G., Stougie, L.: A stochastic method for
  global optimization.
\newblock Mathematical Programming \textbf{22}(1), 125--140 (1982)

\bibitem{byrd2000trust}
Byrd, R.H., Gilbert, J.C., Nocedal, J.: A trust region method based on interior
  point techniques for nonlinear programming.
\newblock Mathematical Programming \textbf{89}(1), 149--185 (2000)

\bibitem{calver2017numerical}
Calver, J., Enright, W.: {Numerical methods} for {computing sensitivities} for
  {ODEs} and {DDEs}.
\newblock Numerical Algorithms \textbf{74}(4), 1101--1117 (2017)

\bibitem{cao2019nonlinear}
Cao, C.J., Hill, T.L., Conn, A.T., Li, B., Gao, X.: Nonlinear dynamics of a
  magnetically coupled dielectric elastomer actuator.
\newblock Physical Review Applied \textbf{12}(4), 044,033 (2019)

\bibitem{chai2013unified}
Chai, Q., Loxton, R., Teo, K.L., Yang, C.: A unified parameter identification
  method for nonlinear time-delay systems.
\newblock Journal of Industrial and Management Optimization (JIMO)
  \textbf{9}(2), 471--486 (2013)

\bibitem{chavez2020numerical}
Ch{\'a}vez, J.P., Zhang, Z., Liu, Y.: A numerical approach for the bifurcation
  analysis of nonsmooth delay equations.
\newblock Communications in Nonlinear Science and Numerical Simulation
  \textbf{83}, 105,095 (2020)

\bibitem{chicone2004asymptotic}
Chicone, C., Liu, W.: Asymptotic phase revisited.
\newblock Journal of Differential Equations \textbf{204}(1), 227--246 (2004)

\bibitem{chong2016numerical}
Chong, A.: Numerical modelling and stability analysis of non-smooth dynamical
  systems vie {ABESPOL}.
\newblock Ph.D. thesis, University of Aberdeen (2016)

\bibitem{crisfield1983arc}
Crisfield, M.: An arc-length method including line searches and accelerations.
\newblock International Journal for Numerical Methods in Engineering
  \textbf{19}(9), 1269--1289 (1983)

\bibitem{dankowicz2011extended}
Dankowicz, H., Schilder, F.: An extended continuation problem for bifurcation
  analysis in the presence of constraints.
\newblock Journal of Computational and Nonlinear Dynamics \textbf{6}(3) (2011)

\bibitem{dankowicz2013recipes}
Dankowicz, H., Schilder, F.: Recipes for continuation.
\newblock SIAM (2013)

\bibitem{dankowicz2011continuation}
Dankowicz, H., Schilder, F., Saghafi, M.: Continuation of connecting orbits
  with {L}in’s method using {COCO}.
\newblock In: Proceedings of the 7th European Nonlinear Dynamics Conference
  (ENOC 2011) (2011)

\bibitem{dankowicz2020multidimensional}
Dankowicz, H., Wang, Y., Schilder, F., Henderson, M.E.: Multidimensional
  manifold continuation for adaptive boundary-value problems.
\newblock Journal of Computational and Nonlinear Dynamics \textbf{15}(5) (2020)

\bibitem{d2010choice}
D'Avino, G., Crescitelli, S., Maffettone, P., Grosso, M.: On the choice of the
  optimal periodic operation for a continuous fermentation process.
\newblock Biotechnology Progress \textbf{26}(6), 1580--1589 (2010)

\bibitem{dellnitz1996computation}
Dellnitz, M., Hohmann, A.: The computation of unstable manifolds using
  subdivision and continuation.
\newblock In: Nonlinear dynamical systems and chaos, pp. 449--459. Springer
  (1996)

\bibitem{5654185}
{Demirt}, A., {Gu}, C., {Roychowdhury}, J.: Phase equations for quasi-periodic
  oscillators.
\newblock In: 2010 IEEE/ACM International Conference on Computer-Aided Design
  (ICCAD), pp. 292--297 (2010)

\bibitem{dercole2005slidecont}
Dercole, F., Kuznetsov, Y.A.: {SlideCont}: {A}n {A}uto97 driver for bifurcation
  analysis of filippov systems.
\newblock ACM Transactions on Mathematical Software (TOMS) \textbf{31}(1),
  95--119 (2005)

\bibitem{dhooge2003matcont}
Dhooge, A., Govaerts, W., Kuznetsov, Y.A.: {MATCONT}: a {MATLAB} package for
  numerical bifurcation analysis of {ODEs}.
\newblock ACM Transactions on Mathematical Software (TOMS) \textbf{29}(2),
  141--164 (2003)

\bibitem{doedel2007lecture}
Doedel, E.J.: Lecture notes on numerical analysis of nonlinear equations.
\newblock In: Numerical continuation methods for dynamical systems, pp. 1--49.
  Springer (2007)

\bibitem{doedel2007auto}
Doedel, E.J., Champneys, A.R., Dercole, F., Fairgrieve, T.F., Kuznetsov, Y.A.,
  Oldeman, B., Paffenroth, R., Sandstede, B., Wang, X., Zhang, C.: {AUTO}-07p:
  Continuation and bifurcation software for ordinary differential equations.
\newblock \url{https://github.com/auto-07p/auto-07p}.
\newblock Accessed: 2021-04-22

\bibitem{doedel2006global}
Doedel, E.J., Krauskopf, B., Osinga, H.M.: Global bifurcations of the {L}orenz
  manifold.
\newblock Nonlinearity \textbf{19}(12), 2947 (2006)

\bibitem{engelborghs2001collocation}
Engelborghs, K., Luzyanina, T., Hout, K.I., Roose, D.: Collocation {methods}
  for the {computation} of {periodic solutions} of {delay differential
  equations}.
\newblock SIAM Journal on Scientific Computing \textbf{22}(5), 1593--1609
  (2001)

\bibitem{engelborghs2002numerical}
Engelborghs, K., Luzyanina, T., Roose, D.: Numerical bifurcation analysis of
  delay differential equations using {DDE-BIFTOOL}.
\newblock ACM Transactions on Mathematical Software (TOMS) \textbf{28}(1),
  1--21 (2002)

\bibitem{england2005computing}
England, J.P., Krauskopf, B., Osinga, H.M.: Computing one-dimensional global
  manifolds of poincar{\'e} maps by continuation.
\newblock SIAM Journal on Applied Dynamical Systems \textbf{4}(4), 1008--1041
  (2005)

\bibitem{ermentrout1996type}
Ermentrout, B.: Type {I} membranes, phase resetting curves, and synchrony.
\newblock Neural Computation \textbf{8}(5), 979--1001 (1996)

\bibitem{ermentrout2002simulating}
Ermentrout, B.: Simulating, analyzing, and animating dynamical systems: a guide
  to {XPPAUT} for researchers and students.
\newblock SIAM (2002)

\bibitem{formica2013coupling}
Formica, G., Arena, A., Lacarbonara, W., Dankowicz, H.: Coupling {FEM} with
  parameter continuation for analysis of bifurcations of periodic responses in
  nonlinear structures.
\newblock Journal of Computational and Nonlinear Dynamics \textbf{8}(2) (2013)

\bibitem{fotsch2016thesis}
Fotsch, E.L.: Bifurcation analysis near the cessation of complete chatter and
  {S}hil'nikov homoclinic trajectories in a pressure relief valve model.
\newblock Master's thesis, University of Illinois at Urbana-Champaign (2016)

\bibitem{gelfand2000calculus}
Gelfand, I.M., Silverman, R.A., et~al.: Calculus of variations.
\newblock Courier Corporation (2000)

\bibitem{glass2010mackey}
Glass, L., Mackey, M.: Mackey-{G}lass equation.
\newblock Scholarpedia \textbf{5}(3), 6908 (2010)

\bibitem{gollmann2009optimal}
G{\"o}llmann, L., Kern, D., Maurer, H.: Optimal {control problems} with
  {delays} in {state} and {control variables subject} to {mixed control--state
  constraints}.
\newblock Optimal Control Applications and Methods \textbf{30}(4), 341--365
  (2009)

\bibitem{gonzalez2017assessing}
Gonzalez-Buelga, A., Lazar, I.F., Jiang, J.Z., Neild, S.A., Inman, D.J.:
  Assessing the effect of nonlinearities on the performance of a tuned inerter
  damper.
\newblock Structural Control and Health Monitoring \textbf{24}(3), e1879 (2017)

\bibitem{govaerts2000numerical}
Govaerts, W.: Numerical bifurcation analysis for {ODEs}.
\newblock Journal of Computational and Applied Mathematics \textbf{125}(1-2),
  57--68 (2000)

\bibitem{govaerts2006computation}
Govaerts, W., Sautois, B.: Computation of the phase response curve: a direct
  numerical approach.
\newblock Neural Computation \textbf{18}(4), 817--847 (2006)

\bibitem{guckenheimer2015invariant}
Guckenheimer, J., Krauskopf, B., Osinga, H.M., Sandstede, B.: Invariant
  manifolds and global bifurcations.
\newblock Chaos: An Interdisciplinary Journal of Nonlinear Science
  \textbf{25}(9), 097,604 (2015)

\bibitem{guddat1990parametric}
Guddat, J., Vazquez, F.G., Jongen, H.T.: Parametric optimization:
  singularities, pathfollowing and jumps.
\newblock Springer (1990)

\bibitem{guglielmi2001implementing}
Guglielmi, N., Hairer, E.: Implementing {R}adau {IIA} methods for stiff delay
  differential equations.
\newblock Computing \textbf{67}(1), 1--12 (2001)

\bibitem{haller2016nonlinear}
Haller, G., Ponsioen, S.: Nonlinear normal modes and spectral submanifolds:
  existence, uniqueness and use in model reduction.
\newblock Nonlinear Dynamics \textbf{86}(3), 1493--1534 (2016)

\bibitem{heitmann2021arrhythmogenic}
Heitmann, S., Shpak, A., Vandenberg, J.I., Hill, A.P.: Arrhythmogenic effects
  of ultra-long and bistable cardiac action potentials.
\newblock PLOS Computational Biology \textbf{17}(2), e1008,683 (2021)

\bibitem{henderson2002multiple}
Henderson, M.E.: Multiple parameter continuation: {C}omputing implicitly
  defined k-manifolds.
\newblock International Journal of Bifurcation and Chaos \textbf{12}(03),
  451--476 (2002)

\bibitem{izhikevich1997weakly}
Izhikevich, E., Hoppensteadt, F.: Weakly connected neural networks.
\newblock New York: Springer-Verlag  (1997)

\bibitem{izhikevich2007dynamical}
Izhikevich, E.M.: Dynamical systems in neuroscience.
\newblock MIT press (2007)

\bibitem{ssmtool2}
Jain, S., Thurnher, T., Li, M.: {SSMTool} 2.0: {C}omputation of invariant
  manifolds \& their reduced dynamics in high-dimensional mechanics problems
  (v1.0.0).
\newblock \url{http://doi.org/10.5281/zenodo.4614202}.
\newblock Accessed: 2021-04-11

\bibitem{kelley1995iterative}
Kelley, C.T.: Iterative methods for linear and nonlinear equations.
\newblock SIAM (1995)

\bibitem{kernevez1987optimization}
Kern{\'e}vez, J., Doedel, E.: Optimization {in} {bifurcation} {problems}
  {using} {a} {continuation} {method}.
\newblock In: Bifurcation: Analysis, Algorithms, Applications, pp. 153--160.
  Springer (1987)

\bibitem{kewlani2012polynomial}
Kewlani, G., Crawford, J., Iagnemma, K.: A polynomial chaos approach to the
  analysis of vehicle dynamics under uncertainty.
\newblock Vehicle System Dynamics \textbf{50}(5), 749--774 (2012)

\bibitem{khasawneh2011multi}
Khasawneh, F.A., Mann, B.P., Butcher, E.A.: A multi-interval {C}hebyshev
  collocation approach for the stability of periodic delay systems with
  discontinuities.
\newblock Communications in Nonlinear Science and Numerical Simulation
  \textbf{16}(11), 4408--4421 (2011)

\bibitem{koh2016optimizing}
Koh, M.H., Sipahi, R.: Optimizing agent coupling strengths in a network
  dynamics with inter-agent delays for achieving fast consensus.
\newblock In: 2016 American Control Conference (ACC), pp. 5358--5363. IEEE
  (2016)

\bibitem{krantz2012implicit}
Krantz, S.G., Parks, H.R.: The implicit function theorem: history, theory, and
  applications.
\newblock Springer Science \& Business Media (2012)

\bibitem{krauskopf2007numerical}
Krauskopf, B., Osinga, H.M., Gal{\'a}n-Vioque, J.: Numerical continuation
  methods for dynamical systems, vol.~2.
\newblock Springer (2007)

\bibitem{krauskopf2008lin}
Krauskopf, B., Rie{\ss}, T.: A {L}in's method approach to finding and
  continuing heteroclinic connections involving periodic orbits.
\newblock Nonlinearity \textbf{21}(8), 1655 (2008)

\bibitem{kuehn2015efficient}
Kuehn, C.: Efficient gluing of numerical continuation and a multiple solution
  method for elliptic {PDE}s.
\newblock Applied Mathematics and Computation \textbf{266}, 656--674 (2015)

\bibitem{kuznetsov2013elements}
Kuznetsov, Y.A.: Elements of applied bifurcation theory, vol. 112.
\newblock Springer Science \& Business Media (2013)

\bibitem{kuznetsov1997content}
Kuznetsov, Y.A., Levitin, V.V.: {CONTENT}: integrated environment for analysis
  of dynamical systems.
\newblock \url{https://webspace.science.uu.nl/~kouzn101/CONTENT/}.
\newblock Accessed: 2021-04-21

\bibitem{langfield2020continuation}
Langfield, P., Krauskopf, B., Osinga, H.M.: A continuation approach to
  computing phase resetting curves.
\newblock In: Advances in Dynamics, Optimization and Computation, pp. 3--30.
  Springer (2020)

\bibitem{mingwu2020thesis}
Li, M.: Dynamics and optimal control of information transmission in complex
  systems.
\newblock Ph.D. thesis, University of Illinois at Urbana-Champaign (2020)

\bibitem{li2020tor}
Li, M.: Tor: a toolbox for the continuation of two-dimensional tori in
  autonomous systems and non-autonomous systems with periodic forcing.
\newblock arXiv preprint arXiv:2012.13256  (2020).
\newblock \url{https://github.com/mingwu-li/torus\_collocation}. Accessed:
  2021-04-01

\bibitem{coco-shoot}
Li, M., Dankowicz, H.: A {COCO}-based shooting toolbox for dynamical systems.
\newblock \url{https://github.com/mingwu-li/forward}.
\newblock Accessed: 2021-04-21

\bibitem{coco-fmincon}
Li, M., Dankowicz, H.: Coupling {COCO} with fmincon for constrained
  optimization of dynamical systems.
\newblock \url{https://github.com/mingwu-li/coco\_fmincon}.
\newblock Accessed: 2021-04-03

\bibitem{li2017staged}
Li, M., Dankowicz, H.: Staged {construction} of {adjoints} for {constrained
  optimization} of {integro-differential boundary-value problems}.
\newblock SIAM Journal on Applied Dynamical Systems \textbf{17}(2), 1117--1151
  (2018)

\bibitem{li2020optimization}
Li, M., Dankowicz, H.: Optimization with equality and inequality constraints
  using parameter continuation.
\newblock Applied Mathematics and Computation \textbf{375}, 125,058 (2020)

\bibitem{liberzon2011calculus}
Liberzon, D.: Calculus of variations and optimal control theory: a concise
  introduction.
\newblock Princeton university press (2011)

\bibitem{liu2017controlling}
Liu, Y., Ch{\'a}vez, J.P.: Controlling multistability in a vibro-impact capsule
  system.
\newblock Nonlinear Dynamics \textbf{88}(2), 1289--1304 (2017)

\bibitem{luzyanina1997computation}
Luzyanina, T., Engelborghs, K., Lust, K., Roose, D.: Computation, continuation
  and bifurcation analysis of periodic solutions of delay differential
  equations.
\newblock International Journal of Bifurcation and Chaos \textbf{7}(11),
  2547--2560 (1997)

\bibitem{munoz2003continuation}
Munoz-Almaraz, F.J., Freire, E., Gal{\'a}n, J., Doedel, E., Vanderbauwhede, A.:
  Continuation of periodic orbits in conservative and {H}amiltonian systems.
\newblock Physica D: Nonlinear Phenomena \textbf{181}(1-2), 1--38 (2003)

\bibitem{novivcenko2012phase}
Novi{\v{c}}enko, V., Pyragas, K.: Phase reduction of weakly perturbed limit
  cycle oscillations in time-delay systems.
\newblock Physica D: Nonlinear Phenomena \textbf{241}(12), 1090--1098 (2012)

\bibitem{osinga2010continuation}
Osinga, H.M., Moehlis, J.: Continuation-based computation of global isochrons.
\newblock SIAM Journal on Applied Dynamical Systems \textbf{9}(4), 1201--1228
  (2010)

\bibitem{otter1996modeling}
Otter, M., Elmqvist, H., Cellier, F.E.: Modeling of multibody systems with the
  object-oriented modeling language {D}ymola.
\newblock Nonlinear Dynamics \textbf{9}(1), 91--112 (1996)

\bibitem{paul2000designing}
Paul, C.A.: Designing efficient software for solving delay differential
  equations.
\newblock Journal of Computational and Applied Mathematics \textbf{125}(1-2),
  287--295 (2000)

\bibitem{ponsioen2018automated}
Ponsioen, S., Pedergnana, T., Haller, G.: Automated computation of autonomous
  spectral submanifolds for nonlinear modal analysis.
\newblock Journal of Sound and Vibration \textbf{420}, 269--295 (2018)

\bibitem{porter2009communities}
Porter, M.A., Onnela, J.P., Mucha, P.J.: Communities in networks.
\newblock Notices of the AMS \textbf{56}(9), 1082--1097 (2009)

\bibitem{renson2016robust}
Renson, L., Gonzalez-Buelga, A., Barton, D., Neild, S.: Robust identification
  of backbone curves using control-based continuation.
\newblock Journal of Sound and Vibration \textbf{367}, 145--158 (2016)

\bibitem{renson2019application}
Renson, L., Shaw, A., Barton, D., Neild, S.: Application of control-based
  continuation to a nonlinear structure with harmonically coupled modes.
\newblock Mechanical Systems and Signal Processing \textbf{120}, 449--464
  (2019)

\bibitem{renson2019numerical}
Renson, L., Sieber, J., Barton, D., Shaw, A., Neild, S.: Numerical continuation
  in nonlinear experiments using local {G}aussian process regression.
\newblock Nonlinear Dynamics \textbf{98}(4), 2811--2826 (2019)

\bibitem{rheinboldt1996manpak}
Rheinboldt, W.C.: {MANPAK}: A set of algorithms for computations on implicitly
  defined manifolds.
\newblock Computers \& Mathematics with Applications \textbf{32}(12), 15--28
  (1996)

\bibitem{roose2007continuation}
Roose, D., Szalai, R.: Continuation and bifurcation analysis of delay
  differential equations.
\newblock In: Numerical continuation methods for dynamical systems, pp.
  359--399. Springer (2007)

\bibitem{samaey2002numerical}
Samaey, G., Engelborghs, K., Roose, D.: Numerical computation of connecting
  orbits in delay differential equations.
\newblock Numerical Algorithms \textbf{30}(3), 335--352 (2002)

\bibitem{schiehlen2013advanced}
Schiehlen, W.: Advanced multibody system dynamics: simulation and software
  tools, vol.~20.
\newblock Springer Science \& Business Media (2013)

\bibitem{schilder2015experimental}
Schilder, F., Bureau, E., Santos, I.F., Thomsen, J.J., Starke, J.: Experimental
  bifurcation analysis—continuation for noise-contaminated zero problems.
\newblock Journal of Sound and Vibration \textbf{358}, 251--266 (2015)

\bibitem{COCO}
Schilder, F., Dankowicz, H., Li, M.: Continuation {Core} and {Toolboxes}
  ({COCO}).
\newblock \url{https://sourceforge.net/projects/cocotools}.
\newblock Accessed: 2021-03-26

\bibitem{schilder2005continuation}
Schilder, F., Osinga, H.M., Vogt, W.: Continuation of quasi-periodic invariant
  tori.
\newblock SIAM Journal on Applied Dynamical Systems \textbf{4}(3), 459--488
  (2005)

\bibitem{seydel2009practical}
Seydel, R.: Practical bifurcation and stability analysis, vol.~5.
\newblock Springer Science \& Business Media (2009)

\bibitem{shampine2009numerical}
Shampine, L.F., Thompson, S.: Numerical solution of delay differential
  equations.
\newblock In: Delay Differential Equations, pp. 1--27. Springer (2009)

\bibitem{shinohara2007boundary}
Shinohara, Y., Fujimori, H., Suzuki, T., Kurihara, M.: On a boundary value
  problem for delay differential equations of population dynamics and
  {C}hebyshev approximation.
\newblock Journal of Computational and Applied mathematics \textbf{201}(2),
  348--355 (2007)

\bibitem{PhysRevE.95.012212}
Shirasaka, S., Kurebayashi, W., Nakao, H.: Phase reduction theory for hybrid
  nonlinear oscillators.
\newblock Phys. Rev. E \textbf{95}, 012,212 (2017)

\bibitem{ddebiftoolmanual}
Sieber, J., Engelborghs, K., Luzyanina, T., Samaey, G., Roose, D.:
  {DDE-BIFTOOL} {M}anual --- {B}ifurcation analysis of delay differential
  equations.
\newblock \url{sourceforge.net/projects/ddebiftool} and
  \url{sourceforge.net/p/ddebiftool/git/ci/master/tree/ddebiftool\_coco}

\bibitem{sieber2008tracking}
Sieber, J., Krauskopf, B.: Tracking {oscillations} in the {presence} of
  {delay-Induced essential instability}.
\newblock Journal of Sound and Vibration \textbf{315}(3), 781--795 (2008)

\bibitem{smithoptimal}
Smith, S.: Optimal control of delay differential equations using evolutionary
  algorithms.
\newblock Complexity International \textbf{12}, 1--10 (2005)

\bibitem{Kunt}
Szalai, R.: Knut: A continuation and bifurcation software for
  delay-differential equations.
\newblock \url{https://rs1909.github.io/knut/}.
\newblock Accessed: 2021-03-26

\bibitem{szalai2019model}
Szalai, R.: Model reduction of non-densely defined piecewise-smooth systems in
  banach spaces.
\newblock Journal of Nonlinear Science \textbf{29}(3), 897--960 (2019)

\bibitem{szalai2020invariant}
Szalai, R.: Invariant spectral foliations with applications to model order
  reduction and synthesis.
\newblock Nonlinear Dynamics \textbf{101}(4), 2645--2669 (2020)

\bibitem{thota2008tc}
Thota, P., Dankowicz, H.: {TC-HAT}: {A} novel toolbox for the continuation of
  periodic trajectories in hybrid dynamical systems.
\newblock SIAM Journal on Applied Dynamical Systems \textbf{7}(4), 1283--1322
  (2008)

\bibitem{toilliez2008optimized}
Toilliez, J.O., Szeri, A.J.: Optimized translation of microbubbles driven by
  acoustic fields.
\newblock The Journal of the Acoustical Society of America \textbf{123}(4),
  1916--1930 (2008)

\bibitem{touze2006nonlinear}
Touz{\'e}, C., Amabili, M.: Nonlinear normal modes for damped geometrically
  nonlinear systems: Application to reduced-order modelling of harmonically
  forced structures.
\newblock Journal of Sound and Vibration \textbf{298}(4-5), 958--981 (2006)

\bibitem{TraversoMagri2019}
Traverso, T., Magri, L.: Data assimilation in a nonlinear time-delayed
  dynamical system with {L}agrangian optimization.
\newblock Lecture Notes in Computer Science (including subseries Lecture Notes
  in Artificial Intelligence and Lecture Notes in Bioinformatics)
  \textbf{11539}, 156--168 (2019)

\bibitem{uecker2014pde2path}
Uecker, H., Wetzel, D., Rademacher, J.D.: pde2path-{A} {M}atlab package for
  continuation and bifurcation in 2{D} elliptic systems.
\newblock Numerical Mathematics: Theory, Methods and Applications
  \textbf{7}(1), 58--106 (2014)

\bibitem{wallace2005adaptive}
Wallace, M., Wagg, D., Neild, S.: An adaptive polynomial based forward
  prediction algorithm for multi-actuator real-time dynamic substructuring.
\newblock Proceedings of the Royal Society A: Mathematical, Physical and
  Engineering Sciences \textbf{461}(2064), 3807--3826 (2005)

\bibitem{yuqing2018thesis}
Wang, Y.: Multidimensional continuation of families of periodic orbits.
\newblock Master's thesis, University of Illinois at Urbana-Champaign (2018)

\bibitem{watson1987algorithm}
Watson, L.T., Billups, S.C., Morgan, A.P.: Algorithm 652: {HOMPACK}: A suite of
  codes for globally convergent homotopy algorithms.
\newblock ACM Transactions on Mathematical Software (TOMS) \textbf{13}(3),
  281--310 (1987)

\bibitem{wyczalkowski2003optimization}
Wyczalkowski, M., Szeri, A.J.: Optimization of acoustic scattering from
  dual-frequency driven microbubbles at the difference frequency.
\newblock The Journal of the Acoustical Society of America \textbf{113}(6),
  3073--3079 (2003)

\bibitem{zhong2020global}
Zhong, J., Ross, S.D.: Global invariant manifolds delineating transition and
  escape dynamics in dissipative systems: an application to snap-through
  buckling.
\newblock Nonlinear Dynamics pp. 1--29 (2021)

\end{thebibliography}

\end{document}